\documentclass[11pt]{book}
\usepackage{amsfonts}
\usepackage{bbm}
\usepackage{mathrsfs}
\usepackage{amssymb}

\def\hang{\hangindent\parindent}
\def\textindent#1{\indent\llap{#1\enspace}\ignorespaces}
\def\re{\par\hang\textindent}

\def\v5{\vskip .5truecm}\def\QED{\hfill{$\Box$}}\def\hang{\hangindent\parindent}
\def\textindent#1{\indent\llap{#1\enspace}\ignorespaces}
\def\item{\par\hang\textindent}
\def \r{\rightarrow}\def\OV#1{\overline {#1}}

\def\mapdown#1{\llap{$\vcenter {\hbox {$\scriptstyle #1$}}$}
                                \Bigg\downarrow}

\def\mapright#1#2{\smash{\mathop{\longrightarrow}\limits^{#1}_{#2}}}
\def\NZ{\mathbb{N}}

\def\LM{{\bf LM}}\def\LT{{\bf
LT}}\def\KX{K\langle X\rangle} \def\KS{K\langle X\rangle}
\def\B{\mathscr{B}} \def\LC{{\bf LC}} \def\G{{\cal G}} 
\def\SUM^#1_#2{\displaystyle{\sum^{#1}_{#2}}}   \def\BE{\B (e)}
\def\PRCVE{\prec_{\varepsilon\hbox{-}gr}}\def\BV{\B (\varepsilon )}\def\PRCEGR{\prec_{e\hbox{\rm -}gr}}

\def\KS{K\langle X\rangle}

\begin{document}

\thispagestyle{empty}


\quad\quad~{\bf\LARGE\textsf Notes  on}\vskip 6pt

\quad\quad~{\bf \LARGE \textsf Gr\"obner Bases and Free}\vskip 6pt

\quad\quad~{\bf\LARGE\textsf Resolutions of Modules over}\vskip 6pt

\quad\quad~{\bf\LARGE\textsf Solvable Polynomial Algebras}

~~~~~~~~\underline{~~~~~~~~~~~~~~~~~~~~~~~~~~~~~~~~~~~~~~~~~~~~~~~~~~~~~~~~~~~~~~~~~~~~~~}
\vskip 11truecm \quad\quad~{\bf\large\textsf Huishi Li}\vskip 6pt
\quad\quad~{\bf\textsf Department of Applied Mathematics}\par

\quad\quad~{\bf\textsf College of Information Science \& 
Technology}\par

\quad\quad~{\bf\textsf Hainan University, Haikou 570228, China}\par 
\vfill

\newpage\thispagestyle{myheadings}\pagenumbering{roman}
\setcounter{page}{1}

\chapter*{Introduction}
\markboth{\rm Introduction}{\rm Introduction} \vskip 
2.5truecm\def\GR{Gr\"obner} 
\def\B{{\cal B}}

\hfill{\emph{Polynomials and power series},~~~~~~~~~~~~}\par

\hfill{~~~\emph{May they forever rule the world}.}\vskip 6pt 
\hfill{--- Shreeram S. Abhyankar}{\parindent=0pt\vskip 
6pt{\parindent=0pt

Since the late 1980s, the  Gr\"obner basis theory for commutative 
polynomial algebras and their modules (cf. [Bu1, 2],  [Sch], [BW], 
[AL2], [Fr\"ob], [KR1, 2]) has been successfully generalized to 
(noncommutative) solvable polynomial algebras and their modules (cf. 
[AL1], [Gal], [K-RW], [Kr2], [LW], [Li1], [Lev]). In particular, 
there have been a noncommutative version of  Buchberger's criterion 
and a noncommutative version of the Buchberger algorithm for 
checking and  computing {\parindent=.8truecm\par

\item{(1)} Gr\"obner bases of (one-sided) ideals of solvable polynomial 
algebras, and\par

\item{(2)} Gr\"obner bases of submodules of free modules over solvable 
polynomial algebras,\par}

while the noncommutative version of Buchberger algorithm has been 
implemented in some well-developed computer algebra systems, such as 
\textsc{Modula-2} [KP] and \textsc{Singular} [DGPS]. By the 
definition (see Section 1 of [K-RW],  or Definition 1.1.3 given in 
our Chapter 1), a solvable polynomial algebra is a polynomial-like 
but generally noncommutative algebra. Nowadays it is well known that 
the class of solvable polynomial algebras covers numerous 
significant noncommutative algebras such as enveloping algebras of 
finite dimensional Lie algebras, Weyl algebras (including  algebras 
of partial differential operators with polynomial coefficients over 
a field of characteristic 0), more generally a lot of operator 
algebras, iterated Ore extensions, and  quantum (quantized) 
algebras. So, comparing with the commutative case, the theory of 
Gr\"obner bases for solvable polynomial algebras and their modules   
has created great possibility of studying certain noncommutative 
algebras and  their modules in a  ``solvable setting" (see [Wik1] 
for an introduction to ``Decision problem" or ``Solvable problem", 
which may help us to better understand why a ``solvable polynomial 
algebra" deserves its name).}}\par

Along the lines in the literature for developing a computational 
(one-sided) ideal theory and more generally, a computational module 
theory over solvable polynomial algebras,  it seems that a rapid but 
relatively systematic and concrete introduction to the subject is 
worthwhile. Thereby,  we wrote these lecture notes so as to provide  
graduates and researchers (who are interested in {\it 
noncommutative} computational algebra) with an accessible reference 
on{\parindent=.5truecm\par

\item{$\bullet$}  an concise introduction to solvable polynomial algebras and the theory of 
Gr\"obner bases for submodules of free nodules over solvable 
polynomial algebras, and  \par

\item{$\bullet$} some details concerning applications of Gr\"obner bases in constructing finite free 
resolutions over an arbitrary solvable polynomial algebras, minimal 
finite $\NZ$-graded free resolutions over an $\NZ$-graded solvable 
polynomial algebra with the  degree-0 homogeneous part  being the 
ground field $K$, and minimal finite $\NZ$-filtered  free 
resolutions  over an $\NZ$-filtered solvable polynomial algebra  
(where the $\NZ$-filtration is determined by a positive-degree 
function). 
\par} \vskip 6pt

The first three chapters of these notes  grew out of a course of 
lectures given to graduate students at Hainan University, and the 
last two chapters are adapted from the author's recent research work 
[Li7] (or [Li5] arXiv:1401.5206v2 [math.RA], [Li6] arXiv:1401.5464 
[math.RA]). Also, at the level of module theory, these notes may be 
viewed as supplements of the sections concerning solvable polynomial 
algebras and their modules in  ([Li1], [Li2]). \vskip 6pt

Throughout the text, $K$ denotes a field, $K^*=K-\{ 0\}$; 
$\mathbb{N}$ denotes the additive monoid of all nonnegative 
integers, and $\mathbb{Z}$ denotes the additive group of all 
integers; all algebras are associative $K$-algebras with the 
multiplicative identity 1, and modules over an algebra are meant 
left unitary modules.\newpage

\chapter*{Contents}\markboth{\rm Contents}{\rm Contents} \vskip 
2.5truecm\def\GR{Gr\"obner} \def\B{{\cal B}}

\vskip 1.5truecm {\parindent=0pt

{\large\bf Introduction}\hfill{i}\vskip 6pt

{\large\bf 1.} {\large\bf Solvable Polynomial Algebras
}\hfill{1}}{\parindent=1.1truecm\par

\re{1.1} Definition and Basic Properties\hfill{}

\re{1.2} Left Gr\"obner Bases of Left Ideals\hfill{}

\re{1.3} The Noetherianess \hfill{}

\re{1.4} A Constructive Characterization\hfill{}}

{\parindent=0pt\vskip 6pt {\large\bf 2.} {\large\bf Left Gr\"obner
Bases for Modules}\hfill{33}}{\parindent=1.1truecm\par

\re{2.1} Left Monomial Orderings on Free Modules \hfill{}

\re{2.2} Left Gr\"obner Bases of Submodules \hfill{}

\re{2.3} The Noncommutative  Buchberger Algorithm \hfill{}}

{\parindent=0pt\vskip 6pt {\large\bf 3.} {\large\bf Finite Free
Resolutions}\hfill{51}}{\parindent=1.1truecm\par

\re{3.1} Computation of Syzygies \hfill{}

\re{3.2} Computation of Finite Free Resolutions\hfill{}

\re{3.3} Global Dimension and Stability \hfill{~}

\re{3.4} Calculation of p.dim$_AM$ \hfill{}}

{\parindent=0pt\vskip 6pt {\large\bf 4.} {\large\bf Minimal Finite
Graded Free Resolutions}\hfill{69}} {\parindent=1.1truecm\par

\re{4.1} $\NZ$-graded Solvable Polynomial Algebras \hfill{}

\re{4.2} $\NZ$-Graded Free Modules\hfill{}

\re{4.3} Computation of Minimal Homogeneous Generating Sets\hfill{}

\re{4.4} Computation of Minimal Finite Graded Free
Resolutions\hfill{}}

{\parindent=0pt\vskip 6pt

{\large\bf 5.} {\large\bf  Minimal Finite Filtered Free
Resolutions}\hfill{89}} {\parindent=1.1truecm\par

\re{5.1} $\NZ$-Filtered Solvable Polynomial Algebras \hfill{}

\re{5.2} $\NZ$-Filtered Free Modules\hfill{}

\re{5.3} Filtered-Graded Transfer of Gr\"obner Bases for
Modules\hfill{}

\re{5.4} F-Bases and Standard Bases with Respect to  Good
Filtration\hfill{}

\re{5.5} Computation of Minimal F-Bases and Minimal Standard
Bases\hfill{}

\re{5.6} Minimal Filtered Free Resolutions and Their
Uniqueness\hfill{}

\re{5.7} Computation of  Minimal Finite Filtered Free
Resolutions\hfill{}}

{\parindent=0pt\vskip 6pt

{\large\bf 6.} {\large\bf  Some Examples of 
Applications}\hfill{125}} {\parindent=1.1truecm\par}

\newpage\pagestyle{myheadings} \pagenumbering{arabic}\setcounter{page}{1}

\chapter*{1. Solvable Polynomial\\ \hskip 1.2truecm  Algebras}
\markboth{\rm Solvable Polynomial Algebras}{\rm Solvable Polynomial
Algebras} \vskip 2.5truecm\def\GR{Gr\"obner} \def\B{{\cal B}}

In this chapter we give a concise but easily accessible introduction
to solvable polynomial algebras and the theory of left Gr\"obner
bases for left ideals of such algebras. In the first section we 
introduce the  definition of solvable polynomial algebras and some 
typical examples; also from the definition we highlight some often 
used properties of solvable polynomial algebras.  In Section 2, we 
introduce left Gr\"obner bases for left ideals of solvable 
polynomial algebras via a left division algorithm, and we discuss 
some basic properties of left Gr\"obner bases. In Section 3, on the 
basis of Dickson's lemma we show that every left ideal of a solvable 
polynomial algebra has a finite Gr\"obner basis, thereby solvable 
polynomial algebras are left Noetherian. Since every solvable 
polynomial algebra has a (two-sided) monomial ordering, a right 
division algorithm and a theory of right Gr\"obner bases for right 
ideals hold true as well, it follows that every solvable polynomial 
algebra is also right Noetherian. Concerning the algorithmic 
approach to computing a finite left  Gr\"obner basis, it will be a 
job of Chapter 2 for modules. In the final  Section 4, by employing 
Gr\"obner bases of two-sided ideals in free algebras we give a 
characterization of solvable polynomial algebras, so that such  
algebras are completely recognizable and constructible in a 
computational way. \par

The main references of this chapter are [AL1], [Gal], [K-RW], [Kr2],
[LW], [Li1], [Li4], [DGPS], [Ber2], [Mor], [Gr], [Uf].

\section*{1.1. Definition and Basic Properties}

Let $K$ be a field and let $A=K[a_1,\ldots ,a_n]$ be a finitely
generated $K$-algebra with the set of generators $\{ a_1,\ldots
,a_n\}$, that is, $A$ is an associative ring with the multiplicative
identity 1, every element $a\in A$ is expressed as a finite sum of
the form $a=\sum_i\lambda_ia_{i_1}^{\alpha_1}\cdots
a_{i_t}^{\alpha_t}$ with $\lambda_i\in K$, $a_{i_j}\in\{ a_1,\ldots
,a_n\}$, $\alpha_j\in\NZ$, $t\ge 1$, and $A$ is also a $K$-vector
space with respect to its additive operation, such that
$\lambda(ab)=(\lambda a)b=a(\lambda b)$ holds for all $\lambda\in K$
and $a,b\in A$. Moreover, we assume that  any proper subset of the
given generating set $\{ a_1,\ldots ,a_n\}$ of $A$ cannot generate
$A$ as a $K$-algebra, i.e., the given set of generators is {\it
minimal}.\par

If, for some permutation $\tau =i_1i_2\cdots i_n$ of $1,2,\ldots
,n$, the set $\B =\{ a^{\alpha}=a_{i_1}^{\alpha_1}\cdots
a_{i_n}^{\alpha_n}~|~\alpha =(\alpha_1,\ldots ,\alpha_n)\in\NZ^n\}
,$ forms a $K$-basis of $A$, then $\B$ is referred to as a {\it PBW
$K$-basis} of $A$ (where the phrase ``PBW $K$-basis" is  abbreviated
from the well-known {\it Poincar\'e-Birkhoff-Witt Theorem}
concerning the standard $K$-basis of the enveloping algebra of a Lie
algebra, e.g., see [Hu], P. 92).  It is clear that if $A$ has a PBW
$K$-basis, then we can always assume that $i_1=1,\ldots ,i_n=n$.
Thus, we make the following convention once for
all.{\parindent=0pt\v5

{\bf Convention} From now on in this paper, if we say that an
algebra $A$ has the PBW $K$-basis $\B$, then it means that
$$\B =\{ a^{\alpha}=a_{1}^{\alpha_1}\cdots
a_{n}^{\alpha_n}~|~\alpha =(\alpha_1,\ldots ,\alpha_n)\in\NZ^n\} .$$
Moreover, adopting the commonly used terminology in computational
algebra, elements of $\B$ are referred to as {\it monomials} of
$A$.\v5

{\bf Monomial ordering and admissible system}\vskip 6pt

Suppose that the $K$-algebra $A=K[a_1,\ldots ,a_n]$ has the PBW
$K$-basis $\B$ as presented above and that  $\prec$ is a total
ordering on $\B$. Then every nonzero element $f\in A$ has a unique
expression
$$\begin{array}{rcl} f&=&\lambda_1a^{\alpha (1)}+\lambda_2a^{\alpha (2)}+\cdots +\lambda_ma^{\alpha (m)},\\
&{~}&\hbox{where}~ \lambda_j\in K^*,~a^{\alpha
(j)}=a_1^{\alpha_{1j}}a_2^{\alpha_{2j}}\cdots a_n^{\alpha_{nj}}\in\B
,~1\le j\le m.\end{array}$$ If  $a^{\alpha (1)}\prec a^{\alpha
(2)}\prec\cdots \prec a^{\alpha (m)}$ in the above representation,
then the {\it leading monomial of $f$} is defined as $\LM
(f)=a^{\alpha (m)}$, the {\it leading coefficient of $F$} is defined
as $\LC (f)=\lambda_m$, and the {\it leading term of $f$} is defined
as $\LT (f)=\lambda_ma^{\alpha (m)}$. {\parindent=0pt\v5

{\bf 1.1.1. Definition}  Suppose that the $K$-algebra
$A=K[a_1,\ldots ,a_n]$ has the PBW $K$-basis $\B$. If $\prec$ is a
total ordering on $\B$ that satisfies the following three
conditions:}{\parindent=1.34truecm\par

\item{(1)} $\prec$ is a well-ordering (i.e. every nonempty subset of $\B$ has a minimal element);\par

\item{(2)} For $a^{\gamma},a^{\alpha},a^{\beta},a^{\eta}\in\B$, if $a^{\gamma}\ne 1$, $a^{\beta}\ne
a^{\gamma}$, and $a^{\gamma}=\LM (a^{\alpha}a^{\beta}a^{\eta})$,
then $a^{\beta}\prec a^{\gamma}$ (thereby $1\prec a^{\gamma}$ for
all $a^{\gamma}\ne 1$),\par

\item{(3)} For $a^{\gamma},a^{\alpha},a^{\beta}, a^{\eta}\in\B$, if
$a^{\alpha}\prec a^{\beta}$, $\LM (a^{\gamma}a^{\alpha}a^{\eta})\ne
0$, and $\LM (a^{\gamma}a^{\beta}a^{\eta})\not\in \{ 0,1\}$, then
$\LM (a^{\gamma}a^{\alpha}a^{\eta})\prec\LM
(a^{\gamma}a^{\beta}a^{\eta})$;\par

}{\parindent=0pt

then $\prec$ is called a {\it monomial ordering} on $\B$ (or a
monomial ordering on $A$). }}\par

If $\prec$ is a monomial ordering on $\B$, then the data $(\B
,\prec)$ is referred to as an {\it admissible system} of $A$.
{\parindent=0pt\v5

{\bf Remark.} (i) Definition 1.1.1 above is  borrowed from the 
theory of Gr\"obner bases  for  general finitely generated 
$K$-algebras, in which the algebras considered may  be 
noncommutative, may have divisors of zero, and the $K$-bases used 
may  not be a PBW basis, but with a (one-sided, two-sided) monomial 
ordering that, theoretically, enables such algebras  have a 
(one-sided, two-sided) Gr\"obner basis theory. For more details on 
this topic, one may referrer to ([Li2], Section 1 of Chapter 3 and 
Section 3 of Chapter 8). Also, to see the essential difference 
between Definition 1.2.1 and the classical definition of a monomial 
ordering in the commutative case, one may refer to (Definition 1.4.1 
and the proof of Theorem 1.4.6 given in [AL2]). \vskip 6pt

(ii) The conditions (2) and (3) listed in Definition 1.1.1 mean that
$\prec$ is {\it two-sided compatible with the multiplication
operation of the algebra $A$}. As one will see soon, that the use of
a (two-sided) monomial ordering $\prec$  on a solvable polynomial
algebra $A$ first  guarantees  that {\it $A$ is a domain}, and
furthermore guarantees an effective (left, right,  two-sided) finite
Gr\"obner basis theory for $A$.    }\v5

Note that if a $K$-algebra $A=K[a_1,\ldots ,a_n]$ has the PBW
$K$-basis $\B =\{ a^{\alpha}=a_1^{\alpha_1}\cdots
a_n^{\alpha_n}~|~\alpha =(\alpha_1,\ldots ,\alpha_n)\in\NZ^n\}$,
then for any given $n$-tuple $(m_1,\ldots ,m_n)\in\mathbb{N}^n$, a
{\it weighted degree function} $d(~)$ is well defined on nonzero
elements of $A$, namely, for each $a^{\alpha}=a_1^{\alpha_1}\cdots
a_n^{\alpha_n}\in\B$, $d(a^{\alpha})=m_1\alpha_1+\cdots
+m_n\alpha_n,$ and for each nonzero
$f=\sum_{i=1}^s\lambda_ia^{\alpha (i)}\in A$ with $\lambda_i\in K^*$
and $a^{\alpha (i)}\in\B$, $$d(f)=\max\{ d(a^{\alpha (i)})~|~1\le
i\le s\}.$$ If $d(a_i)=m_i>0$ for $1\le i\le n$, then $d(~)$ is
referred to as a {\it positive-degree function} on $A$.
{\parindent=0pt\v5

{\bf 1.1.2. Definition} Let  $d(~)$ be a  positive-degree function
on $A$. If $\prec$ is a monomial ordering on $\B$ such that for all
$a^{\alpha},a^{\beta}\in\B$,
$$a^{\alpha}\prec a^{\beta}~\hbox{implies}~d(a^{\alpha})\le d(a^{\beta}),$$
then we call $\prec$ a {\it graded monomial ordering} with respect
to $d(~)$. \v5

{\bf Convention} Unless otherwise stated, from now on in this book
we always use $\prec_{gr}$ to denote a  graded monomial ordering
with respect to a positive-degree function on $A$.}\v5

As one may see from  the literature (or loc. cit) that in both the
commutative and noncommutative computational algebra, the most
popularly used graded monomial orderings on an algebra
$A=K[a_1,\ldots ,a_n]$ with the PBW $K$-basis $\B$ are those graded
(reverse) lexicographic orderings with respect to the degree
function $d(~)$ such that $d(a_i)=1$, $1\le i\le
n$.{\parindent=0pt\v5

{\bf Definition of solvable polynomial algebras and examples}\vskip
6pt

Originally, a (noncommutative) solvable polynomial algebra (or an
algebra of solvable type) $R'$ was defined in [K-RW] by first fixing
a monomial ordering $\prec$ on the standard $K$-basis $\mathscr{B}
=\{ X_1^{\alpha_1}\cdots X_n^{\alpha_n}~|~\alpha_i\in\NZ\}$ of the
commutative polynomial algebra  $R=K[X_1,\ldots ,X_n]$ in $n$
variables $X_1,\ldots ,X_n$ over a field $K$, and then introducing a
new multiplication $*$ on $R$, such that certain axioms ([K-RW],
AXIOMS 1.2) are satisfied. In [LW], a solvable polynomial algebra
was redefined in the formal language of associative $K$-algebras, as
follows. \v5

{\bf 1.1.3. Definition} If a $K$-algebra $A=K[a_1,\ldots ,a_n]$
satisfies the following two conditions:\par

(S1) $A$ has the PBW $K$-basis $\B$;

(S2) There is a monomial ordering $\prec$ on $\B$, i.e., $(\B ,\prec
)$ is an admissible system of $A$,  such that for all
$a^{\alpha}=a_1^{\alpha_1}\cdots a_n^{\alpha_n}$,
$a^{\beta}=a_1^{\beta_1}\cdots a^{\beta_n}_n\in\B$,
$$\begin{array}{rcl} a^{\alpha}a^{\beta}&=&\lambda_{\alpha ,\beta}a^{\alpha +\beta}+f_{\alpha ,\beta},\\
&{~}&\hbox{where}~\lambda_{\alpha ,\beta}\in K^*,~a^{\alpha
+\beta}=a_1^{\alpha_1+\beta_1}\cdots
a_n^{\alpha_n+\beta_n},~\hbox{and}\\
&{~}&f_{\alpha ,\beta}\in K\hbox{-span}\B~\hbox{with}~\LM (f_{\alpha
,\beta})\prec a^{\alpha +\beta}~\hbox{if}~f_{\alpha ,\beta}\ne
0,\end{array}$$ or alternatively, such that for all generators 
$a_1,\ldots ,a_n$ of $A$ and $i<j$, 
$$\begin{array}{rcl} a_ja_i&=&\lambda_{ji}a_ia_j+f_{ji},\\
&{~}&\hbox{where}~\lambda_{ji}\in K^*,~\hbox{and}\\
&{~}&f_{ji}\in K\hbox{-span}\B~\hbox{with}~\LM (f_{ji})\prec 
a_ia_j~\hbox{if}~f_{ji}\ne 0,\end{array}$$                                  
then $A$ is said to be a {\it solvable polynomial algebra}.\v5

{\bf Remark} As one will see later,  that the pure-ring-theoretical 
definition of a solvable polynomial algebra above  may at least 
provide us with the following two advantages.\par

(i) Solvable polynomial algebras can be characterized in a 
constructive way (see subsequent Section 1.4), so that more 
noncommutative algebras (and hence their modules) can be studied in 
a  ``solvable setting" (see [Wik1] for an introduction on ``Decision 
problem" (or ``Solvable problem")).\par

(ii) It is quite helpful in determining whether a solvable 
polynomial algebra $A$ is an $\NZ$-graded algebra as specified in 
(Section 4.1 of Chapter 4), or an $\NZ$-filtered algebra as 
specified in (Section 5.1 of Chapter 5) respectively. \v5

{\bf Example} (1) By Definition 1.1.3, every commutative polynomial
algebra $K[x_1,\ldots ,x_n]$ in variables $x_1,\ldots ,x_n$ over a
field $K$ is trivially a solvable polynomial algebra.}\vskip 6pt

From [KR-W] and [Li1] we recall several most well-known
noncommutative solvable polynomial algebras as
follows.{\parindent=0pt\vskip 6pt

{\bf Example} (2)  The $n$th Weyl algebra\par

The $n$th Weyl algebra $A_n(k)$ over a field $K$ is defined to be
the $K$-algebra generated by $2n$ generators $x_1,...,x_n,y_1,...,
y_n$ subject to the relations:
$$\begin{array}{ll}
x_ix_j=x_jx_i,~y_iy_j=y_jy_i,&1\le i<j\le n,\\
y_jx_i=x_iy_j+\delta_{ij}~\hbox{the~Kronecker~delta} , &1\le i,j\le
n.
\end{array}$$
Historically, the Weyl algebra is the first ``quantum algebra''
([Dir] 1926, [Wey] 1928). It is a well-known fact that if char$K=0$,
then $A_n(k)$ coincides with the algebra of linear partial
differential operators $K[x_1,\ldots ,x_n,\partial_1,\ldots
,\partial_n]$ of the polynomial ring $K[x_1,...,x_n]$ (or the ring
of  polynomial functions of the affine $n$-space {\bf A}$^n_k$),
where $\partial_i=\frac{\partial~}{\partial x_i}$, $1\le i\le
n$.\vskip 6pt

{\bf Example} {(3) The additive analogue of the Weyl algebra}\par

This algebra was introduced in Quantum Physics in ([Kur] 1980) and
studied in ([JBS] 1981), that is, the algebra $A_n(q_1,...,q_n)$
generated over a field $K$ by $x_1,...,x_n,y_1,...,y_n$ subject to
the relations:
$$\begin{array}{ll}
x_ix_j=x_jx_i,~y_iy_j=y_jy_i,~&1\le i<j\le n,\\
y_ix_i=q_i x_iy_i+1,&1\le i\le n,\\
x_jy_i=y_ix_j,&i\ne j, \end{array}$$ where $q_i\in K^*$.  If
$q_i=q\ne 0$, $i=1,...,n$, then this algebra becomes the algebra of
$q$-{\it differential operators}. \vskip 6pt\def\G{{\bf g}}

{\bf Example} {(4) The multiplicative analogue of the Weyl 
algebra}\par

This is the algebra stemming from ([Jat] 1984 and [MP] 1988, where
one may see why this algebra deserves its title), that is, the
algebra ${\cal O}_n(\lambda_{ji})$ generated over a field $K$ by
$x_1,...,x_n$ subject to the relations:
$$x_jx_i=\lambda_{ji}x_ix_j,\quad 1\le i<j\le n,$$
where $\lambda_{ji}\in K^*$.  If $n=2$, then ${\cal 
O}_2(\lambda_{21})$ is the {\it quantum plane} in the sense of Manin 
[Man]. If $\lambda_{ji}=q^{-2}\ne 0$ for some $q\in K^*$ and all 
$i<j$, then ${\cal O}_n(\lambda_{ji})$ becomes the coordinate ring 
of the so called {\it quantum affine $n$-space} (see [Sm]).\vskip 
6pt

{\bf Example} {(5) The enveloping algebra of a finite dimensional 
Lie algebra}\par

Let $\G$ be a finite dimensional vector space over the field $k$
with basis $\{x_1,...,x_n\}$ where $n=$ dim$_k\G$. If there is a
binary operation on $\G$, called the {\it bracket product} and
denoted $[~,~]$, which is bilinear, i.e., for $a,b,c\in\G$,
$\lambda\in k$,
$$\begin{array}{l}
[a+b,c]=[a,c]+[b,c]\\

 [a,b+c]=[a,b]+[a,c]\\

\lambda [a,b]=[\lambda a,b]=[a,\lambda b],
\end{array}$$
 and satisfies:
$$\begin{array}{l}
[a,b]=-[b,a],\quad a, b\in\G ,\\

[[a,b],c]+[[c,a],b]+[[b,c],a]=0 \quad \hbox{Jacobi ~identity},
~a,b,c\in\G ,
\end{array}$$
then $\G$ is called a finite dimensional Lie algebra over $k$. Note
that $[~,~]$ need not satify the associative low. If $[a,b]=[b,a]$
for every $a,b$ in a Lie algebra {\bf g}, {\bf g} is called {\it
abelian}. The  enveloping algebra of $\G$, denoted $U(\G )$, is
defined to be the {\it associative} $k$-algebra generated by $
x_1,...,x_n$ subject to the relations:
$$x_jx_i-x_ix_j=[x_j,x_i],\quad 1\le i<j\le n.$$
For instance, the Heisenberg Lie algebra {\bf h} has the $k$-basis
$\{ x_i,y_j,z~|~ i,j=1,...,n\}$ and the bracket product is given by
$$\begin{array}{ll}
[x_i,y_i]=z,& 1\le i\le n,\\

[x_i,x_j]=[x_i,y_j]=[y_i,y_j]=0,&i\ne j,\\

[z,x_i]=[z,y_i]=0,& 1\le i\le n.
\end{array}$$
\vskip 6pt

{\bf Example} {(6) The $q$-Heisenberg algebra}\par

This is the algebra stemming from ([Ber] 1992, [Ros] 1995) which has
its root in $q$-calculus (e.g., [Wal] 1985), that is, the algebra
${\bf h}_n(q)$ generated over a field $K$ by the set of elements $\{
x_i,y_i,z_i~|~i=1,...,n\}$ subject to the relations:
$$\begin{array}{ll}
x_ix_j=x_jx_i,~y_iy_j=y_jy_i,~z_jz_i=z_iz_j,& 1\le i<j\le n,\\
x_iz_i=qz_ix_i,&1\le i\le n,\\
z_iy_i=qy_iz_i,&1\le i\le n,\\
x_iy_i=q^{-1}y_ix_i+z_i,&1\le i\le n,\\
x_iy_j=y_jx_i,~x_iz_j=z_jx_i,~y_iz_j=z_jy_i,&i\ne j,\\
\end{array}$$
where $q\in K^*$. \vskip 6pt\def\X{{\bf x}}

{\bf Example} {(7) The algebra of $2\times 2$ quantum matrices}\par 
This is the algebra $M_q(2,K)$ introduced in ([Man] 1988), which is  
generated over a field $K$ by $a,b,c,d$ subject to the relations:
$$\begin{array}{lll}
ba=qab,&ca=qac,&dc=qcd,\\
db=qbd,&cb=bc,&da-ad=(q-q^{-1})bc,\end{array}$$ where $q\in K^*$.
\vskip 6pt

{\bf Example} {(8) The algebra {\bf U} in constructing Hayashi 
algebra}\par

In order to get bosonic representations for the types of ${\bf A}_n$
and ${\bf C}_n$ of the Drinfield-Jimbo quantum algebras, Hayashi
introduced in ([Hay] 1990) the $q$-Weyl algebra ${\cal A}^-_q$,
which is constructed as follows (see [Ber1]). Let ${\bf U}$ be the
algebra generated over the field $K=\mathbb{C}$ by the set of
elements $\{x_i,y_i,z_i~|~i=1,...,n\}$ subject to the relations:
$$\begin{array}{ll}
x_jx_i=x_ix_j,~y_jy_i=y_iy_j,~z_jz_i=z_iz_j,&1\le i<j\le n,\\
x_jy_i=q^{-\delta_{ij}}y_ix_j,~z_jx_i=q^{-\delta_{ij}}x_iz_j,&1\le
i,j\le n,\\
z_jy_i=y_iz_j,&i\ne j,\\
z_iy_i-q^2y_iz_i=-q^2x_i^2,&1\le i\le n,\end{array}$$ where $q\in
K^*$. Then ${\cal A}^-_q=S^{-1}{\bf U}$, }\v5

One may also use the technique  given in the subsequent Section 1.4  
to directly verify that the algebras listed above are solvable 
polynomial algebras. {\parindent=0pt\v5

{\bf Basic properties}\vskip 6pt

By Definition 1.1.1 and Definition 1.1.3,  the properties listed in
the next two propositions are straightforward. \v5

{\bf 1.1.4. Proposition}  Let $A=K[a_1,\ldots ,a_n]$ be a solvable
polynomial algebra with   admissible system $(\B ,\prec )$. The
following statements hold.\par

(i) If $f,g\in A$ with $\LM (f)=a^{\alpha}$, $\LM (g)=a^{\beta}$,
then
$$\LM (fg)=\LM (\LM (f)\LM (g))=\LM (a^{\alpha}a^{\beta})=a^{\alpha +\beta}.$$\par

(ii) $A$ is a domain, that is, $A$ has no (left and right) divisors
of zero.\vskip 6pt

{\bf Proof} Exercise.\QED\v5

{\bf 1.1.5. Proposition} Let $A_1=K[a_1,\ldots ,a_n]$ and
$A_2=K[b_1,\ldots ,b_m]$ be solvable polynomial algebras with
admissible systems $(\B_1,\prec_1)$ and $(\B_2,\prec_2)$
respectively. Then $A=A_1\otimes_KA_2$ is a solvable polynomial
algebra with the admissible system $(\B ,\prec)$, where
$\B=\{a^{\alpha}\otimes
b^{\beta}~|~a^{\alpha}\in\B_1,~b^{\beta}\in\B_2\}$, while $\prec$ is
defined on $\B$ subject to the rule: for $a^{\alpha}\otimes
b^{\beta}$, $a^{\gamma}\otimes b^{\eta}\in\B$,
$$a^{\alpha}\otimes
b^{\beta}\prec a^{\gamma}\otimes
b^{\eta}\Leftrightarrow\left\{\begin{array}{l}
a^{\alpha}\prec_1a^{\gamma};\\
\hbox{or}\\
a^{\alpha}=a^{\gamma}~\hbox{and}~b^{\beta}\prec_2b^{\eta}.\end{array}\right.$$\vskip 
6pt

{\bf Proof} Exercise. }\v5

\section*{1.2. Left Gr\"obner Bases of Left Ideals}

Let $A=K[a_1,\ldots ,a_n]$ be a solvable polynomial algebra with
admissible system $(\B ,\prec )$. In this section we introduce left
Gr\"obner bases for left ideals of $A$ via a left division algorithm
in $A$, and we record some basic facts determined by left Gr\"obner
bases. Moreover, minimal left Gr\"obner bases and reduced left
Gr\"obner bases are discussed. {\parindent=0pt\v5

{\bf Left division algorithm}\vskip 6pt

Let $a^{\alpha},a^{\beta}\in\B$. We say that $a^{\alpha}$ {\it
divides $a^{\beta}$ from the left side}, denoted $a^{\alpha}|_{_{\rm
L}}a^{\beta}$, if there exists $a^{\gamma}\in\B$ such that
$$a^{\beta}=\LM (a^{\gamma}a^{\alpha}) .$$
It follows from Proposition 1.1.4(i) that the division defined above
is implementable.}\v5

Let $F$ be a nonempty subset of $A$. Then the division of monomials
defined above yields a subset of $\B$:
$$\mathcal {N}(F)=\{ a^{\alpha}\in\B~|~\LM (f)\not{|_{_{\rm L}}}~a^{\alpha},~f\in F\}.$$
If $F=\{ g\}$ consists of a single element $g$, then we simply write
$\mathcal {N}(g)$ in place of $\mathcal {N}(F)$. Also we write
$K$-span$\mathcal {N}(F)$ for the $K$-subspace of $A$ spanned by
$\mathcal {N}(F)$.{\parindent=0pt\v5

{\bf 1.2.1. Definition} Elements of $\mathcal {N}(F)$ are referred
to as {\it normal monomials} (mod $F$). Elements of
$K$-span$\mathcal {N}(F)$ are referred to as {\it normal elements}
(mod $F$).}\v5

In view of Proposition 1.1.4(i), the left division we defined for
monomials in $\B$ can naturally be used to define a left division
procedure for elements in $A$. More precisely, let $f,g\in A$ with
$\LC (f)=\mu \ne 0$, $\LC (g)=\lambda\ne 0$. If $\LM (g)|_{_{\rm
L}}\LM (f)$, i.e., there exists $a^{\alpha}\in\B$ such that $\LM
(f)=\LM (a^{\alpha}\LM (g))$, then put $f_1=f-\lambda^{-1}\mu
a^{\alpha}g$; otherwise, put $f_1=f-\LT (f)$. Note that in both
cases we have $f_1=0$, or $f_1\ne 0$ and $\LM (f_1)\prec\LM (f)$. At
this stage, let us refer to such a procedure of canceling the 
leading term of $f$  as the {\it left division procedure by $g$}.  
With $f:=f_1\ne 0$, we can repeat the left division procedure by $g$ 
and so on. This returns successively   a descending sequence $$\LM 
(f)\succ\LM (f_1)\succ\LM (f_2)\succ\cdots .$$  Since $\prec$ is a 
well-ordering, it follows that such a division procedure terminates 
after a finite number of repetitions, and consequently $f$ is 
expressed as
$$f=qg+r,$$
where $q,r\in A$ with $r\in K$-span$\mathcal {N}(g)$, i.e., $r$ is
normal (mod $g$), such that either $\LM (f)=\LM (qg)$ or $\LM
(f)=\LM (r)$.\par

Furthermore, the left division procedure demonstrated above can be 
extended to a left division procedure  by a finite subset $G=\{ 
g_1,\ldots ,g_s\}$ in $A$, which yields the following  division 
algorithm:{\parindent=0pt\vskip 6pt

\underline{\bf Algorithm-LDIV~~~~~~
~~~~~~~~~~~~~~~~~~~~~~~~~~~~~~~~~~~~~~~~~~~~~~~~~~~~~~}\vskip 6pt

\textsc{INPUT}: ~$ f,~G=\{g_1,\ldots,g_s\}~\hbox{with}~g_i\ne 
0~(1\le i\le s)$\par \textsc{OUTPUT}: ~$q_1,\ldots ,q_s,r~\hbox{such 
that}~r\in K\hbox{-span}\mathcal {N}(G),f=\sum^s_{i=1}q_ig_i+r,$\par
~~~~~~~~~~~~~~~~~~$\LM (q_ig_i)\preceq\LM (f)~\hbox{for}~q_i\ne 
0,~\LM (r)\preceq\LM (f)~\hbox{if}~r\ne 0$\par 
\textsc{INITIALIZATION}:~$q_1:=0,~q_2:=0,~\cdots
,~q_s:=0;~r:=0;~h:=f$ \par

\textsc{BEGIN}\par ~~~~~\textsc{WHILE}~$h\ne 0$~\textsc{DO}\par 
~~~~~~~~~~\hbox{IF}~\hbox{there exist}~$i$~\hbox{and}~$a^{\alpha 
(i)}\in\B$\par ~~~~~~~~~~~~~~\hbox{such that}~$\LM (h)=\LM 
(a^{\alpha (i)}\LM (g_i))$~\textsc{THEN}\par 
~~~~~~~~~~~~~~\hbox{choose}~$i$~\hbox{least such that}~$\LM (h)=\LM 
(a^{\alpha (i)}\LM (g_i))$\par ~~~~~~~~~~$q_i:=q_i+\LC (g_i)^{-1}\LC 
(h)a^{\alpha (i)}$\par ~~~~~~~~~~$h:=h-\LC (g_i)^{-1}\LC 
(h)a^{\alpha (i)}g_i$\par ~~~~~~~~~~\hbox{ELSE}\par 
~~~~~~~~~~$r:=r+\LT (h)$\par ~~~~~~~~~~$h:=h-\LT (h)$\par 
~~~~~~~~~~\hbox{END}\par ~~~~~\textsc{END}\par                   
\textsc{END}\par\vskip -.15truecm\underline{ 
~~~~~~~~~~~~~~~~~~~~~~~~~~~~~~~~~~~~~~~~~~~~~~~~~~~~~~~~~~~~~~~~~~~~~~~~~~~~~~~~~~~~~~~~~~~~~~} 
\v5

{\bf 1.2.2. Definition} The element $r$ obtained in {\bf
Algorithm-LDIV} is called a {\it remainder} of $f$ on left division
by $G$, and is denoted by $\OV f^{G}$, i.e., $\OV f^G=r$. If $\OV
f^G=0$, then we say that $f$ {\it is reduced to $0$} (mod $G$).\v5

{\bf Remark} Note that the element $r$ obtained in {\bf
Algorithm-LDIV} depends on how the order of elements in $G$ is
arranged. This implies that $r$ may be different if a different
order for elements of $G$ is given (acturally as in the commutative
case, e.g. [AL2], P.31). That is why we use the phrase ``a
remainder" in the above definition.}\v5

Summing up, we have reached the following{\parindent=0pt\v5

{\bf 1.2.3. Theorem} Given a set of nonzero elements $G=\{
g_1,\ldots ,g_s\}$ and $f$ in $A$, the {\bf Algorithm-LDIV} produces
elements $q_1,\ldots ,q_s,~r\in A$ with $r\in K$-span$\mathcal
{N}(G)$, such that $f=\sum^s_{i=1}q_ig_i+r$ and $\LM
(q_ig_i)\preceq\LM (f)$ whenever $q_i\ne 0$, $\LM (r)\preceq\LM (f)$
if $r\ne 0$.\QED\v5

{\bf Left Gr\"obner bases}\vskip 6pt

Theoretically the left division procedure by using a finite subset
$G$ of nonzero elements in $A$ can be extended to a left division
procedure by means of an {\it arbitrary proper subset} $G$ of
nonzero elements, for,  it is a true statement that if $f\in A$ with
$\LM (f)\ne 0$, then either there exists $g\in G$ such that $\LM
(g)|_{_{\rm L}}\LM (f)$ or such $g$ does not exist.  This leads to
the following\v5

{\bf 1.2.4. Proposition} Let $N$ be a left ideal of $A$ and $G$ a
proper subset of nonzero elements in $N$. The following two
statements are equivalent:\par

(i) If $f\in N$ and $f\ne 0$, then there exists $g\in G$ such that
$\LM (g)|_{_{\rm L}}\LM (f)$.\par

(ii) Every nonzero $f\in N$ has a representation
$f=\sum^s_{j=1}q_ig_{i_j}$ with $q_i\in A$ and $g_{i_j}\in G$, such
that $\LM (q_ig_{i_j})\preceq\LM (f)$ whenever $q_i\ne 0$.\QED\v5

{\bf 1.2.5. Definition}  Let $N$ be a left ideal of $A$. With
respect to a given monomial ordering $\prec$ on $\B$, a proper
subset $\G$ of nonzero elements in $N$ is said to be a {\it left
Gr\"obner basis } of $N$ if $\G$ satisfies one of the equivalent
conditions of Proposition 1.2.4.}
\par

If $\G$ is a left Gr\"obner basis of $N$, then the expression
$f=\sum^s_{j=1}q_ig_{i_j}$ appeared in Proposition 1.2.4(ii) is
called a {\it left Gr\"obner representation} of $f$. \v5

Clearly, if $N$ is a left ideal of $A$ and $N$ has a left Gr\"obner
basis $\G$, then $\G$ is certainly a generating set of $N$, i.e.,
$N=\sum_{g\in\G}Ag$. But the converse is not necessarily true. For
instance, in the solvable polynomial algebra
$A=\mathbb{C}[a_1,a_2,a_3]$ generated by $\{ a_1,a_2,a_3\}$ subject
to the relations
$$a_2a_1=3a_1a_2,\quad a_3a_1=a_1a_3,\quad a_3a_2=5a_2a_3,$$
let $g_1=a_1^2a_2-a_3$, $g_2=a_2$, and $N=Ag_1+Ag_2$. Then since
$a_3\in N$, $G=\{ g_1,g_2\}$ is not a left Gr\"obner basis with
respect to any given monomial ordering $\prec$ on $\B$ such that
$\LM (g_1)=a_1^2a_2$, $\LM (g_2)=a_2$.{\parindent=0pt\v5

{\bf The existence of left Gr\"obner bases}\vskip 6pt

Noticing $1\in\B$, if we define on $\B$ the ordering:
$$a^{\alpha}\prec 'a^{\beta}\Leftrightarrow a^{\alpha}|_{_{\rm
L}}a^{\beta},\quad a^{\alpha},~a^{\beta}~\hbox{in}~\B ,$$ then, by
Definition 1.1.1, Proposition 1.1.4(i) and the divisibility  defined
for monomials,  it is an easy exercise to check that $\prec '$ is
reflexive, anti-symmetric, transitive, and moreover,
$$a^{\alpha}\prec 'a^{\beta}~\hbox{implies}~a^{\alpha}\prec a^{\beta}.$$
Since the given monomial ordering $\prec$ is a well-ordering on
$\B$, it follows that every nonempty subset of $\B$ has a minimal
element with respect to the ordering $\prec '$ on $\B$. }\v5

Now, let $N\ne \{ 0\}$ be a left ideal of $A$, $\LM (N)=\{ \LM
(f)~|~f\in N\}$, and $\Omega =\{ a^{\alpha}~|~a^{\alpha}~\hbox{is
minimal in}~\LM (N)~\hbox{w.r.t.}~\prec'\} .${\parindent=0pt\v5

{\bf 1.2.6. Theorem} With the notation as above, the following
statements hold.\par

(i) $\Omega\ne \emptyset$ and $\Omega$ is a proper subset of $\LM
(N)$.\par

(ii) Let $\G =\{ g\in N~|~\LM (g)\in\Omega\}$. Then $\G$ is a left
Gr\"obner basis of $N$.\vskip 6pt

{\bf Proof} (i) That $\Omega\ne\emptyset$ follows from $N\ne \{ 0\}$
and the remark about $\prec '$ we made above. Since $N$ is a left
ideal and $A$ is a domain (Proposition 1.1.4), if $f\in N$ and $f\ne
0$,  then $hf\in N$ for all nonzero $h\in A$.   It follows from
Proposition 1.1.4(i) that $\LM (hf)\in\LM (N)$ and $\LM (f)|_{_{\rm
L}}\LM (hf)$, i.e., $\LM (f)\prec '\LM (hf)$ in $\LM (N)$. It is
clear that if $\LM (h)\ne 1$, then $\LM (hf)\ne \LM (f)$. This shows
that $\Omega$ is a proper subset of $\LM (N)$.}\par

(ii) By the definition of $\prec '$ and the remark we made above the
theorem, if $f\in N$ with $\LM (f)\ne 0$ and $\LM (f)\not\in\Omega$,
then there is some $g\in\G$ such that $\LM (g)|_{_{\rm L}}\LM (f)$.
Hence the selected $\G$ is a left Gr\"obner basis for $N$.\QED\v5

In Section 3 we will show that every nonzero left ideal $N$ of $A$
has a {\it finite} left Gr\"obner basis. {\parindent=0pt\v5

{\bf Basic facts determined by left Gr\"obner bases}\vskip 6pt

By referring to Proposition 1.1.4, Theorem 1.2.3 and Proposition
1.2.4, the foregoing discussion allows us to summarize some basic
facts determined by left Gr\"obner bases, of which the detailed
proof is left as an exercise.\v5

{\bf 1.2.7. Proposition} Let $N$ be a left ideal of $A$, and let
$\G$ be a left Gr\"obner basis of $N$. Then the following statements
hold.\par

(i) If $f\in N$ and $f\ne 0$, then $\LM (f)=\LM (qg)$ for some $q\in
A$ and $g\in\G$.\par

(ii) If $f\in A$ and $f\ne 0$, then $f$ has a unique remainder $\OV
f^{\G}$ on division by $\G$. \par

(iii) If $f\in A$ and $f\ne 0$, then $f\in N$ if and only if $\OV
f^{\G}=0$. Hence the membership problem for left ideals of $A$ can
be solved by using left Gr\"obner bases.\par

(iv) As a $K$-vector space, $A$ has the decomposition
$$A=N\oplus K\hbox{-span}\mathcal {N}(\G ).$$\par

(v) As a $K$-vector space, $A/N$ has the $K$-basis
$$\OV{\mathcal {N}(\G )}=\{ \OV{a^{\alpha}}~|~a^{\alpha}\in\mathcal {N}(\G )\} ,$$
where $\OV{a^{\alpha}}$ denotes the coset represented by
$a^{\alpha}$ in $A/N$.\par

(vi)  $\mathcal {N}(\G)$ is a finite set, or equivalently,
Dim$_KA/N<\infty$, if and only if  for each $i=1, \ldots ,n$, there
exists $g_{j_i}\in\G$ such that $\LM (g_{j_i})=a_i^{m_i}$, where
$m_i\in\NZ$. \QED  \v5

{\bf Minimal and reduced left Gr\"obner bases}\vskip 6pt

{\bf 1.2.8. Definition} Let $N\ne \{ 0\}$ be a left ideal of $A$ and
let $\G$ be a left gr\"obner basis of $N$. If any proper subset of
$\G$ cannot be a left Gr\"obner basis of $N$, then $\G$ is called a
{\it minimal left Gr\"obner basis} of $N$.\v5

{\bf 1.2.9. Proposition} Let $N\ne \{ 0\}$ be a left ideal of $A$. A
left Gr\"obner basis $\G$ of $N$ is minimal if and only if ~$\LM
(g_i){\not |_{_{\rm L}} }\LM (g_j)$ for all $g_i,g_j\in\G$ with
$g_i\ne g_j$.\vskip 6pt

{\bf Proof} Suppose that $\G$ is minimal. If there were $g_i\ne g_j$
in $\G$ such that  $\LM (g_i)|_{_{\rm L}}\LM (g_j)$, then since the
left division is transitive, the proper subset $\G'=\G-\{ g_j\}$ of
$\G$ would form a left Gr\"obner basis of $N$. This contradicts the
minimality of $\G$.}\par

Conversely, if the condition $\LM (g_i){\not |_{_{\rm L}}}\LM (g_j)$
holds for all $g_i,g_j\in\G$ with $g_i\ne g_j$, then  the definition
of a left Gr\"obner basis entails that any proper subset of $\G$
cannot be a left Gr\"obner basis of $N$. This shows that $\G$ is
minimal.\QED{\parindent=0pt\v5

{\bf 1.2.10. Corollary} Let $N\ne \{ 0\}$ be a left ideal of $A$,
and let the notation be as in Theorem 1.2.6.\par

(i) The left Gr\"obner basis $\G$ obtained there is indeed a minimal
left Gr\"obner basis for $N$.\par

(ii) If $G$ is any left Gr\"obner basis of $N$, then
$$\Omega =\LM (\G )=\{ \LM (g)~|~g\in\G\}\subseteq\LM (G)=\{\LM (g')~|~g'\in G\} .$$
Therefore, any two minimal left Gr\"obner bases of $N$ have the same
set of leading monomials $\Omega$.\vskip 6pt

{\bf Proof} This follows immediately from the definition of a left
Gr\"obner basis, the construction  of $\Omega$ and Proposition
1.2.9.\QED}\v5

We next introduce the notion of a reduced left Gr\"obner
basis.{\parindent=0pt\v5

{\bf 1.2.11. Definition} Let $N\ne \{ 0\}$ be a left ideal of $A$
and let $\G$ be a left gr\"obner basis of $N$. If $\G$ satisfies the
following conditions:\par

(1) $\G$ is a minimal left Gr\"obner basis;\par

(2) $\LC (g)=1$ for all $g\in\G$;\par

(3) For every $g\in\G$, $h=g-\LM (g)$ is a normal element (mod
$\G$), i.e., $h\in K$-span$\mathcal {N}(\G )$,\par

then $\G$ is called a {\it reduced left Gr\"obner basis} of $N$.}\v5

{\parindent=0pt\v5

{\bf 1.2.12. Proposition} Every nonzero left ideal $N$ of $A$ has a
unique reduced left Gr\"obner basis. Therefore, two left ideals
$N_1$, $N_2$ have the same reduced left Gr\"obner basis if and only
if $N_1=N_2$.\vskip 6pt

{\bf Proof}  Note that if $\G$ and $\G '$ are reduced left Gr\"obner
bases for $I$, then $\LM (\G )=\LM (\G ')$ by Corollary 1.2.10. If
$\LM (g_i)=\LM (g_j')$ then $g_i-g_j'\in N\cap K$-span$\mathcal
{N}(\G )=N\cap K\hbox{-span}\mathcal {N}(\G ')$. It follows from
Proposition 1.2.7(iv) that  $g_i=g_j'$. Hence $\G =\G '$. By the
uniqueness, the second assertion is clear.\QED}\v5

Since we will see in Section 3 that every left ideal $N$ of $A$ has
a finite left Gr\"obner basis, a minimal left Gr\"obner basis,
thereby the reduced left Gr\"obner basis for $N$, can be obtained in
an algorithmic way. More precisely, the next proposition holds
true.{\parindent=0pt\v5

{\bf 1.2.13. Proposition} Let $N\ne \{ 0\}$ be a left ideal of $A$,
and let $\G =\{ g_1,\ldots ,g_m\}$ be a finite left Gr\"obner basis
of $N$. \par

(i) The subset $\G_0=\{ g_i\in\G~|~\LM (g_i)~\hbox{is minimal
in}~\LM (\G )~\hbox{w.r.t.}\prec '\} $ of $\G$ forms a minimal left
Gr\"obner basis of $N$ (see the definition of $\prec '$ given before
Theorem 1.2.6).  An algorithm written in pseudo-code is omitted
here.\par

(ii) With $\G_0$ as in (i) above, we may assume, without loss of
generality, that $\G_0=\{ g_1,\ldots ,g_s\}$ with $\LC (g_i)=1$ for
$1\le i\le s$. Put $\G_1=\{ g_2,...,g_s\}$ and
$h_1=\OV{g_1}^{\G_1}$. Then $\LM (h_1)=\LM (g_1)$. Put $\G_2=\{
h_1,g_3,...,g_s\}$ and $h_2=\OV{g_2}^{\G_2}$. Then  $\LM (h_2)=\LM
(g_2)$. Put $\G_3=\{ h_1,h_2,g_4,...,g_s\}$ and
$h_3=\OV{g_3}^{\G_3}$, and so on. The last obtained
$\G_{s+1}=\{h_1\ldots ,h_{s-1}\}\cup\{h_s\}$ is then the reduced
left Gr\"obner basis. An algorithm written in pseudo-code is omitted
here.\vskip 6pt

{\bf Proof} This can be verified directly, so we leave it as an
exercise.\QED }\v5

\section*{1.3. The Noetherianess}

Let $A=K[a_1,\ldots ,a_n]$ be a solvable polynomial algebra with
admissible system $(\B ,\prec )$. In this section we show that $A$
is a (left and right) Noetherian $K$-algebra by showing that every
(left, right) ideal of $A$ has a finite (left, right) Gr\"obner
basis. \v5

The following lemma, which will be essential not only in proving our
next theorem but also in proving the existence of finite left
Gr\"obner bases and the termination property of Buchberger's
algorithm for modules over solvable polynomial algebras (Section 3 
of Chapter 2), is usually attributed to the American algebraist L. 
E. Dickson 
(http:$//$en.wikipedia.org$/$wiki$/$Dickson's$_{-}$lemma). For a 
detailed argument on Dickson's lemma in commutative computational 
algebra, we refer to ([BW], 4.4).{\parindent=0pt \v5

{\bf 1.3.1. Lemma} (Dickson's Lemma) Every subset $S$ of $\NZ^n$ has
a finite subset $B$ such that for each $(\alpha_1,\ldots
,\alpha_n)\in S$, there exists $(\gamma_1,\ldots ,\gamma_n)\in B$
with $\gamma_i\le \alpha_i$ for $1\le i\le n$.\QED}\v5

Let $a^{\alpha},a^{\beta}\in\B$ with $\alpha =(\alpha_1,\ldots
,\alpha_n)$, $\beta =(\beta_1,\ldots ,\beta_n)\in\NZ^n$. Recall that
$a^{\alpha}\prec 'a^{\beta}$ if and only if $a^{\alpha}|_{_{\rm
L}}a^{\beta}$, while the latter is defined subject to the property
that  $a^{\beta}=\LM (a^{\gamma}a^{\alpha})$ for some
$a^{\gamma}\in\B$. Thus, by Proposition 1.1.4(i), $a^{\alpha}\prec
'a^{\beta}$ is actually equivalent to the property that $\alpha_i\le
\beta_i$ for $1\le i\le n$. Also note that the correspondence
$\alpha =(\alpha_1,\ldots ,\alpha_n)\longleftrightarrow a^{\alpha}$
gives the bijection between $\NZ^n$ and $\B$. Combining this
observation with Dickson's lemma,  the following result is
obtained.{\parindent=0pt\v5

{\bf 1.3.2. Theorem}  With notations as in Theorem 1.2.6, let $N\ne
\{ 0\}$ be a left ideal of $A$, $\LM (N)=\{ \LM (f)~|~f\in N\}$,
$\Omega =\{ a^{\alpha}\in \LM (N)~|~a^{\alpha}~\hbox{is minimal
w.r.t.}~\prec'\}$, and $\G =\{ g\in N~|~\LM (g)\in\Omega\}$. Then
$\Omega$ is a finite subset of $\LM (N)$, thereby $\G$ is a finite
left Gr\"obner basis of $N$.\QED}\v5

Since every solvable polynomial algebra has a (two-sided) monomial
ordering, a right division algorithm and a theory of right Gr\"obner
bases for right ideals hold true as well, we then have the
following{\parindent=0pt\v5

{\bf 1.3.3. Corollary} Every solvable polynomial algebra $A$ is
(left and right) Noetherian.}\QED\v5

In the next chapter we will see that the Buchberger algorithm, which  
computes a finite Gr\"obner basis for a finitely generated 
commutative polynomial ideal, has a complete noncommutative version 
({\bf Algorithm-LGB} presented in Section 3 of Chapter 2), that 
computes a finite left Gr\"obner basis for a finitely generated 
submodule of a free module over a solvable polynomial algebra.

\section*{1.4. A Constructive Characterization}\def\MB{\mathbb{B}}

From Definition 1.1.3 we see that  the two conditions (S1) and (S2),
which determine a solvable polynomial algebra $A=K[a_1,\ldots
,a_n]$, are mutually {\it independent} factors. In this section we
give a characterization of solvable polynomial algebras  by
employing Gr\"obner bases of ideals in free algebras, so that
solvable polynomial algebras are completely recognizable and
constructible in a computational way. \v5

To make the text self-contained, we start by recalling some basics
on Gr\"obner bases of {\it two-sided ideals} in a free $K$-algebra
$\KS =K\langle X_1,\ldots ,X_n\rangle$ on $X=\{ X_1,\ldots
,X_n\}$.{\parindent=0pt\v5

{\bf Division algorithm in $\KX$}\vskip 6pt

Let $\mathbb{B}=\{ 1,~X_{i_1}\cdots X_{i_s}~|~X_{i_j}\in X,~s\ge
1\}$ be the standard $K$-basis of $\KS$. For convenience, elements
of $\mathbb{B}$ are also referred to as {\it monomials}, and we use
capital letters $U,V,W,S,\ldots$ to denote monomials in
$\mathbb{B}$. Recall that a {\it monomial ordering} $\prec_{_X}$ on
$\mathbb{B}$ (or equivalently on $\KX$) is a well-ordering such that
for any $W,U,V, S\in\mathbb{B}$,{\parindent=1.3truecm\par

\item{(1)} $U\prec_{_X}V~\hbox{implies}~WU\prec_{_X}WV,
~US\prec_{_X}VS,$  or equivalently, $WUS\prec_{_X}WVS;$ \par

\item{(2)} if $U\ne V$, then $V=WUS$ implies $U\prec_{_X}V$ (thereby
$1\prec_{_X}W$ for all $1\ne W\in\mathbb{B}$).\par} \par

If $\prec_{_{\rm X}}$ is a monomial ordering on $\mathbb{B}$, then
the data $(\mathbb{B},\prec_{_{\rm X}})$ is referred to as an {\it
admissible system} of $\KX$.}\v5

Note that for any given $n$-tuple $(m_1,\ldots
,m_n)\in\mathbb{N}^n$, a {\it weighted degree function} $d(~)$ is
well defined on nonzero elements of $\KX$, namely, by assigning each
$X_i$ the degree $d(X_i)=m_i$, $1\le i \le n$, we may define for
each $W=X_{i_1}\cdots X_{i_s}\in\mathbb{B}$ the degree
$d(W)=m_{i_1}+\cdots +m_{i_s}$, and for each nonzero
$f=\sum_{i=1}^s\lambda_iW_i\in \KS$ with $\lambda_i\in K^*$ and
$W_i\in\mathbb{B}$, the degree of $f$ is then defined as
$$d(f)=\max\{ d(W_i)~|~1\le i\le s\}.$$ If $d(X_i)=m_i>0$
for $1\le i\le n$, then $d(~)$ is referred to as a {\it
positive-degree function} on $\KS$. \par

Let  $d(~)$ be a  positive-degree function on $\KX$. If
$\prec_{_{\rm X}}$ is a monomial ordering on $\mathbb{B}$ such that
for all $U,V\in\mathbb{B}$,
$$U\prec_{_{\rm X}} V~\hbox{implies}~d(U)\le d(V),$$
then, as in Section 1,  we call $\prec_{_{\rm X}}$ a {\it graded
monomial ordering} with respect to $d(~)$. For instance, with
respect to any given positive-degree function $d(~)$, the
lexicographic graded ordering $\prec_{grlex}$ on $\mathbb{B}$ can be
defined as follows. For $U,V\in\mathbb{B}$,
$$U\prec_{grlex}V\Leftrightarrow\left\{\begin{array}{l} d(U)<d(V);\\
\hbox{or}\\
d(U)=d(V)~\hbox{and}~U\prec_{lex}V,\end{array}\right.$$ where the
lexicographic ordering $\prec_{lex}$ on $\mathbb{B}$ may be defined
by ordering $X_1,\ldots ,X_n$ arbitrarily, say
$X_{i_1}\prec_{lex}X_{i_2}\prec_{lex}\cdots\prec_{lex}X_{i_n},$
i.e., if $U=X_{\ell_1}\cdots X_{\ell_s},~V=X_{t_1}\cdots
X_{t_m}\in\MB$ and $s\le m$, then
$$U\prec_{lex}V\Leftrightarrow\left\{\begin{array}{l} s<m,~X_{\ell_k}=X_{t_k},~1\le k\le s;\\
\hbox{or}\\
s=m,~\hbox{there exists}~p\le s~\hbox{such that}~
X_{\ell_k}=X_{t_k}~\hbox{for}~k <p\\
~~~~~~~~~~~~~~~~~~~~~~~~~~~~~~~~~~\hbox{and}~
X_{\ell_p}\prec_{lex}X_{t_p}.\end{array}\right.$$\par

At this point we should point out that  though the ordering
$\prec_{lex}$ is a total ordering on $\MB$, {\it it is not a
monomial orderings} on $\mathbb{B}$. For instance, considering the
free algebra $K\langle X_1,X_2\rangle$ with $X_1\prec_{lex} X_2$, we
have
$$X_2\succ_{lex}X_1X_2\succ_{lex}X_1X_1X_2\succ_{lex}X_1X_1X_1X_2\succ_{lex}\cdots .$$
This shows that $\prec_{lex}$ {\it is not a well-ordering on} $\KX$.
\par

One may refer to loc. cit. (e.g. [Gr]) for more monomial orderings
on $\KX$.  \v5

Let $\prec_{_X}$ be a monomial ordering  on $\mathbb{B}$. If $f\in
\KX$ is such that $f=\sum_{i=1}^m\lambda_iW_i$ with $\lambda_i\in
K^*$, $W_i\in\MB$, and $W_1\prec_{_X} W_2\prec_{_X}\cdots\prec_{_X}
W_m$, then we write $\LM (f)=W_m$ for the {\it leading monomial} of
$f$,  $\LC (f)=\lambda_m$ for the {\it leading coefficient} of $f$,
and we write $\LT (f)=\lambda W_m$ for the {\it leading term} of
$f$.\v5

Note that $\MB$ forms a multiplicative monoid with the identity 1
and the multiplication in $\B$ is just simply the concatenation of
monomials,  i.e., if $U,V\in\MB$ with $U=X_{i_1}\cdots X_{i_p}$,
$V=X_{j_1}\cdots X_{i_q}$, then $UV=X_{i_1}\cdots
X_{i_p}X_{j_1}\cdots X_{j_q}\in\MB$.  Thereby if there exist
$W,S\in\MB$ such that $V=WUS$, then we say that $U$ {\it divides}
$V$, denoted $U|V$.\v5

Let $F$ be a nonempty subset of $\KX$. Then the division of
monomials defined above yields a subset of $\MB$:
$$\mathcal {N}(F)=\{ W\in\MB~|~\LM (f){\not |~ W},~f\in F\}.$$
If $F=\{ g\}$ consists of a single element $g$, then we simply write
$\mathcal {N}(g)$ in place of $\mathcal {N}(F)$. Also we write
$K$-span$\mathcal {N}(F)$ for the $K$-subspace of $\KX$ spanned by
$\mathcal {N}(F)$.{\parindent=0pt\v5

{\bf 1.4.1. Definition} Elements of $\mathcal {N}(F)$ are referred
to as {\it normal monomials} (mod $F$). Elements of
$K$-span$\mathcal {N}(F)$ are referred to as {\it normal elements}
(mod $F$).}\v5

Given a monomial ordering $\prec_{_{\rm X}}$ on $\MB$, if $f\in\KX$
and $U,V\in\MB$, then the definition of $\prec_{_{\rm X}}$ entails
that $\LM (UfV)=U\LM (f)V$. So, the division we defined for
monomials in $\MB$ can naturally be used to define a  division
procedure for elements in $\KX$. More precisely, let $f,g\in \KX$
with $\LC (f)=\mu \ne 0$, $\LC (g)=\lambda\ne 0$. If $\LM (g)|\LM
(f)$, i.e., there exists $U,V\in\MB$ such that $\LM (f)=U\LM (g)V$,
then put $f_1=f-\lambda^{-1}\mu UgV$; otherwise, put $f_1=f-\LT
(f)$. Note that in both cases we have $f_1=0$, or $f_1\ne 0$ and
$\LM (f_1)\prec_{_{\rm X}}\LM (f)$. At this stage, let us refer to 
such a procedure of canceling the leading term of $f$  as the {\it 
division procedure by $g$}.  With $f:=f_1\ne 0$, we can repeat the 
division procedure by $g$ and so on. This, in turn, gives rise to
$$\LM (f_{i+1})\prec_{_{\rm X}}\LM (f_i)\prec_{_{\rm
X}}\cdots\prec_{_{\rm X}}\LM (f_1)\prec_{_{\rm X}}\LM (f),\quad i\ge
2.$$  Since $\prec_{_{\rm X}}$ is a well-ordering, it follows that
such a division procedure terminates after a finite number of
repetitions, and consequently $f$ is expressed as
$$f=\sum_{i}^t\lambda_{i}U_igVi+r,$$
where $\lambda_i\in K^*$, $U_i,V_i\in\MB$, and $r\in
K$-span$\mathcal {N}(g)$, i.e., $r$ is normal (mod $g$), such that
either $\LM (f)=\max\{\LM (U_igV_i)~|~1\le i\le t\}$ or $\LM (f)=\LM
(r)$.\par

Furthermore, the division procedure demonstrated above can be
extended to a  division procedure canceling the leading term of $f$
by a finite subset $G=\{ g_1,\ldots ,g_s\}$, which therefore gives
rise to an effective division algorithm in $\KX$, that is, we have
reached the following {\parindent=0pt\v5

{\bf 1.4.2. Theorem} Given a set of nonzero elements $G=\{
g_1,\ldots ,g_s\}$ and $f$ in $A$, the {\bf Algorithm-DIV} given
below produces finitely many  $\lambda_{ij}\in K^*$,
$U_{ij},V_{ij}\in\MB$, and an $r\in K$-span$\mathcal {N}(G)$, such
that $f=\sum_{i,j}\lambda_{ij}U_{ij}g_jV_{ij}+r$ and $\LM
(U_{ij}g_jV_{ij})\preceq\LM (f)$, $\LM (r)\preceq\LM (f)$ if $r\ne
0$.\vskip 6pt

\underline{\bf Algorithm-DIV~~~~~~~~~~~~
~~~~~~~~~~~~~~~~~~~~~~~~~~~~~~~~~~~~~~~~~~~~~~~~~~}\vskip 6pt

\textsc{INPUT}: $f,~G=\{g_1,\ldots,g_s\}~\hbox{with}~g_i\ne 0~(1\le 
i\le s)$\par                                                            
\textsc{OUTPUT}: $\lambda_{ij}\in K^*,~U_{ij},V_{ij}\in\MB,~r\in 
K\hbox{-span}\mathcal {N}(G)$\par \textsc{INITIALIZATION}: 
$i=0;~r:=0;~h:=f$\par 
~~~~~~~~~~~~~~~~~~~~~~~~~~~~~$\Lambda_1:=\emptyset 
,~\cdots,~\Lambda_s=\emptyset ;~Q_1:=\emptyset ,~\cdots 
,Q_s:=\emptyset$\newpage

\textsc{BEGIN}\par ~~~~~\textsc{WHILE}~$h\ne 0~\textsc{DO}$\par 
~~~~~~~~~~\hbox{IF}~\hbox{there 
exist}~$j$~\hbox{and}~$U,V\in\MB$\par                                  
~~~~~~~~~~~~~~\hbox{such that}~$\LM (h)=U\LM 
(g_j)V$~\textsc{THEN}\par 
~~~~~~~~~~~~~~\hbox{choose}~$j$~\hbox{least such that}~$\LM (h)=U\LM 
(g_j)V$\par                                                              
~~~~~~~~~~$i:=i+1,$ $\lambda_{ij}:=\LC(g_j)^{-1}\LC
(h),~U_{ij}:=U,~V_{ij}:=V$\par                                        
~~~~~~~~~~$\Lambda_j :=\Lambda_j\cup\{\lambda_{ij}\} 
,~Q_j:=Q_j\cup\{ U_{ij},V_{ij}\}$\par                                                       
~~~~~~~~~~$h:=h-\LC (g_j)^{-1}\LC (h)U_{ij}g_jV_{ij}$\par  
~~~~~~~~~~\hbox{ELSE}\par                                                    
~~~~~~~~~~$r:=r+\LT (h)$\par                                           
~~~~~~~~~~$h:=h-\LT (h)$\par                                            
~~~~~~~~~~\hbox{END}\par                                                     
~~~~~\textsc{END}\par                                                 
\textsc{END}\par\vskip -.2truecm 
\underline{~~~~~~~~~~~~~~~~~~~~~~~~~~~~~~~~~~~~~~~~~~~~~~~~~~~~~~~~~~~~~~~~~~~~~~~~~~~~~~~~~~~~~~~~~~~~~~~} 
\v5

{\bf 1.4.3. Definition} The element $r$ obtained in {\bf
Algorithm-DIV} is called a {\it remainder} of $f$ on division by
$G$, and is denoted by $\OV f^{G}$, i.e., $\OV f^G=r$. If $\OV
f^G=0$, then we say that $f$ {\it is reduced to $0$} (mod $G$).\v5

{\bf Remark} Actually as in the case of (Section 2, {\bf
Algorithm-LDIV}),  the element $r$ obtained in {\bf Algorithm-DIV}
depends on how the order of elements in $G$ is arranged. This
implies that the $r$ may not be unique. That is why the phrase ``a
remainder" is used in the above definition.\v5

{\bf Gr\"obner bases for two-sided ideals of $\KX$}\vskip 6pt

For the remainder of this section, ideals of $\KX$ are always meant
two-sided ideals; if $M\subset\KX$ is a nonempty subset, then we
write $\LM (M)=\{\LM(f)~|~f\in M\}$ for the set of leading monomials
of $M$, and we write $I=\langle M\rangle$ for the two-sided ideal
$I$ of $\KX$ generated by $M$; moreover, we fix an admissible system
$(\MB ,\prec_{_{\rm X}})$ of $\KX$.\par}

As with a solvable polynomial algebra in Section 2, in principle the
division procedure by using a finite subset $G$ of nonzero elements
in $\KX$ can be extended to a left division procedure by means of an
{\it arbitrary proper subset} $G$ of nonzero elements, for,  it is a
true statement that if $f\in \KX$ with $\LM (f)\ne 0$, then either
there exists $g\in G$ such that $\LM (g)|\LM (f)$ or such $g$ does
not exist.  This leads to the following{\parindent=0pt\v5

{\bf 1.4.4. Proposition} Let $I$ be an ideal of $\KX$ and $G$ a
proper subset of nonzero elements in $I$. The following statements
are equivalent:\par

(i) If $f\in I$ and $f\ne 0$, then there exists $g\in G$ such that
$\LM (g)|\LM (f)$.\par

(ii) Every nonzero $f\in I$ has a finite representation
$f=\sum_{i,j}\lambda_{ij}U_{ij}g_jV_{ij}$ with $\lambda_{ij}\in
K^*$, $U_{ij},V_{ij}\in\MB$ and $g_j\in G$, such that $\LM
(U_{ij}g_jV_{ij})\preceq_{_{\rm X}}\LM (f)$.\par

(iii) $\langle \LM (I)\rangle =\langle\LM (G)\rangle$.\QED\v5

{\bf 1.4.5. Definition}  Let $I$ be an ideal of $\KX$. With respect
to a given monomial ordering $\prec_{_{\rm X}}$ on $\MB$, a proper
subset $\G$ of nonzero elements in $I$ is said to be a {\it
Gr\"obner basis } of $I$ if $\G$ satisfies one of the equivalent
conditions of Proposition 1.4.4.}
\par

If $\G$ is a Gr\"obner basis of $I$, then the expression
$f=\sum_{i,j}\lambda_{ij}U_{ij}g_jV_{ij}$ appeared in Proposition
1.4.4(ii) is called a {\it Gr\"obner representation} of $f$. \v5

Clearly, if $I$ is an ideal of $\KX$ and $I$ has a Gr\"obner basis
$\G$, then $\G$ is certainly a generating set of $I$, i.e.,
$I=\langle \G\rangle$. But the converse is not necessarily true. For
instance, let $g_1=X_1^2X_2-X_3$, $g_2=X_2$, and $I=\langle
g_1,g_2\rangle$. Then since $X_3\in I$, $G=\{ g_1,g_2\}$ is not a
Gr\"obner basis with respect to any given monomial ordering
$\prec_{_{\rm X}}$ on $\MB$ such that $\LM (g_1)=X_1^2X_2$, $\LM
(g_2)=X_2$.{\parindent=0pt\v5

{\bf The existence of Gr\"obner bases}\vskip 6pt

Noticing $1\in\MB$, if we define on $\MB$ the ordering:
$$U\prec V\Leftrightarrow U|V,\quad U,~V~\hbox{in}~\MB ,$$
then, by the divisibility defined for monomials,  it is easy to see
that $\prec $ is reflexive, anti-symmetric, transitive, and
moreover, $$U\prec V~\hbox{implies}~U\prec_{_{\rm X}}V.$$ Since the
given monomial ordering $\prec_{_{\rm X}}$ is a well-ordering on
$\MB$, it follows that every nonempty subset of $\MB$ has a minimal
element with respect to the ordering $\prec $ on $\MB$. }\v5

Now, let $I$ be a nonzero ideal of $\KX$, $\LM (I)=\{ \LM (f)~|~f\in
I\}$, and $\Omega =\{ U\in \LM (I)~|~U~\hbox{is minimal
w.r.t.}~\prec\} .${\parindent=0pt\v5

{\bf 1.4.6. Theorem} With the notation as above, the following
statements hold.\par

(i) $\Omega\ne \emptyset$ and $\Omega$ is a proper subset of $\LM
(I)$.\par

(ii) Let $\G =\{ g\in I~|~\LM (g)\in\Omega\}$. Then $\G$ is a
Gr\"obner basis of $I$.\vskip 6pt

{\bf Proof} (i) That $\Omega\ne\emptyset$ follows from $I\ne \{ 0\}$
and the remark about $\prec$ we made above. Since $I$ is a two-sided
ideal and $\KX$ is a domain, if $f\in I$ then $hfg\in I$, thereby
$\LM (hfg)\in\LM (I)$  for all $h,g\in \KX$. Also note that the
leading monomials of elements in $\KX$ are taken with respect to the
given monomial ordering $\prec_{_{\rm X}}$. It follows that $\LM
(hfg)=\LM (h)\LM (f)\LM (g)$ which implies $\LM (f)|\LM (hfg)$,
thereby $\LM (f)\prec \LM (hfg)$ in $\LM (I)$. It is clear that if
$\LM (h)\ne 1$ or $\LM (g)\ne 1$, then $\LM (hfg)\ne \LM (f)$. This
shows that $\Omega$ is a proper subset of $\LM (I)$.}\par

(ii) By the definition of $\prec$ and the remark we made above the
theorem, if $f\in I$ with $\LM (f)\ne 0$ and $\LM (f)\not\in\Omega$,
then there is some $g\in\G$ such that $\LM (g)|\LM (f)$. Hence the
selected $\G$ is a Gr\"obner basis for $I$.\QED\v5 {\parindent=0pt

{\bf Basic facts determined by Gr\"obner bases}\vskip 6pt

By referring to Proposition 1.4.4, Definition 1.4.5 and the
foregoing discussion,  the next proposition summarizes some basic
facts determined by Gr\"obner bases of ideals in $\KX$, of which the
detailed proof is left as an exercise.\v5

{\bf 1.4.7. Proposition} Let $I$ be an ideal of $\KX$, and let  $\G$
be a Gr\"obner basis of $I$. Then the following statements hold.\par

(i) If $f\in I$ and $f\ne 0$, then $\LM (f)=\LM (UgV)$ for some
$U,V\in\MB$ and $g\in\G$.\par

(ii) If $f\in \KX$ and $f\ne 0$, then $f$ has a unique remainder
$\OV f^{\G}$ on division by $\G$. \par

(iii) If $f\in \KX$ and $f\ne 0$, then $f\in I$ if and only if $\OV
f^{\G}=0$. Hence the membership problem for ideals of $\KX$ can be
solved by using Gr\"obner bases.\par

(iv) As a $K$-vector space, $\KX$ has the decomposition
$$\KX=I\oplus K\hbox{-span}\mathcal {N}(\G ).$$\par

(v) As a $K$-vector space, $\KX /I$ has the $K$-basis
$$\OV{\mathcal {N}(\G )}=\{ \OV{U}~|~U\in\mathcal {N}(\G )\} ,$$
where $\OV{U}$ denotes the coset represented by $U$ in $\KX /I$.\QED
\v5

{\bf Minimal and reduced Gr\"obner bases}\vskip 6pt

{\bf 1.4.8. Definition} Let $I\ne \{ 0\}$ be an ideal of $\KX$ and
let $\G$ be a gr\"obner basis of $I$. If any proper subset of $\G$
cannot be a  Gr\"obner basis of $I$, then $\G$ is called a {\it
minimal  Gr\"obner basis} of $I$.\v5

{\bf 1.4.9. Proposition} Let $I\ne \{ 0\}$ be an ideal of $\KX$. A
Gr\"obner basis $\G$ of $I$ is minimal if and only if $\LM
(g_i){\not |}~\LM (g_j)$ for all $g_i,g_j\in\G$ with $g_i\ne
g_j$.\vskip 6pt

{\bf Proof} Suppose that $\G$ is minimal. If there were $g_i\ne g_j$
in $\G$ such that  $\LM (g_i)|\LM (g_j)$, then since the division of
monomials is transitive, the proper subset $\G'=\G-\{ g_j\}$ of $\G$
would form a Gr\"obner basis of $I$. This contradicts the minimality
of $\G$.}\par

Conversely, if the condition $\LM (g_i){\not |}~\LM (g_j)$ holds for
all $g_i,g_j\in\G$ with $g_i\ne g_j$, then  the definition of a
Gr\"obner basis entails that any proper subset of $\G$ cannot be a
 Gr\"obner basis of $I$. This shows that $\G$ is
minimal.\par\QED{\parindent=0pt\v5

{\bf 1.4.10. Corollary} Let $I\ne 0$ be an ideal $I$ of $\KX$. With
the notation as in Theorem 1.4.6, the Gr\"obner basis $\G$ obtained
there is indeed a minimal  Gr\"obner basis for $I$. Moreover,  if
$G$ is any  Gr\"obner basis of $I$, then
$$\Omega =\LM (\G )=\{ \LM (g)~|~g\in\G\}\subseteq\LM (G)=\{\LM (g')~|~g'\in G\} .$$
Therefore, any two minimal  Gr\"obner bases of $I$ have the same set
of leading monomials $\Omega$.\vskip 6pt

{\bf Proof} This follows immediately from the definition of a
Gr\"obner basis, the construction  of $\Omega$ and Proposition
1.4.9.\QED}\v5

We next introduce the notion of a reduced  Gr\"obner
basis.{\parindent=0pt\v5

{\bf 1.4.11. Definition} Let $I\ne \{ 0\}$ be an ideal of $\KX$ and
let $\G$ be a gr\"obner basis of $I$. If $\G$ satisfies the
following conditions:\par

(1) $\G$ is a minimal  Gr\"obner basis;\par

(2) $\LC (g)=1$ for all $g\in\G$;\par

(3) For every $g\in\G$, $h=g-\LM (g)$ is a normal element (mod
$\G$), i.e., $h\in K$-span$\mathcal {N}(\G )$,\par

then $\G$ is called a {\it reduced Gr\"obner basis} of $N$.\v5

{\bf 1.4.12. Proposition} Every nonzero ideal $I$ of $\KX$ has a
unique reduced Gr\"obner basis. Therefore,two ideals $I$ and $J$
have the same reduced Gr\"obner basis if and only if $I=J$.\vskip
6pt

{\bf Proof}  Note that if $\G$ and $\G '$ are reduced Gr\"obner
bases for $I$, then $\LM (\G )=\LM (\G ')$ by Corollary 1.4.10. If
$\LM (g_i)=\LM (g_j')$ then $g_i-g_j'=I\cap K$-span$\mathcal {N}(\G
)=I\cap K\hbox{-span}\mathcal {N}(\G ')$. It follows from
Proposition 1.4.7(iv) that $g_i=g_j'$. Hence $\G =\G '$. By the
uniqueness, the second assertion is clear.\QED}\v5

 In the case that an ideal $I$ has a finite Gr\"obner basis, a minimal Gr\"obner basis,
thereby the reduced Gr\"obner basis for $I$, can be obtained in an
algorithmic way. More precisely, the next proposition holds
true.{\parindent=0pt\v5

{\bf 1.4.13. Proposition} Let $I\ne \{ 0\}$ be an ideal of $\KX$,
and let $\G =\{ g_1,\ldots ,g_m\}$ be a finite  Gr\"obner basis of
$I$. \par

(i) The subset $\G_0=\{ g_i\in\G~|~\LM (g_i)~\hbox{is minamal
in}~\LM (\G )~\hbox{w.r.t.}\prec \} $ of $\G$ forms a minimal
Gr\"obner basis of $I$ (see the definition of $\prec$ given before
Theorem 1.4.6). An algorithm written in pseud-code is omitted
here.\par

(ii) With $\G_0$ as in (i) above, we may assume, without loss of
generality, that $\G_0=\{ g_1,\ldots ,g_s\}$ with $\LC (g_i)=1$ for
$1\le i\le s$. Put $\G_1=\{ g_2,...,g_s\}$ and
$h_1=\OV{g_1}^{\G_1}$. Then $\LM (h_1)=\LM (g_1)$. Put $\G_2=\{
h_1,g_3,...,g_s\}$ and $h_2=\OV{g_2}^{\G_2}$. Then  $\LM (h_2)=\LM
(g_2)$. Put $\G_3=\{ h_1,h_2,g_4,...,g_s\}$ and
$h_3=\OV{g_3}^{\G_3}$, and so on. The last obtained $\G_s$ is then
the reduced  Gr\"obner basis of $I$. An algorithm written in
pseud-code is omitted here.\vskip 6pt

{\bf Proof} This can be verified directly, so we leave it as an
exercise.\QED \v5\v5

{\bf Construction of Gr\"obner bases}\vskip 6pt

Let the admissible system $(\MB ,\prec_{_{\rm X}})$ of $\KX$ be as
fixed before. A subset $G$ of nonzero elements in $\KX$ is said to
be LM-{\it reduced} if $\LM (g_i){\not |}~\LM (g_j)$ for all $g_i\ne
g_j$ in $G$. For $f,g\in G$, if there are monomials $u,v\in\MB$ such
that {\parindent=.75truecm\par
\item{(1)} $\LM (f)u=v\LM (g)$, and
\item{(2)} $\LM (f){\not |}~v$ and $\LM (g)\not |~u$,}{\parindent=0pt\par

then the element
$$o(f,u;~v,g)=\frac{1}{\LC (f)}(f\cdot u)-
\frac{1}{\LC (g)}(v\cdot g)$$ is referred to as an {\it overlap
element} of $f$ and $g$.}}\par

Obviously $o(f,u;~v,g)$ is generally not unique and it includes also
the case $f=g$. For instance, look at the cases where $\LM
(f)=X_1X_2^2$ and $\LM (g)=X_2^2X_1$; $\LM (f)=X_3^4=\LM (g)$. \v5

Let $I=\langle G\rangle$ be an ideal of $\KX$ generated by a finite
subset $G=\{ g_1,\ldots ,g_t\}$ of nonzero elements. Then, $G$ may
be reduced to an LM-reduced subset $G'$ in an algorithmic way as in
Proposition 1.4.13, such that $I=\langle G'\rangle$. Thus, we may
always assume that $I$ is generated by a reduced finite subset
$G$.\par

The next theorem and the following algorithm are known as the
implementation of Bergman's diamond lemma ([Ber2], [Mor], [Gr]), and
are practically used to check whether $G$ is a Gr\"obner basis of
$I$ or not; if not, the algorithm produces a (possibly infinite)
Gr\"obner basis for $I$. For the detailed proof of the theorem we
refer the reader to loc. cit.{\parindent=0pt\v5

{\bf 1.4.14. Theorem} Let $G =\{ g_1,\ldots ,g_t\}$ be an LM-reduced
subset of nonzero elements in $\KX$. Then $G$ is a Gr\"obner basis
for the ideal $I=\langle G\rangle$ if and only if for each pair
$g_i,g_j\in\G$, including $g_i=g_j$, every overlap element
$o(g_i,u;~v,g_j)$ of $g_i$ and $g_j$ has the property
$\overline{o(g_i,u;~v,g_j)}^{G}=0,$ that is, $o(g_i,u;~v,g_j)$ is
reduced to 0 (mod $G)$.}\par

If $G$ is not a Gr\"obner basis of $I$, then the  algorithm
presented below returns a (finite or countably infinite) Gr\"obner 
basis $\G$ for $I$.{\parindent=0pt\vskip 6pt

\underline{{\bf Algorithm-GB}
~~~~~~~~~~~~~~~~~~~~~~~~~~~~~~~~~~~~~~~~~~~~~~~~~~~~~~~~~~~~~~~~~~~~~~~~}{\vskip 
6pt

\textsc{INPUT:} $G_0=\{ g_1,...,g_t\}$\par

\textsc{OUTPUT:} $\G =\{ g_1,...,g_m,...\}$, a Gr\"obner basis for
$I$\par

\textsc{INITIALIZATION:} $\G :=G_0$, $O :=\{ o(g_i,g_j)~|~g_i,g_j\in
G_0\}-\{0\}$\par

\textsc{BEGIN}\par

~~~~~\textsc{WHILE} $O\ne\emptyset$ \textsc{DO}\par

~~~~~~~~~~~Choose any $o(g_i,g_j)\in O$\par

~~~~~$O :=O-\{o(g_i,g_j)\}$\par

~~~~~~~~~~~\textsc{IF} $\overline{o(g_i,g_j)}^{\G}=r\ne 0$ 
\textsc{THEN}\par

~~~~~~~~~~~$O :=\{o(g,r), ~o(r,g),~o(r,r)~|~g\in\G\}-\{0\}$ \par

~~~~~~~~~~~$\G :=\G\cup\{ r\}$}\par

~~~~~~~~~~~\textsc{END}{\parindent=0pt\par

~~~~~$\vdots$

\underline{~~~~~~~~~~~~~~~~~~~~~~~~~~~~~~~~~~~~~~~~~~~~~~~~~~~~~~~~~~
~~~~~~~~~~~~~~~~~~~~~~~~~~~~~~~~~~~~}}}\par\QED\v5

Since generally the \textsc{WHILE} loop in algorithm {\bf 
Algorithm-GB} does not terminate after a finite number of 
executions, the algorithm is  written without the ending statement  
for the \textsc{WHILE} loop. {\parindent=0pt\v5

{\bf Getting PBW $K$-bases via Gr\"obner bases}\vskip 6pt

Concerning the first condition (S1) that a solvable polynomial
algebra should satisfy, we recall from [Li2] the following result,
which  is a generalization of ([Gr], Proposition 2,14; [Li1],
CH.III, Theorem 1.5).{\parindent=0pt\v5

{\bf 1.4.15. Proposition} ([Li2], Ch 4, Theorem 3.1) Let $I\ne \{
0\}$ be an ideal of the free $K$-algebra $\KX =K\langle X_1,\ldots
,X_n\rangle$, and $A=\KS /I$. Suppose that $I$ contains a subset of
$\frac{n(n-1)}{2}$ elements
$$G=\{ g_{ji}=X_jX_i-F_{ji}~|~F_{ji}\in\KS ,~1\le i<j\le n\}$$
such that with respect to some monomial ordering $\prec_{_X}$ on
$\mathbb{B}$, $\LM (g_{ji})=X_jX_i$ holds for all the $g_{ji}$. The
following two statements are equivalent:\par

(i) $A$ has the  PBW $K$-basis $\mathscr{B}=\{ \OV X_1^{\alpha_1}\OV
X_2^{\alpha_2}\cdots \OV X_n^{\alpha_n}~|~\alpha_j\in\NZ\}$ where
each $\OV X_i$ denotes the coset of $I$ represented by $X_i$ in
$A$.\par

(ii) Any  subset $\G$ of $I$ containing $G$ is a  Gr\"obner basis
for $I$ with respect to $\prec_{_X}$.\QED \v5

{\bf Remark} (i) Obviously, Proposition 1.4.15 holds true if we use
any permutation $\{X_{k_1},\ldots ,X_{k_n}\}$ of $\{ X_1,\ldots
,X_n\}$ (see an example given in the end of this section). So, in
what follows we conventionally use only $\{ X_1,\ldots ,X_n\}$.}\par

(ii) Let the ideal $I$ be as in Proposition 1.4.15.  If $G=\{
g_{ji}=X_jX_i-F_{ji}~|~F_{ji}\in\KS ,~1\le i<j\le n\}$ is a
Gr\"obner basis of $I$ such that $\LM (g_{ji})=X_jX_i$ for all the
$g_{ji}$, then it is not difficult to see that the  reduced
Gr\"obner basis of $I$ is of the form
$$\G =\left\{\left.
g_{ji}=X_jX_i-\sum_q\mu^{ji}_qX_1^{\alpha_{1q}}X_2^{\alpha_{2q}}\cdots
X_n^{\alpha_{nq}}~\right |~\begin{array}{l} \LM (g_{ji})=X_jX_i,~\\
1\le i<j\le n\end{array}\right\}
$$  where $\mu^{ji}_q\in K$ and $(\alpha_{1q},\alpha_{2q},\ldots ,\alpha_{nq})\in\NZ^n$.\v5

{\bf A  characterization of solvable polynomial algebras}\vskip 6pt

Bearing in mind Definition 1.1.3 and the remark made above, we are
now able to present the main result of this section.
{\parindent=0pt\v5

{\bf 1.4.16. Theorem}  ([Li4], Theorem?) Let $A=K[a_1,\ldots ,a_n]$
be a finitely generated algebra over the field $K$, and let $\KS
=K\langle X_1,\ldots ,X_n\rangle$ be the free $K$-algebras with the
standard $K$-basis $\mathbb{B}=\{ 1,~X_{i_1}\cdots
X_{i_s}~|~X_{i_j}\in X,~s\ge 1\}$. With notations as before, the
following two statements are equivalent:\par

(i) $A$ is a solvable polynomial algebra in the sense of Definition
1.1.3.\par

(ii) $A\cong \OV A=\KS /I$ via the $K$-algebra epimorphism $\pi_1$:
$\KS \r A$ with $\pi_1(X_i)=a_i$, $1\le i\le n$, $I=$ Ker$\pi_1$,
satisfying  {\parindent=1.3truecm

\item{(a)} with respect to some monomial ordering $\prec_{_X}$ on $\mathbb{B}$, the ideal $I$ has a
finite Gr\"obner basis $G$ and the reduced Gr\"obner basis of $I$ is
of the form
$$\G =\left\{ g_{ji}=X_jX_i-\lambda_{ji}X_iX_j-F_{ji}
~\left |~\begin{array}{l} \LM (g_{ji})=X_jX_i,\\ 1\le i<j\le
n\end{array}\right. \right\} $$ where $\lambda_{ji}\in K^*$,
$\mu^{ji}_q\in K$, and
$F_{ji}=\sum_q\mu^{ji}_qX_1^{\alpha_{1q}}X_2^{\alpha_{2q}}\cdots
X_n^{\alpha_{nq}}$ with $(\alpha_{1q},\alpha_{2q},\ldots
,\alpha_{nq})\in\NZ^n$, thereby $\mathscr{B}=\{ \OV
X_1^{\alpha_1}\OV X_2^{\alpha_2}\cdots \OV X_n^{\alpha_n}~|~$
$\alpha_j\in\NZ\}$ forms  a PBW $K$-basis for $\OV A$, where each
$\OV X_i$ denotes the coset of $I$ represented by $X_i$ in $\OV A$;
and

\item{(b)} there is a  monomial ordering
$\prec$ on $\mathscr{B}$ such that $\LM (\OV{F}_{ji})\prec \OV
X_i\OV X_j$ whenever $\OV F_{ji}\ne 0$, where $\OV
F_{ji}=\sum_q\mu^{ji}_q\OV X_1^{\alpha_{1i}}\OV
X_2^{\alpha_{2i}}\cdots \OV X_n^{\alpha_{ni}}$, $1\le i<j\le n$.
\vskip 6pt}

{\bf Proof} (i) $\Rightarrow$ (ii) Let $\B=\{
a^{\alpha}=a_1^{\alpha_1}\cdots a_n^{\alpha_n}~|~\alpha
=(\alpha_1,\ldots ,\alpha_n)\in\NZ^n\}$ be the PBW $K$-basis of the
solvable polynomial algebra $A$ and $\prec$ a  monomial ordering on
$\B$. By Definition 1.1.3, the generators of $A$ satisfy the
relations:
$$a_ja_i=\lambda_{ji}a_ia_j+f_{ji},\quad 1\le i<j\le n,\eqno{(*)}$$
where $\lambda_{ji}\in K^*$ and  $f_{ji}=\sum_q\mu^{ji}_qa^{\alpha
(q)}\in K$-span$\B$ with $\LM (f_{ji})\prec a_ia_j$. Consider in the
free $K$-algebra $\KS =K\langle X_1,\ldots ,X_n\rangle$ the subset
$$\G =\{ g_{ji}=X_jX_i-\lambda_{ji}X_iX_j-F_{ji}~|~1\le i<j\le n\} ,$$
where if
$f_{ji}=\sum_q\mu^{ji}_qa_1^{\alpha_{1q}}a_2^{\alpha_{2q}}\cdots
a_n^{\alpha_{nq}}$ then
$F_{ji}=\sum_q\mu^{ji}_qX_1^{\alpha_{1q}}X_2^{\alpha_{2q}}\cdots
X_n^{\alpha_{nq}}$ for $1\le i<j\le n$. We write
$J=\langle\G\rangle$ for the ideal of $\KS$ generated by $\G$ and
put $\OV A=\KS /J$. Let $\pi_1$: $\KS\r A$ be the $K$-algebra
epimorphism with $\pi_1(X_i)=a_i$, $1\le i\le n$, and let $\pi_2$:
$\KS\r\OV A$ be the canonical algebra epimorphism. It follows from
the universal property of the canonical homomorphism that there is
an algebra epimorphism $\varphi$: $\OV A\r A$ defined by $\varphi
(\OV X_i)=a_i$, $1\le i\le n$, such that the following diagram of
algebra homomorphisms is commutative:
$$\begin{array}{cccc}
\KS&\mapright{\pi_2}{}&\OV A&\\
\mapdown{\pi_1}&\swarrow\scriptstyle{\varphi}&&\varphi\circ\pi_2=\pi_1\\
A&&&\end{array}$$  On the other hand, by the definition of each
$g_{ji}$ we see that every element $\OV H\in\OV A$ may be written as
$\OV H=\sum_j\mu_j\OV X_1^{\beta_{1j}}\OV X_2^{\beta_{2j}}\cdots \OV
X_n^{\beta_{nj}}$ with $\mu_j\in K$ and $(\beta_{1j},\ldots
,\beta_{nj})\in\NZ^n$, where each $\OV X_i$ is the coset of $J$
represented by $X_i$ in $\OV A$. Noticing the relations presented in
$(*)$,  it is straightforward to check that the correspondence
$$\begin{array}{cccc} \psi :&A&\mapright{}{}&\OV A\\
&\displaystyle{\sum_i}\lambda_ia_1^{\alpha_{1i}}\cdots
a_n^{\alpha_{ni}}&\mapsto&\displaystyle{\sum_i}\lambda_i\OV
X_1^{\alpha_{1i}}\cdots \OV X_n^{\alpha_{ni}}\end{array}$$ is an
algebra homomorphism such that $\varphi\circ\psi =1_A$ and
$\psi\circ\varphi =1_{\OV A}$, where $1_A$ and $1_{\OV A}$ denote
the identity maps of $A$ and $\OV A$ respectively. This shows that
$A\cong\OV A$, thereby Ker$\pi_1=I=J$; moreover, $\mathscr{B}=\{ \OV
X_1^{\alpha_1}\OV X_2^{\alpha_2}\cdots \OV
X_n^{\alpha_n}~|~\alpha_j\in\NZ\}$ forms  a PBW $K$-basis for $\OV
A$, and $\prec$ is a  monomial ordering on $\mathscr{B}$.}}
\par
We next show that $\G$ forms the reduced Gr\"obner basis for $I$ as
described in (a). To this end, we first show that the monomial
ordering $\prec$ on $\B$ induces a monomial ordering $\prec_{_X}$ on
the standard $K$-basis $\mathbb{B}$ of $\KS$. For convenience,
recall that we have used capital letters $U,V,W,S,\ldots$ to denote
elements (monomials) in $\mathbb{B}$. Also let us fix a graded
lexicographic ordering $\prec_{grlex}$ on $\mathbb{B}$ with respect
to a given positive-degree function $d(~)$ on $\KS$ (see the
definition given before Definition 1.4.1), such that
$$X_1\prec_{lex}X_2\prec_{lex}\cdots\prec_{lex}X_n.$$
Then, for $U,V\in\mathbb{B}$ we define
$$U\prec_{_X}V~\hbox{if}~\left\{\begin{array}{l} \LM (\pi_1 (U))\prec\LM (\pi_1(V)),\\
\hbox{or}\\
\LM (\pi_1 (U))=\LM (\pi_1
(V))~\hbox{and}~U\prec_{grlex}V.\end{array}\right.$$ Since  $A$ is a
domain (Proposition 1.1.4(ii)) and $\pi_1$ is an algebra
homomorphism with $\pi_1(X_i)=a_i$ for $1\le i\le n$, it follows
that $\LM (\pi_1(W))\ne 0$ for all $W\in\mathbb{B}$. We also note
from Proposition 1.1.4(i) that if $f,g\in A$ are nonzero elements,
then $\LM (fg)=\LM (\LM (f)\LM (g))$. Thus, if $U,V,W\in\mathbb{B}$
and $U\prec_{_X}V$ subject to $\LM (\pi_1(U))\prec\LM (\pi_1(V))$,
then
$$\begin{array}{cc}\begin{array}{rcl} \LM
(\pi_1(WU))&=&\LM (\LM (\pi_1(W))\LM (\pi_1(U)))\\
&\prec&\LM (\LM (\pi_1(W))\LM (\pi_1(V)))\\
&=&\LM (\pi_1(WV))\end{array}&\end{array}$$ implies
$WU\prec_{_X}WV$; if $U\prec_{_X}V$ subject to $\LM (\pi_1 (U))=\LM
(\pi_1 (V))$ and $U\prec_{grlex}V$, then
$$\begin{array}{rcl} \LM (\pi_1(WU))&=&\LM (\LM (\pi_1(W))\LM (\pi_1(U)))\\
&=&\LM (\LM (\pi_1(W))\LM (\pi_1(V)))\\
&=&\LM (\pi_1(WV))\end{array}$$ \\
and $WU\prec_{grlex}WV$ implies $WU\prec_{_X}WV$.  Similarly, if
$U\prec_{_X}V$ then $US\prec_{_X}VS$ for all $S\in\mathbb{B}$.
Moreover, if $W,U,V,S\in\mathbb{B}$, $W\ne V$, such that $W=UVS$,
then $\LM (\pi_1(W))=\LM (\pi_1(UVS))$ and clearly
$V\prec_{grlex}W$, thereby $V\prec_{_X}W$. Since  $\prec$ is a
well-ordering on $\B$ and $\prec_{grlex}$ is a well-ordering on
$\mathbb{B}$, the above argument shows that $\prec_{_X}$ is a
 monomial ordering on $\mathbb{B}$. With this monomial
ordering $\prec_{_X}$ in hand,  by the definition of $F_{ji}$ we see
that $\LM (F_{ji})\prec_{_X}X_iX_j$. Furthermore, since $\LM
(\pi_1(X_jX_i))=a_ia_j=\LM (\pi_1(X_iX_j))$ and $X_iX_j\prec_{grlex}
X_jX_i$, we see that $X_iX_j\prec_{_X}X_jX_i$. It follows that $\LM
(g_{ji})=X_jX_i$ for $1\le i<j\le n$. Now, by Proposition 1.4.14 we
conclude that $\G$ forms a Gr\"obner basis for $I$ with respect to
$\prec_{_X}$. Finally, by the definition of $\G$, it is clear that
$\G$ is the reduced Gr\"obner basis of $I$, as desired.\par

(ii) $\Rightarrow$ (i)  Note that (a) $+$ (b) tells us that  the
generators of $\OV A$ satisfy the relations $\OV X_j\OV
X_i=\lambda_{ji}\OV X_i\OV X_j+\OV F_{ji}$, $1\le i<j\le n$, and
that if $\OV F_{ji}\ne 0$ then $\LM (\OV F_{ji})\prec \OV X_i\OV
X_j$ with respect to the given monomial ordering $\prec$ on
$\mathscr{B}$. It follows that  $\OV A$ and hence $A$ is a solvable
polynomial algebra in the sense of Definition 1.1.3.\QED
{\parindent=0pt\v5

{\bf Remark} The monomial ordering $\prec_{_X}$ we defined in the
proof of Theorem 1.4.16 is a modification of the {\it lexicographic
extension} defined in [EPS]. But our definition of $\prec_{_X}$
involves a graded monomial ordering $\prec_{grlex}$ on the standard
$K$-basis $\mathbb{B}$ of the free $K$-algebra $\KS =K\langle
X_1,\ldots ,X_n\rangle$. The reason is that {\it the monomial
ordering $\prec_{_X}$ on $\mathbb{B}$ must be compatible with the
usual rule of division}, namely, $W,U,V,S\in\mathbb{B}$, $W\ne V$,
and $W=UVS$ implies $V\prec_{_X}W$. While it is clear that if we use
any lexicographic ordering $\prec_{lex}$ in the definition of
$\prec_{_X}$, then this rule will not work in general.}\v5

We end this section by an example illustrating Theorem 1.4.16, in
particular, illustrating that the monomial ordering $\prec_{_X}$
used in the condition (a) and the monomial ordering $\prec$ used in
the condition (b) may be mutually independent, namely $\prec$ may
not necessarily be the restriction of $\prec_{_X}$ on $\mathscr{B}$,
and the choice of $\prec$ is indeed quite
flexible.{\parindent=0pt\v5

{\bf Example} (1)  Considering the positive-degree function $d(~)$ 
on the free $K$-algebra $\KS =K\langle X_1,X_2, X_3\rangle$ such 
that $d(X_1)=2$, $d(X_2)=1$, and $d(X_3)=4$, let $I$ be the ideal of 
$\KS$ generated by the elements
$$\begin{array}{l} g_1=X_1X_2- X_2X_1,\\
g_2=X_3X_1-\lambda X_1X_3-\mu X_3X_2^2-f(X_2),\\
g_3=X_3X_2- X_2X_3,\end{array}$$ where $\lambda\in K^*$, $\mu\in K$,
$f(X_2)$ is a polynomial in $X_2$ which has degree $\le 6$, or
$f(X_2)=0$. The following properties hold.   \par

(1) If we use the graded lexicographic ordering
$X_2\prec_{grlex}X_1\prec_{grlex}X_3$ on $\KS$, then the three
generators have the leading monomials $\LM (g_1)=X_1X_2$, $\LM
(g_2)=X_3X_1$, and $\LM (g_3)=X_3X_2$. It is an exercise to verify
that $\G =\{ g_1,g_2,g_3\}$ forms a Gr\"obner basis for $I$ with
respect to $\prec_{grlex}$.
\par

(2) With respect to the fixed $\prec_{grlex}$ in (1), the reduced
Gr\"obner basis $\G '$ of $I$ consists of
$$\begin{array}{l} g_1=X_1X_2- X_2X_1,\\
g_2=X_3X_1- \lambda X_1X_3-\mu X_2^2X_3-f(X_2),\\
g_3=X_3X_2- X_2X_3,\end{array}$$\par

(3) Writing $A=K[a_1,a_2, a_3]$ for the quotient algebra $\KS /I$,
where $a_1$, $a_2$ and $a_3$ denote the cosets $X_1+I$, $X_2+I$ and
$X_3+I$ in $\KS /I$ respectively, it follows that $A$ has the PBW
basis $\B =\{
a^{\alpha}=a_2^{\alpha_2}a_1^{\alpha_1}a_3^{\alpha_3}~|~\alpha
=(\alpha_2,\alpha_1,\alpha_3)\in\NZ^3\}$. Noticing that $
a_2a_1=a_1a_2$, it is clear that $\B ' =\{
a^{\alpha}=a_1^{\alpha_1}a_2^{\alpha_2}a_3^{\alpha_3}~|~\alpha
=(\alpha_1,\alpha_2,\alpha_3)\in\NZ^3\}$ is also a PBW basis for
$A$. Since $a_3a_1=\lambda a_1a_3+\mu a_2^2a_3+f(a_2)$, where
$f(a_2)\in K$-span$\{ 1,a_2,a_2^2,\ldots ,a_2^6\}$, we see that $A$
has the monomial ordering $\prec_{lex}$ on $\B '$ such that
$a_3\prec_{lex}a_2\prec_{lex}a_1$ and $\LM (\mu
a_2^2a_3+f(a_2))\prec_{lex}a_1a_3$, thereby $A$ is turned into a
solvable polynomial algebra with respect to $\prec_{lex}$.}\par

Moreover, it is also an exercise to check that if  we use the
positive-degree function $d(~)$ on $\B '$ such that $d(a_1)=2$,
$d(a_2)=1$, and $d(a_3)=4$, then $A$ has another monomial ordering
on $\B '$, namely the graded lexicographic ordering $\prec_{grlex}$
such that $a_3\prec_{grlex}a_2\prec_{grlex}a_1$ and $\LM (\mu
a_2^2a_3+f(a_2))\prec_{grlex}a_1a_3$, thereby $A$ is turned into a
solvable polynomial algebra with respect to $\prec_{grlex}$. 
\newpage\setcounter{page}{33}

\chapter*{2. Left Gr\"obner Bases for\\ \hskip 1.25truecm  Modules }\par
\vskip 2.5truecm\markboth{\rm Gr\"obner Bases for Modules}{ \rm 
Gr\"obner Bases for Modules}

Based on the theory of left Gr\"obner bases for left ideals of
solvable polynomial algebras presented in Chapter 1, in this chapter
we introduce {\it left Gr\"obner bases} for submodules of {\it free
left modules} over solvable polynomial algebras. More precisely, let
$A=K[a_1,\ldots ,a_n]$ be a solvable polynomial algebra with
admissible system $(\B ,\prec )$ as described in Chapter 1. In the
first section, we define left monomial orderings, especially the
graded left monomial orderings and Schreyer orderings,  on free
$A$-modules. In Section 2, we introduce  left Gr\"obner bases,
minimal left Gr\"obner bases and reduced left Gr\"obner bases for
submodules of free left $A$-modules via a left division algorithm,
and we discuss some basic properties of left Gr\"obner bases, in
particular, we show that  every submodule of a free left $A$-module
$L=\oplus_{I=1}^sAe_i$ has a finite left Gr\"obner basis. In Section
3, we define left S-polynomials for pairs of elements in free
$A$-modules, and we show that a noncommutative version of
Buchberger's criterion holds true for  submodules of free
$A$-modules, and that a noncommutative version of the Buchberger 
algorithm works effectively for computing finite left Gr\"obner 
bases of submodules in free $A$-modules.
 \par

The main references of this chapter are [AL2], [Eis], [KR1], [K-RW],
[Li1], [Lev], [DGPS].\par

Throughout this chapter, modules are meant left $A$-modules, and all
notions and notations used in Chapter 1 are maintained.

\section*{2.1. Left Monomial Orderings on Free Modules}

Let $A=K[a_1,\ldots ,a_n]$ be a solvable polynomial algebra with
admissible system $(\B ,\prec )$ in the sense of Definition 1.1.3,
where $\B =\{ a^{\alpha}=a_1^{\alpha_1}\cdots
a_n^{\alpha_n}~|~\alpha =(\alpha_1,\ldots ,\alpha_n)\in\NZ^n\}$ is
the PBW $K$-basis of $A$ and $\prec$ is a monomial ordering on $\B$,
and let $L=\oplus_{i=1}^sAe_i$ be a free $A$-module with the
$A$-basis $\{ e_1,\ldots ,e_s\}$. Then $L$ has the $K$-basis
$$\BE =\{ a^{\alpha}e_i~|~a^{\alpha}\in\B ,~1\le i\le s\} .$$ For
convenience, elements of $\BE$ are also referred to as {\it
monomials} in $L$.\v5

If $\prec_{e}$ is a total ordering on $\BE$, and if $\xi
=\sum_{j=1}^m\lambda_ja^{\alpha (j)}e_{i_j}\in L$, where
$\lambda_j\in K^*$ and $\alpha (j)=(\alpha_{j_1},\ldots
,\alpha_{j_n})\in\NZ^n$, such that $$a^{\alpha (1)}e_{i_1}\prec_{e}
a^{\alpha (2)}e_{i_2}\prec_{e}\cdots\prec_{e} a^{\alpha
(m)}e_{i_m},$$ then by $\LM (\xi )$ we denote the {\it leading
monomial} $a^{\alpha (m)}e_{i_m}$ of $\xi $, by $\LC (\xi )$ we
denote the {\it leading coefficient} $\lambda_m$ of $\xi $,  and by
$\LT (\xi )$ we denote the {\it leading term} $\lambda_ma^{\alpha
(m)}e_{i_m}$ of $f$.{\parindent=0pt\v5

{\bf 2.1.1. Definition} With respect to the given monomial ordering
$\prec$ on $\B$, a total ordering $\prec_{e}$ on $\BE$ is called a
{\it left monomial ordering} if the following two conditions are
satisfied:
\par

(1) $a^{\alpha}e_i\prec_{e} a^{\beta}e_j$ implies  $\LM
(a^{\gamma}a^{\alpha}e_i)\prec_{e} \LM (a^{\gamma}a^{\beta}e_j)$ for
all $a^{\alpha}e_i$, $a^{\beta}e_j\in\BE$, $a^{\gamma}\in\B$;\par

(2) $a^{\beta}\prec a^{\beta}$ implies $a^{\alpha}e_i\prec_{e}
a^{\beta}e_i$ for all $a^{\alpha},a^{\beta}\in\B$ and $1\le i\le
s$.} \par

If $\prec_e$ is a left monomial ordering on $\BE$, then we also say
that $\prec_e$ is a left monomial ordering on the free module $L$,
and the data $(\BE ,\prec_e)$ is referred to as an {\it left
admissible system} of $L$.\v5

By referring to Proposition 1.1.4, we record two easy but useful
facts on a left monomial ordering  $\prec_{e}$ on $\BE$, as follows.

{\parindent=0pt\v5

{\bf 2.1.2. Lemma} (i) Every left monomial ordering  $\prec_{e}$ on
$\BE$ is a well-ordering, i.e., every nonempty subset of $\BE$ has a
minimal element. \par

(ii) If  $f\in A$ with $\LM (f)=a^{\gamma}$ and $\xi\in L$ with $\LM
(\xi )=a^{\alpha}e_i$, then
$$\LM (f\xi )=\LM (\LM (f)\LM (\xi ))=\LM (a^{\gamma}a^{\alpha}e_i)=a^{\gamma +\alpha}e_i.$$
\par\QED}\v5

Actually as in the commutative case ([AL2], [Eis], [KR1]), any
monomial ordering $\prec$ on $\B$ may induce  two  left monomial
orderings on $\BE$:
$$\begin{array}{l} (\hbox{{\bf TOP} ordering})\quad a^{\alpha}e_i\prec_{e} a^{\beta}e_j\Leftrightarrow
a^{\alpha}\prec a^{\beta},~\hbox{or}~ a^{\alpha}=a^{\beta}~
\hbox{and}~i<j;\\
(\hbox{{\bf POT} ordering})\quad a^{\alpha}e_i\prec_{e}
a^{\beta}e_j\Leftrightarrow i<j,~
\hbox{or}~i=j~\hbox{and}~a^{\alpha}\prec a^{\beta}.\end{array}$$ \v5

Let $d(~)$ be a positive-degree function on $A$ (see Section 1 of
Chapter 1) such that $d(a_i)=m_i>0$, $1\le i\le n$, and let
$(b_1,\ldots b_s)\in\NZ^n$ be any fixed $s$-tuple. Then, by
assigning $e_j$ the degree $b_j$, $1\le j\le s$, every monomial
$a^{\alpha}e_j$ in the $K$-basis $\BE$ of $L$ is endowed with the
degree $d(a^{\alpha})+b_j$. Similar to Definition 1.1.2, if a left
monomial ordering $\prec_{e}$ on $\BE$ satisfies
$$a^{\alpha}e_i\prec_ea^{\beta}e_j~\hbox{implies}~d(a^{\alpha})+b_i\le d(a^{\beta})+b_j,$$
then we call it a {\it graded left monomial ordering} on
$\BE$.\par

We refer the reader to Chapter 4 and Chapter 5 for examples of
graded left monomial orderings, also for the reason that a graded
left monomial ordering is defined in such a way. \v5

Let  $\prec_{e}$ be a  left monomial ordering on $\BE$, and let
$L_1=\oplus_{i=1}^mA\varepsilon_i$ be another free $A$-module with
the $A$-basis $\{ \varepsilon_1,\ldots ,\varepsilon_m\}$. Then, as
in the commutative case ([AL2], [Eis], [KR1]), for any given finite
subset $G=\{ g_1,\ldots ,g_m\}\subset L$,  an ordering on the
$K$-basis $\B (\varepsilon )=\{
a^{\alpha}\varepsilon_i~|~a^{\alpha}\in\B ,~1\le i\le m\}$ of $L_1$
can be defined subject to the rule: for
$a^{\alpha}\varepsilon_i,a^{\beta}\varepsilon_j\in\B (\varepsilon
)$,
$$a^{\alpha}\varepsilon_i\prec_{s\hbox{-}\varepsilon} a^{\beta}\varepsilon_j
\Leftrightarrow\left\{\begin{array}{l}
\LM (a^{\alpha}g_i)\prec_{e}\LM (a^{\beta}g_j),\\
\hbox{or}\\
\LM (a^{\alpha}g_i)=\LM (a^{\beta}g_j)~\hbox{and}~i<j,
\end{array}\right.$$
It is an exercise to check that this ordering is a left monomial
ordering on  $\B (\varepsilon )$.  $\prec_{s\hbox{-}\varepsilon}$ is
usually referred to as the {\it Schreyer ordering} induced by $G$
with respect to $\prec_{e}$.   \v5

\section*{2.2. Left Gr\"obner Bases of Submodules}

Let $A=K[a_1,\ldots ,a_n]$ be a solvable polynomial algebra with
admissible system $(\B ,\prec )$, and let $L=\oplus_{i=1}^sAe_i$ be
a free $A$-module with left admissible system $(\BE ,\prec_e)$. In
this section we introduce left Gr\"obner bases for submodules of $L$
via a left division algorithm in $L$, and we record some basic facts
determined by left Gr\"obner bases. Moreover, minimal left Gr\"obner
bases and reduced left Gr\"obner bases are discussed. Finally, we
show that  every nonzero submodule of $L$ has a finite left
Gr\"obner basis. {\parindent=0pt\v5

{\bf Left division algorithm}\vskip 6pt

Let $a^{\alpha}e_i$, $a^{\beta}e_j\in\BE$, where $\alpha
=(\alpha_1,\ldots ,\alpha_n)$, $\beta =(\beta_1,\ldots
,\beta_n)\in\NZ^n$. Then, the left division of monomials in $\B$ we
introduced in (Section 2 of Chapter 1) gives rise to a left division
for monomials in $\BE$, that is, we say that {\it $a^{\alpha}e_i$
{\it divides $a^{\beta}e_j$} from left side}, denoted
$a^{\alpha}e_i|_{_{\rm L}}a^{\beta}e_j$, if $i=j$ and there is some
$a^{\gamma}\in\B$ such that
$$a^{\beta}e_i=\LM (a^{\gamma}a^{\alpha}e_i).$$
It follows from Lemma 2.1.2(ii) that the division defined above is
implementable.}\v5

Let $\Xi$ be a nonempty subset of $L$. Then the division of
monomials defined above yields a subset of $\BE$:
$$\mathcal {N}(\Xi )=\{ a^{\alpha}e_i\in\BE~|~\LM (\xi)\not{|_{_{\rm L}}}~a^{\alpha}e_i,~\xi\in \Xi\}.$$
If $\Xi =\{ \xi\}$ consists of a single element $\xi$, then we
simply write $\mathcal {N}(\xi )$ in place of $\mathcal {N}(\Xi )$. 
Also we write $K$-span$\mathcal {N}(\Xi )$ for the $K$-subspace of 
$L$ spanned by $\mathcal {N}(\Xi )$.{\parindent=0pt\v5

{\bf 2.2.1. Definition} Elements of $\mathcal {N}(\Xi )$ are
referred to as {\it normal monomials} (mod $\Xi $). Elements of
$K$-span$\mathcal {N}\Xi )$ are referred to as {\it normal elements} 
(mod $\Xi $).}\v5

In view of Lemma 2.1.2(ii), the left division we defined for
monomials in $\BE$ can naturally be used to define a left division
procedure for elements in $L$. More precisely, let $\xi, \zeta\in L$
with $\LC (\xi)=\mu \ne 0$, $\LC (\zeta)=\lambda\ne 0$. If $\LM
(\zeta)|_{_{\rm L}}\LM (\xi)$, i.e., there exists $a^{\alpha}\in\B$
such that $\LM (\xi)=\LM (a^{\alpha}\LM (\zeta))$, then put
$\xi_1=\xi-\lambda^{-1}\mu a^{\alpha}\zeta$; otherwise, put
$\xi_1=\xi-\LT (\xi)$. Note that in both cases we have $\xi_1=0$, or
$\xi_1\ne 0$ and $\LM (\xi_1)\prec_e\LM (\xi )$. At this stage, let
us refer to such a procedure of canceling the leading term of $\xi$  
as the {\it left division procedure by $\zeta$}.  With $\xi 
:=\xi_1\ne 0$, we can repeat the left division procedure by $\zeta$ 
and so on. This returns successively   a descending sequence $$\LM 
(\xi_{i+1} )\prec_e\LM (\xi_i)\prec_e\cdots\prec_e\LM 
(\xi_1)\prec_e\LM (\xi ),\quad i\ge 2.$$  Since $\prec_e$ is a 
well-ordering, it follows that such a division procedure terminates 
after a finite number of repetitions, and consequently $\xi$ is 
expressed as
$$\xi =q\zeta+\eta,$$
where $q\in A$ and $\eta\in K$-span$\mathcal {N}(\zeta )$, i.e.,
$\eta$ is normal (mod $\zeta$), such that either $\LM (\xi )=\LM
(q\zeta )$ or $\LM (\xi )=\LM (\eta )$.\par

Actually as in (Section 2 of Chapter 1), the left division procedure
demonstrated above can be extended to a left division procedure by a 
finite subset $\Xi =\{ \xi_1,\ldots ,\xi_s\}$ in $L$, which yields 
the following left division algorithm:{\parindent=0pt\v5

\underline{\bf Algorithm-DIV-L
~~~~~~~~~~~~~~~~~~~~~~~~~~~~~~~~~~~~~~~~~~~~~~~~~~~~~~~~~~~}\vskip 
6pt

\textsc{INPUT}: $\xi ,~\Xi 
=\{\xi_1,\ldots,\xi_s\}~\hbox{with}~\xi_i\ne 0~(1\le i\le s)$\par  
\textsc{OUTPUT}:~$q_1,\ldots ,q_s\in A,~\eta\in
K\hbox{-span}\mathcal {N}(\Xi ),~\hbox{such that}$\par 
~~~~~~~~~~~~~~~~~$\xi =\sum^s_{i=1}q_i\xi_i+\eta,~\LM
(q_i\xi_i)\preceq_e\LM (\xi )~\hbox{for}~q_i\ne 0,$\par 
~~~~~~~~~~~~~~~~~~~~~~~~~~~~~~~~~~~~~~~~~~~$\LM (\eta )\preceq_e\LM
(\xi )~\hbox{if} ~\eta\ne 0$\par                                         
\textsc{INITIALIZATION}:~$q_1:=0,~q_2:=0,~\cdots 
,~q_s:=0;~\eta:=0;~\omega :=\xi$\par

\textsc{BEGIN}\par ~~~~~\textsc{WHILE}~$\omega \ne 0$ 
\textsc{DO}\par                                                
~~~~~~~~~~~\textsc{IF}~\hbox{there exist}~$i$~\hbox{and}~$a^{\alpha 
(i)}\in\B$~\hbox{such that}\par ~~~~~~~~~~~~~~~$\LM (\omega )=\LM 
(a^{\alpha (i)}\LM (\xi_i))$ \textsc{THEN}\par                         
~~~~~~~~~~~~~~~\hbox{choose}~$i$~\hbox{least such that}~$\LM (\omega 
)=\LM (a^{\alpha (i)}\LM (\xi_i))$\par             
~~~~~~~~~~~$q_i:=q_i+\LC (\xi_i)^{-1}\LC (\omega )a^{\alpha 
(i)}$\par                                                      
~~~~~~~~~~~$\omega :=\omega -\LC (\xi_i)^{-1}\LC (\omega )a^{\alpha 
(i)}\xi_i$\par                                 
~~~~~~~~~~~\textsc{ELSE}\par                                   
~~~~~~~~~~~$\eta:=\eta+\LT (\omega )$\par      
~~~~~~~~~~~$\omega:=\omega-\LT (\omega )$\par 
~~~~~~~~~~\textsc{END}\par ~~~~~\textsc{END}\par                   
\textsc{END}\par\vskip -.2truecm 
\underline{~~~~~~~~~~~~~~~~~~~~~~~~~~~~~~~~~~~~~~~~~~~~~~~~~~~~~~~~~~~~~~~~~~~~~~~~~~~~~~~~~~~~~~~~~~~~~~~} 
\v5

{\bf 2.2.2. Definition} The element $\eta$ obtained in {\bf
Algorithm-DIV-L} is called a {\it remainder} of $\xi$ on left
division by $\Xi$, and is denoted by $\OV{\xi}^{\Xi}$, i.e.,
$\OV{\xi}^{\Xi}=\eta$. If $\OV{\xi}^{\Xi}=0$, then we say that $\xi$ 
{\it is reduced to $0$} (mod $\Xi$).\v5

{\bf Remark} For the reason that the element $\eta$ obtained in {\bf
Algorithm-DIV-L} depends on how an order of elements in $\Xi$ is
arranged,  we used the phrase ``a remainder" in the above 
definition.  In other words,  the element $\eta$ may be different if 
a different order for elements in $\Xi$ is given. }\v5

Summing up, we have reached the following{\parindent=0pt\v5

{\bf 2.2.3. Theorem} Given a set of nonzero elements $\Xi =\{
\xi_1,\ldots ,\xi_s\}$ and $\xi$ in $L$, the {\bf Algorithm-DIV-L}
produces elements $q_1,\ldots ,q_s\in A$ and $\eta\in
K$-span$\mathcal {N}(\Xi )$, such that $\xi
=\sum^s_{i=1}q_i\xi_i+\eta$ and $\LM (q_i\xi_i)\preceq_e\LM (\xi )$
whenever $q_i\ne 0$, $\LM (\eta )\preceq_e\LM (\xi )$ if $\eta\ne
0$.\QED\v5

{\bf Left Gr\"obner bases}\vskip 6pt

It is theoretically correct that the left division procedure by
using a finite subset $\Xi $ of nonzero elements in $L$ can be
extended to a left division procedure by means of an {\it arbitrary
proper subset} $G$ of nonzero elements, for,  it is a true statement
that if $\xi\in L$ with $\LM (\xi )\ne 0$, then either there exists
$g\in G$ such that $\LM (g)|_{_{\rm L}}\LM (\xi )$ or such $g$ does
not exist.  This leads to the following\v5

{\bf 2.2.4. Proposition} Let $N$ be a submodule of $L$ and $G$ a
proper subset of nonzero elements in $N$. The following two
statements are equivalent:\par

(i) If $\xi\in N$ and $\xi\ne 0$, then there exists $g\in G$ such
that $\LM (g)|_{_{\rm L}}\LM (\xi )$.\par

(ii) Every nonzero $\xi\in N$ has a representation $\xi
=\sum^s_{j=1}q_ig_{i_j}$ with $q_i\in A$ and $g_{i_j}\in G$, such
that $\LM (q_ig_{i_j})\preceq_e\LM (\xi )$ whenever $q_i\ne
0$.\QED\v5

{\bf 2.2.5. Definition}  Let $N$ be a submodule of $L$. With respect
to a given left monomial ordering $\prec_e$ on $\BE$, a proper
subset $\G$ of nonzero elements in $N$ is said to be a {\it left
Gr\"obner basis } of $N$ if $\G$ satisfies one of the equivalent
conditions of Proposition 2.2.4.}
\par

If $\G$ is a left Gr\"obner basis of $N$, then the expression $\xi
=\sum^s_{j=1}q_ig_{i_j}$ appeared in Proposition 2.2.4(ii) is called
a {\it left Gr\"obner representation} of $\xi$. \v5

Clearly, if $N$ is a submodule $L$ and $N$ has a left Gr\"obner
basis $\G$, then $\G$ is certainly a generating set of $N$, i.e.,
$N=\sum_{g\in\G}Ag$. But the converse is not necessarily true. For
instance, consider the solvable polynomial algebra
$A=\mathbb{C}[a_1,a_2,a_3]$ generated by $\{ a_1,a_2,a_3\}$ subject
to the relations
$$a_2a_1=3a_1a_2,\quad a_3a_1=a_1a_3,\quad a_3a_2=5a_2a_3,$$
and let $L=Ae_1\oplus Ae_2\oplus Ae_3$ be the free $A$-module with
the $A$-basis $\{ e_1,e_2,e_3\}$.  Then, with
$g_1=a_1^2a_2e_1-a_3e_3$, $g_2=a_2e_1$, and $N=Ag_1+Ag_2$, we have
$a_3e_e\in N$. Hence the set  $G=\{ \xi_1,g_2\}$ is not a left
Gr\"obner basis with respect to any given monomial ordering
$\prec_e$ on $\BE$ such that $\LM (g_1)=a_1^2a_2e_1$, $\LM
(g_2)=a_2e_1$.{\parindent=0pt\v5

{\bf The existence of left Gr\"obner bases}\vskip 6pt

Noticing $1\in\B$, if we define on $\BE$ the ordering:
$$a^{\alpha}e_i\prec_e 'a^{\beta}e_j\Leftrightarrow a^{\alpha}e_i|_{_{\rm
L}}a^{\beta}e_j,\quad a^{\alpha}e_i,~a^{\beta}e_j~\hbox{in}~\BE ,$$
then, by Definition 2.1.1, Lemma 2.1.2(ii) and the divisibility
defined for monomials in $\BE$,  it is an easy exercise to check
that $\prec '$ is reflexive, anti-symmetric, transitive, and
moreover,
$$a^{\alpha}e_i\prec_e 'a^{\beta}e_j~\hbox{implies}~a^{\alpha}e_i\prec_e a^{\beta}e_j.$$
Since the given left monomial ordering $\prec_e$ is a well-ordering
on $\BE$, it follows that every nonempty subset of $\BE$ has a
minimal element with respect to the ordering $\prec_e '$ on $\BE$.
}\v5

Now, let $N\ne \{ 0\}$ be a submodule of $L$, $\LM (N)=\{ \LM (\xi
)~|~\xi \in N\}$, and $\Omega =\{ a^{\alpha}e_i\in \LM
(N)~|~a^{\alpha}e_i~\hbox{is minimal in}~\LM
(N)~\hbox{w.r.t.}\prec_e'\} .${\parindent=0pt\v5

{\bf 2.2.6. Theorem} With the notation as above, the following
statements hold.\par

(i) $\Omega\ne \emptyset$ and $\Omega$ is a proper subset of $\LM
(N)$.\par

(ii) Let $\G =\{ g\in N~|~\LM (g)\in\Omega\}$. Then $\G$ is a left
Gr\"obner basis of $N$.\vskip 6pt

{\bf Proof} (i) That $\Omega\ne\emptyset$ follows from $N\ne \{ 0\}$
and the remark about $\prec_e '$ we made above. Since $N$ is a
submodule of $L$ and $A$ is a domain (Proposition 1.1.4), if $\xi\in
N$ and $\xi\ne 0$, then $h\xi\in N$ for all nonzero $h\in A$.  It
follows from Lemma 2.1.2(ii) that $0\ne \LM (h\xi )\in\LM (N)$ and
$\LM (\xi )|_{_{\rm L}}\LM (h\xi )$, i.e., $\LM (\xi )\prec_e '\LM
(h\xi )$ in $\LM (N)$. It is clear that if $\LM (h)\ne 1$, then $\LM
(h\xi )\ne \LM (\xi )$. This shows that $\Omega$ is a proper subset
of $\LM (N)$.}\par

(ii) By the definition of $\prec_e '$ and the remark we made above
the theorem, if $\xi\in N$ with $\LM (\xi )\ne 0$ and $\LM (\xi
)\not\in\Omega$, then there is some $g\in\G$ such that $\LM
(g)|_{_{\rm L}}\LM (\xi )$. Hence the selected $\G$ is a left
Gr\"obner basis for $N$.\QED\v5

Since every solvable polynomial algebra $A$ is (left and right)
Noetherian by Corollary 1.3.3, the free module
$L=\oplus_{i=1}^sAe_i$ is a Noetherian left $A$-module, thereby
every submodule $N$ of $L$ is finitely generated.  In the next
section we will show that if a finite generating set $\Xi
=\{\xi_1,\ldots ,\xi_t\}$ of $N$ is given, then a {\it finite} left
Gr\"obner basis $\G$ for $N$ can be produced by means of a
noncommutative Buchberger algorithm. {\parindent=0pt\v5

{\bf Basic facts determined by left Gr\"obner bases}\vskip 6pt

By referring to Lemma 2.1.2 Theorem 2.2.3 and Proposition 2.2.4, the
foregoing discussion allows us to summarize some basic facts
determined by left Gr\"obner bases, of which the detailed proof is
left as an exercise.\v5

{\bf 2.2.7. Proposition} Let $N$ be a submodule of the free
$A$-module $L=\oplus_{i=1}^sAe_i$, and let $\G$ be a left Gr\"obner
basis of $N$. Then the following statements hold.\par

(i) If $\xi\in N$ and $\xi\ne 0$, then $\LM (\xi )=\LM (qg)$ for
some $q\in A$ and $g\in\G$.\par

(ii) If $\xi\in L$ and $\xi\ne 0$, then $\xi$ has a unique remainder
$\OV{\xi}^{\G}$ on division by $\G$. \par

(iii) If $\xi\in L$ and $\xi\ne 0$, then $\xi\in N$ if and only if
$\OV{\xi}^{\G}=0$. Hence the membership problem for submodules of
$L$ can be solved by using left Gr\"obner bases.\par

(iv) As a $K$-vector space, $L$ has the decomposition
$$L=N\oplus K\hbox{-span}\mathcal {N}(\G ).$$\par

(v) As a $K$-vector space, $L/N$ has the $K$-basis
$$\OV{\mathcal {N}(\G )}=\{ \OV{a^{\alpha}e_i}~|~a^{\alpha}e_i\in\mathcal {N}(\G )\} ,$$
where $\OV{a^{\alpha}e_i}$ denotes the coset represented by
$a^{\alpha}e_i$ in $L/N$.\par

(vi)  $\mathcal {N}(\G)$ is a finite set, or equivalently,
Dim$_KA/N<\infty$, if and only if  for $j=1,\ldots ,n$ and  each
$i=1, \ldots ,s$, there exists $g_{j_i}\in\G$ such that $\LM
(g_{j_i})=a_{j}^{m_{j}}e_i$, where $m_{j}\in\NZ$. \QED  \v5

{\bf Minimal and reduced left Gr\"obner bases}\vskip 6pt

{\bf 2.2.8. Definition} Let $N\ne \{ 0\}$ be a submodule of $L$ and
let $\G$ be a left gr\"obner basis of $N$. If any proper subset of
$\G$ cannot be a left Gr\"obner basis of $N$, then $\G$ is called a
{\it minimal left Gr\"obner basis} of $N$.\v5

{\bf 2.2.9. Proposition} Let $N\ne \{ 0\}$ be a submodule of $L$. A
left Gr\"obner basis $\G$ of $N$ is minimal if and only if ~$\LM
(g_i){\not |_{_{\rm L}} }\LM (g_j)$ for all $g_i,g_j\in\G$ with
$g_i\ne g_j$.\vskip 6pt

{\bf Proof} Suppose that $\G$ is minimal. If there were $g_i\ne g_j$
in $\G$ such that  $\LM (g_i)|_{_{\rm L}}\LM (g_j)$, then since the
left division is transitive, the proper subset $\G'=\G-\{ g_j\}$ of
$\G$ would form a left Gr\"obner basis of $N$. This contradicts the
minimality of $\G$.}\par

Conversely, if the condition $\LM (g_i){\not |_{_{\rm L}}}\LM (g_j)$
holds for all $g_i,g_j\in\G$ with $g_i\ne g_j$, then  the definition
of a left Gr\"obner basis entails that any proper subset of $\G$
cannot be a left Gr\"obner basis of $N$. This shows that $\G$ is
minimal.\QED{\parindent=0pt\v5

{\bf 2.2.10. Corollary} Let $N\ne \{ 0\}$ be a submodule of $L$, and
let the notation be as in Theorem 2.2.6.\par

(i) The left Gr\"obner basis $\G$ obtained there is indeed a minimal
left Gr\"obner basis for $N$.\par

(ii) If $G$ is any left Gr\"obner basis of $N$, then
$$\Omega =\LM (\G )=\{ \LM (g)~|~g\in\G\}\subseteq\LM (G)=\{\LM (g')~|~g'\in G\} .$$
Therefore, any two minimal left Gr\"obner bases of $N$ have the same
set of leading monomials $\Omega$.\vskip 6pt

{\bf Proof} This follows immediately from the definition of a left
Gr\"obner basis, the construction  of $\Omega$ and Proposition
2.2.9.\QED}\v5

We next introduce the notion of a reduced left Gr\"obner
basis.{\parindent=0pt\v5

{\bf 2.2.11. Definition} Let $N\ne \{ 0\}$ be a submodule of $L$ and
let $\G$ be a left gr\"obner basis of $N$. If $\G$ satisfies the
following conditions:\par

(1) $\G$ is a minimal left Gr\"obner basis;\par

(2) $\LC (g)=1$ for all $g\in\G$;\par

(3) For every $g\in\G$, $\xi =g-\LM (g)$ is a normal element (mod
$\G$), i.e., $\xi\in K$-span$\mathcal {N}(\G )$,\par

then $\G$ is called a {\it reduced left Gr\"obner basis} of $N$.}\v5

{\parindent=0pt\v5

{\bf 2.2.12. Proposition} Every nonzero submodule $N$ of $L$ has a
unique reduced Gr\"obner basis.\vskip 6pt

{\bf Proof}  Note that if $\G$ and $\G '$ are reduced left Gr\"obner
bases for $N$, then $\LM (\G )=\LM (\G ')$ by Corollary 2.2.10. If
$\LM (g_i)=\LM (g_j')$ then $g_i-g_j'\in N\cap K$-span$\mathcal
{N}(\G )=N\cap K\hbox{-span}\mathcal {N}(\G ')$. It follows from
Proposition 2.2.7(iv) that  $g_i=g_j'$. Hence $\G =\G '$. \QED}\v5

The next proposition shows that if a finite left Gr\"obner basis
$\G$ of $N$ is given, then a minimal left Gr\"obner basis, thereby
the reduced left Gr\"obner basis for $N$ can be obtained in an
algorithmic way. {\parindent=0pt\v5

{\bf 2.2.13. Proposition} Let $N\ne \{ 0\}$ be a submodule of $L$,
and let $\G =\{ g_1,\ldots ,g_m\}$ be a finite left Gr\"obner basis
of $N$.
\par

(i) The subset $\G_0=\{ g_i\in\G~|~\LM (g_i)~\hbox{is minimal
in}~\LM (\G )~\hbox{w.r.t.}\prec_e '\} $ of $\G$ forms a minimal
left Gr\"obner basis of $N$ (see the definition of $\prec_e '$ given
before Theorem 2.2.6).  An algorithm written in pseudo-code is
omitted here.\par

(ii) With $\G_0$ as in (i) above, we may assume, without loss of
generality, that $\G_0=\{ g_1,\ldots ,g_s\}$ with $\LC (g_i)=1$ for
$1\le i\le s$. Put $\G_1=\{ g_2,...,g_s\}$ and
$\xi_1=\OV{g_1}^{\G_1}$. Then $\LM (\xi_1)=\LM (g_1)$. Put $\G_2=\{
\xi_1,g_3,...,g_s\}$ and $\xi_2=\OV{g_2}^{\G_2}$. Then  $\LM
(\xi_2)=\LM (g_2)$. Put $\G_3=\{ \xi_1,\xi_2,g_4,...,g_s\}$ and
$\xi_3=\OV{g_3}^{\G_3}$, and so on. The last obtained
$\G_{s+1}=\{\xi_1,\ldots ,\xi_{s-1}\}\cup\{\xi_s\}$ is then the
reduced left Gr\"obner basis. An algorithm written in pseud0-code is
omitted here.\vskip 6pt

{\bf Proof} This can be verified directly, so we leave it as an
exercise.\QED \v5

{\bf The existence of finite left Gr\"obner bases}\par

Finally we show that every nonzero submodule $N$ of the free
$A$-module $L=\oplus_{i=1}^sAe_i$ has a finite left Gr\"obner basis
$\G$ with respect to a given left monomial ordering $\prec_e$.} \v5

Let $a^{\alpha}e_i,a^{\beta}e_j\in\BE$ with $\alpha
=(\alpha_1,\ldots ,\alpha_n)$, $\beta =(\beta_1,\ldots
,\beta_n)\in\NZ^n$. Recall that $a^{\alpha}e_i\prec_e 'a^{\beta}$ if
and only if $i=j$ and  $a^{\alpha}|_{_{\rm L}}a^{\beta}$, while the
latter is defined subject to the property that  $a^{\beta}=\LM
(a^{\gamma}a^{\alpha})$ for some $a^{\gamma}\in\B$. Thus, by Lemma
2.1.2(ii),  $a^{\alpha}e_i\prec_e 'a^{\beta}e_j$ is actually
equivalent to $i=j$ and  $\alpha_i\le \beta_i$ for $1\le i\le n$.
Also note that the correspondence $\alpha =(\alpha_1,\ldots
,\alpha_n)\longleftrightarrow a^{\alpha}$ gives the bijection
between $\NZ^n$ and $\B$. Combining this observation with Dickson's
lemma (Lemma 1.3.1),  we are able to reach the following
result.{\parindent=0pt\v5

{\bf 2.2.14. Theorem}  With notations as in Theorem 2.2.6, let $N\ne
\{ 0\}$ be a submodule of $L$, $\LM (N)=\{ \LM (\xi )~|~\xi\in N\}$,
$\Omega =\{ a^{\alpha}e_i\in \LM (N)~|~a^{\alpha}e_i~\hbox{is
minimal in}~\LM (N)~\hbox{w.r.t.}\prec_e'\}$, and $\G =\{ g\in
N~|~\LM (g)\in\Omega\}$. Then $\Omega$ is a finite subset of $\LM
(N)$, thereby $\G$ is a finite left Gr\"obner basis of $N$.\vskip
6pt

{\bf Proof} For $1\le i\le s$, put $$\begin{array}{l} \BE_i=\{
a^{\alpha}e_i~|~a^{\alpha}\in \B\} ,\\
\Omega_i=\Omega\cap\BE_i.\end{array}$$  Then $\Omega
=\cup_{i=1}^s\Omega_i$. Applying  Dickson's lemma (Lemma 1.3.1) to
$\Omega_i$, it turns out that each $\Omega_i$ has a finite subset
$\Omega_i'=\{ a^{\alpha (1_i)}e_i,\ldots ,a^{\alpha (t_i)}e_i\}$
such that if $a^{\beta}e_i\in\Omega_i$,  then $a^{\alpha
(j_i)}e_i|_{_{\rm L}}a^{\beta}e_i$ for some $a^{\alpha
(j_i)}\in\Omega_i'$. But this implies $a^{\alpha
(j_i)}e_i\prec_e'a^{\beta}e_i$ in $\Omega$. It follows from the
definition of $\Omega$ that $\Omega=\cup_{i=1}^s\Omega_i'$ is a
finite subset of $\LM (N)$, and consequently $\G$ is a finite
(minimal) left Gr\"obner basis of $N$.}\QED\def\S{{\cal S}}

\section*{2.3.  The Noncommutative  Buchberger Algorithm}

Recall from the theory of Gr\"obner bases for  commutative
polynomial algebras ([Bu1], [Bu2], [AL2], [BW], [Fr\"ob], [KR1],
[KR2]) that the celebrated Buchberger algorithm depends on
Buchberger's  criterion which establishes the strategy for computing
Gr\"obner  bases of polynomial ideals.  Let $A=K[a_1,\ldots ,a_n]$
be a solvable polynomial algebra with admissible system $(\B ,\prec
)$, and let $L=\oplus_{i=1}^sAe_i$ be a free $A$-module with left
admissible system $(\BE ,\prec_e)$. Then by Theorem 2.2.14, every
nonzero submodule $N$ has a finite left Gr\"obner basis $\G$. In
this section, by introducing left S-polynomials for elements  $(\xi
,\zeta )\in L\times L$, we show that a noncommutative version of
Buchberger's criterion holds true  for submodules of $L$, and that a
noncommutative version of the Buchberger algorithm works effectively 
for computing finite left Gr\"obner bases of  submodules in $L$.\par

Since the noncommutative version of Buchberger's criterion and the 
noncommutative version of Buchberger's algorithm for modules over a 
solvable polynomial algebra $A$ (Theorem 2.3.3 and {\bf 
Algorithm-LGB} presented below)  look as if working the same way as 
in the commutative case by reducing the S-polynomials, at this stage 
one is asked  to pay more attention to compare the argumentation  
concerning  Buchberger's criterion in the commutative case  (e.g. 
[AL2], P.40-42) and that in the noncommutative case, given in [K-RW] 
and we are going to give respectively, so as to see  how the barrier 
made by the noncommutativity of $A$ (i.e., the product of two 
monomials in $A$ is no longer necessarily a monomial of $A$) can be 
broken down. Also, at this point we remind that  in the proof of 
([K-RW], Theorem 3.11), the argumentation was given in the language 
of abstract rewriting (see [Wik2] for an introduction of this 
topic), namely, it was shown that the left reduction of left 
S-polynomials by $G$ is locally confluent; while in the 
argumentation we are going to give below,   Lemma 2.3.1 will play 
the key role as the breakthrough point, though  our presentation 
looks quite similar to the most popularly known presentation in the 
commutative case  (e.g. [AL2]). \v5

Let $\xi ,\zeta$ be nonzero elements of $L$ with $\LM
(\xi)=a^{\alpha}e_{i}$, $\LM (\zeta )=a^{\beta}e_{j}$, where $\alpha
=(\alpha_{1},\ldots,\alpha_{n})$, $\beta =(\beta_{1},\ldots
,\beta_{n})$. Put $\gamma =(\gamma_1,\ldots ,\gamma_n)$ with
$\gamma_k =\max\{ \alpha_{k},\beta_k)$. The {\it left S-polynomial}
of $\xi$ and $\zeta$ is defined as the element
$$S_{\ell}(\xi,\zeta )=\left\{\begin{array}{ll}
\displaystyle{\frac{1}{\LC (a^{\gamma -\alpha}\xi)}}a^{\gamma
-\alpha}\xi-\displaystyle{\frac{1}{\LC (a^{\gamma
-\beta}\zeta)}}a^{\gamma
-\beta}\zeta,&\hbox{if}~ i=j\\
0,&\hbox{if}~i\ne j.\end{array}\right.$$ {\parindent=0pt\def\S{{\cal
S}}\v5

{\bf Observation} If $S_{\ell}(\xi,\zeta)\ne 0$ then, with respect
to the given $\prec_e$ on $\BE$ we have $\LM (S_{\ell}(\xi
,\zeta))\prec_ea^{\gamma}e_i$.\v5

{\bf 2.3.1. Lemma} Let $\xi ,\zeta$ be as above such that
$S_{\ell}(\xi ,\zeta )\ne 0$, and put $\lambda =\LC (a^{\gamma
-\alpha}\xi)$, $\mu =\LC (a^{\gamma -\beta}\zeta )$. If there exist
$a^{\theta (1)},a^{\theta (2)}\in\B$ such that $\LC (a^{\theta
(1)}\xi )=\lambda_1$, $\LC (a^{\theta (2)}\zeta )=\mu_1$, $\LM
(a^{\theta (1)}\xi )=a^{\rho}e_i=\LM (a^{\theta (2)}\zeta),$ where
$\rho =(\rho_1,\ldots ,\rho_n)$, then  there exists
$a^{\delta}\in\B$ such that
$$S_{\ell}( a^{\theta (1)}\xi ,a^{\theta (2)}\zeta )=b\left (
a^{\delta}S_{\ell}(\xi ,\zeta )-da^{\theta (2)}\zeta -f_1\xi
-f_2\zeta\right ),$$ where $b=\frac{\lambda}{\lambda_{\rho
,\alpha}\lambda_1}$, $d=\frac{\lambda_{\rho
,\alpha}\lambda_1}{\lambda\mu_1}-\frac{\mu_{\rho ,\beta}}{\mu}$,
$\LM (a^{\delta}S_{\ell}(\xi ,\zeta ))\prec_ea^{\rho}e_i,$ $\LM
(f_1\xi )\prec_e a^{\rho}e_i,$ $\hbox{and}~\LM (f_2\zeta )\prec_e
a^{\rho}e_i,$ $\lambda_{\rho ,\alpha},\mu_{\rho ,\beta}\in K^*$,
$f_1,f_2\in A$, which are given in terms of
$$\begin{array}{l} a^{\rho
-\gamma}a^{\gamma -\alpha}=\lambda_{\rho ,\alpha}a^{\theta
(1)}+f_1~\hbox{with}~f_1\in A~\hbox{and}~\LM (f_1)\prec a^{\theta
(1)};\\
a^{\rho -\gamma}a^{\gamma -\beta}=\mu_{\rho ,\beta}a^{\theta
(2)}+f_2~\hbox{with}~f_2\in A~\hbox{and}~\LM (f_2)\prec a^{\theta
(2)}.\end{array}$$ \vskip 6pt

{\bf Proof} By the assumption and Lemma 2.1.2, $\theta (1)+\alpha
=\rho =\theta (2)+\beta$. Since $\gamma =(\gamma_1,\ldots
,\gamma_n)$ with $\gamma_k =\max\{ \alpha_{k},\beta_k)$, we may put
$\delta =(\delta_1,\ldots ,\delta_n)$ with
$\delta_i=\rho_i-\gamma_i$, i.e., $\delta =\rho -\gamma$. Now we see
that
$$\begin{array}{rcl} a^{\delta}S_{\ell}(\xi ,\zeta )&=&\frac{1}{\lambda}a^{\delta}a^{\gamma-\alpha}\xi -
\frac{1}{\mu}a^{\delta}a^{\gamma -\beta}\zeta\\
&=&\frac{1}{\lambda}a^{\rho -\gamma}a^{\gamma -\alpha}\xi
-\frac{1}{\mu}a^{\rho -\gamma}a^{\gamma -\beta}\zeta\\
&=&\left (\frac{\lambda_{\rho ,\alpha}}{\lambda}a^{\theta (1)}\xi
+f_1\xi\right )-\left ( \frac{\mu_{\rho ,\beta}}{\mu}a^{\theta
(2)}\zeta +f_2\zeta\right ),\\
&=&\frac{\lambda_{\rho ,\alpha}\lambda_1}{\lambda}\left
(\frac{1}{\lambda_1}a^{\theta (1)}\xi -\frac{1}{\mu}a^{\theta
(2)}\zeta\right ) \\
&{~}&+\left (\frac{\lambda_{\rho ,\alpha}\lambda_1}{\lambda\mu_1}
-\frac{\mu_{\rho ,\beta}}{\mu}\right )a^{\theta (2)}\zeta +f_1\xi
+f_2\zeta\\
&=&\frac{\lambda_{\rho ,\alpha}\lambda_1}{\lambda}S_{\ell}(a^{\theta
(1)}\xi ,a^{\theta (2)}\zeta )\\
&{~}&+\left (\frac{\lambda_{\rho ,\alpha}\lambda_1}{\lambda\mu_1}
-\frac{\mu_{\rho ,\beta}}{\mu}\right )a^{\theta (2)}\zeta +f_1\xi
+f_2\zeta ,
\end{array}$$ in which $\LM
(a^{\delta}S_{\ell}(\xi ,\zeta ))\prec_ea^{\rho}e_i,~\LM (f_1\xi
)\prec_e a^{\rho}e_i,~\hbox{and}~\LM (f_2\zeta )\prec_e
a^{\rho}e_i$.\QED\v5

{\bf 2.3.2. Lemma } Let $\xi_1,...,\xi_t\in L$ be such that $\LM
(\xi_i)= a^{\rho}e_q$ for all $i=1,...,t$. If, for $c_1,\ldots
,c_t\in K^*$, the element $\xi =\sum^t_{i=1}c_i\xi_i$ satisfies $\LM
(\xi )\prec_e a^{\rho}e_q$, then $\xi$ can be expressed as  a
$K$-linear combination of the form
$$\xi =d_{12}S_{\ell}(\xi_1,\xi_2)+d_{23}S_{\ell}(\xi_2,\xi_3)+\cdots +d_{t-1},t(\xi_{t-1},\xi_t),$$
where $d_{ij}\in K$. \vskip 6pt

{\bf Proof}  If we rewrite each $\xi_i$ as
$\xi_i=\lambda_ia^{\rho}e_q+$ lower terms, where $\lambda_i=\LC
(\xi_i)$, then the assumption yields a cancelation of leading terms
which gives rise to $\sum^s_{i=1}c_i\lambda_i=0$. Since $\LM
(\xi_i)=a^{\rho}e_q=\LM (\xi_j)$, it follows that
$S_{\ell}(\xi_i,\xi_j)=\frac{1}{\lambda_i}\xi_i-\frac{1}{\lambda_j}\xi_j$.
Thus
\begin{eqnarray*}
\xi&=&c_1\xi_1+\cdots +c_s\xi_t\\
&=&c_1\lambda_1\left (\frac{1}{\lambda_1}\xi_1\right ) +\cdots
+c_t\lambda_t\left (\frac{1}{\lambda_t}\xi_t
\right )\\
&=&c_1\lambda_1\left
(\frac{1}{\lambda_1}\xi_1-\frac{1}{\lambda_2}\xi_2\right
)+(c_1\lambda_1+c_2\lambda_2)
\left (\frac{1}{\lambda_2}\xi_2-\frac{1}{\lambda_3}\xi_3\right )\\
&{~}&+\cdots +(c_1\lambda_1+\cdots +c_{t-1}\lambda_{t-1}) \left
(\frac{1}{\lambda_{t-1}}\xi_{t-1}-\frac{1}{\lambda_t}\xi_t\right
)\\
&{~}&+(c_1\lambda_1+\cdots +
c_t\lambda_t)\frac{1}{\lambda_t}\xi_t\\
&=&c_1\lambda_1S_{\ell}(\xi_1,\xi_2)+(c_1\lambda_1+c_2\lambda_2)S_{\ell}(\xi_2,\xi_3)+\cdots\\
&~~&+(c_1\lambda_1+\cdots
+c_{t-1}\lambda_{t-1})S_{\ell}(\xi_{t-1},\xi_t),
\end{eqnarray*}
because $\sum^t_{i=1}c_i\lambda_i=0$. This completes the proof.
\QED\v5}

Considering the submodule $N=\sum_{i=1}^mA\xi_i$ of $L$ generated by
$\Xi =\{\xi_1,\ldots ,\xi_m\}\subset L$, it is clear that
$S_{\ell}(\xi_i,\xi_j)\in N$ for $1\le i<j\le m$. {\parindent=0pt\v5

{\bf 2.3.3. Theorem} (Noncommutative version of Buchberger's
criterion) Let $N=\sum_{i=1}^mA\xi_i$ be a submodule of $L$
generated by the set of nonzero elements $\Xi =\{\xi_1,\ldots
,\xi_s\}$. Then, with respect to the given left monomial ordering
$\prec_e$ on $\BE$, $\Xi$ is a left Gr\"obner basis of $N$ if and
only if every nonzero $S_{\ell}(\xi_i,\xi_j)$ is reduced to 0 (mod
$\Xi$), i.e.,   $\OV{S_{\ell}(\xi_i,\xi_j)}^{\Xi}=0$.\vskip 6pt

{\bf Proof}  Note that $S_{\ell}(\xi_i,\xi_j)\in N$. If $\Xi$ is a
left Gr\"obner basis of $N$, then
$\OV{S_{\ell}(\xi_i,\xi_j)}^{\Xi}=0$ whenever
$S_{\ell}(\xi_i,\xi_j)\ne 0$.}\par

Conversely, suppose that $\OV{S_{\ell}(\xi_i,\xi_j)}^{\Xi}=0$ for
every nonzero $S_{\ell}(\xi_i,\xi_j)$. We will show that $\Xi$
satisfies Proposition 2.2.4(ii), or in other words, that every
nonzero element of $N$ has a left Gr\"obner representation by $\Xi$. 
For this purpose, we argue by contradiction. Suppose that there were 
a $\xi =\sum^s_{i=1}f_i\xi_i\in N$ with nonzero $f_i\in A$, such 
that $\LM (\xi )\prec_e\LM (f_j\xi_j)$ for some $j\le s$. Let $\LM 
(f_j)=a^{\theta (j)}$ and $\LM (\xi_j)=a^{\alpha (j)}e_{i_j}$. 
Comparing the linear expressions of both sides of $\xi 
=\sum^s_{i=1}f_i\xi_i$ in terms of base elements in $\BE$,  we may 
assume, without loss of generality,  that for some $2\le t\le s$,
$$\begin{array}{rcl} a^{\rho}e_q&=&\LM (f_1\xi_1)=\LM (f_2\xi_2)=\cdots =\LM
(f_t\xi_t)\\
&=&\max\{\LM (f_1\xi_1),\ldots ,\LM (f_s\xi_s)\}
,\end{array}\eqno{(1)}$$ and moreover, we assume that the
representation $\xi =\sum^s_{i=1}f_i\xi_i$ is chosen so that
$a^{\rho}e_q$ is the least one. We now proceed to produce a new
representation of $\xi$ by elements of $\Xi$ in which the maximal
monomial is strictly less than $a^{\rho}e_q$. To this end, we write
$f_i=c_ia^{\theta (i)}+$ lower terms, where $c_i=\LC (f_i)$, and let
$\eta_0: =\sum_{i=1}^tc_ia^{\theta (i)}\xi_i$, $\xi^*=\xi -\eta_0$.
Then $\LM (\eta_0)\prec_ea^{\rho}e_q$, $\LM
(\xi^*)\prec_ea^{\rho}e_q$. With $\lambda_i=\LC (a^{\theta
(i)}\xi_i)$, $1\le i\le t$, we may apply Lemma 2.3.2 to $\eta_0$ so
that
$$\begin{array}{rcl} \eta_0 &=&d_{12}S_{\ell}(a^{\theta (1)}\xi_1, a^{\theta (2)}\xi_2)+
d_{23}S_{\ell}(a^{\theta (2)}\xi_2, a^{\theta (3)}\xi_3)+\cdots\\
&{~}& + d_{t_{t-1},t}S_{\ell}(a^{\theta (t-1)}\xi_{t-1}, a^{\theta
(t)}\xi_t),\end{array}\eqno{(2)}$$ where $d_{ij}\in K$ and
$S_{\ell}(a^{\theta (i)}\xi_i, a^{\theta
(j)}\xi_j)=\frac{1}{\lambda_i}a^{\theta
(i)}\xi_i-\frac{1}{\lambda_j}a^{\theta (j)}\xi_j$. We next apply
Lemma 2.3.1 to each $S_{\ell}(a^{\theta (i)}\xi_i, a^{\theta
(j)}\xi_j)$ so that
$$S_{\ell}( a^{\theta (i)}\xi_i ,a^{\theta (j)}\xi_j )=b_{ij}\left (
a^{\delta (ij)}S_{\ell}(\xi_i ,\xi_j )-z_{ij}a^{\theta (j)}\xi_j
-h_i\xi_i -h_j\xi_j\right ),\eqno{(3)}$$ where $b_{ij}$, $z_{ij}\in
K$, $h_i,~h_j\in A$, $\LM (a^{\delta (ij)}S_{\ell}(\xi_i
,\xi_j))\prec_ea^{\rho}e_q,$ $\LM (h_i\xi_i )\prec_e a^{\rho}e_q,$
$\hbox{and}~\LM (h_j\xi_j )\prec_e a^{\rho}e_q,$ thereby $$\LM
(S_{\ell}( a^{\theta (i)}\xi_i ,a^{\theta (j)}\xi_j ))=\LM
(a^{\theta (j)}\xi_j)=a^{\rho}e_q~\hbox{ if}~b_{ij},~z_{ij}\ne
0.\eqno{(4)}$$ Note that $j=i+1>i$ in (3). After substituting (3)
into (2), and $\eta_0$ into $\xi$, we see that $a^{\theta (1)}\xi_1$
does not appear in the new representation of $\xi$. So, if we take
the sub-sum  $\eta_1:=\sum_{j=2}^td_{ij}b_{ij}z_{ij}a^{\theta
(j)}\xi_j$ of $\xi$, then $\eta_1$ clearly satisfies $\LM
(\eta_1)\prec_ea^{\rho}e_q=\LM (a^{\theta (j)}\xi_j)$, $2\le j\le
t$. Hence we may apply Lemma 2.3.1 and Lemma 2.3.2 to $\eta_1$ so
that we obtain a similar result of (3) which gives rise to a sub-sum
$\eta_2:=\sum_{k=3}^tu_ka^{\theta (k)}\xi_k$ of $\xi$ with $u_k\in
K$ and $\LM (\eta_2)\prec_ea^{\rho}e_q=\LM (a^{\theta (k)}\xi_k)$,
$3\le k\le t$. Repeating such a procedure for at most $t-1$ times
and at each time, applying the assumption
$\OV{S_{\ell}(\xi_i,\xi_j)}^{\Xi}=0$ to (3), we eventually obtain a 
representation of $\xi$ by elements of $\Xi$ which yields the 
desired contradiction and finishes the proof.\QED{\parindent=0pt\v5

{\bf 2.3.4. Theorem} (The Noncommutative version of Buchberger 
algorithm) Let $N=\sum_{i=1}^mA\xi_i$ be a submodule of $L$ 
generated by a finite set of nonzero elements $\Xi =\{\xi_1,\ldots 
,\xi_s\}$. Then, with respect to a given left monomial ordering 
$\prec_e$ on $\BE$, the  algorithm presented below returns a finite 
left Gr\"obner basis $\G$ for $N$.{\parindent=0pt\vskip 6pt

\underline{\bf Algorithm-LGB
~~~~~~~~~~~~~~~~~~~~~~~~~~~~~~~~~~~~~~~~~~~~~~~~~~~~~~~~~~~~~~~}\vskip 
6pt

\textsc{INPUT}: $\Xi = \{ \xi_1,...,\xi_m\}$\par                         
\textsc{OUTPUT}:~$\G =\{ g_1,...,g_t\},~\hbox{a left Gr\"obner basis
for}~N=\sum_{i=1}^mA\xi_i$\par          
\textsc{INITIALIZATION}:~$m':=m,~\G :=\{g_1=\xi_1,\ldots 
,g_{m'}=\xi_m\} ,$\par                                                           
\hskip 3.75truecm $\S :=\left\{S_{\ell}(g_i,g_j)~\left 
|~g_i,g_j\in\G ,~i<j\right.\right\}-\{ 0\}$\par 

\textsc{BEGIN}\par ~~~~\textsc{WHILE}~$\S\ne 
\emptyset~\textsc{DO}$\par                      
~~~~~~~~~\hbox{Choose any}~$S_{\ell}(g_i,g_j)~\hbox{from}~\S$\par         
~~~~~$\S :=\S -\{ S_{\ell}(g_i,g_j)\}$\par                                    
~~~~~~~~~\textsc{IF}~$\OV{S_{\ell}(g_i,g_j)}^{\G}=\eta\ne 
0~\hbox{with}~\LM (\eta )= a^{\rho}e_k~\textsc{THEN}$\par      
~~~~~~~~~$m':=m'+1,~g_{m'}:=\eta $\par                                                            
~~~~~~~~~$\S :=\S\cup\{ S_{\ell}(g_j,g_{m'})~|~g_j\in\G,~\LM (g_j 
)=a^{\nu}e_k\} -\{ 0\}$\par                                         
~~~~~~~~~$\G :=\G\cup\{g_{m'}\},$\par         
~~~~~~~~~\textsc{END}\par                                      
~~~~\textsc{END}\par                                   
\textsc{END}\par \vskip -.2truecm 
\underline{~~~~~~~~~~~~~~~~~~~~~~~~~~~~~~~~~~~~~~~~~~~~~~~~~~~~~~~~~~~~~~~~~~~~~~~~~~~~~~~ 
~~~~~~~~~~~~~~~} \vskip 6pt

{\bf Proof} We first prove that the algorithm terminates after a 
finite number of executing the WHILE loop. To this end, let  
$\G_{n+1}$ denote the new set obtained after the $n$-th turn of 
executing the \textsc{WHILE} loop, that is,
$$\G_{n+1}=\G_n\cup\left\{\eta =\OV{S_{\ell}(g_i,g_j)}^{\G_n}\ne 0~\hbox{for some pair}~g_i,g_j\in\G_n\right\} ,
\quad n\in\NZ ,$$ where $\G_0=\Xi$ , and let
$$\G=\bigcup_{n\in\NZ}\G_n,\quad\LM (\G )=\{\LM (g)~|~g\in\G\} .$$
By Dickson's lemma (Lemma 1.3.1), $\LM (\G )$ has a finite subset
$U=\{\LM (g_1),\ldots ,\LM (g_s)\}$, such that if $g\in \G$ then
$\LM (g_i)|\LM (g)$ for some $\LM (g_i)\in U$. Since $U$ is finite,
we may assume that $U\subset\LM (\G_k)$ for some $k$. This shows
that $\OV{S_{\ell}(g_i,g_j)}^{\G_k}=0$ for all $g_i,g_j\in\G_k$,
thereby the algorithm terminates after the $k$-th turn of executing
the \textsc{WHILE} loop.}}\par

Now that the algorithm terminates in a finite number of executions,
we may assume that $\G =\{ g_1,\ldots ,g_t\}$. Then, since
$\G_0=\Xi\subset\G$, the algorithm itself tells us that $\G$
generates the submodule $N$ and $\OV{S_{\ell}(g_i,g_j)}^{\G}=0$ for
all $g_i,g_j\in\G$. It follows from Theorem 2.3.3 that $\G$ is a
left Gr\"obner basis for $N$.\QED\v5

One is referred to the up-to-date computer algebra systems
\textsc{Singular} [DGPS] for the implementation of {\bf
Algorithm-LGB}. Also, nowadays there have been optimized algorithms,
such as the signature-based algorithm for computing Gr\"obner bases
in solvable polynomial algebras [SWMZ], which is based on the
celebrated F5 algorithm [Fau],  may be used to speed-up the
computation of left Gr\"obner bases for modules. {\parindent=0pt\v5

{\bf 2.3.5. Corollary} Let $N=\sum_{i=1}^mA\xi_i$ be a submodule of
$L$ generated by a finite set of nonzero elements $\Xi
=\{\xi_1,\ldots ,\xi_m\}$, and let $\G=\{ g_1,...,g_t\}$ be a left
Gr\"obner basis of $N$ produced by running {\bf Algorithm-LGB}. Then
the following statements hold true.\par

(i) All the properties listed in Proposition 2.2.7 can be recognized
in a computational way.\par

(ii) There is a $t\times m$ matrix
$$V_{t\times m}=\left (\begin{array}{ccc} h_{11}&\cdots&h_{1m}\\
\cdots&\cdots&\cdots \\
h_{t1}&\cdots&h_{tm}\end{array}\right ) ,\quad h_{ij}\in A$$ which
is obtained after $\G$ is computed by running {\bf Algorithm-LGB},
such that
$$\G =\left (\begin{array}{c} g_1\\ \vdots\\ g_t\end{array}\right )
=V_{t\times m}\left (\begin{array}{c} \xi_1\\ \vdots\\
\xi_m\end{array}\right ), ~ 1\le j\le s.$$ (Note that for
convenience, we formally used the matrix expression to demonstrate
an obvious meaning.)\par

(iii) If $\xi\in N$ and $\xi\ne 0$, then a representation of $\xi$
by $\Xi=\{\xi_1,\ldots ,\xi_m\}$, say $\xi =\sum^m_{i=1}f_i\xi_i$,
can be computed.\vskip 6pt

{\bf Proof} (i) Now that a left Gr\"obner basis $\G$ of $N$ has been
computed, all the properties listed in Proposition 2.2.7 can be
recognized by means of the division by $\G$.\par

(ii) Recall from {\bf Algorithm-LGB} that for each $g_j\in\G$,
either $g_j\in\Xi = \{ \xi_1,...,\xi_m\}$ or $g_j$ is a newly added 
member of $\G$ obtained after a certain pass through the 
\textsc{WHILE} loop.  So, starting with the expressions obtained 
from the first pass through the \textsc{WHILE} loop:
$$S_{\ell}(\xi_i,\xi_j)=\sum_kh_k\xi_k+\OV{S_{\ell}(\xi_i,\xi_j)}^{\Xi},~h_k\in A,~1\le i<j\le m,$$
and  keeping track of the linear combinations that give rise to the 
new elements of $\G$, the algorithm eventually yields the desired 
matrix $V_{t\times m}$.\par

(iii) By Proposition 2.2.7(iii), after dividing by $\G$ the element 
$\xi$ has a representation $\xi =\sum^t_{j=1}s_jg_j$ with $s_j\in A$ 
and $g_j\in \G$. Now  the conclusion (ii) yields the desired 
representation of $\xi$ by $\Xi$.}\newpage\setcounter{page}{51}

\chapter*{3. Finite Free Resolutions}
\markboth{\rm Finite Free Resolutions}{\rm Finite Free Resolutions}
\vskip 2.5truecm

Since every solvable polynomial algebra $A$ is (left and right) 
Noetherian (Section 3 of Chapter 1), and every nonzero submodule of 
a free left $A$-module has a finite left Gr\"obner basis which can 
be produced by running {\bf Algorithm-LGB} (Section 3 of Chapter 2), 
starting from this chapter we shall successively present some 
details concerning applications of Gr\"obner bases in constructing 
finite free resolutions over an arbitrary solvable polynomial 
algebra $A$, minimal finite $\NZ$-graded free resolutions over an 
$\NZ$-graded solvable polynomial algebra $A$ with the  degree-0 
homogenous part  being the ground field $K$, and minimal finite 
$\NZ$-filtered  free resolutions  over an $\NZ$-filtered solvable 
polynomial algebra $A$ (where the $\NZ$-filtration of $A$ is 
determined by a positive-degree function). 
\par

Let $A$ be a solvable polynomial algebra as before. The current 
chapter consists of four sections.  In Section 1 we demonstrate, for 
a finitely generated submodule $N$ of a free left $A$-module $L$, 
how to compute a generating set of the syzygy module Syz$(N)$ of $N$ 
via computing a left Gr\"obner basis $\G$ of $N$.  In Section 2 we 
show, in a constructive way,  that a noncommutative version of 
Hilbert's syzygy theorem holds true for $A$, and consequently, that 
a finite free resolution  can be algorithmically constructed for 
every finitely generated $A$-module. In Section 3,  the 
noncommutative version of Hilbert's syzygy theorem is applied to 
highlight two homological properties of $A$, that is, $A$ has finite 
global homological dimension, and every finitely generated 
projective $A$-module is stably free.  Based on Section 1 -- Section 
3, the final Section 4 is devoted to the calculation of projective 
dimension of a finitely generated $A$-module $M$, and meanwhile, the 
proposed algorithmic procedure verifies whether $M$ is a projective 
module or not.\par

The main references of this chapter are [AL2], [AF], [Rot],  [Li1],
[GV], [Lev], [DGPS].\par

Throughout this chapter, modules are meant left modules over
solvable polynomial algebras, and all notions and notations used in
previous chapters are maintained. \par

\section*{3.1. Computation of Syzygies}

Let $A=K[a_1,\ldots ,a_n]$ be a solvable polynomial algebra with
admissible system $(\B ,\prec )$ in the sense of Definition 1.1.3,
where $\B =\{ a^{\alpha}=a_1^{\alpha_1}\cdots
a_n^{\alpha_n}~|~\alpha =(\alpha_1,\ldots ,\alpha_n)\in\NZ^n\}$ is
the PBW $K$-basis of $A$ and $\prec$ is a monomial ordering on $\B$.
Let $L_0=\oplus_{i=1}^sAe_i$ be a free left $A$-module with
$A$-basis $\{ e_1,\ldots ,e_s\}$, and $\prec_{e}$ a  left monomial
ordering on the $K$-basis $\BE=\{ a^{\alpha}e_i~|~a^{\alpha}\in\B
,~1\le i\le s\}$ of $L_0$. As in  (Section 3 of Chapter 2) we write
$S_{\ell}(\xi,\zeta)$ for the left S-polynomial of two elements
$\xi$, $\zeta\in L_0$, that is, if  $\LM (\xi)=a^{\alpha}e_{i}$ with
$\alpha =(\alpha_{1},\ldots,\alpha_{n})$, $\LM (\zeta
)=a^{\beta}e_{j}$ with   $\beta =(\beta_{1},\ldots ,\beta_{n})$,
then
$$S_{\ell}(\xi,\zeta )=\left\{\begin{array}{ll}
\displaystyle{\frac{1}{\LC (a^{\gamma -\alpha}\xi)}}a^{\gamma
-\alpha}\xi-\displaystyle{\frac{1}{\LC (a^{\gamma
-\beta}\zeta)}}a^{\gamma
-\beta}\zeta,&\hbox{if}~ i=j\\
0,&\hbox{if}~i\ne j\end{array}\right.$$ where $\gamma
=(\gamma_1,\ldots ,\gamma_n)$ with $\gamma_k =\max\{
\alpha_{k},\beta_k)$; moreover, if $S_{\ell}(\xi,\zeta)\ne 0$, then
with respect to the given $\prec_e$ on $\BE$ we have $\LM
(S_{\ell}(\xi ,\zeta))\prec_ea^{\gamma}e_i$.\par

Let $N=\sum^m_{i=1}A\xi_i$ be a  submodule of $L_0$ generated by the
set of nonzero elements $ U =\{\xi_1,\ldots ,\xi_m\}$, and let
$L_1=\oplus_{i=1}^mA\omega_i$ be the free $A$-module with $A$-basis
$\{\omega_1,\ldots ,\omega_m\}$. Then the syzygy module of $U$ (or
equivalently the syzygy module of $N$), denoted Syz$(U)$, is the
submodule of $L_1$ defined by
$$\hbox{Syz}(U)=\left\{\sum^m_{i=1}h_i\omega_i\in L_1~\left |~\sum^m_{i=1}h_i\xi_i=0\right.\right\}.$$
Our aim of this section is to show that a generating set of the
syzygy module Syz$( U )$ can be computed by means of a left
Gr\"obner basis of $N$. To this end, let $\G=\{ g_1,\ldots ,g_t\}$
be a left Gr\"obner basis of $N$ with respect to $\prec_{e}$, then
every nonzero left S-polynomial $S_{\ell}(g_i,g_j)$ has a left
Gr\"obner representation $S_{\ell}(g_i,g_j)=\sum_{i=1}^tf_ig_i$ with
$\LM(f_ig_i)\preceq_e\LM (S_{\ell}(g_i,g_j))$ whenever $f_i\ne 0$
(note that such a representation is obtained by using the division
by $\G$ during executing the WHILE loop of {\bf Algorithm-LGB}
presented in  (Chapter 2, Theorem 2.3.4)). Considering the syzygy
module Syz$(\G )$ of $\G$ in the free $A$-module
$L_2=\oplus_{i=1}^tA\varepsilon_i$ with $A$-basis
$\{\varepsilon_1,\ldots ,\varepsilon_t\}$, if we put
$$\begin{array}{rcl} s_{ij}&=&f_1\varepsilon_1+\cdots +\left (f_i-\frac{a^{\gamma-\alpha (i)}}{\LC (a^{\gamma -\alpha (i)}\xi_i)}\right )\varepsilon_i+
\cdots \\
\\
&{~}&~~~~~~~~~~~~~+\left (f_j+\frac{a^{\gamma -\alpha (j)}}{\LC
(a^{\gamma -\alpha (j)}\xi_j)}\right )\varepsilon_j+\cdots
+f_t\varepsilon_t,\\
\\
{\cal S} &=&\{ s_{ij}~|~1\le i<j\le t\},\end{array}$$  then it can
be shown, actually as in the commutative case (cf. [AL2], Theorem
3.7.3), that ${\cal S}$ generates Szy$(\G )$ in $L_2$. However, by
employing the Schreyer ordering $\prec_{s\hbox{-}\varepsilon}$ on
the $K$-basis $\B (\varepsilon )=\{
a^{\alpha}\varepsilon_i~|~a^{\alpha}\in\B ,~1\le i\le m\}$ of $L_2$
induced by $\G$ with respect to $\prec_{e}$ (see Section 1 of
Chapter 2), which is defined subject to the rule: for
$a^{\alpha}\varepsilon_i,a^{\beta}\varepsilon_j\in\B (\varepsilon
)$,
$$a^{\alpha}\varepsilon_i\prec_{s\hbox{-}\varepsilon} a^{\beta}\varepsilon_j
\Leftrightarrow\left\{\begin{array}{l}
\LM (a^{\alpha}g_i)\prec_{e}\LM (a^{\beta}g_j),\\
\hbox{or}\\
\LM (a^{\alpha}g_i)=\LM (a^{\beta}g_j)~\hbox{and}~i<j,
\end{array}\right.$$
there is indeed a much stronger result, namely the noncommutative
analogue of Schreyer theorem [Sch] (cf. Theorem 3.7.13 in [AL2] for 
free modules over commutative polynomial algebras; Theorem 4.8 in 
[Lev] for free modules over solvable polynomial 
algebras):{\parindent=0pt \v5

{\bf 3.1.1. Theorem} (2.3.1) With respect to the left monomial
ordering $\prec_{s\hbox{-}\varepsilon}$ on $\B (\varepsilon )$ as
defined above, the following statements hold.\par

(i) Let $s_{ij}$ be determined by $S_{\ell}(g_i,g_j)$, where $i<j$,
$\LM (g_i)=a^{\alpha (i)}e_s$ with $\alpha (i)=(\alpha_{i_1},\ldots
,\alpha_{i_n})$, and $\LM (g_j)=a^{\alpha (j)}e_s$ with $\alpha
(j)=(\alpha_{j_1},\ldots ,\alpha_{j_n})$. Then $\LM
(s_{ij})=a^{\gamma-\alpha (j)}\varepsilon_j$, where $\gamma
=(\gamma_1,\ldots ,\gamma_n)$ with each
$\gamma_k=\max\{\alpha_{i_k},\alpha_{j_k}\}$. \par

(ii) ${\cal S}$ is a left Gr\"obner basis of Syz$(\G )$, thereby
${\cal S}$ generates Syz$(\G )$.\vskip 6pt

{\bf Proof} (i) By the definition of $S_{\ell}(g_i,g_j)$ we know
that $\LM (a^{\gamma -\alpha (i)}g_i)=a^{\gamma}e_s=\LM (a^{\gamma
-\alpha (j)}g_j)$. Since $i<j$, by the definition of
$\prec_{s\hbox{-}\varepsilon}$ we have $a^{\gamma -\alpha
(i)}\varepsilon_i\prec_{s\hbox{-}\varepsilon}a^{\gamma -\alpha
(j)}\varepsilon_j$. Consequently, it follows from  $\LM
(S_{\ell}(g_i,g_j))\prec_{e}\LM (a^{\gamma -\alpha (j)}g_n)$, the 
Gr\"obner representation $S_{\ell}(g_i,g_j)=\sum_{i=1}^tf_ig_i$, and 
the definition of $\prec_{s\hbox{-}\varepsilon}$ that $\LM 
(s_{ij})=a^{\gamma -\alpha (j)}\varepsilon_j$.}\par

(ii) To show that $\S$ forms a left Gr\"obner basis for Syz$(\G )$
with respect to $\prec_{s\hbox{-}\varepsilon}$, take any nonzero
$\xi \in $Syz$(\G )$. After executing the division algorithm by
$\S$, $\xi$ has an expression $\xi =\sum_{s_{ij}\in
S}f_{ij}s_{ij}+\eta$, where
$\LM(f_{ij}s_{ij})\preceq_{s\hbox{-}\varepsilon}\LM (\xi )$ for each
term $f_{ij}s_{ij}$, and $\eta=\OV\xi^{\S}$ is the remainder of
$\xi$ on division by $\S$. We claim that $\eta =0$, thereby $\S$ is
a left Gr\"obner basis for Syz$(\G )$. Otherwise, assuming the
contrary that $\eta =\sum_{i,j}\lambda_{ij}a^{\theta
(i_j)}\varepsilon_j\ne 0$, where $\lambda_{ij}\in K^*$ and $\theta
(i_j)=(\theta_{i_{j1}},\ldots ,\theta_{i_{jn}})\in\NZ^n$, then $\LM
(\eta )=a^{\theta (l_k)}\varepsilon_k$  for some pair $(l,k)$ and
$$\LM (s_{ij})\not |~a^{\theta (l_k)}\varepsilon_k~\hbox{for
all}~s_{ij}\in S.\eqno{(1)}$$ Let $s_{ij}$ be determined by
$S_{\ell}(g_i,g_j)$, where $i<j$, $\LM (g_i)=a^{\alpha (i)}e_s$ and
$\LM (g_j)=a^{\alpha (j)}e_s$. Then by (i) we have $\LM
(s_{ij})=a^{\gamma -\alpha (j)}\varepsilon_j$. If $j=k$ then by (1),
$$a^{\gamma -\alpha (k)}{\not |}~a^{\theta (l_k)}.\eqno{(2)}$$
Since $\LM (\eta )=a^{\theta (l_k)}\varepsilon_k$, by the definition
of $\prec_{s\hbox{-}\varepsilon}$ we have $\LM (a^{\theta
(i_j)}g_j)\preceq_{e}\LM (a^{\theta (l_k)}g_k),$ and for $(i,j)\ne 
(l,k)$,
$$\LM(a^{\theta (i_j)}g_j)=\LM (a^{\theta (l_k)}g_k)~\hbox{implies}~j<k.\eqno{(3)}$$
Noticing $\eta =\xi -\sum_{s_{ij}\in \S}f_{ij}s_{ij}\in$ Syz$(\G )$,
we have $\sum_{i,j}\lambda_{ij}a^{\theta (i_j)}g_j=0$. Since $A$ is
a domain, if $\LM (a^{\theta (i_j)}g_j)\ne \LM (a^{\theta
(l_k)}g_k)$ for all $(i,j)\ne (l,k)$, then we would have $\eta =0$
(note that all $g_i\ne 0$), which is a contradiction. So, by $(3)$
we may assume  that $\LM (a^{\theta (i_j)}g_j)=\LM (a^{\theta
(l_k)}g_k)$ for some $j< k$. Let $\LM (g_j)=a^{\alpha (j)}e_s$ and
$\LM (g_k)=a^{\alpha (k)}e_s$. Then
$$\begin{array}{l} \LM (a^{\theta (i_j)}g_j)=\LM (a^{\theta (i_j)}\LM (g_j))=\LM (a^{\theta (i_j)}a^{\alpha (j)}e_s)=a^{\theta (i_j)+\alpha (j)}e_s,\\
\LM (a^{\theta (l_k)}g_k)=\LM (a^{\theta (l_k)}\LM (g_k))=\LM
(a^{\theta (l_k)}a^{\alpha (k)}e_s)=a^{\theta (l_k)+\alpha
(k)}e_s,\end{array}$$ thereby $\theta (i_j)+\alpha (j)=\theta
(l_k)+\alpha (k)$. Now, taking $\gamma =(\gamma_1,\ldots ,\gamma_n)$
in which $\gamma_{\ell}=\max\{\alpha_{j_{\ell}},\alpha_{k_{\ell}}\}$
with respect to $\alpha (j)=(\alpha_{j_1},\ldots ,\alpha_{j_n})$ and
$\alpha (k)=(\alpha_{k_1},\ldots ,\alpha_{k_n})$, it follows that
$\theta (l_k)+\alpha (k)=\rho +\gamma$ for some $\rho=(\rho_1,\ldots
,\rho_n)$.  Hence $\theta (l_k)=\rho +(\gamma -\alpha (k))$, and
this gives rise to $a^{\theta (l_k)}=\LM (a^{\rho}a^{\gamma -\alpha
(k)})$, i.e., $a^{\gamma -\alpha (k)}|a^{\theta (l_k)}$,
contradicting (2). Therefore, we must have $\eta =0$, as claimed.
This complets the proof. \QED\v5

To go further, again let $\G=\{ g_1,\ldots ,g_t\}$ be  the left
Gr\"obner basis of $N$ produced by running {\bf Algorithm-LGB}
presented in  (Chapter 2, Theorem 2.3.4) with the initial input data
$ U=\{\xi_1,\ldots ,\xi_m\}$. Using the usual matrix notation for
convenience, we have
$$\left (\begin{array}{l} \xi_1\\ \vdots\\ \xi_m\end{array}\right )=U_{m\times t}\left (\begin{array}{l} g_1
\\ \vdots\\ g_t\end{array}\right ) ,\quad
\left (\begin{array}{l} g_1\\ \vdots\\ g_t\end{array}\right
)=V_{t\times m}\left (\begin{array}{l} \xi_1
\\ \vdots\\ \xi_m\end{array}\right ) ,$$
where the $m\times t$ matrix $U_{m\times t}$ (with entries in $A$)
is obtained by the division by $\G$, and the $t\times m$ matrix
$V_{t\times m}$ (with entries in $A$) is obtained by keeping track
of the reductions during executing the WHILE loop of {\bf
Algorithm-LGB} (Chapter 2, Corollary 2.3.5).  By Theorem 3.1.1, we
may write Syz$(\G )=\sum^r_{i=1}A\mathcal {S}_i$ with
$\mathcal{S}_1,\ldots ,\mathcal{S}_r\in
L_2=\oplus_{i=1}^tA\varepsilon_i$;  and if
$\mathcal{S}_i=\sum^t_{j=1}f_{ij}\varepsilon_j$, then we  write
$\mathcal{S}_i$ as a $1\times t$ row matrix, i.e.,
$\mathcal{S}_i=(f_{i1}~\ldots~ f_{it})$,  whenever matrix notation 
is convenient in the according discussion. At this point, we note 
also that all the $\mathcal{S}_i$ may be written down one by one 
during executing the WHILE loop of {\bf Algorithm-LGB} successively. 
Furthermore, we write $D_{(1)},\ldots ,D_{(m)}$ for the rows of the 
matrix $D_{m\times m}=U_{m\times t}V_{t\times m}-E_{m\times m}$ 
where $E_{m\times m}$ is the $m\times m$ identity matrix. The 
following proposition is a noncommutative analogue of ([AL2],  
Theorem 3.7.6). {\parindent=0pt \v5

{\bf 3.1.2. Theorem}  With notation fixed above, the syzygy module
Syz$( U )$ of $ U =\{ \xi_1,\ldots ,\xi_m\}$ is generated by
$$\{ \mathcal{S}_1V_{t\times m},\ldots ,\mathcal{S}_{r}V_{t\times m},D_{(1)},\ldots ,D_{(m)}\} ,$$
where each $1\times m$ row matrix represents an element of the free
$A$-module $L_1=\oplus_{i=1}^mA\omega_i$. \vskip 6pt

{\bf Proof} Since $$0=\mathcal{S}_i\left (\begin{array}{l} g_1\\
\vdots\\ g_t\end{array}\right )=(f_{i_1}~\ldots ~f_{it})\left
(\begin{array}{l} g_1\\ \vdots\\ g_t\end{array}\right
)=(f_{i_1}~\ldots ~f_{it})V_{t\times m}\left (\begin{array}{l}
\xi_1\\ \vdots\\ \xi_m\end{array}\right ) ,$$ we have
$\mathcal{S}_iV_{t\times m}\in$ Syz$( U )$, $1\le i\le r$. Moreover,
since
$$\begin{array}{rcl} D_{m\times m}\left (\begin{array}{l}
\xi_1\\ \vdots\\ \xi_m\end{array}\right )&=&(U_{m\times t}V_{t\times
m}-E_{m\times m})\left (\begin{array}{l} \xi_1\\ \vdots\\
\xi_m\end{array}\right )\\
\\
&=&U_{m\times t}V_{t\times m}\left (\begin{array}{l} \xi_1\\
\vdots\\ \xi_m\end{array}\right )-\left (\begin{array}{l} \xi_1\\
\vdots\\ \xi_m\end{array}\right )\\
\\
&=&U_{m\times t}\left (\begin{array}{l} g_1\\ \vdots\\
g_t\end{array}\right )-\left (\begin{array}{l} \xi_1\\ \vdots\\
\xi_m\end{array}\right )\\
\\
&=&\left (\begin{array}{l} \xi_1\\ \vdots\\
\xi_m\end{array}\right )-\left (\begin{array}{l} \xi_1\\ \vdots\\
\xi_m\end{array}\right )=0,\end{array}$$ we have $D_{(1)},\ldots
,D_{(r)}\in$ Syz$( U )$.}\par

On the other hand, if $H =(h_1~\ldots~h_m)$ represents the element
$\sum^m_{i=1}h_i\omega_i\in \oplus_{i=1}^mA\omega_i$ such that
$H\left (\begin{array}{l} \xi_1\\ \vdots\\ \xi_m\end{array}\right
)=0$, then $0=H U_{m\times t}\left (\begin{array}{l} g_1\\ \vdots\\
g_t\end{array}\right )$. This means $H U_{m\times t}\in$ Syz$(\G )$.
Hence, $H U_{m\times t}=\sum^r_{i=1}f_i\mathcal{S}_i$ with $f_i\in
A$, and it follows that $H U_{m\times t}V_{t\times
m}=\sum^r_{i=1}f_i\mathcal{S}_iV_{t\times m}$. Therefore,
$$\begin{array}{rcl} H&=&H +H U_{m\times t}V_{t\times m}-H U_{m\times t}V_{t\times m}\\
&=&H (E_m-U_{m\times t}V_{t\times
m})+\sum^r_{i=1}f_i\mathcal{S}_iV_{t\times m}\\
&=&-H D_{m\times m}+\sum^r_{i=1}f_i(\mathcal{S}_iV_{t\times
m}).\end{array}$$ This shows that every element of Syz$( U )$ is
generated by $\{ \mathcal{S}_1V_{t\times m},\ldots ,$
$\mathcal{S}_{r}V_{t\times m},D_{(1)},\ldots ,D_{(m)}\} ,$ as
desired. \v5

\section*{3.2. Computation of Finite Free Resolutions}\par

Let $A=K[a_1,\ldots ,a_n]$ be a solvable polynomial algebra with
admissible system $(\B ,\prec )$ in the sense of Definition 1.1.3,
where $\B =\{ a^{\alpha}=a_1^{\alpha_1}\cdots
a_n^{\alpha_n}~|~\alpha =(\alpha_1,\ldots ,\alpha_n)\in\NZ^n\}$ is
the PBW $K$-basis of $A$ and $\prec$ is a monomial ordering on $\B$.
Let $M$ be a  left $A$-module. Then a {\it  free resolution} of $M$
is an {\it exact sequence} by free $A$-modules $L_i$ and  $A$-module
homomorphisms $\varphi_i$:
$${\cal L}_{\bullet}\quad\quad\cdots~\mapright{\varphi_{i+1}}{}~L_i~\mapright{\varphi_i}{}~\cdots ~
\mapright{\varphi_2}{}~L_1~\mapright{\varphi_1}{}~L_0~\mapright{\varphi_0}{}~M~\mapright{}{}~0$$ 
that is, in the sequence ${\cal L}_{\bullet}$, $\varphi_0$ is 
surjective,  each $L_i$ is a free $A$-module and Ker$\varphi_i=$ 
Im$\varphi_{i+1}=\varphi_{i+1}(L_{i+1})$ for all $i\ge 0$. By 
classical homological algebra (e.g. [Rot]), theoretically such a 
free resolution for $M$ exists. If $M$ is a finitely generated 
$A$-module, each $L_i$ in ${\cal L}_{\bullet}$ is finitely 
generated, and $L_{q+1}=0$ for some $q$, then ${\cal L}_{\bullet}$ 
is called a {\it finite free resolution} of $M$. \par

In this section, we  show, in a constructive way, that  a
noncommutative version of Hilbert's syzygy theorem holds true for
$A$, i.e., every finitely generated $A$-module $M$ has a finite free
resolution, and consequently, an algorithm for computing  a finite
free resolution of $M$ is obtained. Our argumentation concerning a
noncommutative version of Hilbert's syzygy theorem  below is adapted
from   ([Eis], Corollary 15.11) and ([Lev], Section 4.4).\v5

Let $L=\oplus^s_{i=1}Ae_i$ be a free $A$-module with left monomial
ordering $\prec_e$ on its $K$-basis $\B (e)$, and let $\G=\{
g_1,\ldots ,g_t\}$ be a left Gr\"obner basis of the submodule
$N=\sum^t_{i=1}Ag_i\subset L$. Then, after relabeling the members of
$\G$ (if necessary), we may always assume that $\G$
satisfies{\parindent=.7truecm\vskip 6pt

\item{($*$)} if $i<j$ and $\LM (g_i)=a^{\alpha (i)}e_k$ with $\alpha
(i)=(\alpha_{i_1},\ldots ,\alpha_{i_n})$, $\LM (g_j)=a^{\alpha
(j)}e_k$ with $\alpha (j)=(\alpha_{j_1},\ldots ,\alpha_{j_n})$, then
$a^{\alpha (i)}\preceq_{lex}a^{\alpha (j)}$ under the lexicographic
ordering with respect to
$a_n\prec_{lex}a_{n-1}\prec_{lex}\cdots\prec_{lex}a_1$,\par}
{\parindent=0pt \v5

{\bf 3.2.1. Lemma} Given a free $A$-module $L=\oplus_{i=1}^sAe_i$
with left monomial ordering $\prec_{e}$ on its $K$-basis $\B (e)$,
let $\G=\{ g_1,\ldots ,g_t\}$ be a left Gr\"obner basis of the
submodule $N=\sum^t_{i=1}Ag_i\subset L$. With notation as in Theorem
3.1.1, let $\prec_{s\hbox{-}\varepsilon}$ be the Schreyer ordering
on  the $K$-basis $\B (\varepsilon )$ of the free $A$-module
$L_1=\oplus_{i=1}^tA\varepsilon_i$ induced by $\G$ with respect to
$\prec_{e}$. Assume that for some $r\le n$, the generators
$a_1,\ldots ,a_r$ of $A$ do not appear in every $\LM (g_{\ell})$,
then the generators $a_1,\ldots ,a_r,a_{r+1}$ of $A$ do not appear
in $\LM (s_{ij})$ for every $s_{ij}\in {\cal S}$.\vskip 6pt

{\bf Proof} If $s_{ij}\in {\cal S}$ with $s_{ij}\ne 0$, then we have
$i<j$ and $\LM (g_i)=a^{\alpha (i)}e_k$ with $\alpha
(i)=(\alpha_{i_1},\ldots ,\alpha_{i_n})$, $\LM (g_j)=a^{\alpha
(j)}e_k$ with $\alpha (j)=(\alpha_{j_1},\ldots ,\alpha_{j_n})$, and
it follows from  the property $(*)$ mentioned before the lemma that
$\alpha_{i_{\ell}}=0=\alpha_{j_{\ell}}$, $1\le \ell\le r$,
$\alpha_{i_{r+1}}\le\alpha_{j_{r+1}}$. By Theorem 3.1.1, $\LM
(s_{ij})=a^{\gamma -\alpha (j)}\varepsilon_j$ where $\gamma
=(\gamma_1,\ldots ,\gamma_n)$ with
$\gamma_k=\max\{\alpha_{i_k},\alpha_{j_k}\}$, in particular,
$\gamma_{r+1}=\alpha_{j_{r+1}}$. This shows  that the generators
$a_1,\ldots ,a_r,a_{r+1}$ of $A$ do not appear in $\LM (s_{ij})$.
\QED }\v5

Let $M=\sum^s_{i=1}Av_i$ be a finitely generated left $A$-module
with generating set $\{ v_1,\ldots ,v_s\}$. Then since $A$ is
Noetherian, if we consider the free $A$-module
$L_0=\oplus_{i=1}^sAe_i$ and the presentation $M\cong L_0/N_0$ of
$M$ by a submodule $N_0$ of $L_0$, then $N_0$ is a finitely
generated submodule of $L_0$.{\parindent=0pt\v5

{\bf 3.2.2. Theorem} (noncommutative version of Hilbert's syzygy
theorem) Let $A=K[a_1,\ldots ,a_n]$ be a solvable polynomial algebra
with admissible system $(\B ,\prec )$. Then every finitely generated
left $A$-module $M$ has a finite free resolution
$$0~\mapright{}{}~L_q~\mapright{\varphi_q}{}~L_{q-1}~\mapright{\varphi_{q-1}}{}\cdots ~
\mapright{\varphi_1}{}~L_0~\mapright{\varphi_0}{}~M~\mapright{}{}~
0$$ such that $q\le n$.  \vskip 6pt

{\bf Proof} As remarked above we may assume $M=L_0/N_0$,  where
$L_0=\oplus_{i=1}^sAe_i$ is a free $A$-module with $A$-basis $\{
e_1,\ldots ,e_s\}$, and $N_0$ is a finitely generated submodule of
$L_0$. Then $\varphi_0$ is given by the canonical $A$-module
homomorphism $L_0~\mapright{\varphi_0}{}~M$. Fixing a left monomial
ordering $\prec_{e}$ on the $K$-basis $\B (e)$ of $L_0$, let
$\G_0=\{ g_1,\ldots ,g_{s_1}\}$ be a left Gr\"obner basis of $N_0$
as specified before Lemma 3.2.1, such that for some $r\le n$, the
generators $a_1,\ldots ,a_r$ of $A$ do not appear in every $\LM
(g_{\ell})$, $1\le \ell\le s_1$. Let $\prec_{s\hbox{-}\varepsilon}$
be the Schreyer ordering on the $K$-basis $\B (\varepsilon )$ of the
free $A$-module $L_1=\oplus_{i=1}^{s_1}A\varepsilon_i$ induced by
$\G_0$ with respect to $\prec_{e}$. Then by Lemma 3.2.1, the
generators $a_1,\ldots ,a_r,a_{r+1}$ do not appear in $\LM (s_{ij})$
for every $s_{ij}\in {\cal S}$, where ${\cal S}$ is the left
Gr\"obner basis of $N_1=$ Syz$(\G_0)\subset L_1$ obtained in Theorem
3.1.1. Defining $\varphi_1$: $L_1\r L_0$ by
$\varphi_1(\varepsilon_i)=g_i$, $1\le i\le s_1$,  and working with
$N_1$ in place of $N_0$ and so on, we then reach an exact sequence
$$ 0~\r~N_{n-r}~\mapright{}{}~L_{n-r}~\mapright{\varphi_{n-r}}{}~L_{n-r-1}~\mapright{\varphi_{n-r-1}}{}~\cdots~\mapright{\varphi_2}{}~L_1~\mapright{\varphi_1}{}~L_0~\mapright{\varphi_0}{}~M~\r~0$$
where every $L_i$ is a free $A$-module of finite $A$-basis and
$N_{n-r}=$Ker$\varphi_{n-r}$ which has a left Gr\"obner basis
$\G_{n-r}=\{ g_1',\ldots ,g_d'\}$ in $L_{n-r}$, such that all the
generators $a_1,\ldots ,a_n$ of $A$ do not appear in every $\LM
(g_j')$, $1\le j\le d$. If $L_{n-r}=\oplus_{i=1}^tA\omega_i$, this
shows that $\LM (g_{j}')=\omega_{k_j}$, $1\le j\le d$. Without loss
of generality, we may assume that all the $\omega_{k_j}$ are
distinct (for instance, $\G_{n-r}$ is minimal), and that $\LM
(g_j')=\omega_j$, $1\le j\le d$. Thus, since $\G_{n-r}$ is a left
Gr\"obner basis, it follows from (Chapter 2, Proposition 2.2.7(iv))
that in this case Ker$\varphi_{n-r-1}\cong
L_{n-r}/N_{n-r}\cong\oplus_{i=d+1}^tA\omega_i$. Now, with $L_{n-r}$
replaced by Ker$\varphi_{n-r-1}$, the desired (length $\le n$) free
resolution of $M$ is obtained. \QED}\v5

Combining the results of the last section, we are ready to have an
algorithmic procedure for constructing a finite free resolution.
{\parindent=0pt \v5

{\bf 3.2.3. Corollary} Let $M=L_0/N_0$ be a finitely generated
$A$-module as in Theorem 3.2.2, and let $\G_0=\{ g_1,\ldots ,g_t\}$
be a  left Gr\"obner basis of the submodule $N_0=\sum^{t}_{i=1}Ag_i$
with respect to a left monomial ordering $\prec_e$ on $L_0$. Then,
starting with the exact sequence
$$
0\r N_0~\mapright{\iota}{}~L_0~\mapright{\varphi_0}{}~M\r 0,$$ where
$\varphi_0$ is the canonical homomorphism, and $\iota$ is the
inclusion map, the following algorithm returns a finite free
resolution of $M$, which is  of length $q\le n$.\vskip 6pt

\underline{\bf Algorithm-FRES
~~~~~~~~~~~~~~~~~~~~~~~~~~~~~~~~~~~~~~~~~~~~~~~~~~~~~~~~~~~~}\vskip
6pt

\textsc{INPUT} $L_0=\oplus_{i=1}^sAe_i$, $\prec_e$, 
$\G_0=\{g_1,\ldots ,g_t\}$, $L_0~\mapright{\varphi_0}{}~M\r 0$
  \vskip 6pt

\textsc{OUTPUT} ${\cal L}_{\bullet}$\quad $ 0\r
L_q~\mapright{\varphi_q}{}~L_{q-1}~\mapright{\varphi_{q-1}}{}~\cdots~\mapright{\varphi_2}{}~
L_1~\mapright{\varphi_1}{}~L_0~\mapright{\varphi_0}{}~M\r 0$\par

~~~~~~~~~~~~~~~~~~~~~~a finite free resolution of $M$\par

{\textsc{INITIALIZATION} $i:=0$,  $\prec\hskip .1truecm :=\hskip 
.1truecm \prec_e$\par

\textsc{LOOP}\par

\textsc{IF} all the generators $a_1,\ldots ,a_n$ of $A$ do not 
appear in $\LM (g_j)$ for every \par                                          
~~~~~$g_j\in\G_{i}$, \textsc{THEN}\par

\hskip .6truecm by the proof of Theorem 3.2.2, there is some $d$ 
such that\par

\hskip .6truecm $\hbox{Ker}\varphi_{i}\cong 
L_i/\sum^t_{j=1}Ag_j\cong\oplus_{i=1}^dAe_i$ 
\par

${\cal L}_{\bullet}:=$ $(0\r 
\hbox{Ker}\varphi_{i}~\mapright{\iota}{}~L_{i}~\mapright{\varphi_{i}}{}~\cdots~\mapright{\varphi_2}{}~ 
L_1~\mapright{\varphi_1}{}~L_0~\mapright{\varphi_0}{}~M\r 0$)\par 

\textsc{ELSE}\par

$i:=i+1$,~$L_i:=\oplus_{j=1}^tAe_j$,~$\varphi_{i}:=$ $(L_i\r 
L_{i-1}$ with $\varphi (e_j)=g_j$, $1\le j\le t$)\par

\hskip .6truecm run {\bf Algorithm-DIV-L} (with respect to $\prec$) 
to compute a left\par\hskip .6truecm Gr\"obner representation of 
each $S_{\ell}(g_k,g_{\ell})\ne 0$ by $\G_{i-1}$ in 
$L_{i-1}$,\par\hskip .6truecm $1\le k<\ell\le t$, so that a left 
Gr\"obner basis $\S =\{ \S_1,\ldots ,\S_r\}$\par\hskip .6truecm of 
Syz$(\G)\subset L_{i}$ is obtained under the Schreyer ordering 
$\prec_{s\hbox{-}e}$ on\par\hskip .6truecm $L_i$  induced by 
$\G_{i-1}$ with respect to $\prec$ (Theorem 3.1.1)\vskip 6pt

$\prec\hskip .1truecm :=\hskip .1truecm \prec_{s\hbox{-}e}$\par

$\G_i :=\{g_1,\ldots ,g_{t}\}$ with $g_j=\S_j\ne 0$ for $1\le j\le 
t\le r$\par

\textsc{END}\par

\textsc{UNTIL} all the generators $a_1,\ldots ,a_n$ of $A$ do not 
appear in $\LM (g_j)$ for\par\hskip 1.5truecm  every $g_j\in\G_{i}$

\textsc{END}\vskip -.2truecm \underline{
~~~~~~~~~~~~~~~~~~~~~~~~~~~~~~~~~~~~~~~~~~~~~~~~~~~~~~~~~~~~~~~~~~~~~~~~~~~~~~~~~~~~~~~~~~~~~~}
\vskip 6pt

{\bf Proof} By the definition of a free resolution, the sequence
${\cal L}_{\bullet}$ returned by the algorithm is clearly the
desired one for $M$. The fact that the {\bf Algorithm-FRES}
terminates and  returns a sequence ${\cal L}_{\bullet}$ of finite
length $q\le n$ is due to Theorem 3.2.2 (or more precisely its
proof)}}\v5

\section*{3.3. Global Dimension and Stability}

In this section, the foregoing Theorem 3.2.2 is applied to highlight 
that every solvable polynomial algebra $A$ has finite global 
homological dimension,  and that every finitely generated projective 
$A$-module is stably free. \v5

Let $A=K[a_1,\ldots ,a_n]$ be a solvable polynomial algebra with 
admissible system $(\B ,\prec )$ in the sense of Definition 1.1.3, 
where $\B =\{ a^{\alpha}=a_1^{\alpha_1}\cdots 
a_n^{\alpha_n}~|~\alpha =(\alpha_1,\ldots ,\alpha_n)\in\NZ^n\}$ is 
the PBW $K$-basis of $A$ and $\prec$ is a monomial ordering on $\B$.  
Recall from classical homological algebra (e.g. [Rot]) that a left 
$A$-module $P$ is said to be {\it projective} if 
$M~\mapright{\beta}{}~N\r0$ is an exact sequence of $A$-modules and 
$P~\mapright{\alpha}{}~N$ is an $A$-module homomorphism, then there 
exists an $A$-module homomorphism $P~\mapright{\gamma}{}~M$ such 
that $\beta\circ\gamma =\alpha$, or in other words, the following 
diagram commutes:

$$\begin{array}{cccc} &&P&\\
&\scriptstyle{\gamma}\swarrow&\mapdown{\alpha}&\\
M&\mapright{}{\beta}&N&\r 0\end{array}\quad\quad \beta\circ\gamma
=\alpha$$\par

Thus, by the universal property of a free module, any free
$A$-module is projective and, an $A$-module $P$ is projective if and
only if there is a free $A$-module $L$ such that $P$ is isomorphic
to a direct summand of $L$, i.e., $L\cong P\oplus L_1$.\par

Let $M$ be a  left $A$-module. Then a {\it  projective resolution}
of $M$ is an {\it exact sequence} by projective  $A$-modules $P_i$
and  $A$-module homomorphisms $\varphi_i$:
$${\cal P}_{\bullet}\quad\quad\cdots~\mapright{\varphi_{i+1}}{}~P_i~\mapright{\varphi_i}{}~\cdots ~
\mapright{\varphi_2}{}~P_1~\mapright{\varphi_1}{}~P_0~\mapright{\varphi_0}{}~M~\mapright{}{}~0$$
that is, in the sequence ${\cal P}_{\bullet}$, $\varphi_0$ is 
surjective,  each $P_i$ is a projective $A$-module and 
Ker$\varphi_i=$ Im$\varphi_{i+1}=\varphi_{i+1}(P_{i+1})$ for all 
$i\ge 0$. Clearly, any free resolution of $M$ is a projective 
resolution of $M$. The {\it projective dimension} of $M$, denoted 
p.dim$_AM$, is defined to be the shortest length $q$ of a projective 
resolution or $\infty$ if no finite projective resolution exists. 
The long Schanuel lemma (cf. [Rot]) shows that any projective 
resolution of $M$ can be terminated at this length. The {\it left 
global dimension} of $A$ is then defined to be
$$\sup\{\hbox{p.dim}_AM~|~M~\hbox{any left}~A\hbox{-module}\} .$$
In terms of right $A$-modules, the right global dimension of $A$ is
defined in a similar way. It follows from classical homological
algebra (e.g. [Rot]) that the left (right) global dimension is
determined by the projective dimension of cyclic modules; moreover,
for a (left and right) Noetherian ring $A$, the left global
dimension of $A$ is equal to the right global dimension of $A$, and
this  common number is called the {\it global dimension of $A$},
denoted  gl.dim$A$. {\parindent=0pt\v5

{\bf 3.3.1. Theorem} Let $A=K[a_1,\ldots ,a_n]$ be a solvable
polynomial algebra. Then any $A$-module $M$ has projective dimension
p.dim$_AM\le n$, thereby $A$ has global  dimension gl.dim$A\le
n$.\vskip 6pt

{\bf Proof} Noticing the classical results on global dimension
reviewed above, this follows from the fact that $A$ is (left and 
right) Noetherian (Corollary 1.3.3) and the noncommutative version 
of Hilbert's syzygy theorem  (Theorem 3.2.2). \QED}\v5

Furthermore, since $A$ is Noetherian, $A$ has {\it IBN} (invariant
basis number), i.e., for every free $A$-module $L$, every two 
$A$-bases of $L$ have the same cardinal (cf. [Rot], Chapter 3). In 
this case, the {\it rank} of $L$, denoted rank$_AL$, is well defined 
as the cardinal of some $A$-basis of $L$, thereby if $L_1$, $L_2$ 
are free $A$-modules, then $L_1\cong L_2$ if and only if 
rank$_AL_1=$ rank$_AL_2$. With this well-defined rank for free 
modules, the stably free modules are then well defined, that is, an 
$A$-module $P$ is said to be {\it stably free of rank} $t$ if 
$P\oplus L_1\cong L_2$, where $L_1$ is a free module of 
rank$_AL_1=s$ and $L_2$ is a free module of rank$_AL_2=s+t$. 
Obviously, a stably free module is necessarily finitely generated 
and projective. It follows from the literature (e.g., [Rot], Chapter 
4; [MR], Chapter 11) that stably free $A$-modules can be 
characterized by finite free resolutions.  {\parindent=0pt\v5

{\bf 3.3.2. Proposition} A finitely generated projective $A$-module
$P$ is stably free if and only if $P$ has a finite free resolution.
Furthermore, if
$${\cal L}_{\bullet}\quad  0\r
L_q~\mapright{\varphi_q}{}~L_{q-1}~\mapright{\varphi_{q-1}}{}~\cdots~\mapright{\varphi_2}{}~
L_1~\mapright{\varphi_1}{}~L_0~\mapright{\varphi_0}{}~P\r 0$$ is a
finite free resolution of $P$, then
rank$_AP=\sum^q_{i=0}(-1)^i\hbox{rank}L_i$.\par\QED\v5

{\bf 3.3.3. Theorem} Let $A=K[a_1,\ldots ,a_n]$ be a solvable
polynomial algebra. Then every finitely generated projective
$A$-module $P$ is stably free, and moreover, rank$_AP$ is computable
via constructing a finite free resolution by running {\bf
Algorithm-FRES} given in Corollary 3.2.3.\vskip 6pt

{\bf Proof} This follows from the noncommutative version of
Hilbert's syzygy theorem  (Theorem 3.2.2), Theorem 3.2.3, and
Proposition 3.3.2 above. }\v5

\section*{3.4. Calculation of p.dim$_AM$}\vskip 6pt

Let $A=K[a_1,\ldots ,a_n]$ be a solvable polynomial algebra with 
admissible system $(\B ,\prec )$ in the sense of Definition 1.1.3, 
where $\B =\{ a^{\alpha}=a_1^{\alpha_1}\cdots 
a_n^{\alpha_n}~|~\alpha =(\alpha_1,\ldots ,\alpha_n)\in\NZ^n\}$ is 
the PBW $K$-basis of $A$ and $\prec$ is a monomial ordering on $\B$. 
Equipped  with the previously developed theory and techniques, in 
this section  we establish an algorithmic procedure which calculates 
the projective dimension of a finitely generated $A$-module $M$ and,  
at the same time, verifies whether $M$ is projective or not. The 
strategy used in our text was proposed by Gago-Vargas in [GV], 
though the algebras considered in [GV] are restricted to Weyl 
algebras.\v5

Let $L~\mapright{\varphi}{}~M\r 0$ be an $A$-module epimorphism, and
$K=$ Ker$\varphi$. Consider the short exact sequence
$0\r~K~\mapright{\iota}{}~L~\mapright{\varphi}{}~M\r 0$, where
$\iota$ is the inclusion map. Suppose that there exists an
$A$-module homomorphism $M~\mapright{\OV{\varphi}}{}~L$ such that
$\varphi\circ\OV{\varphi} =1_M$, where $1_M$ is the identity map of
$M$ to $M$. Then, every $\xi\in L$ has a representation  as $\xi
=(\xi -\OV{\varphi}(\varphi (\xi )))+ \OV{\varphi}(\varphi (\xi ))$
with $\OV{\varphi}(\varphi (\xi ))\in$ Im$\OV{\varphi}$ and $\xi
-\OV{\varphi}(\varphi (\xi ))\in K$. Since it is clear that $K\cap$
Im$\OV{\varphi}=\{ 0\}$ and $\varphi\circ\OV{\varphi}=1_M$ implies
that $\OV{\varphi}$ is a monomorphism, it turns out that
$$L=K\oplus~ \hbox{Im}\OV{\varphi}=K\oplus\OV{\varphi}(M)\cong
K\oplus M.\eqno{(1)}$$ This preliminary enables us to prove the
following{\parindent=0pt\v5

{\bf 3.4.1. Proposition} Let $0\r
L_1~\mapright{\varphi_1}{}~L_0~\mapright{\varphi_0}{}~M\r 0$ be an
exact sequence of $A$-modules in which $L_0$, $L_1$ are free
$A$-modules.  Then $M$ is projective if and only if there exists an
$A$-module homomorphism $L_0~\mapright{\OV{\varphi}_1}{}~L_1$ such
that $\OV{\varphi}_1\circ\varphi_1=1_{L_1}$, where $1_{L_1}$ is the
identity map of $L_1$ to $L_1$.\vskip 6pt

{\bf Proof} Suppose that $M$ is projective. Then there exists an
$A$-module homomorphism $M~\mapright{\OV{\varphi}_0}{}~L_0$ such
that $\varphi\circ\OV{\varphi}_0=1_M$, where $1_M$ is the identity
map of $M$ to $M$. Hence, by the formula (1) above we have
$L_0=\varphi_1(L_1)\oplus\OV{\varphi}_0(M)$. It follows that if we
define $L_0~\mapright{\OV{\varphi}_1}{}~L_1$  by
$$\OV{\varphi}_1 (\varphi_1(\xi_1)+\OV{\varphi}_0(m))=\xi_1,\quad \xi_1\in L_1,~m\in M, $$
then since $\varphi_1$ is injective, it is easy to see that
$\OV{\varphi}_1$ is an $A$-module homomorphism satisfying
$\OV{\varphi}_1\circ\varphi_1=1_{L_1}$. }\par

Conversely, suppose that there exists an $A$-module homomorphism
$L_0~\mapright{\OV{\varphi}_1}{}~L_1$ such that
$\OV{\varphi}_1\circ\varphi_1=1_{L_1}$. Then the sequence
$$0\r K=\hbox{Ker}\OV{\varphi}_1~\mapright{\iota}{}~L_0~\mapright{\OV{\varphi}_1}{}~L_1\r 0$$
is exact. Since $L_1$ is free (hence projective), it follows from
the formula (1) above that $L_0=K\oplus\varphi_1(L_1)$, thereby
$K\cong L_0/\varphi_1(L_1)=L_0/\hbox{Ker}\varphi_0\cong M$. Note
that as a direct summand of the free module $L_0$, $K$ is
projective. Hence $M$ is projective, as desired.\QED\v5

Let $L_0=\oplus_{j=1}^sAe_j$, $L_1=\oplus_{i=1}^tA\varepsilon_i$ be
free left $A$-modules of rank $s$ and $t$ respectively, and let
$L_1~\mapright{\varphi_1}{}~L_0$ be an $A$-module homomorphism with
$\varphi_1 (\varepsilon_i)=\sum^s_{j=1}f_{ij}e_j$, $1\le i\le t$.
Then $\varphi_1$ is uniquely determined by the $t\times s$ matrix
$$Q_{\varphi_1} =\left(\begin{array}{cccc} f_{11}&f_{12}&\cdots&f_{1s}\\
f_{21}&f_{22}&\cdots&f_{2s}\\
\vdots&\vdots&\cdots&\vdots\\
f_{t1}&f_{t2}&\cdots&f_{ts}\end{array}\right ) ,$$ that is, if $\xi
=\sum^t_{i=1}f_i\varepsilon_i\in L_1$, then, $\varphi_1 (\xi )$ is
given by left multiplication by matrices:
$$\varphi_1 (\xi )=\sum^t_{i=1}f_i\varphi_1 (\varepsilon_i)=(f_1,\ldots ,f_t)Q_{\varphi_1}\left (
\begin{array}{c} e_1\\ \vdots \\ e_s\end{array}\right ).$$
The $t\times s$ matrix $Q_{\varphi_1}$ is usually referred to as the
{\it matrix of $\varphi_1$}. Thus, in the language of matrices,
Proposition 3.4.1 can be restated as follows.{\parindent=0pt\v5

{\bf 3.4.2. Proposition} Let the short exact sequence $0\r
L_1~\mapright{\varphi_1}{}~L_0~\mapright{\varphi_0}{}~M\r 0$ be as
in Proposition 3.4.1, and let $Q_{\varphi_1}$ be the matrix of
$\varphi_1$. Then $M$ is projective if and only if the $t\times s$
matrix $Q_{\varphi_1}$ is right invertible, i.e., there is an
$s\times t$ matrix $Q_{\OV{\varphi}_1}$ with entries in $A$, such
that $Q_{\varphi_1}\cdot Q_{\OV{\varphi}_1}=E_{t\times t}$, where
the latter is the $t\times t$ identity matrix with the $(i,i)$-entry
$1\in A$.\par\QED}\v5

Furthermore, let the $t\times s$ matrix $Q_{\varphi_1}$ of
$\varphi_1$ be as above, and write $Q_{\varphi_1}^j$ for the $j$-th
column of $Q_{\varphi_1}$, $1\le j\le s$. Considering the  free {\it
right} $A$-module $\OV L_1=\oplus^t_{i=1}\varepsilon_iA$ of rank
$t$, let  $N=\sum^s_{j=1}\xi_jA\subseteq \OV L_1$ be the
$A$-submodule generated by $\xi_j=(\varepsilon_1,\ldots
,\varepsilon_t)Q_{\varphi_1}^j=\sum^t_{i=1}\varepsilon_if_{ij}$,
$1\le j\le s$. Then Proposition 3.4.2 tells us that the left 
$A$-module $M$ is projective if and only if  $N=\OV L_1$. Since $A$ 
has also a {\it right Gr\"obner basis theory} for {\it right 
modules}, it follows that the following proposition holds true. 
{\parindent=0pt\v5

{\bf 3.4.3. Proposition} With notation as above, let $\G$ be a right
Gr\"obner basis of the right $A$-submodule
$N=\sum^s_{j=1}\xi_jA\subseteq \OV
L_1=\oplus_{i=1}^t\varepsilon_iA$. Then, the left $A$-module $M$ is
projective if and only if  $\varepsilon_i\in\G$ for $1\le i\le t$.
If it is the case, then the right inverse $Q_{\OV{\varphi}_1}$ of
$Q_{\varphi_1}$ is given by the  $s\times t$ matrix $V_{s\times t}$
such that
$$(\varepsilon_1,\ldots ,\varepsilon_t)=(\xi_1,\ldots ,\xi_s)V_{s\times t},$$
which can be computed via using (Chapter 2, Corollary
2.3.5(ii)).\vskip 6pt

{\bf Proof} Note that $N=\OV L_1$ if and only if all the
$\varepsilon_i \in N$, $1\le i\le t$. If $\G$ is a right Gr\"obner
basis of $N$, then by (Chapter 2, Proposition 2.2.7), $N=\OV L_1$ if
and only if all the $\varepsilon_i\in\G$, $1\le i\le t$. Since
$N=\sum^s_{j=1}\xi_jA$, it follows from (Chapter 2, Corollary
2.3.5(ii)) that the specified $s\times t$ matrix $V_{s\times t}$ can
be computed. Also since $(\xi_1,\ldots ,\xi_s)=(\varepsilon_1,\ldots
,\varepsilon_t)Q_{\varphi_1}$, it turns out that
$$\begin{array}{rcl}(\varepsilon_1,\ldots ,\varepsilon_t)&=&(\xi_1,\ldots ,\xi_s)V_{s\times t}\\
&=&(\varepsilon_1,\ldots ,\varepsilon_t)Q_{\varphi_1}V_{s\times
t}.\end{array}$$ This shows that $V_{s\times t}=Q_{\OV{\varphi}_1}$,
as desired.\QED}\v5

Now, let $M$ be a finitely generated $A$-module and $${\cal
L}_{\bullet}\quad  0\r
L_q~\mapright{\varphi_q}{}~L_{q-1}~\mapright{\varphi_{q-1}}{}~\cdots~\mapright{\varphi_2}{}~
L_1~\mapright{\varphi_1}{}~L_0~\mapright{\varphi_0}{}~M\r 0$$ a
finite free resolution of $M$ by free modules of finite rank.
Suppose that Im$\varphi_i$ is projective for some $i\ge 0$. Then, by
the foregoing discussion, it is not difficult to derive inductively 
that Im$\varphi_{i+k}$ is projective for every $k=1,2,\ldots ,q-i$. 
In particular, Im$\varphi_{q-1}$ is projective. By the definition of 
p.dim$_AM$ and the basic property we recalled before Theorem 3.3.1, 
The next proposition is clear.{\parindent=0pt\v5

{\bf 3.4.4. Proposition} With notation as above, p.dim$_AM=q$ if and
only if Im$\varphi_{q-1}$ is not projective if and only if the
matrix $Q_{\varphi_q}$ of $\varphi_q$ is not right invertible, where
the invertibility of $Q_{\varphi_q}$ can be recognized in a
computational way via using Proposition 3.4.3. \par\QED}\v5

Suppose that Im$\varphi_{q-1}$ is projective. It follows from
Proposition 3.4.1 (or its proof) that
$$L_{q-1}=\varphi_q(L_q)\oplus~\hbox{Ker}~\OV{\varphi}_q~\mapright{\psi}{\cong}~
L_q\oplus~\hbox{Im}~\varphi_{q-1},\eqno{(2)}$$ and the latter
isomorphism $\psi$ can be computed, where $\OV{\varphi}_q$ is the
$A$-module homomorphism $L_{q-1}~\mapright{\OV{\varphi}_q}{}~L_q$
such that $\OV{\varphi}_q\circ\varphi_q=1_{L_q}$. The formula (2)
above enables us to construct another finite free resolution of $M$
$$0\r
L_q~\mapright{\varphi_q'}{}~L_q\oplus
L_{q-1}~\mapright{\varphi_{q-1}'}{}~L_q\oplus
L_{q-2}~\mapright{\varphi_{q-2}'}{}~L_{q-3}
~\mapright{\varphi_{q-3}}{}~\cdots~\mapright{\varphi_1}{}~L_0~\mapright{\varphi_0}{}~M\r
0$$ in which  $$\begin{array}{l}
\varphi_q'(\xi_q)=\varphi_q(\xi_q)~\hbox{for all}~\xi_q\in L_q, \\
\\
\varphi_{q-1}'(\xi_q+\xi_{q-1})=\xi_q+\varphi_{q-1}(\xi_{q-1})~\hbox{for
all}~\xi_q+\xi_{q-1}\in L_q\oplus L_{q-1},\\
\\
\varphi_{q-2}'(\xi_q+\xi_{q-2})=\varphi_{q-2}(\xi_{q-2})~\hbox{for
all}~ \xi_q+\xi_{q-2}\in L_q\oplus L_{q-2}.\end{array}$$  By the
exactness of free resolution and the formula (2) above, we then have
$$L_{q-1}~\mapright{\psi}{\cong}~ L_q\oplus\varphi_{q-1}(L_{q-1})=~
\hbox{Im}\varphi_{q-1}'=\hbox{Ker}\varphi_{q-2}',$$ thereby $M$ has
the following finite free resolution
$$0\r
L_{q-1}~\mapright{\psi}{}~L_q\oplus
L_{q-2}~\mapright{\varphi_{q-2}'}{}~L_{q-3}
~\mapright{\varphi_{q-3}}{}~\cdots~\mapright{\varphi_1}{}~L_0~\mapright{\varphi_0}{}~M\r
0$$ in which the homomorphism  $\psi$ can be computed. By
Proposition 3.4.4,  after the projectiveness of Im$\varphi_{q-2}$ is 
checked we can either have p.dim$_AM=q-1$, or repeat the above 
procedure  again in order to get a finite free resolution of $M$ 
which is of length $q-2$. It is clear that after a  finite number of 
repetitions of the same procedure we will eventually have 
p.dim$_AM=m$ with $m\le q$. Obviously, if $m=0$, then $M$ is 
projective.\v5

Summing up, we have reached the following{\parindent=0pt\v5

{\bf 3.4.5. Theorem} Let $M$ be a finitely generated $A$-module, and
let $${\cal L}_{\bullet}\quad  0\r
L_q~\mapright{\varphi_q}{}~L_{q-1}~\mapright{\varphi_{q-1}}{}~\cdots~\mapright{\varphi_2}{}~
L_1~\mapright{\varphi_1}{}~L_0~\mapright{\varphi_0}{}~M\r 0$$ be a
finite free resolution of $M$ computed by running {\bf
Algorithm-FRES}. Then the next algorithm computes p.dim$_AM$ and,  
meanwhile, checks whether $M$ is projective or not.\vskip 6pt

\underline{\bf Algorithm-p.dim
~~~~~~~~~~~~~~~~~~~~~~~~~~~~~~~~~~~~~~~~~~~~~~~~~~~~~~~~~~~~}\v5

\textsc{INPUT} ${\cal L}_{\bullet}$ the given finite free resolution
; $q$ the length of ${\cal L}_{\bullet}$\par

\textsc{OUTPUT} p.dim$_AM$\par

{\textsc{INITIALIZATION} $i:=q$\par

\textsc{IF} $i=0$ \textsc{THEN}\par

p.dim$_AM:=0$\par

\textsc{ELSE} \par

~~~~\textsc{LOOP}\par

~~~~~~~~use Proposition 3.3.6 to check the invertibility of\par 
~~~~~~~~the matrix $Q_{\varphi_i}$ of $\varphi_i$\par

~~~~\textsc{IF} $Q_{\varphi_i}$ is  not right invertible
\textsc{THEN}\par

~~~~p.dim$_AM:=i$\par

~~~~\textsc{ELSE}\par

~~~~$i:=i-1$\par

~~~~~~~~\textsc{IF} $i=0$ \textsc{THEN}\par

~~~~~~~~p.dim$_AM: =0$\par

~~~~~~~~\textsc{ELSE}\par

~~~~~~~~${\cal L}_{\bullet}:=$ $(0\r
L_{i-1}~\mapright{\varphi_{i-1}'}{}~L_{i-2}'~\mapright{\varphi_{i-2}'}{}~L_{i-3}
~\mapright{\varphi_{i-3}}{}~\cdots~\mapright{\varphi_1}{}~L_0~\mapright{\varphi_0}{}~M\r
0$\par

\hskip 2.3truecm in which $L_{i-2}'=L_i\oplus L_{i-2}$, the
homomorphism  $\varphi_{i-1}'=\psi$\par

\hskip 2.3truecm  is computed by using Proposition 3.3.4 (or its
proof), and\par

\hskip 2.3truecm $\varphi_{i-2}'$  is defined before the
theorem.)\par

~~~~~~~~\textsc{END}\par

~~~~\textsc{END}\par

~~~~\textsc{END}\par

\textsc{END}\vskip -.2truecm

 \underline{
~~~~~~~~~~~~~~~~~~~~~~~~~~~~~~~~~~~~~~~~~~~~~~~~~~~~~~~~~~~~~~~~~~~~~~~~~~~~~~~~~~~~~~~~~~~~~~}}}
\newpage\setcounter{page}{69}

\chapter*{4. Minimal Finite Graded\\ \hskip 1.25truecm Free Resolutions}
\markboth{\rm Minimal Graded Free Resolutions}{\rm Minimal Graded
Free Resolutions} \vskip 2.5truecm

In this chapter we demonstrate how the methods and algorithms,
developed in ([CDNR], [KR2]) for computing minimal homogeneous
generating sets of graded submodules and graded quotient modules of
free modules over commutative polynomial algebras, can be adapted
for computing minimal homogeneous generating sets of graded
submodules and graded quotient modules of free modules over an
$\NZ$-graded solvable polynomial $K$-algebras $A$ with the degree-0
homogeneous part $A_0=K$, and how the algorithmic procedures of
computing minimal graded free resolutions for finitely generated
modules over $A$ can be achieved. \par

In the literature, a finitely generated  $\NZ$-graded $K$-algebra 
$A=\oplus_{p\in\NZ}A_p$ with the degree-0 homogeneous part $A_0=K$ 
is usually referred to as a {\it connected $\NZ$-graded 
$K$-algebra}. Concerning introductions to minimal resolutions of 
graded modules over a (commutative or noncommutative) connected 
$\NZ$-graded $K$-algebra (or more generally an $\NZ$-graded local 
$K$-algebra) and relevant results, one may  refer to ([Eis], Chapter 
19), ([Kr1], Chapter 3), and [Li3].\v5

All notions, notations and conventions introduced before are 
maintained. \newpage

\section*{4.1. $\NZ$-graded Solvable Polynomial Algebras }\par

In this section we specify, by means of positive-degree functions 
(as defined in Section 1.1 of Chapter 1), the structure of 
$\NZ$-graded solvable polynomial $K$-algebras with the degree-0 
homogeneous part  equal to $K$. \v5

For convenience, we first  recall that  the condition (S2) given in 
(Definition 1.1.3 of Chapter 1) is equivalent to 
{\parindent=1.6truecm\vskip 6pt

\item{(S2$'$)} There is a monomial ordering $\prec$ on $\B$, i.e., $(\B 
,\prec )$ is an admissible system of $A$,  such that for all 
generators $a_i, a_j$ of $A$ with $1\le i<j\le n$,  
$$\begin{array}{rcl} a_ja_i&=&\lambda_{ji}a_ia_j+f_{ji}\\
&{~}&\hbox{where}~\lambda_{ji}\in K^*, ~f_{ji}=\sum\mu_ka^{\alpha 
(k)}\in K\hbox{-span}\B\\
&{~}&\hbox{with}~\LM (f_{ji})\prec a_ia_j~\hbox{if}~f_{ji}\ne 
0.\end{array}$$\par}

Now, let $A=K[a_1,\ldots ,a_n]$ be a solvable polynomial $K$-algebra 
with admissible system $(\B ,\prec )$, where $\B =\{ 
a^{\alpha}=a_1^{\alpha_1}\cdots a_n^{\alpha_n}~|~\alpha 
=(\alpha_1,\ldots ,\alpha_n)\in\NZ^n\}$ is the PBW $K$-basis of $A$ 
and $\prec$ is a monomial ordering on $\B$. Suppose that $A$ is an 
$\NZ$-graded algebra with the degree-0 homogeneous part equal to 
$K$, namely $A=\oplus_{p\in\NZ}A_p$, where the degree-$p$ 
homogeneous part $A_p$ is a $K$-subspace of $A$, $A_0=K$, and 
$A_{p_1}A_{p_2}\subseteq A_{p_1+p_2}$ for all $p_1,p_2\in\NZ$. Then, 
since conventionally any generator $a_i$ of $A$ is not contained in 
the ground field $K$, writing $d_{\rm gr}(f)=p$ for the {\it 
graded-degree} (abbreviated to {\it gr-degree}) of a nonzero 
homogeneous element $f\in A_p$, we have 
$$d_{\rm gr}(a_i)=m_i,\quad 1\le i\le n,$$
for some positive integers $m_i$. It turns out that if 
$a^{\alpha}=a_1^{\alpha_1}\cdots a_n^{\alpha_n}\in\B$, then $d_{\rm 
gr}(a^{\alpha})=\sum^n_{i=1}\alpha_im_i$. This shows that  $d_{\rm 
gr}(~)$ gives rise to a positive-degree function on $A$ as defined 
in (Section 1.1 of Chapter 1), such that{\parindent=1.2truecm\vskip 
6pt

\item{(1)} $A_p=K\hbox{-span}\{ a^{\alpha}\in\B~|~d_{\rm
gr}(a^{\alpha})=p\}$, $p\in\NZ$;\par

\item{(2)} for $1\le i<j\le n$, all the relations   $a_ja_i=\lambda_{ji}a_ia_j+f_{ji}$ with
$f_{ji}=\sum\mu_ka^{\alpha (k)}$ presented in (S2$'$) above, satisfy  
$d_{\rm gr}(a^{\alpha (k)})=d_{\rm gr}(a_ia_j)$ whenever $\mu_k\ne 
0$.\par}{\parindent=0pt\vskip 6pt

Conversely, given a positive-degree function $d(~)$ on $A$ (as 
defined in  Section 1.1 of Chapter 1) such that $d(a_i)=m_i>0$, 
$1\le i\le n$, then we know that $A$ has  an {\it $\NZ$-graded 
$K$-module structure}, i.e., $A=\oplus_{p\in\NZ}A_p$ with 
$A_p=K\hbox{-span}\{ a^{\alpha}\in\B~|~d(a^{\alpha})=p\}$, in 
particular, $A_0=K$. It is straightforward to verify that if 
furthermore for $1\le i<j\le n$, all the relations   
$a_ja_i=\lambda_{ji}a_ia_j+f_{ji}$ with $f_{ji}=\sum\mu_ka^{\alpha 
(k)}$ presented in (S2$'$) above, satisfy  $d_{\rm gr}(a^{\alpha 
(k)})=d_{\rm gr}(a_ia_j)$ whenever $\mu_k\ne 0$, then 
$A_{p_1}A_{p_2}\subseteq A_{p_1+p_2}$ holds for all $p_1,p_2\in\NZ$, 
i.e., $A$ is turned into an $\NZ$-graded solvable polynomial algebra 
with the degree-0 homogeneous part $A_0=K$.}\par

Summing up, we have reached the following{\parindent=0pt\v5

{\bf 4.1.1. Proposition} A solvable polynomial algebra 
$A=K[a_1,\ldots, a_n]$ is an $\NZ$-graded algebra with the degree-0 
homogeneous part $A_0=K$ if and only if there is a  positive-degree 
function $d(~)$ on $A$ (as defined in Section 1.1 of Chapter 1) such 
that for $1\le i<j\le n$, all the relations   
$a_ja_i=\lambda_{ji}a_ia_j+f_{ji}$ with $f_{ji}=\sum\mu_ka^{\alpha 
(k)}$ presented in (S2$'$) above, satisfy  $d_{\rm gr}(a^{\alpha 
(k)})=d_{\rm gr}(a_ia_j)$ whenever $\mu_k\ne 0$.\par\QED}\v5

To make the compatibility with the structure of $\NZ$-filtered
solvable polynomial algebras specified in Chapter 4, it is necessary 
to emphasize the role played by a positive-degree function in the 
structure of $\NZ$-graded solvable polynomial algebras we specified 
in this section, that is, from now on in the rest of this paper we 
keep using the following{\parindent=0pt\v5

{\bf Convention}  An $\NZ$-graded solvable polynomial $K$-algebra
$A=\oplus_{p\in\NZ}A_p$ with $A_0=K$ is always referred to as an
{\it $\NZ$-graded solvable polynomial algebra with respect to a
positive-degree function $d(~)$}. \v5

{\bf Remark} Let $A=K[a_1,\ldots ,a_n]=\oplus_{p\in\NZ}A_p$ be an
$\NZ$-graded solvable polynomial algebra with respect to a
positive-degree function $d(~)$.\par

(i) We emphasize that   {\it every $a^{\alpha}\in\B$ is a
homogeneous elements of $A$ and $d(a^{\alpha})=d_{\rm
gr}(a^{\alpha})$}, where $d_{\rm gr}(~)$, as we defined above, is
the gr-degree function on nonzero homogeneous elements of $A$.\par

(ii) Since  $A$ is a domain (Theorem 1.2.3), the gr-degree function 
$d_{\rm gr}(~)$  has the property that for all nonzero homogeneous 
elements $h_1,h_2\in A$, 
$$d_{\rm gr}(h_1h_2)=d_{\rm gr}(h_1)+d_{\rm gr}(h_2).\leqno{(\mathbb{P}1)}$$ 
From now on we shall freely use this property without additional 
indication.}\v5

Typical noncommutative $\NZ$-graded solvable polynomial algebras  
are provided by the multiplicative analogues ${\cal 
O}_n(\lambda_{ji})$ of the Weyl algebra (see Example (4) given in 
Section 1.1 of Chapter 1), where the positive-degree function on 
${\cal O}_n(\lambda_{ji})$  can be defined by setting $d(x_i)=m_i$ 
for any fixed tuple $(m_1,\ldots ,m_n)$ of positive integers.  \par

Another family of noncommutative  $\NZ$-graded solvable polynomial 
algebras are provided by the algebras $M_q(2,K)$ of $2\times 2$ 
quantum matrices (see Example (7) given in Section 1.1 of Chapter 
1), where each generator is assigned the degree 1. More generally, 
let $\Lambda =(\lambda_{ij})$ be a multiplicatively antisymmetric 
$n\times n$ matrix over $K$, and let $\lambda\in K^*$ with 
$\lambda\ne -1$. Considering the multiparameter coordinate ring of 
quantum $n\times n$ matrices over $K$ (see [Good]), namely the 
$K$-algebra ${\cal O}_{\lambda ,\Lambda}(M_n(K))$ generated by $n^2$ 
elements $a_{ij}$ ($1\le i,j\le n$) subject to the relations
$$a_{\ell m}a_{ij}=\left\{\begin{array}{ll} \lambda_{\ell i}\lambda_{jm}a_{ij}a_{\ell m}+
(\lambda -1)\lambda_{\ell i}a_{im}a_{\ell j}&(\ell >i,~m>j)\\
\lambda\lambda_{\ell i}\lambda_{jm}a_{ij}a_{\ell m}&(\ell >i,~m\le
j)\\ \lambda_{jm}a_{ij}a_{\ell m}&(\ell =i,~m>j)\end{array}\right.$$
Then ${\cal O}_{\lambda ,\Lambda}(M_n(K))$ is an $\NZ$-graded
solvable polynomial algebra, where each generator has degree 1.\par

Moreover, by ([LW], [Li1], or the later Chapter 4),  the associated 
graded algebra and the Rees algebra of every $\NZ$-filtered solvable 
polynomial algebra with a graded monomial ordering $\prec_{gr}$  are 
$\NZ$-graded solvable polynomial algebras of the type we specified 
in this section.\v5

Finally, we point out that if a solvable polynomial algebra   
$A=K[a_1,\ldots, a_n]$ has a graded monomial ordering $\prec_{gr}$ 
with respect to some given positive-degree function $d(~)$ (see 
Section 1.1 of Chapter 1),  then, by Proposition 4.1.1, it is easy 
to check whether $A$ is an $\NZ$-graded algebra with respect to 
$d(~)$ or not.\v5

\section*{4.2. $\NZ$-Graded Free Modules}

Let $A=K[a_1,\ldots ,a_n]$ be an $\NZ$-graded solvable polynomial 
algebra with respect to a positive-degree function $d(~)$, and let 
$(\B ,\prec )$ be an admissible system of $A$. In this section we 
demonstrate how to construct $\NZ$-graded free left modules over $A$ 
and,  furthermore, we highlight that if the input data is a finite 
set of nonzero homogeneous elements, then {\bf Algorithm-LGB} 
(presented in Theorem 2.3.4 of Chapter 2) produces a homogeneous 
left Gr\"obner bases for the graded submodule generated by the given 
homogeneous elements.\v5

We start by a little generality of $\NZ$-graded modules. Let $M$ be 
a left $A$-module. If $M=\oplus_{q\in\NZ}M_q$ with each $M_q$ a 
$K$-subspace of $M$, such that $A_pM_q\subseteq M_{p+q}$ for all 
$p,q\in\NZ$, then $M$ is called an {\it $\NZ$-graded $A$-module}.  
For each $q\in\NZ$, nonzero elements in $M_q$ are called {\it 
homogeneous elements of degree $q$}, and accordingly $M_q$ is called 
the {\it degree-$q$ homogeneous part} of $M$.  If $\xi\in M_q$ and 
$\xi\ne 0$, then we write $d_{\rm gr}(\xi )$ for the {\it 
graded-degree} (abbreviated to {\it gr-degree}) of $\xi$ as a 
homogeneous element of $M$, i.e., $d_{\rm gr}(\xi )=q$.\par

Let $M=\oplus_{q\in\NZ}M_q$ be a nonzero $\NZ$-graded $A$-module, 
and $T$ a subset of homogeneous elements of $M$. If $T$ generates 
$M$, i.e., $M=\sum_{\xi\in T}A\xi$, then $T$ is called a {\it 
homogeneous generating set} of $M$. Clearly, if $T=\{\xi_i~|~i\in 
I\}$ is a homogeneous generating set of $M$ with $d_{\rm 
gr}(\xi_i)=b_i$ for $\xi_i\in T$, then 
$M_q=\sum_{p_i+b_i=q}A_{p_i}\xi_i$ for all $q\in\NZ$.\par

If a submodule $N$ of the $\NZ$-graded $A$-module 
$M=\oplus_{q\in\NZ}M_q$ is generated by homogeneous elements, i.e., 
$N$ has a homogeneous generating set, then $N$ is called a {\it 
graded submodule} of $M$. A graded submodule $N$ has the 
$\NZ$-graded structure $N=\oplus_{q\in\NZ}N_q$ with $N_q=N\cap M_q$, 
such that $A_pN_q\subseteq N_{p+q}$ for all $p,q\in\NZ$. \par

With the graded submodule $N$ of $M$ as described above, the 
quotient module $M/N$ is an $\NZ$-graded $A$-module with the 
$\NZ$-graded structure $M/N=\oplus_{q\in\NZ}(M/N)_q$, where for each 
$q\in\NZ$, $(M/N)_q=(M_q+N)/N$. Indeed,  a submodule $N$ of $M$ is a 
graded submodule if and only if the quotient module $M/N$ is an  
$\NZ$-graded $A$-module with the $\NZ$-graded structure 
$M/N=\oplus_{q\in\NZ}(M_q+N)/N$.\v5

Now, let $L=\oplus_{i=1}^sAe_i$ be a free left $A$-module with the 
$A$-basis $\{ e_1,\ldots ,e_s\}$. Then $L$ has the $K$-basis $\B 
(e)=\{ a^{\alpha}e_i~|~a^{\alpha}\in\B ,~1\le i\le n\}$ and, for an 
{\it arbitrarily} fixed $\{b_1,\ldots ,b_s\}\subset\NZ$, one checks 
that $L$ can be turned into an {\it $\NZ$-graded free $A$-module} 
$L=\oplus_{q\in\NZ}L_q$ by setting
$$L_q=\{ 0\}~\hbox{if}~q<\min\{ b_1,\ldots ,b_s\};~\hbox{otherwise}~L_q=
\sum_{p_i+b_i=q}A_{p_i}e_i,\quad q\in\NZ ,$$ or alternatively, for
$q\ge \min\{ b_1,\ldots ,b_s\}$,
$$L_q=K\hbox{-span}\{ a^{\alpha}e_i\in\BE~|~d(a^{\alpha})+b_i=q\},
~q\in \NZ ,$$ such that $d_{\rm gr}(e_i)=b_i$, $1\le i\le s$. \par

As with the gr-degree of homogeneous elements in $A$, noticing that 
$d_{\rm gr}(a^{\alpha}e_i)=d(a^{\alpha})+b_i$ for all 
$a^{\alpha}e_i\in\BE$ and that $A$ is a domain, from now on we shall 
freely use the following property without additional indication: for 
all nonzero homogeneous elements $h\in A$ and all nonzero 
homogeneous elements $\xi\in L$,
$$d_{\rm gr}(h\xi )=d_{\rm gr}(h)+d_{\rm gr}(\xi ).\leqno{(\mathbb{P}2)}$${\parindent=0pt\par

{\bf Remark}  Although we have remarked that $d(a^{\alpha})=d_{\rm
gr}(a^{\alpha})$ for all $a^{\alpha}\in\B$, $d(a^{\alpha})$ is used
in constructing $L_q$ just for highlighting the role of $d(~)$.\v5

{\bf Convention} Unless otherwise stated, from now on throughout the
subsequent texts if we say that $L$ is an $\NZ$-graded free module 
over an $\NZ$-graded solvable polynomial algebra $A$ with respect to 
a positive-degree function $d(~)$, then it always means that  $L$ 
has an $\NZ$-gradation as constructed above.}\v5

Let $L=\oplus_{q\in\NZ}L_q$ be an $\NZ$-graded free $A$-module, and 
$N$ a graded submodule of $L$. A left Gr\"obner basis $\G$ of $N$ is 
called a {\it homogeneous left Gr\"obner basis} if $\G$ consists of 
homogeneous elements.\par

Note that monomials in $\B$ are homogeneous elements of $A$, thereby 
left S-polynomials of homogeneous elements are homogeneous elements, 
and remainders of homogeneous elements on division by homogeneous 
remain homogeneous elements. Thus,  the following assertion is clear 
now.{\parindent=0pt\v5

{\bf 4.2.1. Theorem} With notation as above, if a graded submodule 
$N=\sum_{i=1}^mA\xi_i$ of $L$ is generated by the set of nonzero 
homogeneous elements $\{ \xi_1,\ldots ,\xi_m\}$, then, with  the 
initial input data  $\{ \xi_1,\ldots ,\xi_m\}$, {\bf Algorithm-LGB} 
(presented in Theorem 2.3.4 of Chapter 2) produces a finite 
homogeneous left Gr\"obner basis $\G$ for $N$ with respect to any 
given monomial ordering $\prec_{e}$ on $\BE$. }\v5

\section*{4.3. Computation of Minimal Homogeneous Generating 
Sets}\par

In this section, $A=K[a_1,\ldots ,a_n]$ denotes an $\NZ$-graded
solvable polynomial algebra with respect to a positive-degree
function $d(~)$, $(\B,\prec )$ denotes a fixed admissible system of
$A$, $L=\oplus_{i=1}^sAe_i$ denotes an $\NZ$-graded free $A$-module 
such that $d_{\rm gr}(e_i)=b_i$, $1\le i\le s$, and  $\prec_{e}$ 
denotes a fixed left monomial ordering on the $K$-basis $\BE$ of 
$L$. Moreover,  as before we write $S_{\ell}(\xi_i,\xi_j)$ for the 
left S-polynomial of two elements $\xi_i$, $\xi_j\in L$.\par

Let $N$ be a finitely generated graded submodule $N$ of $L$. With 
the preparation made in the previous two sections, our aim of the 
current  section is to provide an algorithmic way of 
computing{\parindent=.7truecm\par

\item{(1)}  a minimal homogeneous generating set of $N$, and \par

\item{(2)} a minimal homogeneous generating set of the graded quotient module $M=L/N$.\v5}

The argumentation  we are going to present below is similar to the 
commutative case (cf. [KR], Section 4.5, Section 4.7). To better 
understand why the similar argumentation can go through the 
noncommutative case, let us again remind that although  monomials 
from the PBW $K$-basis $\B$ of  $A$ {\it can no longer behave as 
well as monomials in a commutative polynomial algebra} (namely the 
product of two monomials is not necessarily a monomial), every 
monomial from $\B$ is a homogeneous element in the $\NZ$-graded 
structure of $A$ (as we remarked in section 4.1), thereby the 
product of two monomials is a homogeneous element. Bearing in mind 
this fact, one will see that the rule of division, Proposition 
1.1.4(i) (Section 1 of Chapter 1), Lemma 2.1.2(ii) (Section 1 of 
Chapter 2), and the properties $(\mathbb{P}1)$, $(\mathbb{P}2)$ 
mentioned in previous Section 1 and Section 2 respectively, all 
together make the work done.\v5

We start by a detailed discussion on computing $n$-truncated left 
Gr\"obner bases for graded submodules of $L$. Except for helping us 
to compute a minimal homogeneous generating set, from the definition 
and the characterization given below one may see clearly that having 
an $n$-truncated left Gr\"obner basis will be very useful in dealing 
with certain problems involving only degree-n homogeneous elements. 
{\parindent=0pt\v5

{\bf 4.3.1. Definition} Let $G=\{ g_1,\ldots ,g_t\}$ be a subset of
nonzero homogeneous elements of $L$, $N=\sum_{i=1}^tAg_i$ the graded
submodule generated by $G$, and let $n\in\NZ$, $G_{\le n}=\{ g_j\in
G~|~d_{\rm gr}(g_j)\le n\}$. If, for each nonzero homogeneous
element $\xi\in N$ with $d_{\rm gr}(\xi )\le n$, there is some
$g_i\in G_{\le n}$ such that $\LM (g_i)|\LM (\xi )$ with respect to
$\prec_{e}$, then we call $G_{\le n}$ an $n$-{\it truncated left
Gr\"obner basis} of $N$ with respect to $(\BE, \prec_{e})$.}\v5

By the definition above,  the lemma below is straightforward. 
{\parindent=0pt\v5

{\bf 4.3.2. Lemma} Let $\G =\{ g_1,\ldots ,g_t\}$ be a homogeneous 
left Gr\"obner basis for the graded submodule $N=\sum_{i=1}^tAg_i$ 
of $L$ with respect to $(\BE ,\prec_{e})$. For each $n\in\NZ$, put 
$\G_{\le n}=\{ g_j\in\G~|~d_{\rm gr}(g_j)\le n\}$, $N_{\le 
n}=\cup_{q=0}^nN_q$ where each $N_q$ is the degree-$q$ homogeneous 
part of $N$, and let $N(n)=\sum_{\xi\in N_{\le n}}A\xi$ be  the 
graded submodule generated by $N_{\le n}$.  The following statements 
hold.\par

(i) $\G_{\le n}$ is an $n$-truncated left Gr\"obner basis of $N$.
Thus, if $\xi\in L$ is a homogeneous element with $d_{\rm gr}(\xi
)\le n$,  then $\xi\in N$ if and only if $\OV{\xi}^{\G_{\le n}}=0$,
i.e., $\xi$ is reduced to zero on division by $\G_{\le n}$,.\par

(ii) $N(n)=\sum_{g_j\in\G_{\le n}}Ag_j$, and $\G_{\le n}$ is an
$n$-truncated left Gr\"obner basis of $N(n)$.\vskip 6pt

{\bf Proof} Exercise.\QED}\v5

In light of Theorem 2.3.3 and  Theorem 2.3.4 (presented in Section 3 
of Chapter 2), an $n$-truncated left Gr\"obner basis is 
characterized as follows.{\parindent=0pt\v5

{\bf 4.3.3. Proposition} Let $N=\sum_{i=0}^sAg_i$ be the graded
submodule of $L$ generated by a set of homogeneous elements $G=\{
g_1,\ldots ,g_m\}$. For each $n\in\NZ$, put  $G_{\le n}=\{ g_j\in
G~|~d_{\rm gr}(g_j)\le n\}$. The following statements are equivalent
with respect to the given $(\BE, \prec_{e})$.\par

(i) $G_{\le n}$ is an $n$-truncated left Gr\"obner basis of $N$.\par

(ii) Every nonzero left S-polynomial $S_{\ell}(g_i,g_j)$ of $d_{\rm
gr}(S_{\ell}(g_i,g_j))\le n$ is reduced to zero on division by
$G_{\le n}$, i.e, $\OV{S_{\ell}(g_i,g_j)}^{G_{\le n}}=0$.\vskip 6pt

{\bf Proof} Recall that if $g_i,g_j\in G$,  $\LT
(g_i)=\lambda_ia^{\alpha}e_t$ with $\alpha =(\alpha_1,\ldots
,\alpha_n)$,  $\LT (g_j)=\lambda_ja^{\beta}e_t$ with $\beta
=(\beta_1,\ldots ,\beta_n)$, and $\gamma =(\gamma _1,\ldots
\gamma_n)$ with $\gamma _i=\max\{\alpha_i,\beta_i\}$, $1\le i\le n$,
then
$$S_{\ell}(g_i,g_j)=\frac{1}{\LC(a^{\gamma -\alpha}g_i)}a^{\gamma -\alpha}g_i-
\frac{1}{\LC(a^{\gamma -\beta}g_j)}a^{\gamma -\beta}g_j$$ is a
homogeneous element in $N$ with $d_{\rm
gr}(S_{\ell}(g_i,g_j))=d(a^{\gamma})+b_t$ by the foregoing property
($\mathbb{P}4$). If $d_{\rm gr}(S_{\ell}(g_i,g_j))\le n$, then it
follows from (i) that (ii) holds.}
\par

Conversely, suppose that (ii) holds. To see that $G_{\le n}$ is an 
$n$-truncated left Gr\"obner basis of $N$, let us run {\bf 
Algorithm-LGB} (presented in Theorem 2.3.4 of Chapter 2) with the 
initial input data $G$. Without optimizing {\bf Algorithm-LGB} we 
may certainly assume that $G\subseteq\G$, thereby $G_{\le 
n}\subseteq\G_{\le n}$ where $\G$ is the new input set returned by 
each pass through the WHILE loop. On the other hand, by the 
construction of $S_{\ell}(g_i,g_j)$ and the property ($\mathbb{P}2$) 
given in the lat section, we know that if $d_{\rm 
gr}(S_{\ell}(g_i,g_j))\le n$, then $d_{\rm gr}(g_i)\le n$, $d_{\rm 
gr}(g_j)\le n$. Hence, the assumption (ii) implies that {\bf 
Algorithm-LGB} does not append any new element of degree $\le n$ to 
$\G$. Therefore, $G_{\le n}=\G_{\le n}$. By Lemma 4.3.2 we conclude 
that $G_{\le n}$ is an $n$-truncated left Gr\"obner basis of 
$N$.\QED {\parindent=0pt\v5

{\bf 4.3.4. Corollary} Let  $N=\sum^m_{i=1}Ag_i$ be the graded
submodule of $L$ generated by a set of homogeneous elements $G=\{
g_1,\ldots ,g_m\}$. Suppose that $G_{\le n}=\{ g_j\in G~|~d_{\rm
gr}(g_j)\le n\}$ is an $n$-truncated left Gr\"obner basis of $N$
with respect to $(\BE, \prec_{e})$. \par

(i) If $\xi\in L$ is a nonzero homogeneous element of $d_{\rm
gr}(\xi )=n$ such that $\LM (g_i){\not |}~\LM (\xi )$ for all
$g_i\in G_{\le n}$, then $G'=G_{\le n}\cup\{\xi\}$ is an
$n$-truncated left Gr\"obner basis for both the graded submodules
$N'=N+A\xi$ and $N''=\sum_{g_j\in G_{\le n}}Ag_j+A\xi$ of $L$.\par

(ii) If $n\le n_1$ and $\xi\in L$ is a nonzero homogeneous element
of $d_{\rm gr}(\xi )=n_1$ such that $\LM (g_i){\not |}~\LM (\xi )$
for all $g_i\in G_{\le n}$, then $G'=G_{\le n}\cup\{\xi\}$ is an
$n_1$-truncated left Gr\"obner basis for the graded submodule
$N'=\sum_{g_j\in G_{\le n}}Ag_j+A\xi$ of $L$. \vskip 6pt

{\bf Proof} If $\xi\in L$ is a nonzero homogeneous element of
$d_{\rm gr}(\xi )=n_1\ge n$ and $\LM (\xi_i){\not |}~\LM (\xi)$ for
all $\xi_i\in G_{\le n}$, then noticing the property mentioned in 
(Lemma 2.1.2(ii) of Chapter 2)  and the property $(\mathbb{P}2)$ 
mentioned in previous Section 2, we see that every nonzero left 
S-polynomial $S_{\ell}(\xi ,\xi_i)$ with $\xi_i\in G$ has $d_{\rm 
gr}(S_{\ell}(\xi ,\xi_i))>n$. Hence both (i) and (ii) hold by 
Proposition 4.3.3.\QED} \v5

Based on the discussion above, the next proposition tells us that 
{\bf Algorithm-LGB} (presented in Theorem 2.3.4 of Chapter 2) can be 
modified for computing $n$-truncated left Gr\"obner bases. 
{\parindent=0pt\v5

{\bf 4.3.5. Proposition} (Compare with ([KR2], Proposition 4.5.10)) 
Given a finite set of nonzero homogeneous elements $ U =\{ 
\xi_1,\ldots ,\xi_m\}\subset L$ with $d_{\rm gr}(\xi_1)\le d_{\rm 
gr}(\xi_2)\le\cdots d_{\rm gr}(\xi_m)$, and a positive integer 
$n_0\ge  d_{\rm gr}(\xi_1)$, the following algorithm computes an 
$n_0$-truncated left Gr\"obner basis $\G =\{ g_1,...,g_t\}$ for the 
graded submodule $N=\sum^m_{i=1}A\xi_i$ such that $d_{\rm 
gr}(g_1)\le d_{\rm gr}(g_2)\le\cdots d_{\rm gr}(g_t)$. 
{\parindent=0pt\vskip 6pt

\underline{\bf Algorithm-TRUNC 
~~~~~~~~~~~~~~~~~~~~~~~~~~~~~~~~~~~~~~~~~~~~~~~~~~~~~~~~~}\vskip 6pt

\textsc{INPUT}: $U = \{ \xi_1,...,\xi_m\}~\hbox{with}~d_{\rm 
gr}(\xi_1)\le d_{\rm gr}(\xi_2)\le\cdots d_{\rm gr}(\xi_m)$\par 
~~~~~~~~~~~~~~~$n_0,~\hbox{where}~n_0\ge d_{\rm gr}(\xi_1)$\par 
\textsc{OUTPUT}: $\G =\{ g_1,...,g_t\}~\hbox{an}~n_0\hbox{-truncated 
left Gr\"obner basis of}~N$\par                                              
\textsc{INITIALIZATION}: $\S_{\le n_0} :=\emptyset ,~W := U,~\G 
:=\emptyset ,~t':=0$\par \textsc{LOOP}\par                              
$n:=\min\{d_{\rm gr}(\xi_i),~d_{\rm 
gr}(S_{\ell}(g_i,g_j))~|~\xi_i\in W ,~S_{\ell}(g_i,g_j)\in\S_{\le 
n_0}\}$\par                                                                      
$\S_n:=\{ S_{\ell}(g_i,g_j)\in\S_{\le n_0}~|~d_{\rm 
gr}(S_{\ell}(g_i,g_j))=n\}$\par                                              
$W_n:=\{\xi_j\in W~|~d_{\rm gr}(\xi_j)=n\}$\par                                                             
$\S_{\le n_0} :=\S_{\le n_0} -\S_n,~W:=W-W_n$\par 
~~~~\textsc{WHILE}~$\S_n\ne\emptyset~\textsc{DO}$\par 
~~~~~~~~~\hbox{Choose any}~$S_{\ell}(g_i,g_j)\in\S_{n}$\par 
~~~~~$\S_n :=\S_n -\{ S_{\ell}(g_i,g_j)\}$\par 
~~~~~~~~~\textsc{IF}~$\OV{S_{\ell}(g_i,g_j)}^{\G}=\eta\ne 
0$~\hbox{with}~$\LM (\eta )= a^{\rho}e_k\textsc{THEN}$\par 
~~~~~~~~~$t':=t'+1,~g_{t'}:=\eta$\par                                                
~~~~~~~~~~~~~~$\S_{\le n_0} :=\S_{\le n_0}\bigcup\left\{ 
S_{\ell}(g_i,g_{t'})\left |\begin{array}{l} g_i\in\G ,1\le i< t',~ 
\LM (g_i )=a^{\tau}e_k,\\ 0<d_{\rm gr}(S_{\ell}(g_i,g_{t'}))\le 
n_0\end{array}\right.\right\}$\par                                              
~~~~~~~~~$\G :=\G\cup\{ g_{t'}\}$\par                                          
~~~~~~~~~\textsc{END}\par                                          
~~~~\textsc{END}\par

~~~~\textsc{WHILE}~$W_n\ne\emptyset~\textsc{DO}$\par
~~~~~~~~~\hbox{Choose any}~$\xi_j\in W_{n}$\par                                     
~~~~$W_n :=W_n -\{ \xi_j\}$\par 
~~~~~~~~~\textsc{IF}~$\OV{\xi_j}^{\G}=\eta\ne 0~\hbox{with}~\LM 
(\eta )= a^{\rho}e_k~\textsc{THEN}$\par                                     
~~~~~~~~~$t':=t'+1,~g_{t'}:=\eta$\par                                                     
~~~~~~~~~$\S_{\le n_0} :=\S_{\le n_0}\cup\left\{ 
S_{\ell}(g_i,g_{t'})\left |\begin{array}{l} g_i\in\G ,1\le i<t',\LM
(g_i )=a^{\tau}e_k,\\
0<d_{\rm gr}(S_{\ell}(g_i,g_{t'}))\le 
n_0\end{array}\right.\right\}$\par                                                  
~~~~~~~~~$\G :=\G\cup\{ g_{t'}\}$\par                                           
~~~~~~~~~\textsc{END}\par                                                       
~~~~\textsc{END}\par                                                                 
\textsc{UNTIL}~$\S_{\le n_0} =\emptyset$\par                                
\textsc{END}\vskip -.2truecm                                                                                     
\underline{
~~~~~~~~~~~~~~~~~~~~~~~~~~~~~~~~~~~~~~~~~~~~~~~~~~~~~~~~~~~~~~~~~~~~~~~~~~~~~~~~ 
~~~~~~~~~~~~~} } \vskip 6pt

{\bf Proof} First note that both the WHILE loops append new elements 
to $\G$ by taking the nonzero normal remainders on division by $\G$. 
Thus, with a fixed $n$, by the definition of a left S-polynomial and 
the normality of $g_{t'}$ (mod $\G$), it is straightforward to check 
that in both the WHILE loops every new appended 
$S_{\ell}(g_i,g_{t'})$ has $d_{\rm gr}(S_{\ell}(g_i,g_{t'}))>n$. To 
proceed, let us write $N(n)$ for the submodule generated by $\G$ 
which is obtained after $W_n$ is exhausted in the second WHILE loop. 
If $n_1$ is the first number after $n$ such that 
$\S_{n_1}\ne\emptyset$, and for some $S_{\ell}(g_i,g_j)\in 
\S_{n_1}$, $\eta =\OV{S_{\ell}(g_i,g_j)}^{\G}\ne 0$ in a certain 
pass through the first WHILE loop, then we note that this $\eta$ is 
still contained in $N(n)$. Hence, after $\S_{n_1}$ is exhausted in 
the first WHILE loop, the obtained $\G$ generates $N(n)$ and $\G$ is 
an $n_1$-truncated left Gr\"obner basis of $N(n)$. Noticing that the 
algorithm starts with $\S =\emptyset$ and $\G =\emptyset$, 
inductively it follows from Proposition 4.3.3 and Corollary 4.3.4 
that after $W_{n_1}$ is exhausted in the second WHILE loop, the 
obtained $\G$ is an $n_1$-truncated left Gr\"obner basis of 
$N(n_1)$. Since $n_0$ is finite and all the generators of $N$ with 
$d_{\rm gr}(\xi_j)\le n_0$ are processed through the second WHILE 
loop, the algorithm terminates and the eventually obtained $\G$ is 
an $n_0$-truncated left Gr\"obner basis of $N$. Finally, the fact 
that the degrees of elements in $\G$ are non-decreasingly ordered 
follows from the choice of the next $n$ in the algorithm. \QED }\v5

Let the data $(A,\B,\prec )$ and $(L, \BE),\prec_{e})$ be as fixed 
before. Combining the foregoing results, we now proceed to show that 
the algorithm given in ([KR], Theorem 4.6.3)) can be adapted for 
computing  minimal homogeneous generating sets of graded submodules 
in free modules over $A$.\par

Let $N$ be a graded submodule of the $\NZ$-graded free $A$-module
$L$ fixed above. We say that a homogeneous generating set $ U$ of
$N$ is a {\it minimal homogeneous generating set} if any proper
subset of $ U$ cannot be a generating set of $N$. As preparatory
result, we first show that the noncommutative analogue of ([KR],
Proposition 4.6.1, Corollary 4.6.2) holds true for
$N$.{\parindent=0pt\v5

{\bf 4.3.6. Proposition} Let $N=\sum^m_{i=1}A\xi_i$ be the graded
submodule of $L$ generated by a set of homogeneous elements $ U=\{
\xi_1,\ldots ,\xi_m\}$, where $d_{\rm gr}(\xi_1)\le d_{\rm
gr}(\xi_2)\le \cdots\le d_{\rm gr}(\xi_m)$. Put $N_1=\{ 0\}$,
$N_i=\sum^{i-1}_{j=1}A\xi_j$, $2\le i\le m$. The following
statements hold.\par

(i) $ U$ is a minimal homogeneous generating set of $N$ if and only
if $\xi_i\not\in N_i$, $1\le i\le m$.\par

(ii) The set $\OV U=\{ \xi_k~|~\xi_k\in U ,~\xi_k\not\in N_k\}$ is a
minimal homogeneous generating set of $N$.\vskip 6pt

{\bf Proof} (i) If $ U$ is a minimal homogeneous generating set of
$N$, then clearly $\xi_i\not\in N_i$, $1\le i\le m$.}\par

Conversely, suppose $\xi_i\not\in N_i$, $1\le i\le m$. If $ U$ is
not a minimal homogeneous generating set of $N$, then, there is some
$i$ such that $N$ is generated by $\{ \xi_1,\ldots
,\xi_{i-1},\xi_{i+1},\ldots ,\xi_m\}$, thereby $\xi_i=\sum_{j\ne
i}h_j\xi_j$ for some nonzero homogeneous elements $h_j\in A$ such
that $d_{\rm gr}(\xi_i)=d_{\rm gr}(h_jg_j)=d_{\rm gr}(h_j)+d_{\rm
gr}(\xi_j)$, where the second equality follows from the foregoing
property $(\mathbb{P}4)$. Thus $d_{\rm gr}(\xi_j)\le d_{\rm
gr}(\xi_i)$ for all $j\ne i$. If $d_{\rm gr}(\xi_j)<d_{\rm
gr}(\xi_i)$ for all $j\ne i$, then $\xi_i\in\sum^{i-1}_{j=1}A\xi_j$,
which contradicts the assumption. If $d_{\rm gr}(\xi_i)=d_{\rm
gr}(\xi_j)$ for some $j\ne i$, then since $h_j\ne 0$ we have $h_j\in
A_0-\{ 0\}=K^*$. Putting $i'=\max\{i,~j~|~f_j\in K^*\}$, we then
have $\xi_{i'}\in\sum^{i'-1}_{j=1}A\xi_j$, which again contradicts
the assumption. Hence, under the assumption we conclude that $ U$ is
a minimal homogeneous generating set of $N$.\par

(ii) In view of (i), it is sufficient to show that $\OV U$ is a
homogeneous generating set of $N$. Indeed, if $\xi_i\in  U-\OV
 U$, then $\xi_i\in\sum^{i-1}_{j=1}A\xi_j$. By checking
$\xi_{i-1}$ and so on, it follows that $\xi_i\in\sum_{\xi_k\in\OV
U}A\xi_k$, as desired.\par\QED{\parindent=0pt\v5

{\bf 4.3.7. Corollary} Let $ U =\{ \xi_1,\ldots ,\xi_m\}$ be a 
minimal homogeneous generating set of a graded submodule $N$ of $L$, 
where $d_{\rm gr}(\xi_1)\le d_{\rm gr}(\xi_2)\le \cdots\le d_{\rm 
gr}(\xi_m)$, and let $\xi\in L-N$ be a homogeneous element with $ 
d_{\rm gr}(\xi_m)\le d_{\rm gr}(\xi )$. Then $\widehat{ U}=
 U\cup\{\xi\}$ is a minimal homogeneous generating set of the
graded submodule $\widehat{N}=N+A\xi$.\par\QED}\v5

We are ready now to reach the following{\parindent=0pt\v5

{\bf 4.3.8. Theorem} (Compare with ([KR2], Theorem 4.6.3)) Let $ U 
=\{ \xi_1,\ldots ,\xi_m\}\subset L$ be a finite set of nonzero 
homogeneous elements of $L$ with $d_{\rm gr}(\xi_1)\le d_{\rm 
gr}(\xi_2)\le\cdots \le d_{\rm gr}(\xi_m)$. Then the algorithm  
presented below returns a minimal homogeneous generating set $ 
U_{\min}=\{\xi_{j_1},\ldots ,\xi_{j_r}\}\subset  U$ for the graded 
submodule $N=\sum^m_{i=1}A\xi_i$; and meanwhile it returns a 
homogeneous left Gr\"obner basis $\G =\{ g_1,...,g_t\}$ for $N$ such 
that $d_{\rm gr}(g_1)\le d_{\rm gr}(g_2)\le\cdots d_{\rm gr}(g_t)$. 
\newpage

\underline{\bf Algorithm-MINHGS 
~~~~~~~~~~~~~~~~~~~~~~~~~~~~~~~~~~~~~~~~~~~~~~~~~~~~~~~}\vskip 6pt

\textsc{INPUT}: $U = \{ \xi_1,...,\xi_m\}~\hbox{with}~d_{\rm 
gr}(\xi_1)\le d_{\rm gr}(\xi_2)\le\cdots \le d_{\rm gr}(\xi_m)$\par          
\textsc{OUTPUT}: $U_{\min}=\{\xi_{j_1},\ldots ,\xi_{j_r}\}\subset
 U
~\hbox{a minimal homogeneous generating}$\par                         
\hskip 6.4truecm\hbox{set of }~$N$;\par                                                          
\hskip 2.1truecm $\G =\{ g_1,...,g_t\}~\hbox{a homogeneous left 
Gr\"obner basis of}~N$\par                                                         
\textsc{INITIALIZATION}: $\S :=\emptyset ,~W := U,~\G :=\emptyset 
,~t':=0,~ U_{\min}:=\emptyset$\par                                             
\textsc{LOOP}\par                                                                  
$n :=\min\{d_{\rm gr}(\xi_i),~d_{\rm 
gr}(S_{\ell}(g_i,g_j))~|~\xi_i\in W,~S_{\ell}(g_i,g_j)\in\S\}$\par             
$\S_n:=\{ S_{\ell}(g_i,g_j)\in\S~|~d_{\rm gr}(S_{\ell}(g_i,g_j))=n\} 
,~ W_n:=\{\xi_j\in W~|~d_{\rm gr}(\xi_j)=n\}$\par                                
$\S :=\S -\S_n,~W:=W-W_n$\par                                                 
~~~~\textsc{WHILE}~$\S_n\ne\emptyset~\textsc{DO}$\par                             
~~~~~~~~~\hbox{Choose any}~$S_{\ell}(g_i,g_j)\in\S_{n}$\par                        
~~~~~$\S_n :=\S_n -\{ S_{\ell}(g_i,g_j)\}$\par                                   
~~~~~~~~~\textsc{IF}~$\OV{S_{\ell}(g_i,g_j)}^{\G}=\eta\ne 
0~\hbox{with}~\LM (\eta )= a^{\rho}e_k~\textsc{THEN}$\par                            
~~~~~~~~~$t':=t'+1,~g_{t'}:=\eta$\par                                          
~~~~~~~~~$\S :=\S\cup\{ S_{\ell}(g_i,g_{t'})\ne 0~|~g_i\in\G ,~1\le 
i< t',~\LM (g_i )=a^{\tau}e_k\}$\par                                           
~~~~~~~~~$\G :=\G\cup\{ g_{t'}\}$\par                                          
~~~~~~~~~\textsc{END}\par                                         
~~~~\textsc{END}\par                                                 
~~~~\textsc{WHILE}~$W_n\ne\emptyset~\hbox{DO}$\par                                   
~~~~~~~~~\hbox{Choose any}~$\xi_j\in W_{n}$\par                                        
~~~~$W_n :=W_n -\{ \xi_j\}$\par                                            
~~~~~~~~~\textsc{IF}~$\OV{\xi_j}^{\G}=\eta\ne 0~\hbox{with}~\LM 
(\eta )= a^{\rho}e_k~\textsc{THEN}$\par                                                
~~~~~~~~~$U_{\min}:= U_{\min}\cup \{\xi_j\}$\par                             
~~~~~~~~~$t':=t'+1,~g_{t'}:=\eta$\par                                          
~~~~~~~~~$\S :=\S\cup\{ S_{\ell}(g_i,g_{t'})\ne 0~|~g_i\in\G,~1\le 
i<t',~\LM (g_i )=a^{\tau}e_k\}$\par                                                        
~~~~~~~~~$\G :=\G\cup\{ g_{t'}\}$\par                                         
~~~~~~~~~\textsc{END}\par                                            
~~~~\textsc{END}\par                                                             
\textsc{UNTIL}~$\S=\emptyset$\par                                                    
\textsc{END}\vskip -.2truecm                                                               
\underline{~~~~~~~~~~~~~~~~~~~~~~~~~~~~~~~~~~~~~~~~~~~~~~~~~~~~~~~~~~~~~~~~~~~~~~~~~~~~~~~~~~~~~~~~~~~~~~~} 
\vskip 6pt

{\bf Proof} Since this algorithm is clearly a variant of {\bf 
Algorithm-LGB}  and {\bf Algorithm-TRUNC} with a minimization 
procedure which works with the finite set $U$, it terminates after a 
certain integer $n$ is executed, and the eventually obtained $\G$ is 
a homogeneous left Gr\"obner basis for $N$ in which the degrees of 
elements are ordered non-decreasingly. It remains to prove that the 
eventually obtained $ U_{\min}$ is a minimal homogeneous generating 
set of $N$.}\par

As in the proof of Proposition 4.3.5, let us first bear in mind that 
for each $n$, in both the WHILE loops every  new appended 
$S_{\ell}(g_i,g_{t'})$ has $d_{\rm gr}(S_{\ell}(g_i,g_{t'}))>n$. 
Moreover, for  convenience, let us write $\G (n)$ for the $\G$ 
obtained after $\S_n$ is exhausted in the first WHILE loop, and 
write $ U_{\min}[n]$, $\G [n]$ respectively  for the $ U_{\min}$, 
$\G$ obtained after $W_n$ is exhausted in the second WHILE loop. 
Since the algorithm starts with ${\cal O}=\emptyset$ and 
$\G=\emptyset$, if, for a fixed $n$, we check carefully how the 
elements of $ U_{\min}$ are chosen during executing the second WHILE 
loop, and how the new elements are appended to $\G$ after each pass 
through the first or the second WHILE loop, then it follows from 
Proposition 4.3.3 and Corollary 4.3.4 that after $W_n$ is exhausted, 
the obtained $ U_{\min}[n]$ and $\G [n]$ generate the same module, 
denoted $N(n)$, such that  $\G [n]$ is an $n$-truncated left 
Gr\"obner basis of $N(n)$.  We now use induction to show that the 
eventually obtained  $ U_{\min}$ is a minimal homogeneous generating 
set for $N$. If $ U_{\min}=\emptyset$, then it is a minimal 
generating set of the zero module.  To proceed,  we assume that $ 
U_{\min}[n]$ is a minimal homogeneous generating set for $N(n)$ 
after $W_n$ is exhausted in the second WHILE loop. Suppose that 
$n_1$ is the first number after $n$ such that 
$\S_{n_1}\ne\emptyset$. We complete the induction proof below by 
showing  that $ U_{\min}[n_1]$ is a minimal homogeneous generating 
set of $N(n_1)$.\par

If in a certain pass through the first WHILE loop, 
$\OV{S_{\ell}(g_i,g_j)}^{\G}=\eta\ne 0$ for some 
$S_{\ell}(g_i,g_j)\in\S_{n_1}$, then we note that $\eta\in N(n)$. It 
follows that after $\S_{n_1}$ is exhausted in the first WHILE loop, 
$\G (n_1)$ generates $N(n)$ and $\G (n_1)$ is an $n_1$-truncated 
left Gr\"obner basis of $N(n)$. Next, assume that $W_{n_1}=\{ 
\xi_{j_1},\ldots ,\xi_{j_s}\}\ne\emptyset$ and that the elements of 
$W_{n_1}$ are processed in the given order during executing the 
second WHILE loop. Since $\G (n_1)$ is an $n_1$-truncated left 
Gr\"obner basis of $N(n)$, if $\xi_{j_1}\in W_{n_1}$ is such that 
$\OV {\xi_{j_1}}^{\G (n_1)}=\eta_1\ne 0$, then $\xi_{j_1}, \eta_1\in 
L -N(n)$. By Corollary 4.3.4, we conclude  that $\G 
(n_1)\cup\{\eta_1\}$ is an $n_1$-truncated Gr\"obner basis for the 
module $N(n)+A\eta_1$; and by Corollary 4.3.7, we conclude that $ 
U_{\min}[n]\cup\{\xi_{j_1}\}$ is a minimal homogeneous generating 
set of $N(n)+A\eta_1$. Repeating this procedure, if $\xi_{j_2}\in 
W_{n_1}$ is such that $\OV{f_{j_2}}^{\G 
(n_1)\cup\{\eta_1\}}=\eta_2\ne 0$, then $\xi_{j_2}, \eta_2\in L 
-(N(n)+A\eta_1)$.  By Corollary 4.3.4, we conclude  that $\G 
(n_1)\cup\{\eta_1,\eta_2\}$ is an $n_1$-truncated left Gr\"obner 
basis for the module  $N(n)+A\eta_1+A\eta_2$; and by Corollary 
4.3.7, we conclude that $ U_{\min}[n]\cup\{\xi_{j_1},\xi_{j_2}\}$ is 
a minimal homogeneous generating set of $N(n)+A\eta_1+A\eta_2$. 
Continuing this procedure until $W_{n_1}$ is exhausted we assert 
that the returned $\G [n_1]=\G$ and $ U_{\min}[n_1]= U_{\min}$ 
generate the same module $N(n_1)$ and $\G [n_1]$ is an 
$n_1$-truncated left Gr\"obner basis of $N(n_1)$ and $ 
U_{\min}[n_1]$ is a minimal homogeneous generating set of $N(n_1)$, 
as desired. As all elements of $ U$ are eventually processed by the 
second WHILE loop, we conclude that the finally obtained $\G$ and $ 
U_{\rm min}$ have the properties: $\G$ generates the module $N$, 
$\G$ is an $n_0$-truncated left Gr\"obner basis of $N$, and $ 
U_{\min}$ is a minimal homogeneous generating set of $N$.\par\QED 
{\parindent=0pt \v5

{\bf Remark.}  If we are only interested in getting a minimal 
homogeneous generating set for the submodule $N$, then {\bf 
Algorithm-MINHGS} can indeed be speed up. More precisely, with  
$$d_{\rm gr}(\xi_1)\le d_{\rm gr}(\xi_2)\le\cdots \le d_{\rm 
gr}(\xi_m)= n_0,$$ it follows from the proof above that if we  stop 
executing the algorithm after $S_{n_0}$ and $W_{n_0}$ are exhausted, 
then the resulted $ U_{\min}[n_0]$ is already the desired minimal 
homogeneous generating set for $N$, while $\G [n_0]$ is an 
$n_0$-truncated left Gr\"obner basis of $N$. \v5

{\bf 4.3.9. Corollary}  Let $ U =\{ \xi_1,\ldots ,\xi_m\}\subset L$ 
be a finite set of nonzero homogeneous elements of $L$ with $d_{\rm 
gr}(\xi_1)=d_{\rm gr}(\xi_2)=\cdots = d_{\rm gr}(\xi_m)=n_0$
\par

(i) If $ U$ satisfies $\LM (\xi_i)\ne\LM (\xi_j)$ for all $i\ne j$, 
then $ U$ is a minimal homogeneous generating set of the graded 
submodule $N=\sum_{i=1}^mA\xi_i$ of $L$, and meanwhile $ U$ is an 
$n_0$-truncated left Gr\"obner basis for $N$.\par

(ii) If $ U$ is a minimal left Gr\"obner basis of the graded 
submodule $N=\sum^m_{i=1}A\xi_i$ (i.e., $ U$ is a left Grobner basis 
of $N$ satisfying $\LM (\xi_i)\ne\LM (\xi_j)$ for all $i\ne j$), 
then  $ U$ is a minimal homogeneous generating set of $N$.\vskip 6pt

{\bf Proof} By the assumption, it follows from the second WHILE loop 
of {\bf Algorithm-MINHGS} that $ U_{\min}= U$.\par \QED}\v5

Let $N$ be an arbitrary  nonzero graded submodule of the 
$\NZ$-graded free $A$-module $L=\oplus_{i=1}^sAe_i$ with $d_{\rm 
gr}(e_i)=b_i$, $1\le i\le s$, and consider  the graded quotient 
module $M=L/N$ (see previous Section 4.2). Our next goal is to 
compute a minimal homogeneous generating set for $M$. \v5

Since $A$ is Noetherian, $N$ is a finitely generated graded 
submodule of $L_0$. Let $N=\sum^m_{j=1}A\xi_j$ be generated by the 
set of nonzero homogeneous elements $ U =\{\xi_1,\ldots ,\xi_m\}$, 
where $\xi_{\ell}=\sum_{k=1}^sf_{k\ell}e_k$ with $f_{k\ell}\in A$, 
$1\le \ell\le m$. Then, every nonzero $f_{k\ell}$ is a homogeneous 
element of $A$ such that $d_{\rm gr}(\xi_{\ell})=d_{\rm 
gr}(f_{k\ell}e_k)=d_{\rm gr}(f_{k\ell})+b_k$, where $b_k=d_{\rm 
gr}(e_k)$, $1\le k\le s$, $1\le \ell\le m$. {\parindent=0pt\v5

{\bf 4.3.10. Lemma} With  every 
$\xi_{\ell}=\sum_{i=1}^sf_{i\ell}e_i$ as fixed above, $1\le \ell\le 
m$, if the $i$-th coefficient $f_{ij}$ of some $\xi_j$ is a nonzero 
constant, say $f_{ij}=1$ without loss of generality, then for each 
$\ell =1,\ldots ,j-1,j+1,\ldots ,m$, the element 
$\xi_{\ell}'=\xi_{\ell}-f_{i\ell}\xi_j$ does not involve $e_i$. 
Putting $ U '=\{ \xi_1',\ldots ,\xi'_{j-1},\xi_{j+1}',\ldots 
,\xi'_m\}$, there is a graded $A$-module isomorphism $M'=L'/N'\cong 
L/N=M$, where $L'=\oplus_{k\ne i}Ae_k$ and $N'=\sum_{\xi_{\ell}'\in 
U '}A\xi_{\ell}'$. \vskip 6pt

{\bf Proof} Since $f_{ij}=1$ by the assumption, we see that every 
$\xi_{\ell}'=\sum_{k\ne i}(f_{k\ell}-f_{i\ell}f_{kj})e_k$ does not 
involve $e_i$. Let $ U '=\{ \xi_1',\ldots 
,\xi'_{j-1},\xi_{j+1}',\ldots ,\xi'_m\}$ and 
$N'=\sum_{\xi_{\ell}'\in U '}A\xi_{\ell}'$. Then $N'\subset 
L'=\oplus_{k\ne i}Ae_k$. Again since $f_{ij}=1$, we have $d_{\rm 
gr}(\xi_j)=d_{\rm gr}(e_i)=b_i$. It follows from the property 
($\mathbb{P}$4) formulated in Subsection 2.1 that
$$\begin{array}{rcl} d_{\rm gr}(f_{i\ell}f_{kj}e_k)&=&d_{\rm gr}(f_{i\ell})+d_{\rm gr}(f_{kj}e_k)\\
&=&d_{\rm gr}(f_{i\ell})+d_{\rm gr}(\xi_j)\\
&=&d_{\rm gr}(f_{i\ell})+b_i\\
&=&d_{\rm gr}(f_{i\ell}e_i)\\
&=&d_{\rm gr}(\xi_{\ell})\\
&=&d_{\rm gr}(f_{k\ell}e_k).\end{array}$$ Noticing that $d_{\rm 
gr}(f_{i\ell}\xi_j)=d_{\rm gr}(f_{i\ell})+d_{\rm gr}(\xi_j)$, this 
shows that in the representation of $\xi_{\ell}'$ every nonzero term 
$(f_{k\ell}-f_{i\ell}f_{kj})e_k$ is a homogeneous element of degree 
$d_{\rm gr}(\xi_{\ell})=d_{\rm gr}(f_{i\ell}\xi_j)$, thereby 
$M'=L'/N'$ is a graded $A$-module. Note that $N=N'+A\xi_j$ and that  
$\xi_j =e_i+\sum_{k\ne i}f_{kj}e_k$. Without making confusion, if we 
use the same notation $\OV e_k$ to denote the coset represented by 
$e_k$ in $M'$ and $M$ respectively,  it is now clear that the 
desired graded $A$-module isomorphism $M'\mapright{\varphi}{}M$ is 
naturally defined by $\varphi (\OV e_k)=\OV e_k$, $k=1,\ldots 
,i-1,i+1,\ldots ,s$. \QED}\v5

Let $M=L/N$ be as fixed above with $N$ generated by the set of 
nonzero homogeneous elements $ U =\{ \xi_1,\ldots ,\xi_m\}$. Then 
since $A$ is $\NZ$-graded with $A_0=K$, it is well known that the 
homogeneous generating set $\OV E=\{ \OV e_1,\ldots ,\OV e_s\}$ of 
$M$ is a minimal homogeneous generating set if and only if 
$\xi_{\ell}=\sum_{k=1}^sf_{k\ell}e_k$ implies $d_{\rm 
gr}(f_{k\ell})>0$ whenever $f_{k\ell}\ne 0$, $1\le\ell\le m$. 
{\parindent=0pt\v5

{\bf 4.3.11. Proposition} (Compare with ([KR2], Proposition 4.7.24)) 
With notation as fixed above, the algorithm presented below returns 
a subset $\{ e_{i_1},\ldots ,e_{i_{s'}}\}\subset\{ e_1,\ldots 
,e_s\}$ and a subset $V=\{ v_1,\ldots ,v_t\}\subset N\cap L'$ such 
that $M\cong L'/N'$ as graded $A$-modules, where 
$L'=\oplus_{q=1}^{s'}Ae_{i_q}$ with $s'\le s$ and 
$N'=\sum^{t}_{k=1}Av_k$, and such that $\{\OV e_{i_1},\ldots ,\OV 
e_{i_{s'}}\}$ is a minimal homogeneous generating set of $M$. \vskip 
6pt

\underline{\bf Algorithm-MINHGSQ 
~~~~~~~~~~~~~~~~~~~~~~~~~~~~~~~~~~~~~~~~~~~~~~~~~~~~~}\vskip 6pt

\textsc{INPUT}: $E=\{ e_1,\ldots ,e_s\};~~ U = \{ 
\xi_1,...,\xi_m\}$\par                                                    
~~~~~~~~~~~~~\hbox{where}~$\xi_{\ell}=\sum_{k=1}^sf_{k\ell}e_k~
\hbox{with homogeneous}~f_{k\ell}\in A, ~1\le \ell\le m$\par                           
\textsc{OUTPUT}: $E' =\{ e_{i_1},\ldots ,e_{i_{s'}}\}; ~~V=\{ 
v_1,\ldots ,v_t\}\subset N\cap L',~\hbox{such that}$\par 
~~~~~~~~~~~~~~~~~$v_{j}=\sum^{s'}_{q=1}h_{qj}e_{i_q}\in 
L'=\oplus_{q=1}^{s'}Ae_{i_q}~ \hbox{with}~h_{qj}\not\in K^*$\par                           
~~~~~~~~~~~~~~~~~\hbox{whenever}~$h_{qj}\ne 0,~1\le j\le t$\par                                                                        
\textsc{INITIALIZATION}: $t :=m;~ V := U ;  ~s':=s; ~E' :=E$\par                             
\textsc{BEGIN}\par                                                                 
~~~~~\textsc{WHILE}~\hbox{there is 
a}~$v_j=\sum^{s'}_{k=1}f_{kj}e_k\in V~\hbox{satisfying}$\par                                       
~~~~~~~~~~~~~~~~~$f_{kj}\not\in 
K^*~\hbox{for}~k<i~\hbox{and}~f_{ij}\in K^*~\hbox{DO}$\par                                             
~~~~~~~~~~~~~~~~~~~~\hbox{for}~$T=\{1,\ldots ,j-1,j+1,\ldots 
,t\}~\hbox{compute}$\par                                                          
~~~~~~~~~~~~~~~~~~~~~$v_{\ell}' 
=v_{\ell}-\frac{1}{f_{ij}}f_{i\ell}v_j,~\ell\in 
T,~r=\#\{\ell~|~\ell\in T,~v_{\ell}' =0\}$\par                                   
~~~~~$t := t-r-1$\par                                                 
~~~~~$V :=\{ v_{\ell}=v_{\ell}'~|~\ell\in T,~v_{\ell}'\ne 0\}$\par                                                                 
~~~~~~~~~$=\{ v_1,\ldots ,v_t\}~(\hbox{after reordered})$\par                     
~~~~~$s' :=s'-1$\par                                                          
~~~~~$E':=E'-\{ e_i\}=\{ e_1,\ldots ,e_{s'}\}~(\hbox{after 
reordered})$\par                                                                
~~~~~\textsc{END}\par                                                              
\textsc{END}\vskip -.2truecm 
\underline{~~~~~~~~~~~~~~~~~~~~~~~~~~~~~~~~~~~~~~~~~~~~~~~~~~~~~~~~~~~~~~~~~~~~~~~~~~~~~~~~~~~~~~~~~~~~~~~~~~} 
}{\parindent=0pt \vskip 6pt

{\bf Proof} It is clear that the algorithm is finite. The 
correctness of the algorithm follows immediately from Lemma 4.3.10 
and the remark we made before the proposition.} \v5

\section*{4.4. Computation of Minimal Finite Graded Free
Resolutions}

Let $A=K[a_1,\ldots ,a_n]=\oplus_{p\in\NZ}A_p$ be an  $\NZ$-graded 
solvable polynomial algebra with respect to a positive-degree 
function $d(~)$, and $(\B,\prec )$ a fixed admissible system of $A$. 
Let $M=\oplus_{q\in\NZ}M_q$, $M'=\oplus_{q\in\NZ}M'_q$ be 
$\NZ$-graded left $A$-modules, and $M~\mapright{\varphi}{}~M'$ an 
$A$-module homomorphism. If $\varphi (M_q)\subseteq M'_q$ for all 
$q\in\NZ$, then $\varphi$ is called a {\it graded  homomorphism}.   
In the literature, such graded homomorphisms are also referred to as 
graded homomorphisms of degree-0 (cf. [NVO]). By the definition it 
is clear that the identity map of $\NZ$-graded $A$-modules is graded 
homomorphism, and compositions of graded homomorphisms are  graded 
homomorphisms.  Thus, all $\NZ$-graded left $A$-modules form a 
subcategory of the category of left $A$-modules, in which morphisms 
are the graded homomorphisms as defined above. Furthermore, if  
$M~\mapright{\varphi}{}~M'$ is a graded homomorphism, then one 
checks that the kernel Ker$\varphi$ of $\varphi$ is a graded 
submodule of $M$, and the image Im$\varphi$ of $\varphi$ is a graded 
submodule of $M'$ (See previous Section 4.2). Consequently, the 
exactness of a  sequence 
$N~\mapright{\varphi}{}~M~\mapright{\psi}{}~M'$ of graded 
homomorphisms in the category of $\NZ$-graded $A$-modules is defined 
as the same as for a sequence of usual $A$-module homomorphisms, 
i.e., the sequence satisfies Im$\varphi =$ Ker$\psi$. Long exact 
sequence in the category  of $\NZ$-graded $A$-modules may be defined 
in an obvious way.\par

Since $A$ is Noetherian and $A_0=K$, it is theoretically well known 
that up to a graded isomorphism of chain complexes in the category 
of graded $A$-modules, every finitely generated graded $A$-module 
$M$ has a unique minimal graded free resolution (cf. [Eis], Chapter 
19; [Kr1], Chapter 3; [Li3]). Based on previously obtained results, 
in this section we establish  the algorithmic procedures for 
constructing  minimal finite graded free resolutions over $A$. All 
notions, notations and conventions used in previous sections are 
maintained. \v5

Given  a finitely generated $\NZ$-graded $A$-module 
$M=\sum^s_{i=1}Av_i$ with the set of nonzero homogeneous generators 
$\{ v_1,\ldots ,v_s\}$ such that $d_{\rm gr}(v_i)=b_i$ for $1\le 
i\le s$,   consider the $\NZ$-graded free $A$-module 
$L_0=\oplus_{i=1}^sAe_i$ with $d_{\rm gr}(e_i)=b_i$, $1\le i\le s$, 
as constructed in Section 4.2. Then, under the $\NZ$-graded 
epimorphism $L_0~\mapright{\varphi_0}{}~M\r 0$ defined by 
$\varphi_0(e_i)=v_i$, $1\le i\le s$, there is a graded isomorphism 
$M\cong L_0/N$, where $N=$ Ker$\varphi_0$.   Thus, we may identify 
$M$ with the graded quotient module $L_0/N$ and write $M=L_0/N$. 
\par

Recall from the literature that a {\it minimal graded free 
resolution} of $M$ is an exact sequence by free $A$-modules and  
$A$-module homomorphisms
$${\cal L}_{\bullet}\quad\quad\cdots~\mapright{\varphi_{i+1}}{}~L_i~\mapright{\varphi_i}{}~\cdots ~
\mapright{\varphi_2}{}~L_1~\mapright{\varphi_1}{}~L_0~\mapright{\varphi_0}{}~M~\mapright{}{}~0$$ 
in which each $L_i$ is an $\NZ$-graded free $A$-module with a finite 
homogeneous $A$-basis $E_i=\{ e_{i_1},\ldots ,e_{i_{s_i}}\}$, and 
each $\varphi_i$ is a graded homomorphism, such 
that{\parindent=1.2truecm\par

\item{(1)} $\varphi_0(E_0)$ is a minimal homogeneous
generating set of $M$, Ker$\varphi_0=N$, and\par

\item{(2)} for $i\ge 1$,  $\varphi_i(E_i)$ is a minimal
homogeneous generating set of 
Ker$\varphi_{i-1}$.\par}{\parindent=0pt\v5

{\bf 4.4.1. Theorem} With notation as fixed above, suppose that 
$N=\sum^m_{i=1}A\xi_i$ with the set of nonzero homogeneous 
generators $ U =\{\xi_1,\ldots ,\xi_m\}$. Then the graded $A$-module 
$M=L_0/N$ has a minimal graded free resolution of length $d\le n$:
$${\cal L}_{\bullet}\quad\quad 0~\mapright{}{}~L_d~\mapright{\varphi_{q}}{}~
\cdots~\mapright{\varphi_2}{}~L_1~\mapright{\varphi_1}{}~L_0~\mapright{\varphi_0}{}~M~\mapright{}{}~0$$ 
which can be constructed by implementing the following 
procedures:}\par

{\bf Procedure 1.}  Run {\bf Algorithm-MINHGSQ} of Proposition 
4.3.11 with the initial input data $E=\{ e_1,\ldots ,e_s\}$ and $ U 
=\{ \xi_1,\ldots ,\xi_m\}$ to compute a subset $E'=\{ e_{i_1},\ldots 
,e_{i_{s'}}\}\subset\{ e_1,\ldots ,e_s\}$ and a subset $V=\{ 
v_1,\ldots ,v_t\}\subset N\cap L_0'$ such that $M\cong L_0'/N'$ as 
graded $A$-modules, where $L_0'=\oplus_{q=1}^{s'}Ae_{i_q}$ with 
$s'\le s$ and $N'=\sum^{t}_{k=1}Av_k$, and such that $\{\OV 
e_{i_1},\ldots ,\OV e_{i_{s'}}\}$ is a minimal homogeneous 
generating set of $M$.\par

For convenience, after accomplishing Procedure 1 we may assume that 
$E=E'$, $ U =V$ and $N=N'$. Accordingly we have the short exact 
sequence $$0~\mapright{}{}~ N~\mapright{}{}~ 
L_0~\mapright{\varphi_0}{}~M~\mapright{}{}  ~0$$ such that 
$\varphi_0 (E)=\{ \OV e_1,\ldots ,\OV e_s\}$  is a minimal 
homogeneous generating set of $M$. \par

{\bf Procedure 2.} Choose a left monomial ordering $\prec_{e}$ on 
the $K$-basis $\BE$ of $L_0$  and run {\bf Algorithm-MINHGS} of 
Theorem 4.3.8 with the initial input data $ U =\{ \xi_1,\ldots 
,\xi_m\}$ to compute a minimal homogeneous generating set $ 
U_{\min}=\{\xi_{j_1},\ldots ,\xi_{j_{s_1}}\}$ and a left Gr\"obner 
basis $\G$ for $N$; at the same time, by keeping track of the 
reductions during executing the first WHILE loop and the second 
WHILE loop respectively, return the matrices ${\cal S}_{r\times t}$ 
and $V_{t\times m}$ required by (Theorem 3.1.2, Chapter 3).
\par

{\bf Procedure 3.} By using the division by the left Gr\"obner basis 
$\G$ obtained in Procedure 2, compute the matrix $U_{m\times t}$ 
required by (Theorem 3.1.2, Chapter 3). Use the matrices ${\cal 
S}_{r\times t}$, $V_{t\times m}$ obtained in Procedure 2, the matrix 
$U_{m\times t}$ and (Theorem 3.1.2, Chapter 3) to compute a 
homogeneous generating set of $N_1=$ Syz$( U_{\min})$ in the 
$\NZ$-graded free $A$-module $L_1=\oplus_{i=1}^{s_1}A\varepsilon_i$, 
where the gradation of $L_1$ is defined by setting $d_{\rm 
gr}(\varepsilon_k)=d_{\rm gr}(\xi_{j_k})$, $1\le k\le s_1$.\par

{\bf Procedure 4.} Construct the exact sequence
$$0~\mapright{}{}~ N_1~\mapright{}{}~L_1~\mapright{\varphi_1}{}~
L_0~\mapright{\varphi_0}{}~M~\mapright{}{}  ~0$$ where 
$\varphi_1(\varepsilon_k)=\xi_{j_k}$, $1\le k\le s_1$.\par

If $N_1\ne 0$, then repeat Procedure 2 -- Procedure 4 for $N_1$ and 
so on.\par

Noticing that $A$ is $\NZ$-graded with the degree-0 homogeneous part 
$A_0=K$, $A$ is Noetherian,  and that every finitely generated 
$A$-module $M$ has finite projective dimension p.dim$_AM\le n$ by  
(Theorem 3.3.1, Chapter 3), thereby $A$ is an $\NZ$-graded local 
ring of finite global homological dimension. It follows from the 
literature ([Eis], Chapter 19; [Kr1], Chapter3;  [Li3]) that
$$\hbox{p.dim}_AM=\hbox{~ the length of a minimal graded free resolution
of}~M.$$  Hence, the desired minimal finite graded free resolution 
${\cal L}_{\bullet}$ for $M$ is then obtained after finite times of 
processing the above procedures.\newpage\setcounter{page}{89}

\chapter*{5. Minimal Finite Filtered\\ \hskip 1.25truecm Free Resolutions}
\markboth{\rm Minimal Filtered Free Resolutions}{\rm Minimal 
Filtered Free Resolutions} \vskip 2.5truecm\def\T#1{\widetilde #1}

Recall that the $\NZ$-filtered solvable polynomial algebras with 
$\NZ$-filtration determined by the natural length of elements from 
the PBW $K$-basis  (especially the quadric solvable polynomial 
algebras are such $\NZ$-filtered algebras) were studied in ([LW], 
[Li1]). In this chapter, after specifying the $\NZ$-filtered 
structure   determined by a positive-degree function $d(~)$ on a 
solvable polynomial algebra $A=K[a_1,\ldots ,a_n]$, we specify the 
corresponding $\NZ$-filtered structure for free $A$-modules. Then, 
we introduce minimal filtered free resolutions for finitely 
generated modules over such $\NZ$-filtered algebra $A$ by 
introducing minimal F-bases and minimal standard bases for 
$A$-modules and their submodules with respect to good filtration;  
we show that  any two minimal F-bases, respectively any two minimal 
standard bases, have the same number of elements and the same number 
of elements of the same filtered-degree, minimal filtered free 
resolutions are unique up to strict filtered isomorphism of chain 
complexes in the category of filtered $A$-modules, and that minimal 
filtered free resolutions can be algorithmically computed in case 
$A$ has a graded monomial ordering $\prec_{gr}$.\par

Since the standard bases we are going to introduce in terms of good 
filtration are generalization of classical Macaulay bases (see a 
remark given in Subsection 3.3), while a classical Macaulay basis 
$V$ is characterized in  terms of both the leading homogeneous 
elements (degree forms) of $V$ and the homogenized elements of $V$ 
(cf. [KR2], P.38, P.55), accordingly, on the basis of Chapter 4,  
our main idea in reaching the goal of this chapter is to use the 
filtered-graded transfer strategy (as proposed in [Li1]) by 
employing both the associated graded algebra (module) and the Rees 
algebra (module)  of an $\NZ$-filtered solvable polynomial algebra 
(of a filtered module).  \v5

All notions, notations and conventions introduced in previous 
chapters  are maintained.\v5

\section*{5.1. $\NZ$-Filtered Solvable Polynomial Algebras}\par

Comparing with the general theory of $\mathbb{Z}$-filtered rings 
[LVO], in this section we formulate the  $\NZ$-filtered structure of 
solvable polynomial algebras determined by means of positive-degree 
functions. \v5

Let $A=K[a_1,\ldots a_n]$ be a solvable polynomial algebra with 
admissible system $(\B ,\prec )$, where $\B =\{ 
a^{\alpha}=a_1^{\alpha_1}\cdots a_n^{\alpha_n}~|~\alpha 
=(\alpha_1,\ldots ,\alpha_n)\in\NZ^n\}$ is the PBW $K$-basis of $A$ 
and $\prec$ is a monomial ordering on $\B$, and let $d(~)$ be a 
positive-degree function on $A$ such that $d(a_i)=m_i>0$, $1\le i\le 
n$ (see Section 1.1 of Chapter 1). Put
$$F_pA=K\hbox{-span}\{ a^{\alpha}\in\B~|~d(a^{\alpha})\le p\} ,\quad p\in\NZ,$$
then  it is clear that  $F_pA\subseteq F_{p+1}A$ for all $p\in\NZ$, 
$A=\cup_{p\in\NZ}F_pA$, and $1\in F_0A=K$.{\parindent=0pt\v5

{\bf 5.1.1. Definition} With notation as above, if 
$F_pAF_qA\subseteq F_{p+q}A$ holds for all $p,q\in\NZ$, then  we 
call $A$ an {\it $\NZ$-filtered solvable polynomial algebra with 
respect to the positive-degree function $d(~)$}, and accordingly we 
call  $FA=\{ F_pA\}_{p\in\NZ}$ the {\it $\NZ$-filtration of $A$ 
determined by $d(~)$}. }\v5

Note that the $\NZ$-filtration $FA$ constructed above is clearly 
{\it separated} in the sense that if $f$ is a {\it nonzero} element 
of $L$, then either $f\in F_0A=K$ or $f\in F_pA-F_{p-1}A$ for some 
$p>0$. Thus, we may define the {\it filtered-degree} (abbreviated to 
{\it fil-degree}) of a nonzero $f\in A$, denoted $d_{\rm fil}(f)$, 
as follows
$$d_{\rm fil}(f)=\left\{\begin{array}{ll} 0,&\hbox{if}~f\in F_0A=K,\\
p,&\hbox{if}~p\in F_pA-F_{p-1}A~\hbox{for 
some}~p>0.\end{array}\right.$$  

Bearing in mind the definition of $d_{\rm fil}(f)$, the following 
featured  property of $FA$ will very much help us to deal with the  
associated graded structures of $A$ and filtered $A$-modules. 
{\parindent=0pt\v5

{\bf 5.1.2. Lemma}  If $f=\sum_i\lambda_ia^{\alpha (i)}\in A$ with 
$\lambda_i\in K^*$ and $a^{\alpha (i)}\in\B$, then $d_{\rm 
fil}(f)=p$ if and only if $d(a^{\alpha (i')})=p$ for some $i'$ if 
and only if $d(f)=p=d_{\rm fil}(f)$, where $d(~)$ is the given 
positive-degree function on $A$. \vskip 6pt

{\bf Proof} Exercise.\QED}\v5

Given a solvable polynomial algebra $A=K[a_1,\ldots ,a_n]$  and a 
positive-degree function $d(~)$ on $A$, it follows from Definition 
1.1.3 (Section 1, Chapter 1), Definition 5.1.1 and Lemma 5.1.2 above 
that the next proposition is clear now.{\parindent=0pt\v5

{\bf 5.1.3. Proposition} $A$ is an $\NZ$-filtered solvable 
polynomial algebra with respect to $d(~)$ if,  for $1\le i<j\le n$, 
all the relations   $a_ja_i=\lambda_{ji}a_ia_j+f_{ji}$ with 
$f_{ji}=\sum\mu_ka^{\alpha (k)}$ presented in (S2$'$) (Section 4.1, 
Chapter 4), satisfy $d(a^{\alpha (k)})\le d(a_ia_j)$ whenever 
$\mu_k\ne 0$. \par\QED}\v5

With the proposition presented above,  the following examples may be 
better understood.{\parindent=0pt\v5

{\bf Example} (1) If $A=K[a_1,\ldots ,a_n]$ is an $\NZ$-graded 
solvable polynomial algebra with respect to a positive-degree 
function $d(~)$, i.e., $A=\oplus_{p\in\NZ}A_p$ with the degree-$p$ 
homogeneous part $A_p=K$-span$\{a^{\alpha}\in\B~|~d(a^{\alpha})=p\}$ 
(see Section 4.1 of Chapter 4), then, with respect to the same 
positive-degree function $d(~)$ on $A$,  $A$ is turned into an 
$\NZ$-filtered solvable polynomial algebra with the $\NZ$-filtration 
$FA=\{F_pA\}_{p\in\NZ}$ where each $F_pA=\oplus_{q\le p}A_q$. \vskip 
6pt

{\bf Example} (2) Let $A=K[a_1,\ldots a_n]$ be a solvable polynomial 
algebra with the admissible system $(\B ,\prec_{gr})$, where 
$\prec_{gr}$ is a graded monomial ordering on $\B$ with respect to a 
given positive-degree function $d(~)$ on $A$ (see the definition of 
$\prec_{gr}$ given in Section 1.1 of Chapter 1). Then by referring 
to (Definition 1.1.3 of Chapter 1) and the above proposition, one 
easily sees that $A$ is an $\NZ$-filtered solvable polynomial 
algebra with respect to the same $d(~)$. In the case where 
$\prec_{gr}$ respects $d(a_i)=1$ for $1\le i\le n$, (Definition 
1.1.3 of Chapter 1) entails that the generators of $A$ satisfy
$$\begin{array}{rcl} a_ja_i&=&\lambda_{ji}a_ia_j+\sum \lambda^{ji}_{k\ell}a_ka_{\ell}+\sum \lambda^{ji}_ta_t+
\mu_{ji},\\ &{~}&\hbox{where}~~1\le i<j\le n,~\lambda_{ji}\in 
K^*,~\lambda^{ji}_{k\ell},\lambda^{ji}_t, \mu_{ji}\in 
K.\end{array}$$ In [Li1] such $\NZ$-filtered solvable polynomial 
algebras are referred to as {\it quadric solvable polynomial 
algebras} which include numerous significant algebras such as Weyl 
algebras and enveloping algebras of finite dimensional  Lie 
algebras. One is referred to [Li1] for some detailed study on 
quadric solvable polynomial algebras.}\vskip 6pt

The next example  provides $\NZ$-filtered solvable polynomial 
algebras in which some generators may have degree $\ge 
2$.{\parindent=0pt\v5

{\bf Example} (3)  Consider the solvable polynomial algebra 
$K$-algebra $A=K[a_1,a_2,a_3]$ constructed in (Example (1) of 
Section 1.4, Chapter 1). Then, one checks that $A$ is turned into an 
$\NZ$-filtered solvable polynomial algebra with respect to the 
lexicographic ordering $a_3\prec_{lex}a_2\prec_{lex}a_1$ and the 
degree function $d(~)$ such that $d(a_1)=2$, $d(a_2)=1$, $d(a_3)=4$. 
Moreover, one may also check that with respect to the same degree 
function $d(~)$, the graded lexicographic ordering 
$a_3\prec_{grlex}a_2\prec_{grlex}a_1$ is another choice to make $A$ 
into an $\NZ$-filtered solvable polynomial algebra.}\v5

Let $A$ be an  $\NZ$-filtered solvable polynomial algebra with 
respect to a given positive-degree function $d(~)$, and let 
$FA=\{F_pA\}_{p\in\NZ}$ be the $\NZ$-filtration of $A$ determined by 
$d(~)$. Then $A$ has the associated $\NZ$-graded $K$-algebra 
$G(A)=\oplus_{p\in\NZ}G(A)_p$ with $G(A)_0=F_0A=K$ and 
$G(A)_p=F_pA/F_{p-1}A$ for $p\ge 1$, where for $\OV f=f+F_{p-1}A\in 
G(A)_p$, $\OV g=g+F_{q-1}A$, the multiplication is given by $\OV 
f\OV g=fg+F_{p+q-1}A\in G(A)_{p+q}$. Another $\NZ$-graded 
$K$-algebra determined by $FA$ is the Rees algebra $\T A$ of $A$, 
which is defined as $\T A=\oplus_{p\in\NZ}\T A_p$ with $\T 
A_p=F_pA$, where the multiplication of $\T A$ is induced by 
$F_pAF_qA\subseteq F_{p+q}A$, $p, q\in\NZ$. For convenience, we fix 
the following notations once for all:{\parindent=.5truecm\par

\item{$\bullet$} If $h\in G(A)_p$ and $h\ne 0$, then
we write $d_{\rm gr}(h)$ for the gr-degree of $h$ as a homogeneous 
element of $G(A)$, i.e., $d_{\rm gr}(h)=p$.\par

\item{$\bullet$} If $H\in \T A_p$ and $H\ne 0$, then
we write $d_{\rm gr}(H)$ for the gr-degree of $H$ as a homogeneous 
element of $\T A$, i.e., $d_{\rm gr}(H)=p$.\v5}

Concerning the $\NZ$-graded structure of $G(A)$, if $f\in A$ with 
$d_{\rm fil}(f)=p$, then by Lemma 5.1.2,  the coset $f+F_{p-1}A$ 
represented by $f$ in $G(A)_p$ is a nonzero homogeneous element of 
degree $p$. If we denote this homogeneous element by $\sigma (f)$ 
(in the literature it is referred to as the principal symbol of 
$f$), then $d_{\rm fil}(f)=p=d_{\rm gr}(\sigma (f))$. However, 
considering the Rees algebra $\T A$ of $A$, we note that a nonzero 
$f\in F_qA$ represents a homogeneous element of degree $q$ in $\T 
A_q$ on one hand, and on the other hand it represents a homogeneous 
element of degree $q_1$ in $\T A_{q_1}$, where $q_1=d_{\rm 
fil}(f)\le q$. So, for a nonzero $f\in F_pA$, we denote the 
corresponding homogeneous element of degree $p$ in $\T A_p$ by 
$h_p(f)$, while we use $\T f$ to denote the homogeneous element 
represented by $f$ in $\T A_{p_1}$ with $p_1=d_{\rm fil}(f)\le p$. 
Thus,  $d_{\rm gr}(\T f)=d_{\rm fil}(f)$, and we see that $h_p(f)=\T 
f$ if and only if $d_{\rm fil}(f)=p$.
\par

Furthermore, if we write $Z$ for the homogeneous element of degree 1 
in $\T A_1$ represented by the multiplicative identity element 1, 
then $Z$ is a central regular element of $\T A$, i.e., $Z$ is not a 
divisor of zero and is contained in the center of $\T A$. Bringing 
this homogeneous element $Z$ into play,  the homogeneous elements of 
$\T A$ are featured as follows:{\parindent=1truecm\vskip6pt

\item{$\bullet$} If $f\in A$ with $d_{\rm fil}(f)=p_1$ then for all $p\ge p_1$,
$h_p(f)=Z^{p-p_1}\T f$. In other words, if $H\in\T A_p$ is a nonzero 
homogeneous element of degree $p$, then there is some $f\in F_pA$ 
such that $H=Z^{p-d_{\rm fil}(f)}\T f=\T f+(Z^{p-d_{\rm 
fil}(f)}-1)\T f$. 
\par}{\parindent=0pt\vskip 6pt

It follows that by sending $H$ to $f+F_{p-1}A$ and sending $H$ to 
$f$ respectively, $G(A)\cong\T A/\langle Z\rangle$ as $\NZ$-graded 
$K$-algebras and  $A\cong \T A/\langle 1-Z\rangle$ as $K$-algebras 
(cf. [LVO]).} \v5

Since a solvable polynomial algebra $A$ is necessarily a domain 
(Proposition 1.1.4, Chapter 1), we summarize two useful properties 
concerning the multiplication of $G(A)$ and $\T A$ respectively  
into the following lemma. {\parindent=0pt\v5

{\bf 5.1.4. Lemma} With notation as before, let $f,g$ be nonzero 
elements of $A$ with $d_{\rm fil}(f)=p_1$, $d_{\rm fil}(g)=p_2$. 
Then\par

(i)  $\sigma (f)\sigma (g)=\sigma (fg)$;\par

(ii)  $\T f\T g=\widetilde{fg}$. If $p_1+p_2\le p$, then 
$h_p(fg)=Z^{p-p_1-p_2}\T f\T g$.\vskip 6pt

{\bf Proof} Exercise. \QED}\v5

The results given in the next theorem,  which are analogues of those 
concerning quadric solvable polynomial algebras in ([LW], Section 3; 
[Li1], CH.IV), may be derived in a similar way as in loc. cit. (With 
the preparation made above, one is also invited to give the detailed 
proof as an exercise). {\parindent=0pt\v5

{\bf 5.1.5. Theorem} Let $A=K[a_1,\ldots ,a_n]$ be a solvable 
polynomial algebra with the admissible system $(\B ,\prec_{gr})$, 
where $\prec_{gr}$ is a graded monomial ordering on $\B$ with 
respect to a given positive-degree function $d(~)$ on $A$, thereby 
$A$ is an $\NZ$-filtered solvable polynomial algebra with respect to 
the same $d(~)$ by the foregoing Example (2), and let $FA=\{ 
F_pA\}_{p\in\NZ}$ be the corresponding $\NZ$-filtration of $A$. 
Considering the associated graded algebra $G(A)$ as well as the Rees 
algebra $\T A$ of $A$, the following statements hold.\par

(i) $G(A)=K[\sigma (a_1),\ldots ,\sigma (a_n)]$,  $G(A)$ has the PBW 
$K$-basis
$$\sigma (\B)=\{
\sigma (a)^{\alpha}=\sigma (a_1)^{\alpha_1}\cdots \sigma 
(a_n)^{\alpha_n}~|~\alpha =(\alpha_1,\ldots ,\alpha_n)\in\NZ^n\} ,$$ 
and, by referring to (Definition 1.1.3, Chapter 1), for  $\sigma 
(a)^{\alpha}$, $\sigma (a)^{\beta}\in\sigma (\B )$ such that 
$a^{\alpha}a^{\beta}=\lambda_{\alpha ,\beta}a^{\alpha 
+\beta}+f_{\alpha ,\beta}$, where $\lambda_{\alpha ,\beta}\in K^*$, 
if $f_{\alpha ,\beta}=0$ then
$$\sigma (a)^{\alpha}\sigma(a)^{\beta}=
\lambda_{\alpha ,\beta}\sigma (a)^{\alpha 
+\beta},~\hbox{where}~\sigma (a)^{\alpha +\beta}=\sigma 
(a_1)^{\alpha_1+\beta_1}\cdots \sigma (a_n)^{\alpha_n+\beta_n};$$ 
and in the case where $f_{\alpha ,\beta}=\sum_j\mu^{\alpha 
,\beta}_ja^{\alpha (j)}\ne 0$ with $\mu^{\alpha ,\beta}_j\in K$,
$$\sigma (a)^{\alpha}\sigma(a)^{\beta}= \lambda_{\alpha ,\beta}\sigma
(a)^{\alpha +\beta}+\displaystyle{\sum_{d(a^{\alpha 
(k)})=d(a^{\alpha +\beta})}}\mu^{\alpha ,\beta}_j\sigma (a)^{\alpha 
(k)}.$$ Moreover, the ordering $\prec_{_{G(A)}}$ defined on $\sigma 
(\B )$ subject to the rule:
$$\sigma (a)^{\alpha}\prec_{_{G(A)}} \sigma (a)^{\beta}\Longleftrightarrow a^{\alpha}\prec_{gr}a^{\beta},
\quad a^{\alpha},a^{\beta}\in\B,$$  is a graded monomial ordering 
with respect to the positive-degree function $d(~)$ on $G(A)$ such 
that $d(\sigma (a_i))=d(a_i)$ for $1\le i\le n$, that turns $G(A)$ 
into an $\NZ$-graded  solvable polynomial algebra.
\par

(ii) $\T A=K[\T a_1,\ldots ,\T a_n, Z]$ where $Z$ is the central 
regular element of degree 1 in $\T A_1$ represented by 1,  $\T A$ 
has the PBW $K$-basis $$\T{\B}=\{\T a^{\alpha}Z^m=\T 
a_1^{\alpha_1}\cdots \T a_n^{\alpha_n}Z^m~|~\alpha =(\alpha_1,\ldots 
,\alpha_n )\in\NZ^n,m\in\NZ\} ,$$ and, by referring to (Definition 
1.1.3, Chapter 1), for  $\T a^{\alpha}Z^s$, $\cdot\T 
a^{\beta}Z^t\in\T{\B}$ such that 
$a^{\alpha}a^{\beta}=\lambda_{\alpha ,\beta}a^{\alpha 
+\beta}+f_{\alpha ,\beta}$, where $\lambda_{\alpha ,\beta}\in K^*$, 
if $f_{\alpha ,\beta}=0$ then
$$\T a^{\alpha}Z^s\cdot\T a^{\beta}Z^t=
\lambda_{\alpha ,\beta}\T a^{\alpha +\beta}Z^{s+t},~\hbox{where}~\T 
a^{\alpha +\beta}=\T a_1^{\alpha_1+\beta_1}\cdots \T 
a_n^{\alpha_n+\beta_n};$$ and in the case where $f_{\alpha 
,\beta}=\sum_j\mu^{\alpha ,\beta}_ja^{\alpha (j)}\ne 0$ with 
$\mu^{\alpha ,\beta}_j\in K$,
$$\begin{array}{rcl} \T a^{\alpha}Z^s\cdot\T a^{\beta}Z^t&=&
\lambda_{\alpha ,\beta}\T a^{\alpha +\beta}Z^{s+t}+\sum_j\mu^{\alpha
,\beta}_j\T a^{\alpha (j)}Z^{q-m_j},\\
&{~}&\hbox{where}~q=d(a^{\alpha +\beta})+s+t,~m_j=d(a^{\alpha 
(j)}).\end{array}$$ Moreover, the ordering $\prec_{_{\T A}}$ defined 
on $\T{\B}$ subject to the rule:
$$\T a^{\alpha}Z^s\prec_{_{\T A}}\T a^{\beta}Z^t\Longleftrightarrow a^{\alpha}\prec_{gr}a^{\beta},
\hbox{or}~a^{\alpha}=a^{\beta}~\hbox{and}~s<t,\quad 
a^{\alpha},a^{\beta}\in\B,$$ is a monomial ordering on $\T{\B}$ 
(which is not necessarily a graded monomial ordering), that turns 
$\T A$ into an $\NZ$-graded solvable polynomial algebra with respect 
to the positive-degree function $d(~)$ on $\T A$ such that $d(Z)=1$ 
and $d(\T{a_i})=d (a_i)$ for $1\le i\le n$.\par\QED }\v5

The corollary presented below  will be very often used in discussing 
left Gr\"obner bases and standard bases for submodules of filtered 
free $A$-modules and their associated graded free $G(A)$-modules as 
well as the graded free $\T A$-modules (Section 5.3, Section 5.4). 
{\parindent=0pt \v5

{\bf 5.1.6. Corollary}  With the assumption and notations as in 
Theorem 3.1.4, if $f=\lambda a^{\alpha}+\sum_j\mu_ja^{\alpha (j)}$ 
with $d(f)=p$ and $\LM (f)=a^{\alpha}$, then $p=d_{\rm 
fil}(f)=d_{\rm gr}(\sigma (f))=d_{\rm gr}(\T f)$, and
$$\begin{array}{l} \sigma (f)=\lambda\sigma (a)^{\alpha}+\sum_{d(a^{\alpha (j_k)})=p}\mu_{j_k}\sigma (a)^{\alpha (j_k)};\\
\LM(\sigma (f))=\sigma (a)^{\alpha}=\sigma (\LM (f));\\
\T f=\lambda\T a^{\alpha}+\sum_j\mu_j\T a^{\alpha
(j)}Z^{p-d(a^{\alpha (j)})};\\
\LM (\T f)=\T a^{\alpha}=\widetilde{\LM (f)},\end{array}$$ where 
$\LM (f)$, $\LM (\sigma (f))$ and $\LM (\T f)$ are taken with 
respect to $\prec_{gr}$, $\prec_{_{G(A)}}$ and $\prec_{_{\T A}}$ 
respectively.\vskip 6pt

{\bf Proof} By referring to Lemma 5.1.2 and  Lemma 5.1.4, this is a 
straightforward exercise.\QED}\v5

\section*{5.2. $\NZ$-Filtered Free Modules}\par

Let $A=K[a_1,\ldots ,a_n]$ be an  $\NZ$-filtered solvable polynomial 
algebra with the filtration $FA=\{F_pA\}_{p\in\NZ}$ determined by a 
positive-degree function $d(~)$ on $A$, and let $(\B ,\prec )$ be a 
fixed admissible system of $A$. Consider a free $A$-module 
$L=\oplus_{i=1}^sAe_i$ with the $A$-basis $\{ e_1,\ldots ,e_s\}$. 
Then $L$ has the $K$-basis $\BE =\{ a^{\alpha}e_i~|~a^{\alpha}\in\B 
,~1\le i\le s\}$. If  $\{ b_1,\ldots ,b_s\}$ is an {\it arbitrarily} 
fixed subset of $\NZ$, then, with $FL=\{ F_qL\}_{q\in\NZ}$ defined 
by putting
$$F_qL=\{ 0\}~\hbox{if}~q<\min\{ b_1,\ldots ,b_s\};~\hbox{otherwise}~F_qL=
\sum^s_{i=1}\left (\sum_{p_i+b_i\le q}F_{p_i}A\right )e_i,$$ or 
alternatively, for $q\ge\min\{ b_1,\ldots, b_s\}$,
$$F_qL=K\hbox{-span}\{
a^{\alpha}e_i\in\BE~|~d(a^{\alpha})+b_i\le q\},$$ $L$ forms an 
$\NZ$-{\it filtered free $A$-module} with respect to the 
$\NZ$-filtered structure of $A$, that is, every $F_qL$ is a 
$K$-subspace of $L$, $F_qL\subseteq F_{q+1}L$ for all $q\in\NZ$, 
$L=\cup_{q\in\NZ}F_qL$, $F_pAF_{q}L\subseteq F_{p+q}L$ for all 
$p,q\in\NZ$, and for each $i=1,\ldots ,s$,
$$e_i\in F_0L~\hbox{if}~b_i=0;~\hbox{otherwise}~
e_i\in F_{b_i}L-F_{b_i-1}L.$$ {\parindent=0pt\par

{\bf Convention.} Let $A$ be an $\NZ$-filtered solvable polynomial 
algebra with respect to a positive-degree function $d(~)$. Unless 
otherwise stated, from now on in the subsequent sections  if we say 
that $L=\oplus_{i=1}^sAe_i$ is a filtered free $A$-module with the 
filtration $FL=\{ F_qL\}_{q\in\NZ}$, then $FL$ is always meant the 
type as constructed above.}\v5

Let $L=\oplus_{i=1}^sAe_i$ be a filtered free $A$-module with the 
filtration $FL=\{ F_qL\}_{q\in\NZ}$, which is constructed with 
respect to a given subset $\{ b_1,\ldots ,b_s\}\subset\NZ$. Then 
$FL$ is {\it separated} in the sense that if $\xi$ is a {\it 
nonzero} element of $L$, then either $\xi\in F_0L$ or $\xi\in 
F_qL-F_{q-1}L$ for some $q>0$. Thus, to make the discussion on $FL$ 
consistent with that on $FA$ in Section 5.1, we define the {\it 
filtered-degree} (abbreviated to {\it fil-degree}) of a nonzero 
$\xi\in L$, denoted $d_{\rm fil}(\xi )$, as follows
$$d_{\rm fil}(\xi )=\left\{\begin{array}{ll} 0,&\hbox{if}~\xi\in F_0L,\\
q,&\hbox{if}~\xi\in F_qL-F_{q-1}L~\hbox{for 
some}~q>0.\end{array}\right.$$                                                     
For instance, we have $d_{\rm fil}(e_i)=b_i$, $1\le i\le s$. 
Comparing with Lemma 5.1.2 we first note the 
following{\parindent=0pt\v5

{\bf 5.2.1. Lemma} Let $\xi\in L-\{ 0\}$. Then $d_{\rm fil}(\xi )=q$ 
if and only if $\xi =\sum_{i,j}\lambda_{ij}a^{\alpha (i_j)}e_j$, 
where $\lambda_{ij}\in K^*$ and $a^{\alpha (i_j)}\in\B$ with $\alpha 
(i_j)=(\alpha_{i_{j1}},\ldots ,\alpha_{i_{jn}})\in\NZ^n$, in which 
some monomial $a^{\alpha (i_j)}e_j$ satisfy $d(a^{\alpha 
(i_j)})+b_j=q$.\vskip 6pt

{\bf Proof} Exercise.\QED}\v5

Let $L=\oplus_{i=1}^sAe_i$ be a filtered free $A$-module with the 
filtration $FL=\{ F_qL\}_{q\in\NZ}$ such that $d_{\rm 
fil}(e_i)=b_i$, $1\le i\le s$. Considering the the associated 
$\NZ$-graded algebra $G(A)$ of $A$, the filtered free $A$ module $L$ 
has  the {\it associated $\NZ$-graded $G(A)$-module} 
$G(L)=\oplus_{q\in\NZ}G(L)_q$ with $G(L)_q=F_qL/F_{q-1}L$, where for 
$\OV f=f+F_{p-1}A\in G(A)_p$, $\OV \xi=\xi+F_{q-1}L\in G(L)_q$, the 
module action is given by $\OV f\cdot \OV \xi =f\xi+F_{p+q-1}L\in 
G(L)_{p+q}$. As with homogeneous elements in $G(A)$, if $h\in 
G(L)_q$ and $h\ne 0$, then  we write $d_{\rm gr}(h)$ for the degree 
of $h$ as a homogeneous element of $G(L)$, i.e., $d_{\rm gr}(h)=q$. 
If $\xi\in L$ with $d_{\rm fil}(\xi )=q$, then the coset 
$\xi+F_{q-1}L$  represented by $\xi$ in $G(L)_q$ is a nonzero 
homogeneous element of degree $q$, and if we denote this homogeneous 
element by $\sigma (\xi )$ (in the literature it is referred to as 
the principal symbol of $\xi$) then $d_{\rm gr}(\sigma (\xi 
))=q=d_{\rm fil}(\xi )$.\par

Furthermore, considering the Rees algebra $\T A$ of $A$, the 
filtration $FL=\{ F_qL\}_{q\in\NZ}$ of $L$ also defines the {\it 
Rees module} $\T L$ of $L$, which is the $\NZ$-graded $\T A$-module 
$\T L=\oplus_{q\in\NZ}\T L_q$, where $\T L_q=F_qL$ and the module 
action is induced by $F_pAF_qL\subseteq F_{p+q}L$. As with 
homogeneous elements in $\T A$, if $H\in \T L_q$ and $H\ne 0$, then 
we write $d_{\rm gr}(H)$ for the degree of $H$ as a homogeneous 
element of $\T L$, i.e., $d_{\rm gr}(H)=q$. Note that any nonzero 
$\xi\in F_qL$ represents a homogeneous element of degree $q$ in $\T 
L_q$ on one hand, and on the other hand it represents a homogeneous 
element of degree $q_1$ in $\T L_{q_1}$, where $q_1=d_{\rm fil}(\xi 
)\le q$. So, for  a nonzero $\xi\in F_qL$ we  denote the 
corresponding homogeneous element of degree $q$ in $\T L_q$ by 
$h_q(\xi )$,  while we use $\T{\xi}$ to denote the homogeneous 
element  represented by $\xi$ in $\T L_{q_1}$ with $q_1=d_{\rm 
fil}(\xi )\le q$. Thus,  $d_{\rm gr}(\T{\xi})=d_{\rm fil}(\xi)$, and 
we see that $h_q(\xi )=\T{\xi}$ if and only if $d_{\rm fil}(\xi)=q$.
\par

We also note that if  $Z$ denotes the homogeneous element of degree 
1 in $\T A_1$ represented by the multiplicative identity element 1, 
then, similar to the discussion given in the last section, one 
checks that there are $A$-module isomorphism  $L\cong \T L/(1-Z)\T 
L$ and graded $G(A)$-module isomorphism $G(L)\cong \T L/Z\T L$. \par

For  $f\in A$, $\xi\in L$, the next  lemma records several 
convenient facts about $d_{\rm fil}(f\xi )$, $\sigma (f\xi )$, 
$h_{\ell}(f\xi )$ and $\widetilde{f\xi}$, respectively.   
{\parindent=0pt\v5

{\bf 5.2.2. Lemma} With notation as above, the following statements 
hold.\par

(i)  $d_{\rm fil}(f\xi )=d(f)+d_{\rm fil}(\xi )$ holds for all 
nonzero $f\in A$ and nonzero $\xi\in L$.\par

(ii) $\sigma (f)\sigma (\xi )=\sigma (f\xi )$ holds for all nonzero 
$f\in A$ and nonzero $\xi\in L$.\par

(iii) If $\xi\in L$ with $d_{\rm fil}(\xi )=q\le \ell$, then 
$h_{\ell}(\xi )=Z^{\ell -q}\T \xi$. Furthermore, let $f\in A$ with 
$d_{\rm fil}(f)=p$, $\xi\in L$ with $d_{\rm fil}(\xi )=q$. Then $\T 
f\T \xi =\widetilde{f\xi}$; if $p+q\le \ell$, then $h_{\ell}(f\xi 
)=Z^{\ell -p-q}\T f\T \xi$. \vskip 6pt

{\bf Proof} Since $A$ is a solvable polynomial algebra,  $G(A)$ and 
$\T A$ are $\NZ$-graded solvable polynomial algebras by Theorem 
5.1.5, thereby they are necessarily domains (Proposition 1.1.4(ii), 
Chapter 1).  Noticing the definition of $d_{\rm fil}(f)$ and $d_{\rm 
fil}(\xi )$, by Lemma 5.2.1, the verification of (i) -- (iii) are 
then straightforward.\QED\v5

{\bf 5.2.3. Proposition}  With notation as fixed before, let 
$L=\oplus_{i=1}^sAe_i$ be a filtered free $A$-module with the 
filtration $FL=\{ F_qL\}_{q\in\NZ}$ such that $d_{\rm 
fil}(e_i)=b_i$, $1\le i\le s$. The following two statements 
hold.\par

(i) $G(L)$ is an $\NZ$-graded free $G(A)$-module with the 
homogeneous $G(A)$-basis $\{ \sigma (e_1),\ldots ,\sigma (e_s)\}$, 
that is,  $G(L)=\oplus_{i=1}^sG(A)\sigma 
(e_i)=\oplus_{q\in\NZ}G(L)_q$ with
$$G(L)_q=\sum_{p_i+b_i=q}G(A)_{p_i}\sigma (e_i)\quad q\in\NZ.
$$
Moreover, $\sigma (\B (e))=\{\sigma (a^{\alpha}e_i)=\sigma 
(a)^{\alpha}\sigma (e_i)~|~a^{\alpha}e_i\in\B (e)\}$ forms a 
$K$-basis for $G(L)$.\par

(ii) $\T L$ is an $\NZ$-graded free $\T A$-module with the 
homogeneous  $\T A$-basis $\{ \T e_1,\ldots ,\T e_s\}$, that is, $\T 
L=\oplus_{i=1}^s\T A\T e_i=\oplus_{q\in\NZ}\T L_q$ with
$$\T L_q=\sum_{p_i+b_i=q}\T A_{p_i}\T e_i, \quad q\in\NZ .$$
Moreover, $\widetilde{\B (e)}=\{\T{a}^{\alpha}Z^m\T e_i~|~\T 
a^{\alpha}Z^m\in\T{\B},~1\le i\le s\}$ forms a $K$-basis for $\T L$, 
where $\T{\B}$ is the PBW $K$-basis of $\T A$ determined in Theorem 
5.1.5(ii).\vskip 6pt

{\bf Proof} Since $d_{\rm fil}(e_i)=b_i$, $1\le i\le s$, if  $\xi 
=\sum^s_{i=1}f_ie_i\in F_qL=\sum^s_{i=1}\left (\sum_{p_i+b_i\le 
q}F_{p_i}A\right )e_i$, then $d_{\rm fil}(\xi )\le q$. By Lemma 
5.2.2,
$$\begin{array}{l} \sigma (\xi )=\sum_{d(f_i)+b_i=q}\sigma (f_i)\sigma (e_i)
\in\sum_{i=1}^sG(A)_{q-b_i}\sigma (e_i)\\
\\
h_q(\xi )=\sum^s_{i=1}Z^{q-d(f_i)-b_i}\T f_i\T e_i\in\sum^s_{i=1}\T 
A_{q-b_i}\T e_i.\end{array}$$ This shows that $\{ \sigma 
(e_1),\ldots ,\sigma (e_s)\}$ and $\{ \T e_1,\ldots ,\T e_s\}$ 
generate the $G(A)$-module $G(L)$ and the $\T A$-module $\T L$, 
respectively. Next, since each $\sigma (e_i)$ is a homogeneous 
element of degree $b_i$, if a degree-$q$ homogeneous element 
$\sum_{i=1}^s\sigma (f_i)\sigma (e_i)=0$, where $f_i\in A$, $d_{\rm 
fil}(f_i)+b_i=q$, $1\le i\le s$, then $\sum_{i=1}^sf_ie_i\in 
F_{q-1}L$ and hence each $f_i\in F_{q-1-b_i}A$ by Lemma 5.2.1, a 
contradiction. It follows that $\{ \sigma (e_1),\ldots ,\sigma 
(e_s)\}$ is linearly independent over $G(A)$. Concerning the linear 
independence of $\{ \T e_1,\ldots ,\T e_s\}$ over $\T A$, since each 
$\T e_i$ is a homogeneous element of degree $b_i$, if a degree-$q$ 
homogeneous element $\sum_{i=1}^sh_{p_i}(f_i)\T e_i=0$, where 
$f_i\in F_{p_i}A$ and $p_i+b_i=q$, $1\le i\le s$, then 
$\sum_{i=1}^sf_ie_i=0$ in $F_qL$ and consequently all $f_i=0$, 
thereby $h_{p_i}(f_i)=0$ as desired. Finally, if $\xi\in F_qL$ with 
$d_{\rm fil}(\xi )=q$, then by Lemma 5.2.1, $\xi 
=\sum_{i,j}\lambda_{ij}a^{\alpha (i_j)}e_j$ with $\lambda_{ij}\in 
K^*$ and $d(a^{\alpha (i_j)})+b_j=\ell_{ij}\le q$. It follows from 
Lemma 5.2.2 that
$$\begin{array}{l} \sigma (\xi )=\sum_{\ell_{ik}=q}
\lambda_{ik}\sigma (a)^{\alpha (i_k)}\sigma (e_k),\\
\T{\xi}=\sum_{i,j}\lambda_{ij}Z^{q-\ell_{ij}}\T a^{\alpha (i_j)}\T 
e_j.\end{array}$$ Therefore, a further application of Lemma 5.2.1  
and Lemma 5.2.2 shows that $\sigma (\B (e))$ and $\T{\BE}$ are 
$K$-bases for $G(L)$ and $\T L$ respectively.} \v5

\section*{5.3. Filtered-Graded Transfer of Gr\"obner Bases for 
Modules}\par\def\PRCEGR{\prec_{e\hbox{\rm -}gr}} 

Throughout this section, we let $A=K[a_1,\ldots ,a_n]$ be a solvable 
polynomial algebra with the admissible system $(\B ,\prec_{gr})$, 
where $\B=\{ a^{\alpha}=a_1^{\alpha_1}\cdots a_n^{\alpha_n}~|~\alpha 
=(\alpha_1,\ldots ,\alpha_n)\in\NZ^n\}$ is the PBW $K$-basis of $A$ 
and $\prec_{gr}$ is a graded monomial ordering  with respect to some 
given positive-degree function $d(~)$ on $A$ (see Section  1). 
Thereby $A$ is turned into an $\NZ$-filtered solvable polynomial 
algebra with  the filtration $FA=\{ F_pA\}_{p\in\NZ}$ constructed 
with respect to the same $d(~)$ (see Example (2) of Section 5.1). In 
order to compute minimal standard bases by employing both 
inhomogeneous and homogenous left Gr\"obner bases in the subsequent  
Section 5.5, our aim of the current section is to establish the 
relations between left Gr\"obner bases in a filtered free (left)  
$A$-module $L$ and homogeneous left Gr\"obner bases in $G(L)$ as 
well as homogeneous left Gr\"obner bases in $\T L$, which are just 
module theory analogues of the results on filtered-graded transfer 
of Gr\"obner bases for left ideals given in ([LW], [Li1]).  All 
notions, notations and conventions introduced in previous sections  
are maintained.\v5

Let $L=\oplus_{i=1}^sAe_i$ be a filtered free  $A$-module with the 
filtration $FL=\{ F_qL\}_{q\in\NZ}$ such that $d_{\rm 
fil}(e_i)=b_i$, $1\le i\le s$. Bearing in mind Lemma 5.2.1, we say 
that a left monomial ordering on $\BE$ is a {\it graded left 
monomial ordering}, denoted by $\PRCEGR$, if for 
$a^{\alpha}e_i,a^{\beta}e_j\in\BE$,
$$a^{\alpha}e_i\PRCEGR a^{\beta}e_j~\hbox{implies}~d_{\rm fil}(a^{\alpha}e_i)=
d(a^{\alpha})+b_i\le d(a^{\beta}) +b_j=d_{\rm fil}(a^{\beta}e_j).$$ 
For instance,  with respect to the given graded monomial ordering 
$\prec_{gr}$ on $\B$ and  the $\NZ$-filtration $FA$ of $A$,   if $\{ 
f_1,\ldots ,f_s\}\subset A$ is a finite subset such that 
$d(f_i)=b_i=d_{\rm fil}(e_i)$, $1\le i\le s$, then it is 
straightforward to check that the Schreyer ordering  (see Section 
2.1 of Chapter 2) induced by $\{ f_1,\ldots ,f_s\}$ subject to the 
rule: for $a^{\alpha}e_i,a^{\beta}e_j\in\BE$,
$$a^{\alpha}e_i\prec_{s\hbox{-}gr} a^{\beta}e_j\Longleftrightarrow\left\{\begin{array}{l}
\LM (a^{\alpha}f_i)\prec_{gr}\LM (a^{\beta}f_j),\\
\hbox{or}\\
\LM (a^{\alpha}f_i)=\LM (a^{\beta}f_j)~\hbox{and}~i<j.\\
\end{array}\right.$$
is a graded left monomial ordering on $\BE$. \par

More generally, let $\{ \xi_1,\ldots ,\xi_m\}\subset L$ be a finite 
subset, where $d_{\rm fil}(\xi_i)=q_i$, $1\le i\le m$, and let 
$L_1=\oplus_{i=1}^mA\varepsilon_i$ be the filtered free $A$-module 
with the filtration $FL_1=\{ F_qL_1\}_{q\in \NZ}$ such that  $d_{\rm 
fil}(\varepsilon_i)=q_i$, $1\le i\le m$.  Then, given {\it any} 
graded left monomial ordering $\PRCEGR$ on $\BE$, the Schreyer 
ordering $\prec_{s\hbox{-}gr}$ defined on the $K$-basis $\B 
(\varepsilon )=\{ a^{\alpha}\varepsilon_i~|~a^{\alpha}\in\B ,~1\le 
i\le m\}$ of $L_1$ subject to the rule: for 
$a^{\alpha}\varepsilon_i,a^{\beta}\varepsilon_j\in\B (\varepsilon 
)$,
$$a^{\alpha}\varepsilon_i\prec_{s\hbox{-}gr} a^{\beta}\varepsilon_j
\Longleftrightarrow\left\{\begin{array}{l}
\LM (a^{\alpha}\xi_i)\PRCEGR\LM (a^{\beta}\xi_j),\\
\hbox{or}\\
\LM (a^{\alpha}\xi_i)=\LM (a^{\beta}\xi_j)~\hbox{and}~i<j,\\
\end{array}\right.$$
is a graded left monomial ordering on $\B (\varepsilon )$.\par

Comparing with Lemma 5.1.2 and Lemma 5.2.1, the lemma given below 
reveals the intrinsic property of a graded left monomial ordering 
employed by a filtered free $A$-module.{\parindent=0pt\v5

\def\PRCVE{\prec_{\varepsilon\hbox{-}gr}}\def\BV{\B (\varepsilon )}

{\bf 5.3.1. Lemma}  Let $L=\oplus_{i=1}^sAe_i$ be a filtered free 
$A$-module with the filtration $FL=\{ F_qL\}_{q\in\NZ}$ such that 
$d_{\rm fil}(e_i)=b_i$, $1\le i\le s$, and let $\prec_{e\hbox{-}gr}$ 
be a graded left monomial ordering on $\BE$. Then $\PRCEGR$ is 
compatible with the filtration $FL$ of $L$ in the sense that $\xi\in 
F_qL-F_{q-1}L$, i.e. $d_{\rm fil}(\xi )=q$, if and only if $\LM (\xi 
)=a^{\alpha}e_i$ with $d_{\rm 
fil}(a^{\alpha}e_i)=d(a^{\alpha})+b_i=q$. \vskip 6pt

{\bf Proof} Let $\xi =\sum_{i,j}\lambda_{ij}a^{\alpha (i_j)}e_j\in 
F_qL-F_{q-1}L$. Then  by Lemma 5.2.1, there is some $a^{\alpha 
(i_{\ell})}e_{\ell}$ such that $d(a^{\alpha 
(i_{\ell})})+b_{\ell}=q$. If $\LM (\xi )=a^{\alpha (i_t)}e_t$ with 
respect to $\PRCEGR$, then $a^{\alpha (i_k)}e_k\PRCEGR a^{\alpha 
(i_t)}e_t$ for all $a^{\alpha (i_k)}e_k$ with $k\ne t$. If $\ell 
=t$, then $d(a^{\alpha (i_t)})+b_t=q$; otherwise, since $\PRCEGR$ is 
a graded left monomial ordering, we have $d(a^{\alpha (i_k)})+b_k\le 
d(a^{\alpha (i_t)})+b_t$, in particular, $q=d(a^{\alpha 
(i_{\ell})})+b_{\ell}\le d(a^{\alpha (i_t)})+b_t\le q$. Hence 
$d_{\rm fil}(a^{\alpha (i_t)}e_t)=d(a^{\alpha (i_t)})+b_t=q$, as 
desired.}\par

Conversely, for $\xi =\sum_{i,j}\lambda_{ij}a^{\alpha (i_j)}e_j\in 
L$, if,  with respect to $\PRCEGR$, $\LM (\xi )=a^{\alpha (i_t)}e_t$ 
with $d_{\rm fil}(a^{\alpha (i_t)}e_t)=d(a^{\alpha (i_t) })+b_t=q$, 
then $a^{\alpha (i_k)}e_k\PRCEGR a^{\alpha (i_t)}e_t$ for all $k\ne 
t$. Since $\PRCEGR$ is a graded left monomial ordering, we have 
$d(a^{\alpha (i_k)})+b_k\le d(a^{\alpha (i_t)})+b_t=q$. It follows 
from Lemma 5.2.1 that $d_{\rm fil}(\xi )=q$, i.e.,  $\xi\in 
F_qL-F_{q-1}L$. \QED \v5

Let  $L=\oplus_{i=1}^sAe_i$ be a filtered free $A$-module with the 
filtration $FL=\{ F_qL\}_{q\in\NZ}$ such that $d_{\rm 
fil}(e_i)=b_i$, $1\le i\le s$.  Then, by Proposition 5.2.3 we know 
that the associated graded $G(A)$-module $G(L)$ of $L$ is an 
$\NZ$-graded free module, i.e., $G(L)=\oplus_{i=1}^sG(A)\sigma 
(e_i)$ with the homogeneous $G(A)$-basis $\{ \sigma (e_1),\ldots 
,\sigma (e_s)\}$, and that $G(L)$ has the  $K$-basis $\sigma(\B 
(e))=\{ \sigma (a^{\alpha}e_i)=\sigma (a)^{\alpha}\sigma 
(e_i)~|~a^{\alpha}e_i\in\B (e)\}$. Furthermore, let 
$\prec_{e\hbox{-}gr}$ be a graded left monomial ordering on $\BE$ as 
defined in the beginning of this section. Then we may define an 
ordering $\prec_{\sigma (e)\hbox{-}gr}$ on $\sigma(\B (e))$ subject 
to the rule:
$$\sigma (a)^{\alpha}\sigma (e_i)\prec_{\sigma (e)\hbox{-}gr}\sigma (a)^{\beta}\sigma (e_j)
\Longleftrightarrow 
a^{\alpha}e_i\prec_{e\hbox{-}gr}a^{\beta}e_j,\quad 
a^{\alpha}e_i,a^{\beta}e_j\in\B (e).$$ {\parindent=0pt\v5
\def\PRCEGR{\prec_{e\hbox{-}gr}}

{\bf 5.3.2. Lemma}  With the ordering $\prec_{\sigma (e)\hbox{-}gr}$ 
defined above, the following statements hold.\par

(i) $\prec_{\sigma (e)\hbox{-}gr}$ is a graded left monomial 
ordering on $\sigma(\B (e))$.\par

(ii) (Compare with Corollary 5.1.6.) $\LM (\sigma (\xi ))=\sigma 
(\LM (\xi ))$ holds for all nonzero $\xi\in L$, , where the monomial 
orderings used for taking $\LM (\sigma (\xi ))$ and $\LM (\xi )$ are 
$\prec_{\sigma (e)\hbox{-}gr}$ and $\PRCEGR$ respectively. \vskip 
6pt

{\bf Proof} (i) Noticing that the given monomial ordering 
$\prec_{gr}$ on $A$ is a graded monomial ordering with respect to a 
positive-degree function $d(~)$ on $A$, it follows from  Theorem 
5.1.5(i) that  $G(A)$ is turned into an $\NZ$-graded solvable 
polynomial algebra by using the graded monomial ordering 
$\prec_{_{G(A)}}$ defined on $\sigma (\B )$ subject to the rule: 
$\sigma (a)^{\alpha}\prec_{_{G(A)}}\sigma (a)^{\beta}$ 
$\Longleftrightarrow$ $a^{\alpha}\prec_{gr}a^{\beta}$, where the 
positive-degree function on $G(A)$ is given by $d(\sigma 
(a_i))=d(a_i)$, $1\le i\le n$. Moreover, since $\sigma (e_i)$ is a 
homogeneous element of degree $b_i$ in $G(L)$, $1\le i\le s$, by 
Lemma 5.2.2, it is then straightforward to verify that 
$\prec_{\sigma (e)\hbox{-}gr}$ is a graded left monomial ordering on 
$\sigma(\B (e))$.}\par

(ii) Let $\xi =\sum_{i,j}\lambda_{ij}a^{\alpha (i_j)}e_j$, where 
$\lambda_{ij}\in K^*$ and $a^{\alpha (i_j)}\in\B$ with $\alpha 
(i_j)=(\alpha_{i_{j1}},\ldots ,\alpha_{i_{jn}})\in\NZ^n$. If $d_{\rm 
fil}(\xi )=q$, i.e., $\xi\in F_qL-F_{q-1}L$, then by Lemma 5.3.1, 
$\LM (\xi )=a^{\alpha (i_t)}e_t$ for some $t$ such that $d_{\rm 
fil}(a^{\alpha (i_t)}e_t)=d(a^{\alpha (i_t)})+b_t=q$. Since 
$\prec_{e\hbox{-}gr}$ is a left graded monomial ordering on $\BE$, 
by Lemma 5.2.2 we have  $\sigma (\xi )=\lambda_{it}\sigma 
(a)^{\alpha (i_t)}\sigma (e_{t})+\sum_{d(a^{\alpha 
(i_k)})+b_k=q}\lambda_{ik}\sigma (a)^{\alpha (i_k)}\sigma (e_k)$. It 
follows from the definition of $\prec_{\sigma (e)\hbox{-}gr}$ that 
$\LM (\sigma (\xi ))=\sigma (a)^{\alpha (i_t)}\sigma (e_t)=\sigma 
(\LM (\xi ))$, as desired.\QED{\parindent=0pt\v5

{\bf 5.3.3. Theorem}  Let $N$ be a submodule of the filtered free 
$A$-module $L=\oplus_{i=1}^sAe_i$, where $L$ is equipped with the 
filtration $FL=\{ F_qL\}_{q\in\NZ}$ such that $d_{\rm 
fil}(e_i)=b_i$, $1\le i\le s$, and let $\prec_{e\hbox{-}gr}$ be a 
graded left monomial ordering on $\BE$. For a subset $\G =\{ 
g_1,\ldots ,g_m\}$ of $N$, the following two statements are 
equivalent.\par

(i) $\G$ is a left Gr\"obner basis of $N$ with respect to 
$\PRCEGR$.\par

(ii) Putting $\sigma (\G )=\{ \sigma (g_1),\ldots ,\sigma (g_m)\}$ 
and considering the filtration $FN=\{ F_qN=F_qL\cap N\}_{q\in\NZ}$ 
of $N$ induced by $FL$, $\sigma (\G )$ is a left Gr\"obner basis for 
the associated graded $G(A)$-module $G(N)$  of $N$ with respect to 
the graded left monomial ordering $\prec_{\sigma (e)\hbox{-}gr}$ 
defined above. \vskip 6pt

{\bf Proof} (i) $\Rightarrow$ (ii) Note that any nonzero homogeneous 
element of $G(N)$ is of the form $\sigma (\xi )$ with $\xi \in N$. 
If $\G$ is a left Gr\"obner basis of $N$, then there exists some 
$g_i\in\G$ such that $\LM (g_i)|\LM (\xi )$, i.e., there is a 
monomial $a^{\alpha}\in\B$ such that $\LM (\xi )=\LM (a^{\alpha}\LM 
(g_i))$. Since the given left monomial ordering $\PRCEGR$ on $\BE$ 
is a graded left monomial ordering, it follows from Lemma 5.2.2 and 
Lemma 5.3.2 that
$$\begin{array}{rcl} \LM (\sigma (\xi )))&=&\sigma (\LM (\xi ))\\
&=&\sigma (\LM (a^{\alpha}\LM (g_i)))\\
&=&\LM (\sigma (a^{\alpha}\LM (g_i)))\\
&=&\LM (\sigma (a)^{\alpha}\sigma (\LM (g_i)))\\
&=&\LM (\sigma (a)^{\alpha}\LM (\sigma (g_i))).\end{array}$$ This 
shows that $\LM (\sigma (g_i))|\LM (\sigma (\xi ) )$, thereby 
$\sigma (\G )$ is a left Gr\"obner basis for $G(N)$.}\par

(ii) $\Rightarrow$ (i) Suppose that $\sigma (\G )$ is a left 
Gr\"obner basis of $G(N)$ with respect to $\prec_{\sigma 
(e)\hbox{-}gr}$. If $\xi\in N$ and $\xi\ne 0$, then $\sigma (\xi 
)\ne 0$, and there exists a $\sigma (g_i)\in\sigma (\G )$ such that 
$\LM (\sigma (g_i))|\LM (\sigma (\xi ))$, i.e., there is a monomial 
$\sigma (a)^{\alpha}\in\sigma (\B )$ such that $\LM (\sigma (\xi 
))=\LM (\sigma (a)^{\alpha}\LM (\sigma (g_i)))$. Again as $\PRCEGR$ 
is a left graded monomial ordering on $\BE$, by Lemma 5.2.2 and 
Lemma 5.3.2 we have
$$\begin{array}{rcl} \sigma (\LM (\xi ))&=&\LM (\sigma (\xi ))\\
&=&\LM (\sigma (a)^{\alpha}\LM (\sigma (g_i)))\\
&=&\LM (\sigma (a)^{\alpha}\sigma (\LM (g_i)))\\
&=&\LM (\sigma (a^{\alpha}\LM (g_i)))\\
&=&\sigma (\LM (a^{\alpha}\LM (g_i))).\end{array}$$ This shows that 
$d_{\rm fil}(\LM (\xi ))=d_{\rm fil}(\LM (a^{\alpha}\LM (g_i)))$. 
Since both $\LM (\xi )$ and $\LM (a^{\alpha}\LM (g_i))$ are 
monomials in $\BE$, it follows from  the construction of $FL$ and 
Lemma 5.3.1  that $\LM (\xi )=\LM (a^{\alpha}\LM (g_i))$, i.e., $\LM 
(g_i)|\LM (\xi )$. This shows that $\G $ is a left Gr\"obner basis 
for $N$.\QED\v5

Similarly, in light of Proposition 5.2.3 we may define an ordering 
$\prec_{\widetilde{e}}$ on the $K$-basis $\widetilde{\B (e)}=\{ 
Z^m\T a^{\alpha}\T e_i~|~Z^m\T a^{\alpha}\in\T\B ,~1\le i\le s\}$ of 
the $\NZ$-graded free $\T A$-module $\T L=\oplus_{i=1}^s\T A\T e_i$ 
subject to the rule: for $Z^s\T a^{\alpha}\T e_i,Z^t\T a^{\beta}\T 
e_j\in\widetilde{\B (e)}$,
$$Z^s\T a^{\alpha}\T e_i\prec_{\widetilde{e}}Z^t\T a^{\beta}\T e_j
\Longleftrightarrow 
a^{\alpha}e_i\prec_{e\hbox{-}gr}a^{\beta}e_j,~\hbox{or}~a^{\alpha}e_i=a^{\beta}e_j~\hbox{and}~s<t,
$$ where $\prec_{e\hbox{-}gr}$ is a given graded left monomial ordering on $\B (e)$. {\parindent=0pt\v5
\def\PRCEGR{\prec_{e\hbox{-}gr}}

{\bf 5.3.4. Lemma}  With the ordering $\prec_{\widetilde{e}}$ 
defined above, the following statements hold.\par

(i)  $\prec_{\widetilde{e}}$  is a  left monomial ordering on 
$\widetilde{\B (e)}$ (but not necessarily a graded left monomial 
ordering).\par

(ii) (Compare with Corollary 5.1.6.) $\LM (\T\xi )=\widetilde{\LM 
(\xi )}$ holds for all nonzero $\xi\in L$, , where the monomial 
orderings used for taking $\LM (\T\xi )$ and $\LM (\xi )$ are 
$\prec_{\widetilde{e}}$ and $\PRCEGR$ respectively. \vskip 6pt

{\bf Proof} (i) Noticing that the given monomial ordering 
$\prec_{gr}$ for $A$ is a graded monomial ordering with respect to a 
positive-degree function $d(~)$ on $A$, it follows from  Theorem 
5.1.5(ii) that  $\T A$ is turned into an $\NZ$-graded solvable 
polynomial algebra by using the  monomial ordering $\prec_{_{\T A}}$ 
defined on $\T\B$ subject to the rule: $$\T a^{\alpha}Z^s\prec_{_{\T 
A}}\T a^{\beta}Z^t\Longleftrightarrow a^{\alpha}\prec_{gr}a^{\beta}, 
\hbox{or}~a^{\alpha}=a^{\beta}~\hbox{and}~s<t,\quad 
a^{\alpha},a^{\beta}\in\B,$$ where the positive-degree function on 
$\T A$ is given by  $d(\T{a_i})=d (a_i)$ for $1\le i\le n$, and 
$d(Z)=1$. Moreover, since $\T e_i$ is a homogeneous element of 
degree $b_i$ in $\T A$, $1\le i\le s$, by Lemma 5.2.2, it is then 
straightforward to verify that $\prec_{\widetilde{e}}$ is a left 
monomial ordering on $\widetilde{B (e)}$.}\par

(ii) Let $\xi =\sum_{i,j}\lambda_{ij}a^{\alpha (i_j)}e_j$, where 
$\lambda_{ij}\in K^*$ and $a^{\alpha (i_j)}\in\B$ with $\alpha 
(i_j)=(\alpha_{i_{j1}},\ldots ,\alpha_{i_{jn}})\in\NZ^n$. If $d_{\rm 
fil}(\xi )=q$, i.e., $\xi\in F_qL-F_{q-1}L$, then by Lemma 5.3.1, 
$\LM (\xi )=a^{\alpha (i_t)}e_t$ for some $t$ such that $d_{\rm 
fil}(a^{\alpha (i_t)}e_t)=d(a^{\alpha (i_t)})+b_t=q$. Since 
$\prec_{e\hbox{-}gr}$ is a left graded monomial ordering on $\BE$, 
by Lemma 5.2.2 we have $\T\xi =\lambda_{it}\T a^{\alpha (i_t)}\T 
e_t+\sum_{j\ne t}\lambda_{ij}Z^{q-\ell_{ij}}\T a^{\alpha (i_j)}\T 
e_j$, where $\ell_{ij}=d_{\rm fil}(a^{\alpha (i_j)}e_j)=d(a^{\alpha 
(i_j)})+d_j$. It follows from the definition of 
$\prec_{\widetilde{e}}$  that $\LM (\T\xi )=\T a^{\alpha (i_t)}\T 
e_t=\widetilde{\LM (\xi )}$, as desired.\QED{\parindent=0pt\v5

{\bf 5.3.5. Theorem}  Let $N$ be a submodule of the filtered free 
$A$-module $L=\oplus_{i=1}^sAe_i$, where $L$ is equipped with the 
filtration $FL=\{ F_qL\}_{q\in\NZ}$ such that $d_{\rm 
fil}(e_i)=b_i$, $1\le i\le s$, and let $\prec_{e\hbox{-}gr}$ be a 
graded left monomial ordering on $\BE$. For a subset $\G =\{ 
g_1,\ldots ,g_m\}$ of $N$, the following two statements are 
equivalent.\par

(i) $\G$ is a left Gr\"obner basis of $N$ with respect to 
$\PRCEGR$.\par

(ii) Putting $\tau (\G )=\{ \T g_1,\ldots ,\T g_m\}$ and considering 
the filtration $FN=\{ F_qN=F_qL\cap N\}_{q\in\NZ}$ of $N$ induced by 
$FL$, $\tau (\G )$ is a left Gr\"obner basis for the Rees module $\T 
N$  of $N$ with respect to the  left monomial ordering 
$\prec_{\widetilde{e}}$ defined before Lemma 5.3.4. \vskip 6pt

{\bf Proof} (i) $\Rightarrow$ (ii) Note that any nonzero homogeneous 
element of $\T N$ is of the form $h_q(\xi )$ for some $\xi\in F_qN$ 
with $d_{\rm fil}(\xi )=q_1\le q$. By Lemma 5.2.2, $h_q(\xi 
)=Z^{q-q_1}\T \xi$. If $\G$ is a left Gr\"obner basis of $N$, then 
there exists some $g_i\in\G$ such that $\LM (g_i)|\LM (\xi )$, i.e., 
there is a monomial $a^{\alpha}\in\B$ such that $\LM (\xi )=\LM 
(a^{\alpha}\LM (g_i))$. It follows from Lemma 5.2.2 and Lemma 5.3.4 
that
$$\begin{array}{rcl} \LM (\T\xi )&=&\widetilde{\LM (\xi )}\\
&=&(\LM (a^{\alpha}\LM (g_i)))\widetilde{~}\\
&=&\LM ((a^{\alpha}\LM (g_i))\widetilde{~})\\
&=&\LM (\T a^{\alpha}\widetilde{\LM (g_i)})\\
&=&\LM (\T a^{\alpha}\LM (\T g_i)).\end{array}$$ Hence, noticing the 
definition of $\prec_{\widetilde{e}}$ we have
$$\begin{array}{rcl} \LM (h_q(\xi ))&=&\LM (Z^{q-q_1}\T\xi )\\
&=&Z^{q-q_1}\LM (\T\xi )\\
&=&Z^{q-q_1}\LM (\T a^{\alpha}\LM (\T g_i))\\
&=&\LM (z^{q-q_1}\T a^{\alpha}\LM (\T g_i)).\end{array}$$ This shows 
that $\LM (\T g_i)|\LM (h_q(\xi ))$, thereby $\tau (\G )$ is a left 
Gr\"obner basis of $\T N$.}\par

(ii) $\Rightarrow$ (i)  If $\xi\in N$ and $\xi\ne 0$,  then 
$\T\xi\ne 0$ and $\LM (\T\xi )=\widetilde{\LM (\xi )}$ by Lemma 
5.3.4. Suppose that $\tau (\G )$ is a left Gr\"obner basis of $\T N$ 
with respect to  $\prec_{\widetilde{e}}$. Then there exists some $\T 
g_i\in\tau (\G )$ such that  $\LM (\T g_i)|\LM (\T\xi )$, i.e., 
there is a monomial $Z^{m}\T a^{\gamma}\in\widetilde{\B}$ such that 
$\LM (\T\xi )=\LM (Z^{m}\T a^{\gamma}\LM (\T g_i))$. Since the given 
left monomial ordering $\PRCEGR$ on $\BE$ is a graded left monomial 
ordering, it follows from Lemma 5.2.2, the definition of 
$\prec_{\widetilde{e}}$ and  Lemma 5.3.2 that
$$\begin{array}{rcl} \widetilde{a^{\alpha}e_j}=\widetilde{\LM (\xi
)}=\LM (\T\xi )&=&\LM (Z^{m}\T a^{\gamma}\LM (\T g_i))\\
&=&Z^{m}(\LM ((a^{\gamma}\LM (g_i))\widetilde{~}))\\
&=&Z^{m}(\LM (a^{\gamma}\LM (g_i)))\widetilde{~}.\end{array}$$ 
Noticing the discussion on  $\T L$ and the role played by $Z$ given 
before Lemma 5.2.2,   we must have $m=0$, thereby $\LM (\xi )=\LM 
(a^{\gamma}\LM (g_i))$. This shows that $\G$ is a left Gr\"obner 
basis for $N$.{\parindent=0pt \v5

{\bf Remark.} It is known that Gr\"obner bases for ungraded ideals 
in both a commutative polynomial algebra and a noncommutative free 
algebra can be obtained via computing homogeneous Gr\"obner bases 
for graded ideals in the corresponding homogenized (graded) algebras 
(cf. [Fr\"ob], [Li2]). In our case here for an $\NZ$-filtered 
solvable polynomial algebra $A$ with respect to a positive-degree 
function $d(~)$, by using a  (de)homogenization-like trick  with 
respect to the central regular element $Z$ in $\T A$, the discussion 
on $\T A$ and $\T L$ presented in previous Section 5.1 indeed 
enables us to have a whole theory and strategy similar to that given 
in (Section 3.6 and Section 3.7 of [Li2]), so that left Gr\"obner 
bases of submodules (left ideals) in $L$ (in $A$) can be obtained 
via computing homogeneous left Gr\"obner bases of graded submodules 
(graded left ideals) in $\T L$ (in $\T A$). Since such a topic is 
beyond the scope of this text,  we omit the detailed discussion 
here. }\v5

\section*{5.4. F-Bases and Standard Bases with Respect to Good
Filtration}

Let $A=K[a_1,\ldots ,a_n]$ be an $\NZ$-filtered solvable polynomial 
algebra with admissible system $(\B ,\prec )$ and the 
$\NZ$-filtration  $FA=\{F_pA\}_{p\in\NZ}$ constructed  with respect 
to a given positive-degree function $d(~)$ on  $A$ (see Section 
5.1). In this section, we introduce F-bases and standard bases 
respectively for $\NZ$-filtered left  $A$-modules and their 
submodules with respect to good filtration, and we show that any two 
minimal F-bases, respectively any two minimal standard bases have 
the same number of elements and the same number of elements of the 
same fil-degree.  Moreover, we show that a standard basis for a 
submodule $N$ of a filtered free $A$-module $L$ can be obtained via 
computing a left Gr\"obner basis  of $N$ with respect to a graded 
left monomial ordering. All notions, notations and conventions used 
before are maintained. \v5

Let $M$ be a left $A$-module. We say that $M$ is an $\NZ$-{\it 
filtered $A$-module} if $M$ has a filtration $FM=\{ 
F_qM\}_{q\in\NZ}$, where each $F_qM$ is a $K$-subspace of $M$, such 
that $M=\cup_{q\in\NZ}F_qM$, $F_qM\subseteq F_{q+1}M$ for all 
$q\in\NZ$, and $F_pAF_qM\subseteq F_{p+q}M$ for all $p,q\in\NZ$. 
{\parindent=0pt\v5

{\bf Convention.}  Unless otherwise stated, from now on in the 
subsequent sections  a filtered $A$-module $M$  is always meant an 
$\NZ$-filtered module with a filtration of the type $FM=\{ 
F_qM\}_{q\in\NZ}$ as described above.}\v5

Let $G(A)$ be the associated graded algebra of $A$, $\T A$ the Rees 
algebra of $A$, and $Z$ the homogeneous element of degree 1 in $\T 
A_1$ represented by the multiplicative identity 1 of $A$ (see 
Section 5.1). If $M$ is a filtered $A$-module with the filtration 
$FM=\{ F_qM\}_{q\in\NZ}$, then, actually as with a  filtered free 
$A$-module in Section 5.2,  $M$  has the  associated  graded 
$G(A)$-module $G(M)=\oplus_{q\in\NZ}G(M)_q$ with $G(M)_0=F_0M$ and 
$G(M)_q=F_qM/F_{q-1}M$ for $q\ge 1$, and the Rees module of $M$ is 
defined as the graded $\T A$-module $\T M=\oplus_{q\in\NZ}\T M_q$ 
with each $\T M_q=F_qM$. Also, we may define the {\it fil-degree} of 
a {\it nonzero} $\xi\in M$, that is, $d_{\rm fil}(\xi )=0$ if 
$\xi\in F_0M$, and if $\xi\in F_qM-F_{q-1}M$ for some $q>0$, then 
$d_{\rm fil}(\xi )=q$. For a nonzero $\xi\in M$ with $d_{\rm 
fil}(\xi )=q\ge 0$, if we write $\sigma (\xi )$ for the nonzero 
homogeneous element of degree $q$ represented by $\xi$ in $G(M)_q$, 
$\T{\xi}$ for the degree-$q$ homogeneous element represented by 
$\xi$ in $\T M_q$, and $h_{q'}(\xi )$ for the degree-$q'$ 
homogeneous element represented by $\xi$ in $\T M_{q'}$ with $q<q'$, 
then $d_{\rm fil}(\xi )=q=d_{\rm gr}(\sigma (\xi ))=d_{\rm gr}(\T\xi 
)$, and $d_{\rm gr}(h_{q'}(\xi ))=q'$. Finally, as for a  filtered 
free $A$-module in Section 5.2, we have  $\T M/Z\T M\cong G(M)$ as 
graded $G(A)$-modules, and $\T M/(1-Z)\T M\cong M$ as $A$-modules.  
\par

With notation as fixed above, the lemma presented below is a version 
of ([LVO], Ch.I, Lemma 5.4, Theorem 5.7) for $\NZ$-filtered modules. 
{\parindent=0pt\v5

{\bf 5.4.1. Lemma}  Let $M$ be a filtered $A$-module with the 
filtration $FM=\{ F_qM\}_{q\in\NZ}$, and $V=\{ v_1,\ldots ,v_m\}$ a 
finite subset of nonzero elements in $M$. The following statements 
are equivalent:\par

(i)  There is a subset $S=\{ n_1,\ldots ,n_m\}\subset\NZ$ such that
$$F_qM=\sum^m_{i=1}\left (\sum_{p_i+n_i\le q}F_{p_i}A\right )v_i,\quad q\in\NZ ;$$\par

(ii) $G(M)=\sum^m_{i=1}G(A)\sigma (v_i)$;\par

(iii) $\T M=\sum^m_{i=1}\T A\T{v_i}$.\par\QED}\v5

{\parindent=0pt\v5

{\bf 5.4.2. Definition}  Let $M$ be a filtered  $A$-module with the 
filtration $FM=\{ F_qM\}_{q\in\NZ}$, and let $V =\{v_1,\ldots, 
v_m\}\subset M$ be a finite subset of nonzero elements. If $V$ 
satisfies one of the equivalent conditions of Lemma 5.4.1, then we 
call $V$ an {\it F-basis} of $M$ with respect to $FM$. }\v5

Let $M$ be a filtered  $A$-module with the filtration $FM=\{ 
F_qM\}_{q\in\NZ}$. If $V$ is an F-basis of $M$ with respect to $FM$, 
then it is necessary to note that{\parindent=1.35truecm\par

\item{(1)} since $M=\cup_{q\in\NZ}F_qM$, it is clear that $V$ is
certainly a generating set of the $A$-module $M$, i.e., 
$M=\sum^m_{i=1}Av_i$;\par

\item{(2)} due to Lemma 5.4.1(i), the  filtration $FM$ is usually
referred to as a {\it good filtration} of $M$ in the literature 
concerning filtered module theory (cf. [LVO]). \par}{\parindent=0pt

Indeed, if $M$ is a finitely generated $A$-module, then {\it any}  
finite generating set $U =\{u_1,\ldots ,u_t\}$ of $M$ can be turned 
into an F-basis with respect to some good filtration  $FM$. More 
precisely, let $S=\{n_1,\ldots ,n_t\}$ be an arbitrarily chosen 
subset of $\NZ$, then the required good filtration 
$FM=\{F_qM\}_{q\in\NZ}$ can be defined by setting
$$\begin{array}{l} F_qM=\{ 0\}~\hbox{if}~q<\min\{ n_1,\ldots ,n_m\} ;\\
F_qM=\sum^t_{i=1}\left (\sum_{p_i+n_i\le q}F_{p_i}A\right 
)u_i~\hbox{otherwise}.\end{array}\quad q\in\NZ .$$  In particular, 
if  $L=\oplus_{i=1}^sAe_i$ is a filtered free $A$-module with the 
filtration $FL=\{ F_qL\}_{q\in\NZ}$ as constructed in Subsection 5.2  
such that $d_{\rm fil}(e_i)=b_i$, $1\le i\le s$, then  $\{ 
e_1,\ldots ,e_s\}$ is an F-basis of $L$ with respect to the good 
filtration $FL$. }{\parindent=0pt\v5

{\bf 5.4.3. Definition}  Let $M$ be a filtered $A$-module with the 
filtration $FM=\{ F_qM\}_{q\in\NZ}$, and suppose that $M$ has an 
F-basis $V =\{v_1,\ldots, v_m\}$ with respect to $FM$. If any proper 
subset of $V$ cannot be an F-basis of $M$ with respect to $FM$, then 
we say that $V$ is a {\it minimal F-basis} of $M$ with respect to 
$FM$. }\v5

Note that $A$ is an $\NZ$-filtered $K$-algebra such that 
$G(A)=\oplus_{p\in\NZ}G(A)_p$ with $G(A)_0=K$, $\T 
A=\oplus_{p\in\NZ}\T A_p$ with $\T A_0=K$, while $K$ is a field.  By 
Lemma 5.4.1 and the well-known result on graded modules over an 
$\NZ$-graded algebra with the degree-0 homogeneous part a field (cf. 
[Eis], Chapter 19; [Kr1], Chapter 3; [Li3]),  we have immediately 
the following{\parindent=0pt\v5

{\bf 5.4.4. Proposition}  Let $M$ be a filtered $A$-module with the 
filtration $FM=\{ F_qM\}_{q\in\NZ}$, and $V=\{ v_1,\ldots 
,v_m\}\subset M$ a subset of nonzero elements.  Then $V$ is  a 
minimal F-basis of $M$ with respect to $FM$ if and only if $\sigma 
(V)=\{\sigma (v_1),\ldots ,\sigma (v_m)\}$ is a minimal homogeneous 
generating set of $G(M)$ if and only if $\tau (V)=\{\T v_1,\ldots 
,\T v_m\}$ is a minimal homogeneous generating set of $\T M$. Hence, 
any two minimal F-bases of $M$ with respect to  $FM$ have the same 
number of elements and the same number of elements of the same 
fil-degree.\par\QED}\v5

Let $M$ be an $\NZ$-filtered $A$-module with the filtration $FM=\{ 
F_qM\}_{q\in\NZ}$, and let $N$ be a submodule of $M$ with the 
filtration $FN=\{ F_qN=N\cap F_qM\}_{q\in\NZ}$ induced by $FM$. 
Then, as with a filtered free $A$-module in Section 5.2, the 
associated graded $G(A)$-module $G(N)=\oplus_{q\in\NZ}G(N)_q$ of $N$ 
with $G(N)_q=F_qN/F_{q-1}N$ is a graded submodule of $G(M)$, and the 
Rees module $\T N=\oplus_{q\in\NZ}\T N_q$ of $N$ with $\T N_q=F_qN$ 
is a graded submodule of $\T M$. {\parindent=0pt\v5

{\bf 5.4.5. Definition}   Let $M$ be a filtered  $A$-module with the 
filtration $FM=\{ F_qM\}_{q\in\NZ}$, and let $N$ be a submodule of 
$M$. Consider the filtration $FN=\{ F_qN=N\cap F_qM\}_{q\in\NZ}$ of 
$N$ induced by $FM$. If $W=\{ \xi_1,\ldots ,\xi_s\}\subset N$ is an 
F-basis with respect to $FN$ in the sense of Definition 5.4.2, then 
we call $W$ a {\it standard basis} of $N$. \v5

{\bf Remark.} By referring to Lemma 5.4.1, one may check that our 
definition 5.4.5 of a standard basis coincides with the classical 
Macaulay basis provided $A=K[x_1,\ldots ,x_n]$ is the commutative 
polynomial $K$-algebra (cf. [KR2], Definition 4.2.13, Theorem 
4.3.19), for, taking the $\NZ$-filtration $FA$ with respect to an 
arbitrarily chosen positive-degree function $d(~)$ on $A$, there are 
graded algebra isomorphisms $G(A)\cong A$ and $\T A\cong 
K[x_0,x_1,\ldots ,x_n]$, where $d(x_0)=1$ and $x_0$ plays the role 
that the central regular element $Z$ of degree 1 does in $\T A$. 
Moreover, if two-sided ideals of an $\NZ$-filtered solvable 
polynomial algebra $A$ are considered, then one may see that our 
definition 5.4.5 of a standard basis coincides with the standard 
basis defined in [Gol].\v5

{\bf 5.4.6. Definition}  Let $M$ be a filtered $A$-module with the 
filtration $FM=\{ F_qM\}_{q\in\NZ}$, and  $N$ a submodule of $M$ 
with the filtration $FN=\{ F_qN=N\cap F_qM\}_{q\in\NZ}$ induced by 
$FM$. Suppose that $N$ has a standard basis  $W=\{\xi_1,\ldots 
,\xi_m\}$ with respect to $FN$. If any proper subset of $W$ cannot 
be a standard basis for $N$ with respect to $FN$, then we call $W$ a 
{\it minimal standard basis} of $N$ with respect to $FN$.}\v5

If $N$ is a submodule of some filtered $A$-module $M$ with 
filtration $FM$, then since a standard basis of $N$ is defined as an 
F-basis of $N$ with respect to the filtration $FN$ induced by $FM$, 
the next proposition follows from Proposition 
5.4.4.{\parindent=0pt\v5

{\bf 5.4.7. Proposition}  Let $M$ be a filtered $A$-module with the 
filtration $FM=\{ F_qM\}_{q\in\NZ}$, and $N$ a submodule of $M$ with 
the induced filtration $FN=\{ F_qN=N\cap F_qM\}_{q\in\NZ}$. A finite 
subset of nonzero elements  $W=\{ \xi_1,\ldots ,\xi_s\}\subset N$ is 
a minimal standard basis of $N$ with respect to $FN$ if and only if 
$\sigma (W)=\{\sigma (\xi_1),\ldots ,\sigma (\xi_m)\}$ is a minimal 
homogeneous generating set of $G(N)$ if and only if $\tau (W)=\{\T 
\xi_1,\ldots ,\T\xi_m\}$ is a minimal homogeneous generating set of 
$\T N$. Hence, any two minimal standard bases of $N$  have the same 
number of elements and the same number of elements of the same 
fil-degree.\par\QED}\v5

Since $A$, $G(A)$ and $\T A$ are all Noetherian domains (Proposition 
1.1.4(ii) of Chapter 1, Theorem 5.1.5 of the current chapter), if a 
filtered $A$-module $M$ has an F-basis $V$ with respect to a given 
filtration $FM$, then  the existence of a standard basis for a 
submodule $N$ of $M$ follows immediately from Lemma 5.4.1. Our next 
theorem shows that a standard basis for a submodule $N$ of a 
filtered free $A$-module $L$ can be obtained via computing a left 
Gr\"obner basis  of $N$ with respect to a {\it graded left monomial 
ordering}.{\parindent=0pt\v5
\def\PRCVE{\prec_{\varepsilon\hbox{-}gr}}\def\BV{\B (\varepsilon )}

{\bf 5.4.8. Theorem}  Let  $L=\oplus_{i=1}^sAe_i$ be a filtered free 
$A$-module with the filtration $FL=\{ F_qL\}_{q\in\NZ}$ such that 
$d_{\rm fil}(e_i)=b_i$, $1\le i\le s$, and let $\prec_{e\hbox{-}gr}$ 
be a graded left monomial ordering on $\BE$ (see Section 5.3). If 
$\G =\{ g_1,\ldots ,g_m\}\subset L$ is a left Gr\"obner basis for 
the submodule $N=\sum_{i=1}^mAg_i$ of $L$ with respect to $\PRCEGR$, 
then $\G$ is a standard basis for $N$ in the sense of Definition 
5.4.5.\vskip 6pt
\def\PRCEGR{\prec_{e\hbox{-}gr}}

{\bf Proof} If $\xi\in F_qN=F_qL\cap N$ and $\xi\ne 0$, then $d_{\rm 
fil}(\xi )\le q$ and $\xi $ has a left Gr\"obner representation by 
$\G$, that is, $\xi =\sum_{i,j}\lambda_{ij}a^{\alpha (i_j)}g_j$, 
where $\lambda_{ij}\in K^*$, $a^{\alpha (i_j)}\in\B$ with $\alpha 
(i_j)=(\alpha_{i_{j1}},\ldots ,\alpha_{i_{jn}}) \in\NZ^n$, 
satisfying $\LM (a^{\alpha (i_j)}g_j)\preceq_{e\hbox{-}gr}\LM (\xi 
)$. Suppose $d_{\rm fil}(g_j)=n_j$, $1\le j\le m$. Since $\PRCEGR$ 
is a graded left monomial ordering on $\BE$, by Lemma 5.3.1 we may 
assume that $\LM (g_j)=a^{\beta (j)}e_{t_j}$ with $\beta 
(j)=(\beta_{j_1},\ldots ,\beta_{j_n})\in\NZ^n$ and $1\le t_j\le s$, 
such that $d(a^{\beta (j)})+b_{t_j}=n_j$, where $d(~)$ is the given 
positive-degree function on $A$. Furthermore, by (Lemma 2.1.2(ii) of 
Chapter 2), we have
$$\LM (a^{\alpha (i_j)}g_j)=\LM (a^{\alpha (i_j)}a^{\beta
(j)}e_{t_j})=a^{\alpha (i_j)+\beta (j)}e_{t_j},$$ and it follows 
from Lemma 5.1.2, Lemma 5.2.2 and Lemma 5.3.1 that $d(a^{\alpha 
(i_j)})+n_j=d(a^{\alpha (i_j)})+d(a^{\beta (j)})+b_{t_j}=d(a^{\alpha 
(i_j)+\beta (j)})+b_{t_j}\le q$. Hence $\xi\in \sum_{j=1}^m\left 
(\sum_{p_j+n_j\le q}F_{p_j}A\right )g_j$. This shows that 
$F_qN=\sum_{j=1}^m\left (\sum_{p_j+n_j\le q}F_{p_j}A\right )g_j$, 
i.e., $\G$ is a standard basis for $N$.}\v5

\section*{5.5. Computation of Minimal F-Bases and Minimal Standard
Bases}

Let $A=K[a_1,\ldots ,a_n]$ be an $\NZ$-filtered solvable polynomial 
algebra with  admissible system $(\B ,\prec )$ and the 
$\NZ$-filtration  $FA=\{F_pA\}_{p\in\NZ}$ constructed  with respect 
to a positive-degree function $d(~)$ on  $A$ (see Section 5.1). In 
this section we show how to algorithmically compute minimal F-bases 
for quotient modules of a filtered free left $A$-module $L$, and how 
to algorithmically compute minimal standard bases for submodules of 
$L$ in the case where a graded left monomial ordering $\PRCEGR$ on 
$L$ is employed.  All notions, notations and conventions used before 
are maintained.\v5

We start by a little more preparation.  Let $M$ and $M'$ be 
$\NZ$-filtered left $A$-modules with the filtration $FM=\{ 
F_qM\}_{q\in\NZ}$ and $FM'=\{ F_qM'\}_{q\in\NZ}$ respectively, and 
$M~\mapright{\varphi}{}~M'$ an $A$-module homomorphism. If $\varphi 
(F_qM)\subseteq F_qM'$ for all $q\in\NZ$, then we call $\varphi$ a 
{\it filtered homomorphism}. In the literature, such filtered 
homomorphisms are also referred to as filtered homomorphism of 
degree-0 (cf. [NVO], [LVO]). By the definition it is clear that the 
identity map of $\NZ$-filtered $A$-modules is filtered homomorphism, 
and compositions of filtered homomorphisms are  filtered  
homomorphisms.  Thus, all $\NZ$-filtered left $A$-modules form a 
subcategory of the category of left $A$-modules, in which morphisms 
are the filtered homomorphisms as defined above. Furthermore, if 
$M~\mapright{\varphi}{}~M'$ is a filtered homomorphism with kernel 
Ker$\varphi =N$, then $N$ is an $\NZ$-filtered submodule of $M$ with 
the induced filtration $FN=\{ F_qN=N\cap F_qM\}_{q\in\NZ}$, and the 
image $\varphi (M)$ of $\varphi$ is an $\NZ$-filtered submodule of 
$M'$ with the filtration $F\varphi (M)=\{F_q\varphi (M)=\varphi 
(F_qM)\}_{q\in\NZ}$. Consequently, the exactness of a  sequence 
$N~\mapright{\varphi}{}~M~\mapright{\psi}{}~M'$ of filtered 
homomorphisms in the category of $\NZ$-filtered $A$-modules is 
defined as the same as for a sequence of usual $A$-module 
homomorphisms, i.e., the sequence satisfies Im$\varphi =$ Ker$\psi$. 
Long exact sequence in the category  of $\NZ$-filtered $A$-modules 
may be defined in an obvious way. \v5

Let $G(A)$ be the associated $\NZ$-graded algebra of $A$ and $\T A$ 
the Rees algebra of $A$. Then naturally, a filtered homomorphism  
$M~\mapright{\varphi}{}~M'$ induces a graded $G(A)$-module 
homomorphism $G(M)~\mapright{G(\varphi )}{}~G(M')$, where if $\xi\in 
F_qM$ and $\OV\xi =\xi+F_{q-1}M$ is the coset represented by $\xi$ 
in $G(M)_q=F_qM/F_{q-1}M$, then $G(\varphi )(\OV\xi )=\varphi (\xi 
)+F_{q-1}M'\in G(M')_q=F_qM'/F_{q-1}M'$,  and $\varphi$ induces a 
graded $\T A$-module homomorphism  $\T 
M~\mapright{\T\varphi}{}~\T{M'}$, where if $\xi\in F_qM$ and 
$h_q(\xi )$ is the homogeneous element of degree $q$ in $\T 
M_q=F_qM$, then $\T\varphi (h_q(\xi ))=h_q(\varphi (\xi 
))\in\T{M'}_q=F_qM'$. Moreover, if 
$M~\mapright{\varphi}{}~M'~\mapright{\psi}{}~M''$ is a sequence of 
filtered homomorphisms, then $G(\psi )\circ G(\varphi 
)=G(\psi\circ\varphi )$ and $\T\psi\circ\T\varphi 
=\widetilde{\psi\circ\varphi}$.
\par

Furthermore, recall that a filtered  homomorphism 
$M~\mapright{\varphi}{}~M'$ is called a {\it strict filtered 
homomorphism} if $\varphi (F_qM)=\varphi (M)\cap F_qM'$ for all 
$q\in\NZ$.  Note that if $N$ is a submodule of $M$ and $\OV M=M/N$, 
then, considering the induced filtration $FN=\{ F_qN=N\cap 
F_qM\}_{q\in\NZ}$ of $N$ and the induced filtration $F(\OV M)=\{ 
F_q\OV M=(F_qM+N)/N\}_{q\in\NZ}$ of $\OV M$, the inclusion map 
$N\hookrightarrow M$ and the canonical map $M\r \OV M$ are strict 
filtered  homomorphisms. Concerning strict filtered homomorphisms 
and the induced graded homomorphisms, the next  proposition is 
quoted from ([LVO], CH.I, Section 4).{\parindent=0pt\v5

{\bf 5.5.1. Proposition}  Given a sequence of filtered homomorphisms 
$$N~\mapright{\varphi}{}~M~\mapright{\psi}{}~M',\leqno{(*)}$$
such that $\psi\circ\varphi =0$, the following statements are 
equivalent.\par

(i) The sequence $(*)$ is exact and $\varphi$, $\psi$ are strict 
filtered homomorphisms.\par

(ii) The induced sequence $G(N)~\mapright{G(\varphi 
)}{}~G(M)~\mapright{G(\psi ) )}{}~G(M')$ is exact.\par

(iii) The induced sequence $\T N~\mapright{\T\varphi}{}~\T 
M~\mapright{\T\psi}{}~\T{M'}$ is exact.\par\QED}\v5

Let $L=\oplus_{i=1}^mAe_i$ be a filtered free $A$-module with the 
filtration $FL=\{ F_qL\}_{q\in\NZ}$ such that $d_{\rm 
fil}(e_i)=b_i$, $1\le i\le m$. Then as we have noted in Section 5.4, 
$\{ e_1,\ldots ,e_m\}$ is an F-basis of $L$ with respect to the good 
filtration $FL$. Let $N$ be a submodule of $L$,  and let the 
quotient module $M=L/N$ be equipped with the filtration 
$FM=\{F_qM=(F_qL+N)/N\}_{q\in\NZ}$ induced by $FL$. Without loss of 
generality, we assume that $\OV e_i\ne 0$ for $1\le i\le m$, where 
each $\OV e_i$ is the coset represented by $e_i$ in $M$. Then we see 
that $\{ \OV e_1,\ldots ,\OV e_m\}$ is an F-basis of $M$ with 
respect to $FM$. {\parindent=0pt\v5

{\bf 5.5.2. Lemma}  Let $M=L/N$ be as fixed above, and let 
$N=\sum^s_{j=1}A\xi_j$ be generated by the set of nonzero elements $ 
U =\{\xi_1,\ldots ,\xi_s\}$, where 
$\xi_{\ell}=\sum_{k=1}^sf_{k\ell}e_k$ with $f_{k\ell}\in A$ and 
$d_{\rm fil}( \xi_{\ell})=q_{\ell}$, $1\le \ell\le s$. The following 
statements hold.\par

(i) If for some $j$, $\xi_j$ has a nonzero term $f_{ij}e_i$ such 
that $d_{\rm fil}(f_{ij}e_i)=d_{\rm fil}(\xi_j)=q_j$ and the 
coefficient $f_{ij}$ is a nonzero constant, say $f_{ij}=1$ without 
loss of generality, then for each $\ell =1,\ldots ,j-1,j+1,\ldots 
,s$, the element $\xi_{\ell}'=\xi_{\ell}-f_{i\ell}\xi_j$ does not 
involve $e_i$. Putting $ U '=\{ \xi_1',\ldots 
,\xi'_{j-1},\xi_{j+1}',\ldots ,\xi'_s\}$, $N'=\sum_{\xi_{\ell}'\in U 
'}A\xi_{\ell}'$, and considering the filtered free $A$-module 
$L'=\oplus_{k\ne i}Ae_k$ with the filtration $FL'=\{ 
F_qL'\}_{q\in\NZ}$  in which each $e_k$ has the same filtered degree 
as it is in $L$, i.e., $d_{\rm fil}(e_k)=b_k$, if the quotient 
module $M'=L'/N'$ is equipped with the  filtration $FM'=\{ 
F_qM'=(F_qL'+N')/N'\}_{q\in\NZ}$ induced by $FL'$, then there is a 
strict filtered isomorphism $\varphi$: $M'\cong M$, i.e., $\varphi$ 
is an $A$-module isomorphism such that $\varphi (F_qM')=F_qM$ for 
all $q\in\NZ$. \par

(ii) With the assumptions and notations as in (i), if $ U 
=\{\xi_1,\ldots ,\xi_s\}$ is a standard basis of $N$ with respect to 
the filtration $FN$ induced by $FL$, then $ U '=\{ \xi_1',\ldots 
,\xi'_{j-1},\xi_{j+1}',\ldots ,\xi'_s\}$ is a standard basis of $N'$ 
with respect to the filtration $FN'$ induced by $FL'$. \vskip 6pt

{\bf Proof} (i) Since $f_{ij}=1$ by the assumption, we see that 
every $\xi_{\ell}'=\xi_{\ell}-f_{i\ell}\xi_j$ does not involve 
$e_i$. Let $ U '=\{ \xi_1',\ldots ,\xi'_{j-1},\xi_{j+1}',\ldots 
,\xi'_s\}$ and  $N'=\sum_{\xi_{\ell}'\in U '}A\xi_{\ell}'$. Then 
$N'\subset L'=\oplus_{k\ne i}Ae_k\subset L$ and $N=N'+A\xi_j$. Since 
$\xi_j =e_i+\sum_{k\ne i}f_{kj}e_k$, the naturally defined 
$A$-module homomorphism $M'=L'/N'\mapright{\varphi}{}L/N=M$ with 
$\varphi (\OV e_k)=\OV e_k$, $k=1,\ldots ,i-1,i+1,\ldots ,m$, is an 
isomorphism, where, without confusion, $\OV e_k$ denotes the coset 
represented by $e_k$ in $M'$ and $M$ respectively. It remains to see 
that $\varphi$ is a strict filtered isomorphism. Note that $\{ 
e_1,\ldots ,e_m\}$ is an F-basis of $L$ with respect to $FL$ such 
that $d_{\rm fil}(e_i)=b_i$, $1\le i\le m$, i.e.,
$$F_qL=\sum_{i=1}^m\left (\sum_{p_i+b_i\le q}F_{p_i}A\right )e_i, \quad q\in\NZ ,$$
that $\{ e_1,\ldots ,e_{i-1},e_{i+1},\ldots ,e_m\}$ is an F-basis of 
$L'$ with respect to $FL'$ such that $d_{\rm fil}(e_k)=b_k$, where 
$k\ne i$, i.e.,
$$F_qL'=\sum_{k\ne i}\left (\sum_{p_i+b_k\le q}F_{p_i}A\right )e_k, \quad q\in\NZ ,$$
and that $\xi_j =e_i+\sum_{k\ne i}f_{kj}e_k$ with $q_j=d_{\rm 
fil}(\xi_j)=d_{\rm fil}(e_i)=b_i$ such that $d_{\rm 
fil}(f_{kj})+b_k\le q_j$ for all $f_{kj}\ne 0$. It follows that 
$\sum_{k\ne i}f_{kj}\OV e_k\in F_{q_j}M'$ and $\varphi (\sum_{k\ne 
i}f_{kj}\OV e_k)=\OV e_i\in F_{q_j}M$, thereby $\varphi 
(F_qM')=F_qM$ for all $q\in\NZ$, as desired.}\par

(ii) Note that $\xi_{\ell}'=\xi_{\ell}-f_{i\ell}\xi_j$. By the 
assumption on $\xi_j$, if $f_{i\ell}\ne 0$ and $d_{\rm 
fil}(f_{i\ell}e_i)=d_{\rm fil}(\xi_{\ell})=q_{\ell}$, then since 
$d_{\rm fil}(\xi_j)=d_{\rm fil}(e_i)$ we have $d_{\rm 
fil}(f_{i\ell}\xi_j)=d_{\rm fil}(\xi_{\ell})=q_{\ell}$. It follows 
that if we equip $N$ with the filtration $FN=\{ F_qN=N\cap 
F_qL\}_{q\in\NZ}$ induced by $FL$ and consider the associated graded 
module $G(N)$ of $N$, then $d_{\rm gr}(\sigma (\xi_{\ell}))=d_{\rm 
gr}(\sigma (f_{i\ell}\xi_j))=d_{\rm gr}(\sigma (f_{i\ell})\sigma 
(\xi_j))$ in $G(N)$, i.e., $\sigma (\xi_{\ell})-\sigma 
(f_{i\ell})\sigma (\xi_j)\in G(N)_{q_{\ell}}$. So, if $\sigma 
(\xi_{\ell})-\sigma (f_{i\ell})\sigma (\xi_j)\ne 0$ then $d_{\rm 
fil}(\xi_{\ell}')=q_{\ell}$ and thus $$\sigma (\xi_{\ell}')=\sigma 
(\xi_{\ell}-f_{i\ell}\xi_j)=\sigma (\xi_{\ell})-\sigma 
(f_{i\ell})\sigma (\xi_j).\eqno{(1)}$$ If $f_{i\ell}\ne 0$ and 
$d_{\rm fil}(f_{i\ell}e_i)<d_{\rm fil}(\xi_{\ell})=q_{\ell}$, then 
$\sigma (\xi_{\ell})=\sum_{d(f_{k\ell})+b_k=q_{\ell}}\sigma 
(f_{k\ell})\sigma (e_k)$ does not involve $\sigma (e_i)$. Also since 
$d_{\rm fil}(\xi_j)=d_{\rm fil}(e_i)$, we have $d_{\rm 
fil}(f_{i\ell}\xi_j)<d_{\rm fil}(\xi_{\ell})=q_{\ell}$. Hence 
$d_{\rm fil}(\xi_{\ell}')=d_{\rm fil}(\xi_{\ell})=q_{\ell}$ and thus
$$\sigma (\xi_{\ell}')=\sigma (\xi_{\ell}-f_{i\ell}\xi_j)=\sigma
(\xi_{\ell}).\eqno{(2)}$$ If $f_{i\ell}=0$, then the equality (1) is 
the same  as equality (2). Now, if $ U =\{\xi_1,\ldots ,\xi_s\}$ is 
a standard basis of $N$ with respect to the induced filtration $FN$, 
then $G(N)=\sum_{\ell =1}^sG(A)\sigma (\xi_{\ell})$ by Lemma 5.4.1. 
But since we have also $G(N)=\sum_{\xi_{\ell}'\in U '}G(A)\sigma 
(\xi_{\ell}')+G(A)\sigma (\xi_j)$ where the $\sigma (\xi_{\ell}')$ 
are those nonzero homogeneous elements obtained according to the 
above equalities (1) and (2), it follows from Lemma 5.4.1 that $ 
U'\cup\{\xi_j\}$ is a standard basis of $N$ with respect to the 
induced filtration $FN$.
\par

We next prove that $ U '$ is a standard basis of 
$N'=\sum_{\xi_{\ell}'\in U '}A\xi_{\ell}'$ with respect to the 
filtration $FN'=\{ F_qN'=N'\cap F_qL'\}_{q\in\NZ}$ induced by $FL'$. 
Since $\{ e_1,\ldots ,e_{i-1},e_{i+1},\ldots ,e_m\}$ is an F-basis 
of $L'$ with respect to the filtration $FL'$ such that each $e_k$ 
has the same filtered degree as it is in $L$, i.e., $d_{\rm 
fil}(e_k)=b_k$, it is clear that $F_qL'=L'\cap F_qL$, $q\in\NZ$, 
i.e., the filtration $FL'$ is the one induced by $FL$. Considering 
the filtration $FN'$ of $N'$ induced by $FL'$, it turns out that
$$F_qN'=N'\cap F_qL'=N'\cap F_qL\subseteq N\cap F_qL=F_qN,\quad q\in\NZ .\eqno{(3)}$$
If $\xi\in F_qN'$, then since $ U '\cup\{\xi_j\}$ is a standard 
basis of $N$ with respect to the induced filtration $FN$, the 
formula (3) entails that
$$\xi =\sum_{\xi_{\ell}'\in U '}f_{\ell}\xi_{\ell}'+f_j\xi_j~\hbox{with}~f_{\ell},f_j\in A,
~d(f_{\ell})+d_{\rm fil}(\xi_{\ell}')\le q,~d(f_j)+d_{\rm 
fil}(\xi_j)\le q.\eqno{(4)}$$ Note that every $\xi_{\ell}'$ does not 
involve $e_i$ and consequently $\xi$ does not involve $e_i$. Hence 
$f_j=0$ in (4) by the assumption on $\xi_j$, and thus 
$\xi\in\sum_{\xi_{\ell}'\in U '}\left (\sum_{p_i+q_{\ell}\le 
q}F_{p_i}A\right )\xi_{\ell}'$. Therefor we conclude that $ U '$ is 
a standard basis for $N'$ with respect to the induced filtration 
$FN'$, as desired. \QED \v5

Combining  Proposition 5.4.4, we now show that for quotient modules 
of filtered free $A$-modules, a result similar to  Proposition 
4.3.11 holds true.  {\parindent=0pt\v5

{\bf 5.5.3. Proposition}  Let $L=\oplus_{i=1}^mAe_i$, $M=L/N$ be as 
in Lemma 5.5.2, and suppose that $ U = \{ \xi_1,...,\xi_s\}$ is now 
a standard basis of $N$ with respect to the filtration $FN$ induced 
by $FL$. The algorithm presented below computes a subset $\{ 
e_{i_1},\ldots ,e_{i_{m'}}\}\subset\{ e_1,\ldots ,e_m\}$ and a 
subset $V=\{ v_1,\ldots ,v_t\}\subset N\cap L'$, where $m'\le m$ and 
$L'=\oplus_{q=1}^{m'}Ae_{i_q}$ such that\par

(i) there is a strict filtered isomorphism $L'/N'=M'\cong M$, where 
$N'=\sum^{t}_{\ell =1}Av_{\ell}$, and $M'$ has the filtration 
$FM'=\{ F_qM'=(F_qL'+N')/N'\}_{q\in\NZ}$ induced by the good 
filtration $FL'$ determined by the F-basis $\{ e_{i_1},\ldots 
,e_{i_{m'}}\}$ of $L'$;\par

(ii) $V=\{ v_1,\ldots ,v_t\}$ is a standard basis of 
$N'=\sum^{t}_{\ell =1}Av_{\ell}$ with respect to the filtration 
$FN'$ induced by $FL'$, such that each 
$v_{\ell}=\sum^{m'}_{k=1}h_{k\ell}e_{i_k}$ has the property that 
$h_{k\ell}\in K^*$ implies $d_{\rm fil}(e_{i_k})=b_{i_k}<d_{\rm 
fil}(v_{\ell})$;\par

(iii) $\{\OV e_{i_1},\ldots ,\OV e_{i_{m'}}\}$ is a minimal F-basis 
of $M$ with respect to the filtration $FM$. \par

\underline{\bf Algorithm-MINFB 
~~~~~~~~~~~~~~~~~~~~~~~~~~~~~~~~~~~~~~~~~~~~~~~~~~~~~~~~~}\vskip 6pt

\textsc{INPUT}: $E=\{ e_1,\ldots ,e_m\};~~ U = \{
\xi_1,...,\xi_s\},$\par
~~~~~~~~~~~~~~~~~~~~~~~~~~~~~~~~~~~~~~~\hbox{where}~$\xi_{\ell}=\sum_{k=1}^mf_{k\ell}e_k~
\hbox{with}~f_{k\ell}\in A,$\par
~~~~~~~~~~~~~~~~~~~~~~~~~~~~~~~~~~~~~~~\hbox{and}~$d(f_{k\ell})+b_k\le 
q_{\ell}=d_{\rm fil}(\xi_{\ell}), ~1\le \ell\le s$\par                                 
\textsc{OUTPUT}: ~$E' =\{ e_{i_1},\ldots ,e_{i_{m'}}\}; ~~V=\{
v_1,\ldots ,v_t\}\subset N\cap L',$\par                                        
~~~~~~~~~~~~~~~~~~~~~~~~~~~~~~~~~~~~~~~~~~~~~~~~~~~~~~\hbox{where}~
$L'=\oplus_{k=1}^{m'}Ae_{i_k}$\par                                                   
\textsc{INITIALIZATION}: $t :=s;~ V := U ;  ~m':=m; ~E' :=E$ \par
\textsc{BEGIN}\par                                                       
~~~~~\textsc{WHILE}~$\hbox{there is 
a}~v_j=\sum^{m'}_{k=1}f_{kj}e_k\in V~\hbox{and}~i~ \hbox{is the 
least index}$\par                                                          
~~~~~~~~~~~~~~~~~\hbox{such that}~$f_{ij}\in K^*~
\hbox{with}~d(f_{ij})+b_i=d_{\rm fil}(v_j)~\hbox{DO}$\par                                
~~~~~~~~~~~~~~~~~\hbox{set}~$T=\{1,\ldots ,j-1,j+1,\ldots
,t\}~\hbox{and compute}$\par                                            
~~~~~~~~~~~~~~~~~$v_{\ell}'
=v_{\ell}-\frac{1}{f_{ij}}f_{i\ell}v_j,~\ell\in T,$\par                                  
~~~~~~~~~~~~~~~~~$r=\#\{\ell~|~\ell\in T,~v_{\ell} =0\}$\par                             
~~~~~~~~~~~$t := t-r-1$\par                                                   
~~~~~~~~~~$V :=\{ v_{\ell}=v_{\ell}'~|~\ell\in T,~v_{\ell}'\ne 
0\}$\par                                                                   
~~~~~~~~~~~~~~$=\{ v_1,\ldots ,v_t\}~(\hbox{after reordered})$\par                           
~~~~~~~~~~$m' :=m'-1$\par                                                      
~~~~~~~~~~$E':=E'-\{ e_i\}=\{ e_1,\ldots ,e_{m'}\}~(\hbox{after 
reordered})$\par                                                          
~~~~~\textsc{END}\par                                             
\textsc{END}\vskip -.2truecm
\underline{~~~~~~~~~~~~~~~~~~~~~~~~~~~~~~~~~~~~~~~~~~~~~~~~~~~~~~~~~~~~~~~~~~~~~~~~~~~~~~~~~~~~~~~~~~~~~~~} 
\vskip 6pt

{\bf Proof} First note that for each $\xi_{\ell}\in U$, $d_{\rm 
fil}(\xi_{\ell})$ is determined by Lemma 5.2.1. Since the algorithm 
is clearly finite, the conclusions (i) and (ii) follow from Lemma 
5.5.2.}\par

To prove the conclusion (iii), by the strict filtered isomorphism 
$M'=L'/N'\cong M$ (or the proof of Lemma 5.5.2(i)) it is sufficient 
to show that $\{\OV e_{i_1},\ldots ,\OV e_{i_{m'}}\}$ is a minimal 
F-basis of $M'$ with respect to the filtration $FM'$. By the 
conclusion (ii), $V=\{ v_1,\ldots ,v_t\}$ is a standard basis of the 
submodule $N'=\sum^{t}_{\ell =1}Av_{\ell}$ of $L'$ with respect to 
the filtration $FN'$ induced by $FL'$ such that each 
$v_{\ell}=\sum^{m'}_{k=1}h_{k\ell}e_{i_k}$ has the property that 
$h_{k\ell}\in K^*$ implies $d_{\rm fil}(e_{i_k})=b_{i_k}<d_{\rm 
fil}(v_{\ell})$. It follows from Lemma 5.4.1 that 
$G(N')=\sum^t_{k=1}G(A)\sigma (v_k)$ in which each $\sigma 
(v_k)=\sum_{d(h_{k\ell})+b_{i_k}=d_{\rm fil}(v_k)}\sigma 
(h_{k\ell})\sigma (e_{i_k})$ and all the coefficients $\sigma 
(h_{k\ell})$ satisfy $d_{\rm gr}(\sigma (h_{k\ell}))>0$ (please note 
the discussion in Section 5.2). Since $G(A)_0=K$, 
$G(L')=\oplus_{k=1}^{m'}G(A)\sigma (e_{i_k})$ (Proposition 5.2.3) 
and $G(N')$ is the graded syzygy module of the graded quotient 
module $G(L')/G(N')$, by the well-known result on finitely generated 
graded modules over $\NZ$-graded algebras with the degree-0 
homogeneous part a field (cf. [Eis], Chapter 19; [Kr1], Chapter 3; 
[Li3]),  we conclude  that $\{ \OV{\sigma (e_{i_1})},\ldots 
,\OV{\sigma (e_{i_{m'}})}\}$ is a minimal homogeneous generating set 
of $G(L')/G(N')$. On the other hand, considering the naturally 
formed exact sequence of  strict filtered homomorphisms
$$0~\mapright{}{}~N'~\mapright{}{}~L'~\mapright{}{}~M'=L'/N'~\mapright{}{}~0,$$
by Proposition 5.5.1 we have the $\NZ$-graded $G(A)$-module 
isomorphism $G(L')/G(N')$ $\cong G(L'/N')=G(M')$ with $\OV{\sigma 
(e_{i_k})}\mapsto \sigma (\OV e_{i_k})$, $1\le k\le m'$. Now, by 
means of Proposition 5.4.4, we conclude that $\{\OV e_{i_1},\ldots 
,\OV e_{i_{m'}}\}$ is a minimal F-basis of $M'$ with respect to the 
filtration $FM'$, as desired. \QED \v5

Finally, let $L=\oplus_{i=1}^sAe_i$ be a filtered free $A$-module 
with the filtration $FL=\{ F_qL\}_{q\in\NZ}$ such that $d_{\rm 
fil}(e_i)=b_i$, $1\le i\le s$, and let $\prec_{e\hbox{-}gr}$ be a 
graded left monomial ordering on the $K$-basis $\BE =\{ 
a^{\alpha}e_i~|~a^{\alpha}\in\B ,~1\le i\le s\}$ of $L$ (see Section 
5.3). Combining Theorem 5.3.3 and (Theorem 4.3.8 of Chapter 4), the 
next theorem shows how to algorithmically compute a minimal standard 
basis. {\parindent=0pt\v5

{\bf 5.5.4. Theorem} Let $N=\sum^c_{i=1}A\theta_i$ be a submodule of 
$L$ generated by the subset of nonzero elements $\Theta =\{ 
\theta_1,\ldots ,\theta_c\}$, and let $FN=\{ F_qN=F_qL\cap 
N\}_{q\in\NZ}$  be the filtration of $N$ induced by $FL$. Then a 
minimal standard basis of  $N$ with respect to $FN$ can be obtained 
by implementing the following procedures:}\par

{\bf Procudure 1.} With the initial input data $\Theta =\{ 
\theta_1,\ldots ,\theta_c\}$, run {\bf Algorithm-LGB} presented in 
(Section 2.3 of Chapter 2) to compute a left Gr\"obner basis $ U =\{ 
\xi_1,\ldots ,\xi_m\}$ for $N$ with respect to $\prec_{e\hbox{-}gr}$ 
on $\BE$.
\par

{\bf Procudure 2.} Let $G(N)$ be the associated graded $G(A)$-module 
of $N$ determined by the induced filtration $FN$. Then $G(N)$ is a 
graded submodule of the associated grade free $G(A)$-module $G(L)$, 
and it  follows from  Theorem 5.3.3 that $\sigma ( U )=\{\sigma 
(\xi_1),\ldots ,\sigma (\xi_m)\}$ is a homogeneous left Gr\"obner 
basis of $G(N)$ with respect to $\prec_{\sigma (e)\hbox{-}gr}$ on 
$\sigma (\BE )$. With the initial input data $\sigma ( U )$, run 
{\bf Algorithm-MINHGS} presented in (Theorem 4.3.8 of Chapter 4) to 
compute a minimal homogeneous generating set $\{ \sigma 
(\xi_{j_1}),\ldots ,\sigma (\xi_{j_t})\}\subseteq\sigma ( U )$ for 
$G(N)$.\par

{\bf Procudure 3.} Write down $W=\{ \xi_{j_1},\ldots ,\xi_{j_t}\}$. 
Then $W$ is a Minimal standard basis for $N$ by Proposition 5.4.7.
\par\QED{\parindent=0pt\v5

{\bf Remark.} (i) Note that the initial input data $\sigma (U)$ in 
{\bf Procedure 2} above is already a left Gr\"obner basis for 
$G(N)$. By (Theorem 4.3.8 of Chapter 4), the finally returned left 
Gr\"obner basis by {\bf Algorithm-MINHGS} is indeed a minimal left 
Gr\"obner basis for $G(N)$.\par

(ii) By Theorem 5.3.5 and Proposition 5.4.7 it is clear that we can 
also obtain a minimal standard basis of the submodule $N$ via 
computing a minimal homogeneous generating set for the Rees module 
$\T N$ of $N$, which is a graded submodule of the Rees module $\T L$ 
of $L$. However, noticing the structure of $\T A$ and $\T L$ (see 
Theorem 5.1.5, Proposition 5.2.3) it is equally clear that the cost 
of working on $\T A$ will be much higher than working on $G(N)$.}\v5

\section*{5.6. Minimal Filtered Free Resolutions and Their 
Uniqueness}

Let $A=K[a_1,\ldots ,a_n]$ be an $\NZ$-filtered solvable polynomial 
algebra with admissible system $(\B ,\prec )$ and the 
$\NZ$-filtration $FA=\{F_pA\}_{p\in\NZ}$ constructed  with respect 
to a positive-degree function $d(~)$ on  $A$ (see Section 5.1). In 
this section, by using minimal F-bases  and minimal standard bases 
in the sense of Definition 5.4.3 and Definition 5.4.6, we define 
minimal filtered free resolutions for finitely generated left 
$A$-modules, and we show that such minimal resolutions are  unique 
up to strict filtered isomorphism of chain complexes in the category 
of filtered $A$-modules. All notions, notations and conventions used 
before are maintained.\v5

Let $M=\sum^m_{i=1}A\xi_i$ be an arbitrary finitely generated 
$A$-module. Then,  as we have noted in section 5.4,  $M$ may be 
endowed with a good filtration $FM=\{ F_qM\}_{q\in\NZ}$ with respect 
to an arbitrarily chosen subset $\{ n_1,\ldots ,n_m\}\subset\NZ$, 
where
$$\begin{array}{l} F_qM=\{ 0\}~\hbox{if}~q<\min\{ n_1,\ldots ,n_m\} ;\\
F_qM=\sum^t_{i=1}\left (\sum_{p_i+n_i\le q}F_{p_i}A\right 
)\xi_i~\hbox{otherwise},\end{array}\quad q\in\NZ .$$ 
{\parindent=0pt\v5

{\bf 5.6.1. Proposition} Let  $M$ and $FM$ be as above. Consider the 
filtered free $A$-module $L=\oplus_{i=1}^mAe_i$ with the good 
filtration $FL=\{ F_qL\}_{q\in\NZ}$ such that $d_{\rm 
fil}(e_i)=n_i$, $1\le i\le m$, and the exact sequence of $A$-module 
homomorphisms
$$(*)\quad\quad 0~\mapright{}{}~N~\mapright{\iota}{}~L~\mapright{\varphi}{}~M~\mapright{}{}~0,$$
in which $\varphi (e_i)=\xi_i$, $1\le i\le m$, $N=$ Ker$\varphi$ 
with the induced filtration $FN=\{ F_qN=N\cap F_qL\}_{q\in\NZ}$, and 
$\iota$ is the inclusion map. The following statements hold.\par 
 
(i) The $A$-module homomorphisms $\iota$ and $\varphi$ are strict 
filtered homomorphisms.  Equipping $\OV L=L/N$ with the induced 
filtration $F\OV L=\{ F_q\OV L=(F_qL+N)/N\}_{q\in\NZ}$,  the induced 
$A$-module isomorphism $\OV L~\mapright{\OV\varphi}{}~M$ is a strict 
filtered  isomorphism, that is, $\OV L\cong M$ and $\OV\varphi$ 
satisfies $\OV\varphi(F_q\OV L)=F_qM$ for all $q\in\NZ$.
\par

(ii) The induced sequence
$$G(*)\quad\quad 0~\mapright{}{}~G(N)~\mapright{G(\iota )}{}~G(L)~\mapright{G(\varphi )}{}~G(M)~\mapright{}{}~0$$
is an exact sequence of graded $G(A)$-module homomorphisms, thereby 
$G(L)/G(N)\cong G(M)\cong G(\OV L)= G(L/N)$ as graded 
$G(A)$-modules. 
\par

(iii) The induced sequence
$$\widetilde{(*)}\quad\quad 0~\mapright{}{}~\T N~\mapright{\T\iota}{}~\T{L}~\mapright{\T{\varphi}}{}~\T M~\mapright{}{}~0$$
is an exact sequence of graded $\T A$-module homomorphisms, thereby 
$\T{L}/\T N\cong \T M\cong \widetilde{\OV L}=\widetilde{L/N}$ as 
graded $\T A$-modules.\vskip 6pt

{\bf Proof}  By the construction of $FL$ (see section 5.2) and 
Proposition 5.5.1, the proof of all assertions is a straightforward 
exercise.\QED}\v5

Proposition 5.6.1(i) enables us to make the 
following{\parindent=0pt\v5

{\bf Convention.} In what follows we shall always assume that a 
finitely generated $A$-module $M$ is of the form as presented in 
Proposition 5.6.1(i), i.e., $M=L/N$, and $M$  has the good 
filtration $$FM=\{F_qM=(F_qL+N)/N\}_{q\in\NZ} .$$}\v5

Comparing with the classical minimal graded free resolutions defined 
for finitely generated graded modules over finitely generated  
$\NZ$-graded  algebras with the degree-0 homogeneous part a field 
(cf. [Eis], Chapter 19; [Kr1], Chapter 3; [Li3]), the results 
obtained in previous sections and the preliminary we made above 
naturally motivate the following{\parindent=0pt\v5

{\bf 5.6.2. Definition}  Let $L_0=\oplus_{i=1}^mAe_i$ be a filtered 
free $A$-module with the filtration $FL_0=\{ F_qL_0\}_{q\in\NZ}$ 
such that $d_{\rm fil}(e_i)=b_i$, $1\le i\le m$, let $N$ be a 
submodule of $L_0$,  and let the $A$-module $M=L_0/N$ be equipped 
with the filtration $FM=\{F_qM=(F_qL_0+N)/N\}_{q\in\NZ}$. A {\it 
minimal filtered free resolution} of $M$ is an exact sequence of 
filtered $A$-modules and filtered homomorphisms
$${\cal L}_{\bullet}\quad\quad\cdots~\mapright{\varphi_{i+1}}{}~L_i~\mapright{\varphi_i}{}~\cdots ~
\mapright{\varphi_2}{}~L_1~\mapright{\varphi_1}{}~L_0~\mapright{\varphi_0}{}~M~\r~0$$ 
satisfying{\parindent=1.3truecm\par

\item{(1)} $\varphi_0$ is the canonical epimorphism, i.e.,
$\varphi_0(e_i)=\OV e_i$ for $e_i\in \mathscr{E}_0=\{ e_1,\ldots 
,e_m\}$ (where each $\OV e_i$ denotes the coset represented by $e_i$ 
in $M$), such that $\varphi_0(\mathscr{E}_0)$ is a minimal F-basis 
of $M$ with respect to  $FM$ (in the sense of Definition 5.4.3);
\par

\item{(2)} for $i\ge 1$,  each $L_i$ is a filtered free
$A$-module with finite $A$-basis $\mathscr{E}_i$,  and each 
$\varphi_i$ is a strict filtered homomorphism, such that 
$\varphi_i(\mathscr{E}_i)$ is a minimal standard basis of 
Ker$\varphi_{i-1}$ with respect to the filtration induced by 
$FL_{i-1}$ (in the sense of Definition 5.4.6). \v5}}

To see that the minimal filtered free resolution introduced above is 
an appropriate definition of ``minimal free resolution" for finitely 
generated modules over the $\NZ$-filtered solvable polynomial 
algebras with filtration determined by positive-degree functions, we 
now show that such a resolution is unique up to a strict filtered 
isomorphism of chain complexes in the category of $\NZ$-filtered 
$A$-modules. {\parindent=0pt\v5

{\bf 5.6.3. Theorem}  Let ${\cal L}_{\bullet}$ be a minimal filtered 
free resolution of $M$ as presented in Definition 5.6.2. The 
following statements hold.\par

(i) The induced sequence of graded $G(A)$-modules and graded 
$G(A)$-module homomorphisms
$$G({\cal L}_{\bullet})\quad\quad\cdots~\mapright{G(\varphi_{i+1})}{}~G(L_i)~\mapright{G(\varphi_i)}{}~\cdots ~
\mapright{G(\varphi_2)}{}~G(L_1)~\mapright{G(\varphi_1)}{}~G(L_0)~\mapright{G(\varphi_0)}{}~G(M)~\r~0$$ 
is a minimal graded free resolution of the finitely generated graded 
$G(A)$-module $G(M)$.\par

(ii) The induced sequence of graded $\T A$-modules and graded $\T 
A$-module homomorphisms
$$\widetilde{\cal L}_{\bullet}\quad\quad\cdots~\mapright{\T\varphi_{i+1}}{}~\T L_i~\mapright{\T\varphi_i}{}
~\cdots ~ \mapright{\T\varphi_2}{}~\T 
L_1~\mapright{\T\varphi_1}{}~\T L_0~\mapright{\T\varphi_0}{}~\T 
M~\r~0$$ is a minimal graded free resolution of the finitely 
generated graded $\T A$-module $\T M$.

(iii) ${\cal L}_{\bullet}$ is uniquely determined by $M$ in the 
sense that if $M$ has another minimal filtered free resolution
$${\cal L'}_{\bullet}\quad\quad\cdots~\mapright{\varphi'_{i+1}}{}~L'_i~\mapright{\varphi'_i}{}
~\cdots ~ 
\mapright{\varphi'_2}{}~L_1'~\mapright{\varphi'_1}{}~L_0~\mapright{\varphi_0}{}~M~\r~0$$ 
then for each $i\ge 1$, there is a strict filtered $A$-module 
isomorphism  $\chi_i$: $L_i\r L_i'$ such that the diagram
$$\begin{array}{ccc} L_i&\mapright{\varphi_i}{}&L_{i-1}\\
\mapdown{\chi_i}\scriptstyle{\cong}&&\mapdown{\chi_{i-1}}\scriptstyle{\cong}\\
L_i'&\mapright{\varphi_i'}{}&L_{i-1}'\end{array}$$ is commutative, 
thereby $\{\chi_i~|~i\ge 1\}$ gives rise to a strict filtered 
isomorphism of chain complexes  ${\cal L}_{\bullet}\cong {\cal 
L}_{\bullet}'$ in the category of $\NZ$-filtered $A$-modules. \vskip 
6pt

{\bf Proof} (i) and (ii) follow from Proposition 5.4.4, Proposition 
5.4.7, and Proposition 5.6.1.}\par

To prove (iii), let the sequence ${\cal L}_{\bullet}'$ be as 
presented above. By the assertion (i), $G({\cal L}_{\bullet}')$ is 
another minimal graded free resolution of $G(M)$. It follows from 
the well-known result on minimal graded free resolutions ([Eis], 
Chapter 19; [Kr1], Chapter 3; [Li3]) that there is a graded 
isomorphism of chain complexes $G({\cal L}_{\bullet})\cong G({\cal 
L}_{\bullet}')$ in the category of graded $G(A)$-modules, i.e., for 
each $i\ge 1$, there is a graded $G(A)$-modules isomorphism  
$\psi_i$: $G(L_i)\r G(L_i')$ such that the diagram
$$\begin{array}{ccc} G(L_i)&\mapright{G(\varphi_i)}{}&G(L_{i-1})\\
\mapdown{\psi_i}\scriptstyle{\cong}&&\mapdown{\psi_{i-1}}\scriptstyle{\cong}\\
G(L_i')&\mapright{G(\varphi_i')}{}&G(L_{i-1}')\end{array}$$ is 
commutative. Our aim below is to construct the desired strict 
filtered isomorphisms $\chi_i$ by using the graded isomorphisms 
$\psi_i$ carefully. So, starting with $L_0$, we assume that we have 
constructed the strict filtered isomorphisms 
$L_j~\mapright{\chi_j}{}~L_j'$, such that $G(\chi_j)=\psi_j$ and 
$\chi_{j-1}\varphi_j=\varphi_j'\chi_j$, $1\le j\le i-1$. Let 
$L_i=\oplus_{j=1}^{s_i}Ae_{i_j}$. Since each $\psi_i$ is a graded 
isomorphism, we have $\psi_i(\sigma (e_{i_j}))=\sigma (\xi_j')$ for 
some $\xi_j'\in L_i'$ satisfying $d_{\rm fil}(\xi_j')=d_{\rm 
gr}(\sigma (\xi_j'))=d_{\rm gr}(\sigma (e_{i_j}))=d_{\rm 
fil}(e_{i_j})$. It follows that if we construct the  filtered 
$A$-module homomorphism $L_i~\mapright{\tau_i}{}~L_i'$ by setting 
$\tau_i(e_{i_j})=\xi_j'$, $1\le j\le s_i$, then $G(\tau_i)=\psi_i$. 
Hence, $\tau_i$ is a strict filtered isomorphism by Proposition 
5.5.1. Since $\psi_{i-1}=G(\chi_{i-1})$, $\psi_i=G(\tau_i)$, and 
thus
$$\begin{array}{l}
\psi_{i-1}G(\varphi_i)=G(\chi_{i-1})G(\varphi_i)=G(\chi_{i-1}\varphi_i)\\
G(\varphi_i')\psi_i=G(\varphi_i')G(\tau_i)=G(\varphi_i'\tau_i),\end{array}$$ 
for each $q\in\NZ$, by the strictness of the $\varphi_j$, 
$\varphi_j'$, $\chi_j$ and $\tau_i$, we have
$$\begin{array}{rcl} G(\chi_{i-1}\varphi_i)(G(L_i)_q)&=&((\chi_{i-1}\varphi_i)(F_qL_i)+F_{q-1}L_{i-1}')/F_{q-1}L_{i-1}'\\
&=&(\chi_{i-1}(\varphi_i(L_i)\cap
F_qL_{i-1})+F_{q-1}L_{i-1}')/F_{q-1}L_{i-1}'\\
&\subseteq&((\chi_{i-1}\varphi_i)(L-i)\cap\chi_{i-1}(F_qL_{i-1})+F_{q-1}L_{i-1}')/F_{q-1}L_{i-1}'\\
&=&((\chi_{i-1}\varphi_i)(L_i)\cap 
F_qL_{i-1}'+F_{q-1}L_{i-1}')/F_{q-1}L_{i-1}',\end{array}$$

$$\begin{array}{rcl} G(\varphi_i'\tau_i)(G(L_i)_q)&=&((\varphi_i'\tau_i)(F_qL_i)+F_{q-1}L_{i-1}')/F_{q-1}L_{i-1}'\\
&=&(\varphi_i'(F_qL_i')+F_{q-1}L_{i-1}')/F_{q-1}L_{i-1}'\\
&=&(\varphi_i'(L_i')\cap
F_qL_{i-1}'+F_{q-1}L_{i-1}')/F_{q-1}L_{i-1}'\\
&=&((\varphi_i'\tau_i)(L_i)\cap 
F_qL_{i-1}'+F_{q-1}L_{i-1}')/F_{q-1}L_{i-1}'.\end{array}$$ Note that 
by the exactness of ${\cal L}_{\bullet}$ and ${\cal L}_{\bullet}'$, 
as well as the commutativity 
$\chi_{i-2}\varphi_{i-1}=\varphi_{i-1}'\chi_{i-1}$, we have 
$(\chi_{i-1}\varphi_i)(L_i)\subseteq 
\varphi_i'(L_i')=(\varphi_i'\tau_i)(L_i)$. Considering the 
filtration of the submodules $(\chi_{i-1}\varphi_i)(L_i)$ and 
$(\varphi_i'\tau_i)(L_i)$ induced by the filtration $FL_{i-1}'$ of 
$L_{i-1}'$, the commutativity 
$\psi_{i-1}G(\varphi_i)=G(\varphi_i')\psi_i$ and the formulas 
derived above show that both submodules have the same associated 
graded module, i.e., 
$G((\chi_{i-1}\varphi_i)(L_i))=G((\varphi_i'\tau_i)(L_i))$. It 
follows from a similar proof of  ([LVO], Theorem 5.7 on P.49) that
$$(\chi_{i-1}\varphi_i)(L_i)=(\varphi_i'\tau_i)(L_i).\eqno{(1)}$$ Clearly, the equality $(1)$
does not necessarily mean the commutativity of the diagram
$$\begin{array}{ccc} L_i&\mapright{\varphi_i}{}&L_{i-1}\\
\mapdown{\tau_i}\scriptstyle{\cong}&&\mapdown{\chi_{i-1}}\scriptstyle{\cong}\\
L_i'&\mapright{\varphi_i'}{}&L_{i-1}'\end{array}$$ To remedy this 
problem, we need to further modify the filtered isomorphism 
$\tau_i$. Since $G(\chi_{i-1}\varphi_i)(\sigma 
(e_{i_j}))=G(\varphi_i'\tau_i)(\sigma (e_{i_j}))$, $1\le j\le s_i$, 
if $d_{\rm fil}(e_{i_j})=b_j$, then by the equality (1) and the 
strictness of $\tau_i$ and $\varphi_i$ we have
$$\begin{array}{rcl} (\chi_{i-1}\varphi_i)(e_{i_j})-(\varphi_i'\tau_i)(e_{i_j})&\in&
(\varphi_i'\tau_i)(L_i)\cap F_{b_j-1}L_{i-1}'\\
&=&\varphi_i'(L_i')\cap F_{b_j-1}L_{i-1}'\\
&=&\varphi_i'(F_{b_j-1}L_i')\\
&=&(\varphi_i'\tau_i)(F_{b_j-1}L_i),\end{array}\eqno{(2)}$$ and 
furthermore from (2) we have a $\xi_j\in F_{b_j-1}L_i$ such that 
$d_{\rm fil}(e_{i_j}-\xi_j)=b_j$ and
$$(\chi_{i-1}\varphi_i)(e_{i_j})=(\varphi_i'\tau_i)(e_{i_j}-\xi_j), ~~~
1\le j\le s_i.\eqno{(3)}$$ Now, if we construct the filtered 
homomorphism $L_i~\mapright{\chi_i}{}~L_i'$ by setting 
$\chi_i(e_{i_j})=\tau_i(e_{i_j}-\xi_j)$, $1\le j\le s_i$, then since 
$\tau_i(\xi_j)\in F_{b_j-1}L_i'$, it turns out that
$$G(\chi_i)(\sigma (e_{i_j}))=G(\tau_i)(\sigma (e_{i_j}-\xi_j))=G(\tau_i)(\sigma (e_{i_j}))=\psi_i(\sigma (e_{i_j})),
~~~1\le j\le s_i, $$ thereby $G(\chi_i)=\psi_i$. Hence, $\chi_i$ is 
a strict filtered isomorphism by Proposition 5.5.1, and moreover, it 
follows from (3) that we have reached the following diagram
$$\begin{array}{ccccc} \cdots~\mapright{\varphi_{i+1}}{}&L_i&\mapright{\varphi_i}{}
&L_{i-1}&\mapright{\varphi_{i-1}}{}~\cdots\\
&\mapdown{\chi_i}\scriptstyle{\cong}&&\mapdown{\chi_{i-1}}\scriptstyle{\cong}&\\
\cdots~\mapright{\varphi_{i+1}'}{}&L_i'&\mapright{\varphi_i'}{}&L_{i-1}'&\mapright{\varphi_{i-1}'}{}~\cdots\end{array}$$ 
in which $\chi_{i-1}\varphi_i=\varphi_i'\chi_i$. Repeating the same 
process to getting the desired $\chi_{i+1}$ and so on, the proof is 
thus finished. \v5

\section*{5.7. Computation of  Minimal Finite Filtered Free
Resolutions}

Let  $A=K[a_1,\ldots ,a_n]$ be a solvable polynomial algebra with 
the  admissible system $(\B ,\prec_{gr})$ in which $\prec_{gr}$ is a 
graded monomial ordering with respect to some given positive-degree 
function $d(~)$ on  $A$ (see section 1.1 of Chapter 1). Thereby $A$ 
is turned into an $\NZ$-filtered solvable polynomial algebra with  
the filtration $FA=\{ F_pA\}_{p\in\NZ}$ constructed with respect to 
the same $d(~)$ (see Example (2) given in Section 5.1). Note that 
Theorem 3.1.1, Theorem 3.1.2, and Theorem 3.2.2 given in Chapter 3  
hold true for {\it any} solvable polynomial algebra.  Combining the 
results of Chapter 4 and previous sections of the current chapter, 
we are now able to work out the algorithmic procedures for computing  
minimal finite filtered free resolutions over $A$ (in the sense of 
Definition 5.6.2) with respect to any graded left monomial ordering 
on free left modules. All notions, notations and conventions used 
before are maintained. {\parindent=0pt\v5

{\bf 5.7.1. Theorem}  Let $L_0=\oplus_{i=1}^mAe_i$ be a filtered 
free $A$-module with the filtration $FL_0=\{ F_qL_0\}_{q\in\NZ}$ 
such that $d_{\rm fil}(e_i)=b_i$, $1\le i\le m$. If 
$N=\sum^s_{i=1}A\xi_i$ is a finitely generated submodule of $L_0$ 
and the quotient module $M=L_0/N$ is equipped with the filtration 
$FM=\{F_qM=(F_qL_0+N)/N\}_{q\in\NZ}$, then $M$ has a minimal 
filtered free resolution of length $d\le n$:
$${\cal L}_{\bullet}\quad\quad 0~\mapright{}{}~L_d~\mapright{\varphi_{q}}{}~
\cdots~\mapright{\varphi_2}{}~L_1~\mapright{\varphi_1}{}~L_0~\mapright{\varphi_0}{}~M~\mapright{}{}~0$$ 
which can be constructed by implementing the following 
procedures:}\par

{\bf Procedure 1.} Fix a graded left monomial ordering 
$\prec_{e\hbox{-}gr}$ on the $K$-basis $\BE$ of $L_0$ (see Section 
5.3), and run {\bf Algorithm-LGB} (presented in Section 2.3 of 
Chapter 2) with the initial input data $ U =\{ \xi_1,\ldots 
,\xi_s\}$ to compute a left Gr\"obner basis $\G=\{ g_1,\ldots 
,g_z\}$ for $N$, so that $N$ has the standard basis $\G$ with 
respect to the induced filtration $FN=\{ F_qN=N\cap 
F_qL_0\}_{q\in\NZ}$ (Theorem 5.4.8).\par

{\bf Procedure 2.} Run {\bf Algorithm-MINFB} (presented in  
Proposition 5.5.3) with the initial input data $E=\{ e_1,\ldots 
,e_m\}$ and $\G=\{ g_1,\ldots ,g_z\}$ to compute a subset 
$\mathscr{E}_0'=\{ e_{i_1},\ldots ,e_{i_{m'}}\}\subset 
\mathscr{E}_0=\{ e_1,\ldots ,e_m\}$ and a subset $V=\{ v_1,\ldots 
,v_t\}\subset N\cap L_0'$ such that there is a strict filtered 
isomorphism $L_0'/N'=M'\cong M$, where 
$L_0'=\oplus_{q=1}^{m'}Ae_{i_q}$ with $m'\le m$ and 
$N'=\sum^{t}_{k=1}Av_k$, and such that $\{\OV e_{i_1},\ldots ,\OV 
e_{i_{m'}}\}$ is a minimal F-basis of $M$ with respect to the 
filtration $FM$.\par

For convenience, after accomplishing Procedure 2 we may assume that 
$\mathscr{E}_0=\mathscr{E}_0'$, $ U =V$ and $N=N'$. Accordingly we 
have the short exact sequence $$0~\mapright{}{}~ 
N~\mapright{\iota}{}~ L_0~\mapright{\varphi_0}{}~M~\mapright{}{}  
~0$$ such that $\varphi_0 (\mathscr{E}_0)=\{ \OV e_1,\ldots ,\OV 
e_m\}$  is a minimal  F-basis of $M$ with respect to the filtration 
$FM$, where $\iota$ is the inclusion map. \par

{\bf Procedure 3.} With the initial input data $ U =V$, implements 
the procedures presented in Theorem 5.5.4 to compute a minimal 
standard basis $W=\{ \xi_{j_1},\ldots ,\xi_{j_{m_1}}\}$ for $N$ with 
respect to the induced filtration $FN$.
\par

{\bf Procedure 4.} Computes a generating set $U_1  =\{\eta_1,\ldots 
,\eta_{s_1}\}$ of $N_1=$ Syz$(W)$ in the free $A$-module 
$L_1=\oplus_{i=1}^{m_1}A\varepsilon_i$ by running {\bf 
Algorithm-LGB} with the initial input data $W$ and using Theorem 
3.1.2.\par

{\bf Procedure 5.} Construct the strict filtered exact sequence
$$0~\mapright{}{}~ N_1~\mapright{}{}~L_1~\mapright{\varphi_1}{}~
L_0~\mapright{\varphi_0}{}~M~\mapright{}{}  ~0$$ where the 
filtration $FL_1$ of $L_1$ is constructed by setting $d_{\rm 
fil}(\varepsilon_k)=d_{\rm fil}(\xi_{j_k})$, $1\le k\le m_1$, and 
$\varphi_1$ is defined by setting 
$\varphi_1(\varepsilon_k)=\xi_{j_k}$, $1\le k\le m_1$. If $N_1\ne 
0$, then, with the initial input data $ U =U_1$,  repeat Procedure 3 
-- Procedure 5 for $N_1$ and so on.\par

By Theorem 5.6.3, a minimal filtered free resolution ${\cal 
L}_{\bullet}$ of $M$ gives rise to a minimal graded free resolution 
$G({\cal L}_{\bullet})$ of $G(M)$. Since $G(A)=K[\sigma (a_1),\ldots 
,\sigma (a_n)]$ is a solvable polynomial algebra by Theorem 5.1.5, 
it follows from  Theorem 4.4.1 that $G({\cal L}_{\bullet})$ 
terminates at a certain step, i.e., Ker$G(\varphi_d)=0$ for some 
$d\le n$. But Ker$G(\varphi_d)=G(\hbox{Ker}\varphi_d)$ by 
Proposition 5.5.1, where Ker$\varphi_d$ has the filtration induced 
by $FL_d$, thereby $G(\hbox{Ker}\varphi_d)=0$. Consequently 
Ker$\varphi_d=0$ since all filtration we are dealing with are 
separated, thereby a minimal finite filtered free resolution of 
length $d\le n$ is achieved for $M$. \v5

\v5 \centerline{References}
\parindent=1.6truecm

\item{[AF]} F.W. Anderson and K.R. Fuller, {\it Rings and categories of Modules}. Springer-Verlag, 1974.

\item{[AL1]} J. Apel and W. Lassner, An extension of Buchberger's
algorithm and calculations in enveloping fields of Lie algebras.
{\it J. Symbolic Comput}., 6(1988), 361--370.

\item{[AL2]} W. W. Adams and P. Loustaunau, {\it An Introduction to Gr\"obner Bases}.
Graduate Studies in Mathematics, Vol. 3. American Mathematical
Society, 1994.

\item{[AVV]} M. J. Asensio,  et al, 
A new algebraic approach to microlocalization of filtered rings. 
Trans. Amer. math. Soc, 2(316)(1989), 537-555.

\item{[Ber1]} R.~Berger, The quantum Poincar\'e-Birkhoff-Witt theorem.
{\it Comm. Math. Physics}, 143(1992), 215--234.

\item{[Ber2]} G. Bergman, The diamond lemma for ring theory. {\it Adv.
Math}., 29(1978), 178--218.

\item{[Bu1]} B. Buchberger, {\it Ein Algorithmus zum Auffinden der
Basiselemente des Restklassenringes nach einem nulldimensionalen
polynomideal}. PhD thesis, University of Innsbruck, 1965.

\item{[Bu2]} B. Buchberger, Gr\"obner bases: ~An algorithmic method
in polynomial ideal theory. In: {\it Multidimensional Systems
Theory} (Bose, N.K., ed.), Reidel Dordrecht, 1985, 184--232.

\item{[BW]} T.~Becker and V.~Weispfenning, {\it Gr\"obner Bases}.
Springer-Verlag, 1993.

\item{[CDNR]} A. Capani,  et al, 
Computing minimal finite free resolutions. {\it Journal of Pure and 
Applied Algebra}, (117\& 118)(1997), 105 -- 117.

\item{[DGPS]} W. Decker,  et al, \newblock {\sc Singular} {3-1-3}
--- {A} computer algebra system for polynomial computations. Available at \par 
\newblock {http://www.singular.uni-kl.de}(2011). 

\item{[Dir]} P. A. M.~Dirac, On quantum algebra. {\it Proc. Camb.
Phil. Soc}., 23(1926), 412--418.

\item{[Eis]} D. Eisenbud, {\it Commutative Algebra with a View
Toward to Algebraic Geometry}, GTM 150. Springer, New York, 1995.

\item{[EPS]} D. Eisenbud,  et al,  {\it Non-commutative
Gr\"obner bases for commutative algebras}. Proc. Amer. Math. Soc., 
126, 1998, pp. 687-691.\par

\item{[Fau]} J.-C. Faug\'ere. A new efficient algorithm for computing Gr\"oobner bases without
reduction to zero (F5). In: {\it proc. ISSAC'02}, ACM Press, New
York, USA, 75-82, 2002.

\item{[Fr\"ob]} R. Fr\"oberg, {\it An Introduction to Gr\"obner Bases}. Wiley, 1997.

\item{[Gal]} A. Galligo, Some algorithmic questions on ideals of
differential operators. {\it Proc. EUROCAL'85}, LNCS 204, 1985,
413--421.

\item{[Gol]} E. S. Golod, Standard bases and homology, in: {\it Some
Current Trends in Algebra}. (Varna, 1986), Lecture Notes in
Mathematics, 1352, Springer-Verlag, 1988, 88-95.

\item{[Gr]~} E. Green, {\it Noncommutative} {\it Gr\"obner Bases},
{\it A Computational} and  {\it Theoretical Tool}. Lecture notes,
Dec. 15, 1996. Available at
www.math.unl.edu/~shermiller2/hs/green2.ps

\item{[GV]} J. Gago-Vargas, Bases for projective modules in $A_n(k)$.
{\it Journal of Symbolic Computation}, 36(2003), 845¨C853.

\item{[Hay]} T. Hayashi, $Q$-analogues of Clifford and Weyl
algebras-Spinor and oscillator representations of quantum enveloping
algebras. {\it Comm. Math. Phys}., 127(1990), 129--144.

\item{[Hu]} J. E. Humphreys, {\it Introduction to Lie Algebras and Representation Theory}. Springer, 1972.

\item{[Jat]} V.A. Jategaonkar, A multiplicative analogue of the Weyl
algebra. {\it Comm. Alg}., 12(1984), 1669--1688.

\item{[JBS]} A. Jannussis, et al,  Remarks on
the $q$-quantization. {\it Lett. Nuovo Cimento}, 30(1981), 123--127.

\item{[Kr1]} U. Kr\"ahmer,  Notes on Koszul algebras. 2010. Available at \\
 http://www.maths.gla.ac.uk/~ukraehmer/connected.pdf

\item{[Kr2]} H. Kredel, {\it Solvable Polynomial Rings}. Shaker-Verlag, 1993.

\item{[KP]} H. Kredel and M.  Pesch,  MAS, modula-2 algebra system. 1998.
 Available at http://krum.rz.uni-mannheim.de/mas/

\item{[KR1]} M. Kreuzer, L. Robbiano, {\it Computational Commutative Algebra 1}. Springer, 2000.

\item{[KR2]} M. Kreuzer, L. Robbiano, {\it Computational Commutative Algebra 2}. Springer, 2005.

\item{[K-RW]} A.~Kandri-Rody and V.~Weispfenning, Non-commutative
Gr\"obner bases in algebras of solvable type. {\it J. Symbolic
Comput.}, 9(1990), 1--26.

\item{[Kur]} M.V. Kuryshkin, Op\'erateurs quantiques
g\'en\'eralis\'es de cr\'eation et d'annihilation. {\it Ann. Fond.
L. de Broglie}, 5(1980), 111--125.

\item{[Lev]} V. Levandovskyy, {\it Non-commutative Computer Algebra for
Polynomial Algebra}: {\it Gr\"obner Bases, Applications and
Implementation}. Ph.D. Thesis, TU Kaiserslautern, 2005.

\item{[Li1]} H. Li, {\it Noncommutative Gr\"obner Bases and
Filtered-Graded Transfer}. Lecture Note in Mathematics, Vol. 1795,
Springer-Verlag, 2002.

\item{[Li2]} H. Li, {\it Gr¡§obner Bases in Ring Theory}. World Scientific
Publishing Co., 2011.

\item{[Li3]} H. Li, On monoid graded local rings. {\it Journal of Pure and Applied Algebra},
216(2012), 2697 -- 2708.

\item{[Li4]} H. Li, A Note on Solvable Polynomial Algebras.
{\it Computer Science Journal of Moldova}, 1(64)(2014), 99--109. 
Available at\par arXiv:1212.5988 [math.RA].

\item{[Li5]} H. Li, Computation of minimal graded free resolutions over
$\NZ$-graded solvable polynomial algebras. Available at \par
arXiv:1401.5206 [math.RA]

\item{[Li6]} H. Li, Computation of minimal filtered free resolutions over
$\NZ$-filtered solvable polynomial algebras. Available at \par  
arXiv:1401.5464 [math.RA]

\item{[Li7]} H. Li, On Computation of Minimal Free Resolutions over Solvable
Polynomial Algebras. (pp66) accepted by {\it Commentationes 
Mathematicae Universitatis Carolinae}, to appear.

\item{[LVO]} H. Li and F. Van Oystaeyen, {\it Zariskian
Filtrations}. $K$-Monograph in Mathematics, Vol.2. Kluwer Academic
Publishers, 1996.

\item{[LW]} H. Li and  Y. Wu, ~Filtered-graded transfer of
Gr\"obner basis computation in solvable polynomial algebras. {\it
Communications in Algebra}, 1(28)(2000), 15--32.

\item{[Man]} Yu.I.~Manin, {\it Quantum Groups and Noncommutative
Geometry}. Les Publ. du Centre de R\'echerches Math., Universite de
Montreal, 1988.

\item{[Mor]} T. Mora, An introduction to commutative and
noncommutative Gr\"obner Bases, {\it Theoretic Computer Science},
134(1994), 131--173.

\item{[MP]} J.C. McConnell and J.J. Pettit, Crossed products and
multiplicative analogues of Weyl algebras. {\it J. London Math.
Soc}., (2)38(1988), 47--55.

\item{[MR]} J.C. McConnell and  J.C. Robson,  {\it Noncommutative Noetherian Rings}. Wiley-Interscience
Publication, 1987.

\item{[NVO]} C.~N$\check{\rm a}$st$\check{\rm a}$sescu and F.~Van 
Oystaeyen, {\it Graded ring theoey}, Math. Library 28, North 
Holland, Amsterdam, 1982.

\item{[Ros]} A.L.~Rosenberg, {\it Noncommutative Algebraic Geometry
and Representations of Quantized Algebras}. Kluwer Academic
Publishers, 1995.

\item{[Rot]} J.J. Rotman, {\it An introduction to homological algebra}. Academic Press, 1979.

\item{[Sch]} F.O. Schreyer, {\it Die Berechnung von Syzygien mit
dem verallgemeinerten Weierstrasschen Divisionsatz}. Diplomarbeit,
Hamburg, 1980.

\item{[Sm]} S.P.~Smith, Quantum groups: an introduction and survey
for ring theorists. in: {\it Noncommutative rings}, S.~Montgomery
and L.~Small eds., MSRI Publ. 24(1992), Springer-Verlag, New York,
131--178.

\item{[SWMZ]} Y. Sun, et al, A signature-based algorithm for computing Gr\"obner bases in
solvable polynomial algebras. In: {\it Proc. ISSAC'12}, ACM Press,
 351-358, 2012.

\item{[Uf]} V. Ufnarovski, Introduction to noncommutative
Gr\"obner basis theory. in: {\it Gr\"obner Bases and Applications}
(Linz, 1998), London Math. Soc. Lecture Notes Ser., 251, Cambridge
Univ. Press, Cambridge, 1998, 259--280.

\item{[Wal]} R. Wallisser, Rationale approximation der $q$-analogues
der exponentialfunktion und Irrationalit\"atsaussagen f\"ur diese
Funktion. {\it Arch. Math}., 44(1985), 59--64.

\item{[Wik1]} Decision problem (Solvable problem). Available at\par

https:\hbox{//}en.wikipedia.org\hbox{/}wiki\hbox{/}Solvable$_{-}$problem

\item{[Wik2]} Confluence (abstract rewriting). Available at\par

https:\hbox{//}en.wikipedia.org\hbox{/}wiki\hbox{/}Confluence$_{-}$(term$_{-}$rewriting)

\end{document}